\documentclass[11pt,leqno,twoside]{amsart}

\usepackage{enumerate}

\usepackage{enumitem}

\usepackage{amssymb}
\usepackage[english]{babel}
\usepackage{amssymb,amsthm,amsmath,eucal,mathrsfs}
\usepackage{bm}

\usepackage{amsmath}
\usepackage{amsfonts}
\usepackage{amsmath}
\usepackage{amssymb}
\usepackage{graphicx}
\usepackage{lscape}
\usepackage{amstext}
\usepackage{amsthm}
\usepackage{color}
\usepackage{float}
\usepackage{mathrsfs}
\usepackage{epsfig}
\usepackage{url}
\usepackage{fancyhdr}
\usepackage{pspicture}
\usepackage{graphicx}
\usepackage{tikz}
\usepackage{cite}

\usepackage{amsmath,verbatim}

\usepackage{amsthm}

\usepackage{amssymb}

\usepackage{amsfonts}

\usepackage{dsfont}

\usepackage{hyperref}

\newcommand{\Inc}{\mathrm{Inc}}

\newtheorem{theorem}{Theorem}[section]

\newtheorem{proposition}[theorem]{Proposition}

\newtheorem{corollary}[theorem]{Corollary}

\newtheorem{assumption}[theorem]{Assumptions}
\newtheorem{notation}[theorem]{Notations}

\theoremstyle{definition}

\newtheorem{definition}[theorem]{Definition}

\newtheorem{remark}[theorem]{Remark}

\numberwithin{equation}{section}

\newcommand{\diam}{\mathrm{diam}}      % diameter

\renewcommand{\Re}{{\ensuremath{\mathrm{Re\,}}}} %Realteil nicht als fraktur

\renewcommand{\div}{\mathrm{div}\,}    %div anstatt geteilt

\renewcommand{\O}{\mathcal{O}\,}

\newcommand{\R}{\mathbb{R}}

\newcommand{\C}{\mathbb{C}}

\newcommand{\LL}{\mathbb{L}}

\newcommand{\PP}{\mathbb{P}}

\newcommand\restr[2]{{% we make the whole thing an ordinary symbol
  \left.\kern-\nulldelimiterspace % automatically resize the bar with \right
  #1 % the function
  \vphantom{\big|} % pretend it's a little taller at normal size
  \right|_{#2} % this is the delimiter
  }}

\usepackage{eucal}

 \newcommand{\mk}[1]{{\color{black}{#1}}}

\title[Electromagnetic Scattering by a cluster of dimers]{Electromagnetic Scattering by a Cluster of Hybrid Dielectric--Plasmonic Dimers}

\author[Cao, Ghandriche and Sini]{Xinlin Cao $^*$ Ahcene Ghandriche $^{**}$ and Mourad Sini$^{\ddag}$}
\thanks{$^*$ Department of Applied Mathematics, The Hong Kong Polytechnic University, Hong Kong SAR. Email: xinlin.cao@polyu.edu.hk. This author is partially supported by the National Natural Science Foundation of China (NSFC): P0055789 and PolyU Internal Research Fund: P0050307.} 
\thanks{$^{**}$ Department of Mathematics, Kuwait University, Kuwait. Email: gh.hsen.math@gmail.com.}
\thanks{$^{\ddag}$ RICAM, Austrian Academy of Sciences, Altenbergerstrasse 69, A-4040, Linz, Austria. Email: mourad.sini@oeaw.ac.at. This author is partially supported by the Austrian Science Fund (FWF): P 30756-NBL and P 32660.}

\begin{document}

% \date{\today}

\allowdisplaybreaks

\begin{abstract}
We consider time-harmonic electromagnetic scattering by a cluster of hybrid dielectric--plasmonic
dimers in $\mathbb{R}^3$. Each dimer consists of a high-contrast dielectric nanoparticle and a
moderately contrasting plasmonic nanoparticle separated by a subwavelength distance. The cluster is
assumed to contain many such dimers whose size $a$ is small compared to the wavelength, with
intra-dimer and inter-dimer distances scaling like $a^{t_1}$ and $a^{t_2}$, and the frequency is tuned
near suitable electric and magnetic resonances of the associated Newtonian and magnetization
operators on the reference shapes.

Under these geometric, contrast and spectral assumptions, we derive a Foldy--Lax type approximation
for the Maxwell system. We show that the scattered field and its far field admit asymptotic
expansions in terms of four moments $(Q_{m_1},R_{m_1},Q_{m_2},R_{m_2})$ attached to each dimer,
which solve an explicit finite-dimensional linear system whose coefficients are expressed through
the free-space Green's functions and the polarization tensors of the reference particles. We prove
invertibility of this system under quantitative smallness conditions on the contrast and the dimer
density, and we obtain error estimates uniform in the number of dimers.

By extracting the dominant components, we further show that each hybrid dimer behaves, at leading
order, as a co-located electric and magnetic dipole driven by the local fields, and we identify the
corresponding $6\times 6$ polarizability matrix. This provides a discrete model for clusters of hybrid
dimers that is suitable for fast forward simulations, inverse schemes, and as input for effective-medium descriptions. In particular, it suggests parameter regimes where clusters of hybrid dimers can generate (double) negative effective permittivity and permeability and bi-anisotropic constitutive
laws and eventually hyperbolic media. 
\end{abstract}

%\keywords{\blue{Electromagnetic scattering; Foldy-Lax type approximation; Subwavelength resonanced; Hybrid dimer; Potential Theory}}


\maketitle

\tableofcontents

\section{Introduction and the main results}\label{prelimilary}

\subsection{Background and motivation}

The analysis of time-harmonic electromagnetic scattering by collections of small inclusions has a long
history, both in classical scattering theory and in the more recent theory of metamaterials and
effective media. On the one hand, the rigorous well-posedness of Maxwell boundary value problems
in non-smooth domains, and the mapping properties of the associated potential operators, are by now
well understood, see for instance the monographs \cite{colton2019inverse, Dautry-Lions, nedelec} and the works
\cite{amrouche-bernardi-dauge-girault,mitrea-mitrea-pipher,friedman-pasciak,raevskii}. On the
other hand, when the scatterers are small compared to the wavelength, the scattered field can often
be represented, up to a controlled error, by a finite system of multipoles whose strengths are
determined by a Foldy--Lax type algebraic system. This philosophy underlies a large part of the
mathematical theory of wave propagation in heterogeneous media and of effective medium theory.

In the context of electromagnetism, small plasmonic nanoparticles (with negative permittivity and
moderate permeability) are known to support strong electric resonances \cite{AmmariEtAl_Maxwell, AmmariEtAl_SPR_ARMA}, while dielectric particles
with high refractive index can support strong magnetic-type Mie resonances \cite{Ammai-Li-Zou-2023, cao-ghandriche-sini-highindex}. These resonant
inclusions can be arranged in clusters to generate effective media with non-trivial effective
permittivity $\varepsilon_{\rm eff}$ and permeability $\mu_{\rm eff}$, including regimes where one of
these tensors has negative real part. A rigorous analysis of such phenomena for clusters of
moderately contrasting plasmonic nanoparticles was carried out in \cite{cao-sini-effective}, where
the effective permittivity and permeability generated by a cloud of nanoparticles were derived,
together with precise error estimates. In \cite{cao-ghandriche-sini-highindex}, clusters of
dielectric nanoparticles with high refractive indices were analyzed and shown to generate strong
magnetic responses. These works show that appropriately designed clusters of small inclusions, even
in the quasi-static or mesoscopic regime, can exhibit effective properties that are markedly
different from those of the background medium.

Beyond the above works, there is by now a substantial mathematical literature on subwavelength resonators and the generation of effective negative material parameters. For plasmonic nanoparticles, Ammari and co-authors have developed a systematic small-volume and spectral theory based on polarization and moment tensors and Neumann–Poincaré type operators; see, for instance, the monograph on polarization tensors \cite{AmmariKang_PMT}, the analysis of scalar plasmonic resonances and surface plasmon resonance \cite{AmmariEtAl_SPR_ARMA}, and the full Maxwell treatment in \cite{AmmariEtAl_Maxwell}. These works show how, near plasmonic resonances, the effective permittivity can become negative and how the geometry and arrangement of nanoparticles control the resonant spectrum and the associated field enhancement. In a more directly metamaterial-oriented direction, the paper \cite{AmmariWuYu_ChiralDoubleNeg} gives formal mathematical arguments of double-negative electromagnetic metamaterials in chiral media, proving that a single type of plasmonic dielectric resonator, embedded in a suitably chiral background, can yield simultaneously negative effective permittivity and permeability near certain resonant frequencies. An acoustic analogue of this double-negative behavior, based on Minnaert-type bubble resonances, is developed in \cite{AmmariEtAl_DoubleNegAcoustic}.

Other approaches to negative-index and double-negative media have been developed using homogenization and spectral methods for high-contrast periodic composites. In \cite{ChenLipton_DoubleNeg}, Chen and Lipton introduce a generic class of metamaterials and show that the interplay between the Dirichlet spectrum of the high-dielectric phase and generalized electrostatic spectra of the complement yields frequency intervals with either double-negative or double-positive effective properties, together with corresponding Bloch waves in subwavelength structures. Bouchitt\'e and co-authors have analyzed artificial magnetism and negative refraction in periodic photonic crystals and split-ring geometries, rigorously deriving effective Maxwell systems with negative or strongly dispersive effective parameters \cite{BouchitteFelbacq_ArtificialMagnetism,BouchitteSchweizer_SplitRing}. Compared to these continuum homogenization frameworks, the present work is based on an explicit Foldy–Lax reduction to a finite-dimensional dipole model for a cluster of hybrid dimers. It provides a complementary route to double negativity and bi-anisotropy, where the effective response is encoded directly in the block structure of the discrete polarizability matrices associated with each dimer and their mutual interactions.

Hybrid structures, consisting of a plasmonic and a dielectric particle placed at a subwavelength
distance, have attracted considerable interest because they combine, in a single meta-atom, an
electric plasmonic resonance and a magnetic dielectric resonance. In \cite{cao-ghandriche-sini-dimer}
the scattering by a \emph{single} hybrid dielectric--plasmonic dimer was analyzed rigorously at the
level of the full Maxwell system. It was shown there that, by tuning the contrast parameters and the
frequency, one can create strong interactions between the electric and magnetic components of the
field, and that the dimer can be modeled at leading order by coupled electric and magnetic dipoles
whose strengths are controlled by suitable polarization tensors associated with the reference
shapes. In a different direction, the work \cite{GS} considered photo-acoustic
imaging in the presence of plasmonic contrast agents, again at the level of the full Maxwell system,
and showed how resonant nanoparticles can enhance the generated acoustic signals and thus improve
the stability of the inversion.

The present work is motivated by the following question: \emph{what is the effective behavior of a
\emph{cluster} of hybrid dielectric--plasmonic dimers when the size of each dimer is small compared
to the wavelength, but the number of dimers is large and their mutual interactions cannot be
neglected?} From the point of view of metamaterials, a cluster of hybrid dimers is a natural
candidate for generating effective media with simultaneous electric and magnetic resonances inside a
single unit cell, and hence for realizing double negativity (negative effective permittivity and
permeability) as well as more general bi-anisotropic responses. From the mathematical point of view,
this setting requires combining several ingredients:
\begin{itemize}
  \item a detailed understanding of the single-dimer problem, including the structure of the
        relevant polarization tensors and the location of the electric and magnetic resonances;
  \item a precise control of multiple scattering between many dimers, quantified in terms of the
        scaling of the size, the minimal distance, and the number of dimers;
  \item a systematic reduction of the full Maxwell system to a discrete Foldy--Lax type model,
        together with error estimates uniform in the number of dimers.
\end{itemize}

Our earlier works \cite{cao-sini-effective,cao-ghandriche-sini-highindex,cao-ghandriche-sini-dimer}
provide the first two ingredients in specific configurations: clusters of purely plasmonic
nanoparticles, clusters of purely dielectric nanoparticles, and a single hybrid dimer. The present
paper can be viewed as a next step in this program, in which we consider a cluster of hybrid dimers
and derive an explicit finite-dimensional model that captures, at leading order, the interactions
between the electric and magnetic components generated by each dimer and by the cluster as a whole.
More precisely, under suitable scaling and spectral assumptions, we show that the scattered field
and its far field can be expressed, up to an explicit remainder, in terms of a system of coupled
electric and magnetic dipole moments attached to each dimer, whose amplitudes are determined by a
linear algebraic system of Foldy--Lax type. The coefficients of this system are written explicitly
in terms of the free-space Green's functions and of the polarization tensors associated with the
reference shapes.

This discrete model has several consequences. First, it provides an efficient and quantitatively
accurate tool for computing the scattered and radiated fields generated by large clusters of hybrid
dimers without solving the full Maxwell boundary value problem. Second, it reveals in a transparent
way how the combination of plasmonic and dielectric components at the level of each dimer, together
with the geometry of the cluster, can give rise to non-trivial effective responses, including
negative $\varepsilon_{\rm eff}$, negative $\mu_{\rm eff}$ and non-zero magneto--electric coupling
tensors. Finally, it offers a natural starting point for a rigorous homogenization theory of
clusters of hybrid dimers, and for the analysis of collective resonances and band-structure
phenomena in periodic or quasi-periodic arrangements of such meta-atoms.

The main objective of this paper is therefore to establish a Foldy--Lax approximation for the
electromagnetic waves generated by a cluster of hybrid dielectric--plasmonic dimers, to quantify
its accuracy in terms of the relevant small parameters, and to analyze the resulting discrete system
with a view towards effective medium and metamaterial interpretations. In a subsequent work, we plan
to use the discrete model derived here to develop a homogenized description of such clusters and to
characterize rigorously the parameter regimes in which double negativity, strong anisotropy and
bi-anisotropic behavior can be achieved.

% Recently, and after closing this work, we got aware of the work \cite{Ammari-Li-Zou-2} where the authors derived similar approximation formulas in the case of a single dielectric nanoparticle. 
\subsection{Main results}
%%%%%%%%%%%%%%%%%%%%%%Main Model%%%%%%%%%%%%%%%%%%

Let $D \, := \, \overset{\aleph}{\underset{m=1}{\cup}} D_{m}$ be a set of coupled Lipschitz domains, representing for a cluster of dimers such that each $D_m:=D_{m_1}\cup D_{m_2}$ for $m=1, 2, \cdots, \aleph$, with $D_{m_1}$ and $D_{m_2}$ being respectively a dieletric nanoparticle with non-magnetic material, and a plasmonic nanoparticle with the moderately contrasting permittivity and permeability.
%
%bounded and Lipschitz-regular domain in $\mathbb{R}^3$, standing for a dimer composed by two nanoparticles $D_1$ and $D_2$. We assume that $D_1$ is a dieletric nanoparticle with non-magnetic material, namely by denoting the relative permittivity and permeability in $D_1$ as $\epsilon_r^{(1)}:=\frac{\epsilon^{(1)}}{\epsilon_0}$ and $\mu_r:=\frac{\mu^{(1)}}{\mu_0}$, we have $\epsilon_r^{(1)}=1$ outside $D_1$ while $\mu_r^{(1)}=1$ in the whole space $\mathbb{R}^3$. $D_2$ is a plasmonic nanoparticle with the moderately contrasting permittivity and permeability satisfying $\epsilon_{r}^{(2)}\sim 1$ ($\epsilon_{r}^{(2)}\neq 1$) and $\mu_2^{(2)}=1$ in $\mathbb{R}^3$. 

The electromagnetic scattering of time-harmonic plane waves from the cluster of dimers $D$ states as follows.
\begin{equation}\label{model-m}
\begin{cases}
\mathrm{Curl} \, E^T-i k \mu_rH^T=0\quad\mbox{in}\ \mathbb{R}^3,\\	
\mathrm{Curl} \, H^T+i k \epsilon_r E^T=0\quad\mbox{in}\ \mathbb{R}^3,\\
E^T=E^{in}+E^s,\quad H^T=H^{in}+H^s,\\
%	\nu\times E|_{+}=\nu\times E|_{-}, \quad \nu\times H|_{+}=\nu\times H|_{-},\quad\mbox{on} \ \partial D,\\
\sqrt{\mu_0\epsilon_0^{-1}}H^s\times\frac{x}{|x|}-E^s=\mathcal{O}({\frac{1}{|x|^2}}),\quad \mbox{as} \ |x|\rightarrow\infty,
\end{cases}
\end{equation}
where $\epsilon_0$ and $\mu_0$ are respectively the electric permittivity and the magnetic permeability of the vacuum outside $D$, $\epsilon_r \, := \, \dfrac{\epsilon}{\epsilon_0}$ and $\mu_r \, := \, \dfrac{\mu}{\mu_0}$ are the relative permittivity and permeability of $D$. Besides, the wave number $k$ is given by $k \, = \, \omega \, \sqrt{\epsilon_{0} \, \mu_{0}}$, where $\omega$ is the used frequency parameter of the incident plane wave $\left( E^{Inc}(x,\theta, \mathrm{p}, \omega); \, H^{Inc}(x,\theta, \mathrm{p}, \omega) \right)$ given by  
\begin{equation}\label{EincHinc}
E^{Inc}(x,\theta, \mathrm{p}, \omega) \, = \, \mathrm{p} \, e^{i \, k \, \theta \cdot x}  \quad \text{and} \quad  H^{Inc}(x,\theta, \mathrm{p}, \omega) \, = \, \left( \theta \times \mathrm{p} \right) \, e^{i \, k \, \theta \cdot x}, \quad \text{for} \; \theta, \mathrm{p} \in \mathbb{S}^{2} \; \text{and} \; x \in \mathbb{R}^{3},
\end{equation}
with $\mathbb{S}^2$ being the unit sphere, $\theta$ is the incident direction vector and $\mathrm{p}$ is the polarization vector such that $\theta \cdot \mathrm{p} \, = \, 0$. The incident field $E^{Inc}(\cdot,\theta, \mathrm{p}, \omega)$ is solution of 
\begin{equation*}
    \underset{x}{\nabla} \times \underset{x}{\nabla} \times \left( E^{Inc}(x,\theta, \mathrm{p}, \omega) \right) \, - \, k^{2} \, E^{Inc}(x,\theta, \mathrm{p}, \omega) \, = \, 0, \, \quad x \in \mathbb{R}^{3},  
\end{equation*}
The problem $(\ref{model-m})$ is well-posed in appropriate Sobolev spaces, see \cite{colton2019inverse, mitrea-mitrea-pipher}, and the scattered wave $(E^s, H^s)$ possesses the following asymptotic behaviors
\begin{equation}\label{far-def}
E^s(x)=\frac{e^{i k|x|}}{|x|}\left(E^\infty(\hat{x}, \theta, p)+O(|x|^{-1})\right),\quad \mbox{as}\quad |x|\rightarrow \infty,
\end{equation}
and
\begin{equation}\notag
H^s(x)=\frac{e^{i k|x|}}{|x|}\left(H^\infty(\hat{x}, \theta, p)+O(|x|^{-1})\right),\quad \mbox{as}\quad |x|\rightarrow \infty,
\end{equation}
where $(E^\infty(\hat{x}, , \theta, p), H^\infty(\hat{x}, \theta, p))$ is the corresponding electromagnetic far field pattern of \eqref{model-m} in the propagation direction $\hat{x} \, := \, \dfrac{x}{|x|}$.
%\begin{equation}\label{U}
%\left \{
%\begin{array}{llrr}
%Curl(E^{T}) - i \, k \, \mu_r\,  H^{T} = 0 \, \quad & \mbox{in } \mathbb{R}^3 \\
%\\
%Curl(H^{T}) + i \, k \, \varepsilon_{r} \, E^{T} = 0 \, \quad & \mbox{in } \mathbb{R}^3,
%\end{array} 
%\right.
%\end{equation}
%where the total field $(E^T, H^{T})$ is of the form $(E^T:=E^{Inc} +E^s,H^T:=H^{Inc} +H^s)$ and the incident plane wave $(E^{Inc},H^{Inc})$ is of the form
%\begin{equation}\notag
%E^{Inc}(x,\theta) = \theta^{\perp} \exp\left(i \, k \; \theta \cdot x \right)\quad\mbox{and}\quad
%H^{Inc}(x,\theta) = \big( \theta^{\perp} \times \theta \big) \; \exp\big(i \, k \; \theta \cdot x \big),
%\end{equation}
%and the scattered field $(E^s, H^s)$ satisfies the Silver-M\"{u}ller radiation condition (SMRC) at infinity:
%
%\begin{equation}\notag
%	\sqrt{\mu_0\epsilon_0^{-1}}H^{s}(x)\times \frac{x}{|x|}-E^s(x)=O(\frac{1}{|x|^2}).
%\end{equation}

\bigskip

%Following the same approach as in \cite{Costabel.D.K}, we write the integral equation solution of our problem. 

Next, we present the necessary assumptions on the model \eqref{model-m} to derive the main results.

\bigskip

%%%%%%%%%%%%%%%%%%%%%%%%%Basic assumption%%%

\begin{assumption}\label{ASass}
\phantom{}
\begin{enumerate}
    \item Assumptions on the cluster of dimers\label{\romannumeral1}. 
    Suppose that each dimer $D_m:=D_{m_1}\cup D_{m_2}$ of $D$ can be represented by $D_{m_1}=a {B}_{m_1}+{z}_{m_1}$, $D_{m_2}=aB_{m_2}+z_{m_2}$ with the parameter $a>0$ and the locations ${z}_{m_\ell}$, $\ell=1,2$, $m=1, 2, \cdots, \aleph$, where
\begin{equation}\label{def-d}
	a \, := \, \max_{1\leq m\leq \aleph}\{\mathrm{diam}(D_{m_1}), \mathrm{diam}(D_{m_2})\}, 
\end{equation} 
with $B_{m_1}$ and $B_{m_2}$ are such that 
\begin{equation*}
	 \max_{1\leq m\leq \aleph}\{\mathrm{diam}(B_{m_1}), \mathrm{diam}(B_{m_2})\} \, = \, 1. 
\end{equation*}
Denote by $d_{\rm in}$ the minimal distance between the nano-particles of each dimer, i.e.,
\begin{equation}\label{dis-in}
	d_{\rm in}=\min_{1\leq m\leq \aleph} \mathrm{dist} (D_{m_1}, D_{m_2}),
\end{equation}
and by $d_{\rm out}$ the minimal distance between the used dimers, i.e., 
\begin{equation}\label{dis-out}
d_{\rm out}:=\min_{1\leq m, j\leq \aleph\atop m\neq j} \mathrm{dist} (D_m, D_j)= \min_{1\leq m, j\leq \aleph\atop m\neq j} |z_{m_0}-z_{j_0}|.
\end{equation}
where $z_{m_0}$ is the intermediate point of $D_m$ between $D_{m_1}$ and $D_{m_2}$ such that 
\begin{equation}\label{def-intermediate point}
z_{m_0} \, := \, z_{m_1} \, + \, \frac{\left(z_{m_2} \, - \, z_{m_1} \right)}{2}, \quad \text{for} \quad  m=1, 2, \cdots, \aleph,    
\end{equation}
see \textbf{Figure \ref{fig:intermediate-point}} for a schematic illustration. 
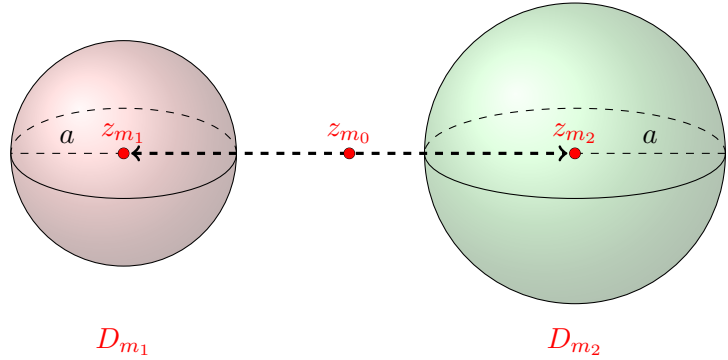
\begin{figure}[H]
	%\centering
	%\includegraphics[width=0.3\linewidth]{"Intermediate point"}
    \begin{center}
\begin{tikzpicture}
%%%%%%%%%%%%%%%%%%%%%%%%%%%%%%%%%%%%%%%%%%%%%%%%%%%%%%%%%
  \shade[ball color = red!40, opacity = 0.4] (0,0) circle (1.5cm);
  \draw (0,0) circle (1.5cm);
  \draw (-1.5,0) arc (180:360:1.5 and 0.6);
  \draw[dashed] (1.5,0) arc (0:180:1.5 and 0.6);
  \fill[fill=black] (0,0) circle (1pt);
  \draw[dashed] (0,0 ) -- node[above]{$a$} (-1.5,0);
%%%%%%%%%%%%%%%%%%%%%%%%%%%%%%%%%%%%%%%%%%%%%%%%%%%%%%%%%
  \shade[ball color = green!40, opacity = 0.4] (6,0) circle (2cm);
  \draw (6,0) circle (2cm);
  \draw (4,0) arc (180:360:2 and 0.6);
  \draw[dashed] (8,0) arc (0:180:2 and 0.6);
  \fill[fill=black] (6,0) circle (1pt);
  \draw[dashed] (6,0 ) -- node[above]{$a$} (8,0);
%%%%%%%%%%%%%%%%%%%%%%%%%%%%%%%%%%%%%%%%%%%%%%%%%%%%%%%%%
\draw[dashed, very thick,<->] (0.1,0) -- (5.9,0);
\draw (0,0) circle[radius=2pt]; 
\fill[fill=red] (0,0) circle[radius=2pt];
\draw (6,0) circle[radius=2pt]; 
\fill[fill=red] (6,0) circle[radius=2pt];
\draw (3,0) circle[radius=2pt]; 
\fill[fill=red] (3,0) circle[radius=2pt];
%%%%%%%%%%%%%%%%%%%%%%%%%%%%%%%%%%%%%%%%%%%%%%%%%%%%%%%
\draw (0,0) node[above] {\textcolor{red}{$z_{m_{1}}$}};
\draw (3,0) node[above] {\textcolor{red}{$z_{m_{0}}$}};
\draw (6,0) node[above] {\textcolor{red}{$z_{m_{2}}$}};
%%%%%%%%%%%%%%%%%%%%%%%%%%%%%%%%%%%%%%%%%%%%%%%%%%%%%%%
\draw node at (0,-2.5)  {\textcolor{red}{$D_{m_{1}}$}};
\draw node at (6,-2.5)  {\textcolor{red}{$D_{m_{2}}$}};
%%%%%%%%%%%%%%%%%%%%%%%%%%%%%%%%%%%%%%%%%%%%%%%%%%%%%%%
\end{tikzpicture}
\end{center}
	\caption{A schematic illustration for the Dimer $D_{m}$ and its intermediate point $z_{m_{0}}$.}
	\label{fig:intermediate-point}
\end{figure}
The parameters $d_{\rm in}, d_{\rm out}$, and $\aleph$ admit the following behaviors with respect to the parameter $a$ as 
	\begin{equation}\label{d-a}
d_{\rm in}=\alpha_{0}\,a^{t_1},\qquad d_{\rm out}=\beta_{0}\,a^{t_2},\qquad \aleph\sim\bigl[d_{\rm out}^{-3}\bigr],
\qquad 0<t_{2}\le t_{1}<1.
\end{equation} 
where $[\cdot]$ stands for the entire part function and  $\alpha_0, \beta_0$ are two positive constants independent on the parameter $a$. 
\item[]
\item Assumptions on the permittivity and permeability of each dimer.\label{\romannumeral3} Regarding the permittivity, we assume that for $m=1, 2, \cdots, \aleph$,	
\begin{equation}\label{def-eta}
	\eta(x) \, := \, \begin{cases}
	\eta_{1} \, := \, \epsilon_{r}^{(1)} \, - \, 1 \, = \, \eta_0 \; a^{-2} & \text{if $x \in D_{m_1}$}\\
	& \\
	\eta_{2} \, := \, \epsilon_{r}^{(2)} \, - \, 1 \, \sim 1 \quad \text{with} \quad \Re \left( \epsilon_{r}^{(2)} \right) \, < \, 0& \text{if $x \in D_{m_2}$}
	\end{cases},
	\end{equation}
	where $\eta_0$ is a constant in the complex plane independent of $a$, with $\Re(\eta_0)>0$.
Thus $\Re(\epsilon_{r}^{(1)})>0$ and $\Re(\epsilon_{r}^{(2)})<0$.
Moreover, we assume $\mu_{r_m}^{(\ell)}=1$ for $\ell=1,2$ and $m=1,2,\dots,\aleph$.
    \item[] 
    \item Assumption on the shape of $B_{m_\ell}$, $\ell=1, 2$.\label{\romannumeral2} For simplicity reason, we assume that the shapes of $B_{m_1}$'s and $B_{m_2}$'s are the same respectively and denote $B_1:=B_{m_1}$, $B_2:=B_{m_2}$ for $m=1, 2, \cdots, \aleph$. Define the vector Newtonian potential operator $N_{B_\ell}(\cdot)$ and the Magnetization operator $\nabla M_{B_\ell} (\cdot)$ as 
\begin{equation}\label{NM0}
N_{B_\ell}(F)(x):=\int_{B_\ell} \Phi_{0}(x,y) \, F(y)\,dy \quad \text{and} \quad \nabla M_{B_\ell}(F)(x):=\underset{x}{\nabla}\int_{B_\ell}\underset{y}{\nabla}\left(\Phi_{0}(x,y)\right)\cdot F(y)\,dy,
\end{equation}
where $\Phi_{0}(\cdot,\cdot)$ is given by 
\begin{equation*}
    \Phi_{0}(x,y) \, := \, \frac{1}{4 \, \pi \, \left\vert x - y \right\vert}, \quad x \neq y, 
\end{equation*}
Under the Helmholtz decomposition of the $\mathbb{L}^2$-space  given by \eqref{L2-decomposition}, denote $\left\{ e_{n, B_\ell}^{(1)} \right\}_{n \in \mathbb{N}}$ as the corresponding eigenfunctions of the operator $N_{B_\ell}(\cdot)$ over the subspace $\mathbb{H}_{0}(\div=0)$. Since 
\begin{equation}\label{Hdiv-curl}
    \mathbb{H}_{0}\left( \div = 0 \right) \equiv Curl \left( \mathbb{H}_{0}\left( Curl \right) \cap \mathbb{H}\left( \div = 0 \right) \right),
\end{equation}
	see for instance \cite{amrouche-bernardi-dauge-girault}, \mk{then there exists $\phi_{n, B_\ell} \in  \mathbb{H}_{0}\left( Curl \right) \cap \mathbb{H}\left( \div = 0 \right)$, such that }  
\begin{equation*}\label{pre-cond}
e_{n, B_\ell}^{(1)}= Curl (\phi_{n, B_\ell}) \,\, \mbox{ with } \,\, \nu\times \phi_{n, B_\ell}=0 \,\, \text{and} \,\, \div(\phi_{n, B_\ell})=0.
\end{equation*}	
We assume that in $B_1$, 
%\label{def-phi-n0}
\begin{equation}\label{DefP011}
{\bf P}_{0, 1}^{(1)} \, := \, \int_{B_1} \phi_{n_{0}, B_1}(y) \, dy  \, \otimes \, \int_{B_1} \phi_{n_{0}, B_1}(y) \, dy \, \neq \, 0 \quad \text{for certain} \; n_{0} \in \mathbb{N},
\end{equation} 
and
\begin{equation}\label{DefP012}
   {\bf P}_{0, 1}^{(2)} \, := \, \sum_{n} \, \frac{1}{\lambda_{n}^{(3)}(B_{1})} \, \int_{B_{1}} e_{n,B_{1}}^{(3)}(x) \, dx \, \otimes \, \int_{B_{1}} e_{n,B_{1}}^{(3)}(x) \, dx \, \neq \, 0, 
\end{equation}
where $\left( \lambda_n^{(3)}(B_1), e_{n, B_1}^{(3)} \right)$ is the eigensystem of $\nabla M_{B_1}(\cdot)$ over $\nabla\mathcal{H}armonic (B_1)$. In $B_2$, it is assumed
\begin{equation}\label{DefP021}
	{\bf P}_{0, 2}^{(1)} \, := \, \sum_{n} \int_{B_2}  \phi_{n, B_2}(y) \, dy \otimes \int_{B_2} \phi_{n, B_2}(y) \, dy \, \neq \, 0,
\end{equation}
and
\begin{equation}\label{DefP022}
    {\bf P}_{0, 2}^{(2)} \, := \, \int_{B_{2}} e_{n_{\star},B_{2}}^{(3)}(x) \, dx \, \otimes \, \int_{B_{2}} e_{n_{\star},B_{2}}^{(3)}(x) \, dx \, \neq \, 0, \quad \text{for certain} \; n_{\star} \in \mathbb{N}, 
\end{equation}
where $e_{n_*,B_{2}}^{(3)}(\cdot) \in \nabla \mathcal{H}armonic(B_{2})$ is an eigenfunction, corresponding to the eigenvalue $\lambda_{n}^{(3)}(B_{2})$, related to the Magnetization operator $\nabla M_{B_{2}}\left( \cdot \right)$.
\item[]
\item Assumption on the used incident frequency $k$. \label{\romannumeral4} Define the vector Magnetization operator $\nabla {M}_{B_2}(\cdot)$ as \eqref{NM0}. Under the Helmholtz decomposition of $\mathbb{L}^2$ space given by \eqref{L2-decomposition},
denote $\left(\lambda_{n}^{(3)}(B_2), e_{n, B_2}^{(3)}\right)$ as the corresponding eigen-system of $\nabla M_{B_2}(\cdot)$ over the subspace $\nabla \mathcal{H}armonic$. There exists positive constants $c_0$ and $d_0$ such that
\begin{equation}\label{condition-on-k}
			1 \, - \, k^2 \, \eta_1 \, a^2 \, \lambda_{n_{0}}^{(1)} (B_1) \, = \, \pm \; c_0\; a^h \quad \text{and} \quad  1+\eta_2\lambda_{n_*}^{(3)}(B_2)=\; \pm \; d_0 \; a^h,\; ~~ a \ll 1,
		\end{equation}
		where $\lambda_{n_0}^{(1)}(B_1)$ is the eigenvalue corresponding to $e_{n_0, B_1}^{(1)}(\cdot)$ for the Newtonian potential operator $N_{B_1}(\cdot)$ in $\mathbb{H}_{0}(\div = 0)$, and $\lambda_{n_*}^{(3)}(B_{2})$ is the eigenvalue corresponding to $e_{n_*, B_2}^{(3)}(\cdot)$ for $\nabla M_{B_2}(\cdot)$ in $\nabla \mathcal{H}armonic$.
        \end{enumerate}
\end{assumption}
%\textcolor{blue}{
%To clarify the assumed assumptions, we set the following remark before presenting the main result of this work.
%\begin{remark}
%    Two remarks are in order. 
%    \begin{enumerate}
 %       \item Assumption $(\ref{\romannumeral1})$ is the are already used for the case of a single Dimer, see \cite[Section 1.2]{cao-ghandriche-sini-dimer}. 
 %       \item[]
 %       \item Assumptions $(\ref{\romannumeral2})$ and  $(\ref{\romannumeral4})$ are already used for the case of a single Dimer, see \cite[Section 1.2]{cao-ghandriche-sini-dimer}.
 %   \end{enumerate}
%\end{remark}}
\medskip
\textit{To avoid cluttering the presentation of the main result of this work, see \textbf{Theorem \ref{main-1}}, we suggest using these notations before presenting it.}
\begin{notation}\label{notbthm}
Before giving the necessary notations, we recall that 
$\Upsilon_k(\cdot, \cdot)$ is the Dyadic Green's kernel given by
  \begin{equation}\label{dyadicG}
  \Upsilon_{k}(\cdot, \cdot) \, := \, \frac{1}{k^{2}} \, \nabla \nabla \Phi_{k}(\cdot, \cdot) \, + \, \Phi_{k}(\cdot, \cdot) \, I_{3}, \quad\mbox{with}\quad \Upsilon_0(\cdot, \cdot):=\nabla\nabla\Phi_0(\cdot, \cdot),
  \end{equation}
where $I_{3} \in \mathbb{R}^{3 \times 3}$ is the identity matrix, $k$ is the wave number, and $\Phi_{k}(\cdot, \cdot)$ being the fundamental solution of the Helmholtz equation given by 
  \begin{equation}\label{Helmholtz-kernel}
      \Phi_{k}(x, y) \, := \, \frac{e^{i \, k \, \left\vert x \, - \, y \right\vert}}{4 \, \pi \, \left\vert x \, - \, y \right\vert}, \quad x \neq y.
  \end{equation}
  % In particular, for $k=0$, we denote $\Upsilon_0(\cdot, \cdot)=\nabla\nabla\Phi_0(\cdot, \cdot)$.
Besides, we recall that ${\bf{P}}_{0, i}^{(j)}$, for $i,j = 1,2$,  are the polarization tensors defined by $(\ref{DefP011})-(\ref{DefP022})$. In what follows, the notations mentioned below will be utilized.
    \begin{enumerate}
        \item[] 
        \item \label{DefBm} A matrix $\mathfrak{B_{m}} \in \mathbb{C}^{12 \times 12}$  given by  
        \begin{equation}\label{detail-B_m}
        \mathfrak{B_{m}} \, := \,  \begin{pmatrix}
    I_{3} & 0 & - \, \mathcal{B}_{13} & - \, \mathcal{B}_{14} \\
    0 & I_{3} & - \, \mathcal{B}_{23}  & -\, \mathcal{B}_{24} \\
    - \, \mathcal{B}_{31} & - \, \mathcal{B}_{32} & I_{3} & 0 \\  
    - \, \mathcal{B}_{41} & - \, \mathcal{B}_{42} & 0 & I_{3}
    \end{pmatrix}_{m},
        \end{equation}
        with 
        \[
\small
\begin{aligned}
\mathcal{B}_{13} &:= \frac{k^{4}\,\eta_{0}}{\pm c_{0}}\,a^{3-h}\,{\bf P}_{0, 1}^{(1)}\cdot\Upsilon_{k}(z_{m_1}, z_{m_2}),\\
\mathcal{B}_{14} &:= \frac{k^{2}\,\eta_{0}}{\pm c_{0}}\,a^{3-h}\,{\bf P}_{0, 1}^{(1)}\cdot\bigl(\nabla \Phi_{k}(z_{m_1}, z_{m_2})\times\bigr),\\
\mathcal{B}_{23} &:= k^{2}\,a^{3}\,{\bf P}_{0, 1}^{(2)}\cdot\bigl(\nabla \Phi_{k}(z_{m_1}, z_{m_2})\times\bigr),\\
\mathcal{B}_{24} &:= k^{2}\,a^{3}\,{\bf P}_{0, 1}^{(2)}\cdot\Upsilon_{k}(z_{m_1}, z_{m_2}),\\
\mathcal{B}_{31} &:= k^{4}\,\eta_{2}\,a^{5}\,{\bf P}_{0, 2}^{(1)}\cdot\Upsilon_{k}(z_{m_2}, z_{m_1}),\\
\mathcal{B}_{32} &:= k^{4}\,\eta_{2}\,a^{5}\,{\bf P}_{0, 2}^{(1)}\cdot\bigl(\nabla \Phi_{k}(z_{m_2}, z_{m_1})\times\bigr),\\
\mathcal{B}_{41} &:= \frac{k^{2}\,\eta_{2}}{\pm d_{0}}\,a^{3-h}\,{\bf P}_{0, 2}^{(2)}\cdot\bigl(\nabla \Phi_{k}(z_{m_2}, z_{m_1})\times\bigr),\\
\mathcal{B}_{42} &:= \frac{k^{2}\,\eta_{2}}{\pm d_{0}}\,a^{3-h}\,{\bf P}_{0, 2}^{(2)}\cdot\Upsilon_{k}(z_{m_2}, z_{m_1}).
\end{aligned}
\]

         \item[]
         \item[] 
         \item A matrix $\Psi_{mj}  \in \mathbb{C}^{12 \times 12}$ given by 
        \begin{equation}\label{HADefPsimj}
            \Psi_{mj} \, := \, \begin{pmatrix}
        \mathcal{C}_{11} & \mathcal{C}_{12} & \mathcal{C}_{13} & \mathcal{C}_{14}\\
        \mathcal{C}_{21} & \mathcal{C}_{22} & \mathcal{C}_{23} & \mathcal{C}_{24} \\
        \mathcal{C}_{31} & \mathcal{C}_{32} & \mathcal{C}_{33} & \mathcal{C}_{34} \\
        \mathcal{C}_{41} & \mathcal{C}_{42} & \mathcal{C}_{43} & \mathcal{C}_{44} \\
    \end{pmatrix}_{mj},
        \end{equation}
        with 
        \[
\small
\begin{aligned}
\mathcal{C}_{11} &:= \frac{k^{4}\,\eta_{0}}{\pm c_{0}}\,a^{3-h}\,{\bf P}_{0, 1}^{(1)}\cdot\Upsilon_{k}(z_{m_1}, z_{j_1}),\\
\mathcal{C}_{12} &:= \frac{k^{2}\,\eta_{0}}{\pm c_{0}}\,a^{3-h}\,{\bf P}_{0, 1}^{(1)}\cdot\bigl(\nabla \Phi_{k}(z_{m_1}, z_{j_1})\times\bigr),\\
\mathcal{C}_{13} &:= \frac{k^{4}\,\eta_{0}}{\pm c_{0}}\,a^{3-h}\,{\bf P}_{0, 1}^{(1)}\cdot\Upsilon_{k}(z_{m_1}, z_{j_2}),\\
\mathcal{C}_{14} &:= \frac{k^{2}\,\eta_{0}}{\pm c_{0}}\,a^{3-h}\,{\bf P}_{0, 1}^{(1)}\cdot\bigl(\nabla \Phi_{k}(z_{m_1}, z_{j_2})\times\bigr),\\
\mathcal{C}_{21} &:= k^{2}\,a^{3}\,{\bf P}_{0, 1}^{(2)}\cdot\bigl(\nabla \Phi_{k}(z_{m_1}, z_{j_1})\times\bigr),\\
\mathcal{C}_{22} &:= k^{2}\,a^{3}\,{\bf P}_{0, 1}^{(2)}\cdot\Upsilon_{k}(z_{m_1}, z_{j_1}),\\
\mathcal{C}_{23} &:= k^{2}\,a^{3}\,{\bf P}_{0, 1}^{(2)}\cdot\bigl(\nabla \Phi_{k}(z_{m_1}, z_{j_2})\times\bigr),\\
\mathcal{C}_{24} &:= k^{2}\,a^{3}\,{\bf P}_{0, 1}^{(2)}\cdot\Upsilon_{k}(z_{m_1}, z_{j_2}),\\
\mathcal{C}_{31} &:= k^{4}\,\eta_{2}\,a^{5}\,{\bf P}_{0, 2}^{(1)}\cdot\Upsilon_{k}(z_{m_2}, z_{j_1}),\\
\mathcal{C}_{32} &:= k^{2}\,\eta_{2}\,a^{5}\,{\bf P}_{0, 2}^{(1)}\cdot\bigl(\nabla \Phi_{k}(z_{m_2}, z_{j_1})\times\bigr),\\
\mathcal{C}_{33} &:= k^{4}\,\eta_{2}\,a^{5}\,{\bf P}_{0, 2}^{(1)}\cdot\Upsilon_{k}(z_{m_2}, z_{j_2}),\\
\mathcal{C}_{34} &:= k^{2}\,\eta_{2}\,a^{5}\,{\bf P}_{0, 2}^{(1)}\cdot\bigl(\nabla \Phi_{k}(z_{m_2}, z_{j_2})\times\bigr),\\
\mathcal{C}_{41} &:= \frac{k^{2}\,\eta_{2}}{\pm d_{0}}\,a^{3-h}\,{\bf P}_{0, 2}^{(2)}\cdot\bigl(\nabla \Phi_{k}(z_{m_2}, z_{j_1})\times\bigr),\\
\mathcal{C}_{42} &:= \frac{k^{2}\,\eta_{2}}{\pm d_{0}}\,a^{3-h}\,{\bf P}_{0, 2}^{(2)}\cdot\Upsilon_{k}(z_{m_2}, z_{j_1}),\\
\mathcal{C}_{43} &:= \frac{k^{2}\,\eta_{2}}{\pm d_{0}}\,a^{3-h}\,{\bf P}_{0, 2}^{(2)}\cdot\bigl(\nabla \Phi_{k}(z_{m_2}, z_{j_2})\times\bigr),\\
\mathcal{C}_{44} &:= \frac{k^{2}\,\eta_{2}}{\pm d_{0}}\,a^{3-h}\,{\bf P}_{0, 2}^{(2)}\cdot\Upsilon_{k}(z_{m_2}, z_{j_2}).
\end{aligned}
\]

\item[]
\item A vector $\mathrm{S}_{m} \, \in \, \mathbb{C}^{1 \times 12}$ given by 
\begin{equation*}
    \mathrm{S}_{m} \, := \, 
\begin{pmatrix}
     \frac{i \, k \, \eta_{0}}{\pm \, c_{0}} \, a^{3-h} \, {\bf P}_{0, 1}^{(1)} \cdot H^{Inc}(z_{m_1}) \\
     a^{3} \, {\bf P}_{0, 1}^{(2)} \cdot E^{Inc}(z_{m_1}) \\
     i \, k \, a^{5} \, {\bf P}_{0, 2}^{(1)} \cdot H^{Inc}(z_{m_2}) \\
     \frac{\eta_{2}}{\pm \, d_{0}} \, a^{3-h} \, {\bf P}_{0, 2}^{(2)} \cdot E^{Inc}(z_{m_2})
    \end{pmatrix}.
    \end{equation*}
    \end{enumerate}
\end{notation}
\medskip

\begin{notation} 
To make a cross reference to the supplementary material document, we have added the prefix $SM-$ before the element that we want to make a reference to. For example, Subsection $SM-1.1$ refers to Subsection 1.1 in the supplementary material document. 
\end{notation}

% \bigskip

The main outcome of this work is now available for presentation. 
%The main result can be stated as follows.
\begin{theorem}\label{main-1}
	Let the \textbf{Assumptions \ref{ASass}}
    %$\ref{\romannumeral1}, \ref{\romannumeral2}, \ref{\romannumeral3},$ and $ \ref{\romannumeral4}$ 
    on the electromagnetic scattering problem \eqref{model-m}, which is generated  by a cluster of dimers $D \, := \, \overset{\aleph}{\underset{m=1}{\cup}} D_{m}   \, = \, \overset{\aleph}{\underset{m=1}{\cup}} \left( D_{m_{1}} \cup D_{m_{2}} \right)$, be satisfied. For $0<t_2\leq t_1<1$ and $h \in \mathbb{R}^+$ with the constants $c_0, d_0$ such that 
\begin{equation}\label{conditions-t-h}
\frac{k^4}{(4\pi)^4\, c_0^2 \, d_0^2}<1\quad \mbox{and}\quad \frac{9}{5}<h<\min\{2, 5-8t_1\},
\end{equation} 
then the scattered wave admits the following expansion
\begin{equation}\label{scattered field expansion-thm}
    E^{s}({x}, \theta, p, \omega) \, 
    = - k^2 \sum_{m=1}^\aleph \mathcal{G}_m\cdot \mathfrak{U}_{m}+\O\left(a^{10-3h-\frac{23}{2}t_2}\right)\quad\mbox{with}\quad \mathcal{G}_m:=\begin{pmatrix}
            \underset{y}{\nabla}\Phi_k(x, z_{m_1})\times\\
            -\Upsilon_k(x, z_{m_1})\\
            \underset{y}{\nabla}\Phi_k(x, z_{m_2})\\
            -\Upsilon_k(x, z_{m_2})
        \end{pmatrix}^\mathsf{T},
\end{equation}
and the corresponding far field admits the following approximation expansion 
\begin{equation}\label{foldy-lax-final}
    E^{\infty}(\hat{x}, \theta, p) \, =  \sum_{m=1}^\aleph \mathfrak{e}_m\cdot \mathfrak{U}_{m}+\O\left(a^{10-3h-\frac{23}{2}t_2}\right) \quad\mbox{with}\quad \mathfrak{e}_m:= \left(
     \begin{array}{c}
	-ik e^{-ik \hat{x}\cdot z_{m_1}} \hat{x}\times\\
	e^{-i k \hat{x}\cdot z_{m_1}}(I-\hat{x}\otimes\hat{x}) \\
	-ik e^{-ik \hat{x}\cdot z_{m_2}} \hat{x}\times\\
	e^{-i k \hat{x}\cdot z_{m_2}}(I-\hat{x}\otimes\hat{x})
\end{array}\right)^\mathsf{T},
\end{equation}
% \begin{eqnarray}\label{approximation-E}
% \nonumber
% E^{\infty}(\hat{x}, \theta, p) \, &=& \, \left( I_{3} \, - \, \hat{x}\otimes\hat{x} \right) \cdot \sum_{m=1}^{\aleph} e^{i k \hat{x}\cdot z_{m_0}} \,  \left( 0 \, I_{3} \, 0 \, I_{3} \right) \cdot \mathfrak{U}_{m} \, - \, i \, k \,  \sum_{m=1}^{\aleph} e^{i k \hat{x}\cdot z_{m_0}} \,   \hat{x} \times \left( I_{3} \, 0 \, I_{3} \, 0 \right) \cdot \mathfrak{U}_{m} \\ \nonumber &+& \mathcal{O}\left(...\right) \\
% && \qquad \qquad + \mathcal{O}\left( a^{t_{1}-3t_{2}+3-h} \right),
% \end{eqnarray} 
%  with 
%  \begin{equation*}
%      \mathfrak{e}_m:= \left(
%      \begin{array}{c}
% 	-ik e^{-ik \hat{x}\cdot z_{m_1}} \hat{x}\times\\
% 	e^{-i k \hat{x}\cdot z_{m_1}}(I-\hat{x}\otimes\hat{x}) \\
% 	-ik e^{-ik \hat{x}\cdot z_{m_2}} \hat{x}\times\\
% 	e^{-i k \hat{x}\cdot z_{m_2}}(I-\hat{x}\otimes\hat{x})
% \end{array}\right)^\mathsf{T},
%  \end{equation*}
 where
 $\left(\mathfrak{U}_{m} \, := \, \left( Q_{m_1}, R_{m_1}, Q_{m_2}, R_{m_2} \right) \right)_{m=1}^{\aleph}$ is the vector solution to the following algebraic system
        \begin{equation}\label{eq-al-D1-}
            \begin{pmatrix}
                 \mathfrak{B_{1}} & - \, \Psi_{12} & \cdots & - \, \Psi_{1\aleph}  \\
                 - \, \Psi_{21} & \mathfrak{B_{2}} & \cdots & - \, \Psi_{2\aleph} \\
                 \vdots & \vdots & \ddots & \vdots \\
                 - \, \Psi_{\aleph1} & - \, \Psi_{\aleph2} & \cdots &  \mathfrak{B_{\aleph}} 
            \end{pmatrix} 
            \cdot
            \begin{pmatrix}
                \mathfrak{U}_{1} \\
                \mathfrak{U}_{2} \\
                \vdots \\
                \mathfrak{U}_{\aleph} 
            \end{pmatrix}
            = 
            \begin{pmatrix}
                \mathrm{S}_{1} \\
                \mathrm{S}_{2} \\
                \vdots \\
                \mathrm{S}_{\aleph}  
            \end{pmatrix},
        \end{equation}
        with $\left\{ \left\{\mathfrak{B_{m}} \right\}_{m=1}^{\aleph}; \left\{ \Psi_{mj}\right\}_{m,j=1}^{\aleph}; \left\{ \mathrm{S}_{m} \right\}_{m=1}^{\aleph} \right\}$ given in $\textbf{Notations \ref{notbthm}}$. 
In particular, the algebraic system $(\ref{eq-al-D1-})$ is invertible under the conditions that
\begin{equation}\label{condi-invertibility}
      \frac{k^2}{2}\left( \frac{k^2 \eta_0}{c_0}|{\bf P}_{0, 1}^{(1)}|+\frac{\eta_2}{d_0}|{\bf P}_{0, 2}^{(2)}|\right)a^{3-h-3t_1}<1\, \mbox{ and }\,  \frac{k^2 \left( \frac{\max\{1, k^2\}\eta_0}{c_0}\left| {\bf P}_{0, 1}^{(1)} \right| + \frac{\eta_2}{d_0} \left| {\bf P}_{0, 2}^{(2)} \right| \right)}{1- \frac{k^2}{2} \left( \frac{k^2 \eta_0}{c_0}\left|{\bf P}_{0, 1}^{(1)} \right| +\frac{\eta_2}{d_0}\left| {\bf P}_{0, 2}^{(2)} \right| \right) a^{3-h-3t_1} } a^{3-h-3t_2} <1,
\end{equation}
% under the conditions 
% \begin{equation}\label{condi-invertibility-1}
%         \frac{2\beta_0}{k \alpha_0} \min\{ 1+\frac{c_0 \eta_2 |{\bf P}_{0, 2}^{(2)}|}{k^2 d_0 \eta_0|{\bf P}_{0, 1}^{(1)}|}, 1+\frac{k^2 d_0 \eta_0 |{\bf P}_{0, 1}^{(1)}|}{c_0 \eta_2 |{\bf P}_{0, 2}^{(2)}|}\} a^{t_1 - t_2}<1,
% \end{equation}
% and 
% \begin{equation}\label{condi-invertibility-2}
%      \frac{k^2}{2}\left( \frac{k^2 \eta_0}{c_0}|{\bf P}_{0, 1}^{(1)}|+\frac{\eta_2}{d_0}|{\bf P}_{0, 2}^{(2)}|\right)a^{3-h-3t_1}<1,
% \end{equation}
with $\alpha_0, \beta_0$ given in \eqref{d-a}.
% \begin{equation}\label{condi-invertibility}
%     \frac{k}{\max\{\alpha,\beta\}}\left(\frac{k^2 \eta_0}{c_0}|{\bf P}_{0, 1}^{(1)}|+\frac{\eta_2}{d_0}|{\bf P}_{0, 2}^{(2)}|\right)a^{3-h-3t_2}<1\quad\mbox{and}\quad 
%     \frac{k^2}{2}\left( \frac{k^2 \eta_0}{c_0}|{\bf P}_{0, 1}^{(1)}|+\frac{\eta_2}{d_0}|{\bf P}_{0, 2}^{(2)}|\right)a^{3-h-3t_1}<1,
% \end{equation}
% where 
% \begin{equation*}
%     \alpha= \frac{k^4 \eta_0}{2c_0} |{\bf P}_{0, 1}^{(1)}|a^{3-h}d_{\rm in}^{-3} ,\quad \beta=\frac{k^2 \eta_2}{2d_0} |{\bf P}_{0, 2}^{(2)}|.
% \end{equation*}
\end{theorem}

\begin{remark}
    The first condition of \eqref{conditions-t-h} can be satisfied not only for appropriate incident frequency $k$, but also certain values of the constants $c_0, d_0$ and $\eta_0$. Indeed, from assumption \eqref{condition-on-k} and \eqref{def-eta}, it is obvious that
    \begin{equation*}
        k^2 = \frac{1\mp c_0 a^h}{\eta_0 \lambda_{n_0}^{(1)}}.
    \end{equation*}
    Hence, to meet \eqref{conditions-t-h}, we have freedom to choose proper $c_0, d_0$ and $\eta_0$ such that 
    \begin{equation*}
        (4\pi)^4 c_0^2 d_0^2 \eta_0^2 (\lambda_{n_0}^{(1)})^2 >1.
    \end{equation*}
\end{remark}

\begin{remark}
    The scattered field approximation expansion form \eqref{scattered field expansion-thm} can be written equivalently as 
    \begin{equation}\label{scattered field-equivalent form}
        E^{s}(x,\theta, p, \omega) \, = \, - \, k^{2} \, \sum_{m=1}^{\aleph} \sum_{\ell=1}^{2} \left[ \underset{y}{\nabla} \Phi_{k}(x,z_{m_{\ell}}) \times {Q}_{m_{\ell}} \, - \, \Upsilon_{k}(x,z_{m_{\ell}}) \cdot {R}_{m_{\ell}} \right] \, + \, \mathcal{O}\left( a^{10-3h - \frac{23}{2}t_2}\right), 
    \end{equation} 
    % and 
    % \begin{equation}\notag
    %     E^{\infty}(\hat{x},\theta, p, \omega)=\sum_{m=1}^\aleph \sum_{i=1}^2 \left(e^{i k \hat{x} \cdot z_{m_i}} (I-\hat{x}\otimes\hat{x})\tilde{R}_{m_i}- i k e^{i k \hat{x}\cdot z_{m_i}} \hat{x}\times \tilde{Q}_{m_i}\right)+\O(a^{7-2h-\frac{15}{2}t_2}).
    % \end{equation}
    which straightforward indicates that the electric field generated by a cluster of dimers is composed by the polarization of the incident electric field associated with the plasmonic particles as monopole and the polarization of the incident magnetic field associated with the dieletric particles as dipole.
\end{remark}

\section{Discussion of the results and possible applications}

In this section, we discuss the information encoded in Theorem \ref{main-1} and outline several
directions in which the result can be used. We work under the assumptions of Theorem \ref{main-1} and
use the notations introduced therein.

\subsection*{Discrete Foldy--Lax model for a hybrid dimer cluster}

Theorem \ref{main-1} states that, under the scaling and spectral assumptions \eqref{def-d}--\eqref{d-a}, \eqref{DefP011}--\eqref{DefP022}, the
scattered field and its far field admit the representations
\begin{equation}\label{eq:FL-main}
E^{s}(x,\theta,p,\omega)
= - k^2
\sum_{m=1}^{\aleph} \mathcal{G}_m(x)\cdot \mathfrak{U}_m
+ O\big(a^{10-3h-\frac{23}{2}t_2}\big),
\end{equation}
\begin{equation}\label{eq:FL-main-ff}
E^{\infty}(\hat x,\theta,p)
=
\sum_{m=1}^{\aleph} \mathfrak{e}_m(\hat x)\cdot \mathfrak{U}_m
+ O\big(a^{10-3h-\frac{23}{2}t_2}\big),
\end{equation}
where, for each dimer $D_m=D_{m_1}\cup D_{m_2}$,
\[
\mathfrak{U}_m := (Q_{m_1},R_{m_1},Q_{m_2},R_{m_2})\in\C^{12}
\]
collects the four moments associated with the projections
$\overset{1}{\PP}(E_{m_\ell})$ and $\overset{3}{\PP}(E_{m_\ell})$ of the total field on $D_{m_\ell}$, $\ell=1,2$, and
$\mathcal{G}_m$ and $\mathfrak{e}_m$ are given explicitly by \eqref{scattered field expansion-thm}--\eqref{foldy-lax-final}.

The vector $\mathfrak{U}=(\mathfrak{U}_1,\dots,\mathfrak{U}_{\mathcal N})\in(\C^{12})^{\aleph}$ is the unique solution to the
finite-dimensional system
\begin{equation}\label{eq:FL-algebraic}
\mathcal B\,\mathfrak{U} = \mathcal S,
\quad
\mathcal B:=
\begin{pmatrix}
\mathfrak{B}_1 & -\Psi_{12} & \cdots & -\Psi_{1\mathcal N}\\
-\Psi_{21} & \mathfrak{B}_2 & \cdots & -\Psi_{2\mathcal N}\\
\vdots & \vdots & \ddots & \vdots\\
-\Psi_{\mathcal N1} & -\Psi_{\mathcal N2} & \cdots & \mathfrak{B}_{\mathcal N}
\end{pmatrix}\quad\mbox{and}\quad \mathcal{S}:= \begin{pmatrix}
                \mathrm{S}_{1} \\
                \mathrm{S}_{2} \\
                \vdots \\
                \mathrm{S}_{\aleph}  
            \end{pmatrix},
\end{equation}
with $\mathfrak{B}_m,\Psi_{mj},\mathrm{S}_m$ given in \textbf{Notations}~1.2. The invertibility conditions (\ref{condi-invertibility}) guarantee that
$\mathcal B$ is invertible and that the error term in
\eqref{eq:FL-main}--\eqref{eq:FL-main-ff} is uniform with respect to~$\aleph$ as long as the
scaling $d_{\rm out}\sim a^{t_2}$ is respected.

Formally, \eqref{eq:FL-main}--\eqref{eq:FL-algebraic} constitute a generalized Foldy--Lax model
for a cluster of hybrid dimers:
\begin{itemize}
  \item each dimer contributes four ``internal'' degrees of freedom
        $(Q_{m_1},R_{m_1},Q_{m_2},R_{m_2})$;
  \item the matrices $\mathfrak{B}_m$ encode self-interaction inside the dimer, through the intra-dimer Green
        kernels and the polarization tensors $\mathbf P^{(j)}_{0,\ell}$, $j, \ell=1,2$;
  \item the matrices $\Psi_{mj}$ encode long-range interactions between different dimers via
        $\Phi_k$ and $\Upsilon_k$.
\end{itemize}
The full Maxwell scattering problem is thus reduced, up to a controlled error, to a finite-dimensional
linear system of size $12\aleph$ with an explicit kernel.

\subsection*{Stability and regime of validity}

The proof of Theorem \ref{main-1} shows that the system \eqref{eq:FL-algebraic} arises by perturbing the
``exact'' system \eqref{eq-al-D1}, satisfied by $\tilde{\mathfrak{U}}_m$, and that the difference $\mathfrak{U}_m-\tilde{\mathfrak{U}}_m$ is
controlled by the error term $\mathfrak{R}_m$ in \eqref{err-prop-la1}. More precisely, using \eqref{err-prop-la1} and the invertibility
conditions \eqref{condi-invertibility}, one obtains the estimate \eqref{far-error-diff}
\[
\Big(\sum_{m=1}^{\aleph} |\mathfrak{U}_m-\tilde{\mathfrak{U}}_m|^2\Big)^{1/2}
\lesssim
a^{10-3h-10t_2} + a^{10-2h-4t_1-6t_2},
\]
which is then propagated to the scattered field, yielding the remainder
$O\big(a^{10-3h-\frac{23}{2}t_2}\big)$ in
\eqref{eq:FL-main}--\eqref{eq:FL-main-ff}. The conditions on $(h,t_1,t_2)$ (in particular
$9/5<h<\min\{2,5-8t_1\}$ and $0<t_2\le t_1<1$) precisely ensure that the exponents in these
powers of~$a$ are positive, so that the approximation improves as $a\to0$ while the number of
dimers $\aleph\sim d_{\rm out}^{-3}$ increases.

From a modelling viewpoint, the inequalities in \eqref{condi-invertibility} have a natural interpretation: they require
that the effective ``single-dimer strength'' (controlled by the norms of $\mathbf P^{(1)}_{0,1}$ and
$\mathbf P^{(2)}_{0,2}$ and by the contrast parameters $c_0,d_0,\eta_0,\eta_2$) and the average number
of neighbours (controlled by $d_{\rm out}$) are small enough so that the multiple-scattering
operator remains a contraction. In particular, the theorem identifies a quantitative regime in which
multiple scattering among many strongly resonant dimers is still well controlled and the
Foldy--Lax reduction remains valid.

\subsection*{Numerical and modelling uses}

The reduction \eqref{eq:FL-main}--\eqref{eq:FL-algebraic} is directly usable in several ways.

\medskip
\noindent\emph{(a) Fast forward simulations.}
Instead of solving the full Maxwell boundary value problem for a complicated cluster, one may:
\begin{enumerate}
  \item compute (or precompute) the polarization tensors $\mathbf P^{(j)}_{0,\ell}$, $j, \ell=1, 2$, associated with the
        reference shapes $B_1,B_2$;
  \item assemble the matrices $\mathfrak{B}_m$ and $\Psi_{mj}$ using \eqref{Helmholtz-kernel}--\eqref{detail-B_m};
  \item solve the $12\aleph \times 12\aleph$ linear system \eqref{eq:FL-algebraic} for~$\mathfrak{U}$;
  \item reconstruct the scattered field or the far field from
        \eqref{eq:FL-main}--\eqref{eq:FL-main-ff}.
\end{enumerate}
For large $\aleph$, this provides a significant reduction in complexity, and the structure of the
kernel (dependence on $\Phi_k$ and $\Upsilon_k$) is well suited to fast algorithms such as fast multipole or hierarchical solvers.

%\medskip\noindent
%\textbf{(b) Inverse scattering and imaging of dimers.}
%In matrix form, \eqref{eq:FL-main-ff} can be rewritten as
%\[
%E_{\infty}(\hat x,\theta,p) \approx \sum_{m=1}^{\mathcal N} e_m(\hat x)\cdot U_m,
%\]
%where the vectors $e_m(\hat x)$ depend only on the geometry and the observation/illumination
%directions. For multi-static measurements (several $(\hat x,\theta,p)$), the far-field data define a
%multi-static response matrix whose columns are approximately linear combinations of the $U_m$.
%This structure can be exploited for:
%\begin{itemize}
%  \item reconstruction of the dimer locations $z_{m\ell}$ (e.g.\ by MUSIC- or factorization-type
 %       methods applied to the far-field operator);
 % \item estimation of effective polarizabilities $(Q_{m\ell},R_{m\ell})$, and hence of local material
 %       parameters, from measured data;
 % \item detection or classification of different types of dimers inside a heterogeneous cluster.
%\end{itemize}
%Thus Theorem \ref{main-1} provides a rigorous justification for using dipole-based sampling or
%factorization methods in imaging problems involving hybrid dimers.

\medskip

\noindent\emph{(b) Analysis of collective resonances.}
The spectrum of the matrix $\mathcal B$ in \eqref{eq:FL-algebraic} encodes not only the single-dimer
resonances, via the entries involving $\mathbf P^{(1)}_{0,1}$ and $\mathbf P^{(2)}_{0,2}$, but also the collective
(resolved) resonances of the cluster, induced by the interaction matrices $\Psi_{mj}$. Frequencies
at which $\mathcal B$ develops small singular values correspond to:
\begin{itemize}
  \item enhanced scattering and strong localisation of the internal fields;
  \item possible opening of band gaps in periodic arrangements of dimers;
  \item shifts and splitting of single-dimer resonances due to coupling within the cluster.
\end{itemize}
The explicit dependence of $\mathcal B$ on the geometry and material parameters makes it possible
to study such collective resonances in a purely algebraic setting.

\medskip

\noindent\emph{(c) Input for effective-medium and metamaterial models.} \,
%%%%%%%%%%%%%%%%%%%%%%%%%%%%%%%%%%%%%%%%%%%%%%%%Corollary%%%%%%%%%%%%%%
If we further extract the very dominant term $(Q_{m_1}, R_{m_2})$ in $\mathfrak{U}_m$ from the solution to the algebraic system \eqref{eq-al-D1-}, for any $m$ such that $1\leq m\leq \aleph$, it leads to the following corollary.

\begin{corollary}\label{corollary-main}
    Under \textbf{Assumption \ref{ASass}} and the invertibility condition \eqref{condi-invertibility}, for $t_1, h \in \mathbb{R}^+$ such that \eqref{conditions-t-h} is fulfilled,
% \begin{equation}\label{conditions-t-h}
% \frac{k^4}{(4\pi)^4\, c_0^2 \, d_0^2}<1\quad \mbox{and}\quad \frac{9}{5}<h<\min\{2, 5-8t_1\},
% \end{equation} 
% with $c_0, d_0$ being positive constants, 
there holds the following expansion for the scattered field and  far field, respectively as
% \begin{equation}\label{far-field-corollary}
%     E^\infty(\hat{x}, \theta, p)=\sum_{m=1}^\aleph e^{-i k \hat{x}\cdot z_{m_0}}\left[ (I-\hat{x}\otimes\hat{x})  R_{m_2} - i k \hat{x}\times Q_{m_1} \right] + \O\left( a^{4-h-\frac{17}{2}t_2} \right),
% \end{equation}
\begin{equation}\notag
      E^s ({x}, \theta, p, \omega) = - k^2 \sum_{m=1}^\aleph  \left[ \underset{y}{\nabla}\Phi_k(x, z_{m_0}) \times \mathring{Q}_{m_1} - \Upsilon_k(x, z_{m_0})\cdot \mathring{R}_{m_2} \right] + \O(a^{10-3h - \frac{23}{2}t_2}),
\end{equation}
and
\begin{equation}\label{far-field-corollary}
    E^\infty (\hat{x}, \theta, p) = - \sum_{m=1}^\aleph e^{- i k\hat{x}\cdot z_{m_0}} \left[ i k \hat{x}\times \mathring{Q}_{m_1}  - (I-\hat{x}\otimes\hat{x})\mathring{R}_{m_2}\right] + \O(a^{10-3h - \frac{23}{2}t_2}),
\end{equation}
where $z_{m_0}$ is the intermediate point between $z_{m_1}$ and $z_{m_2}$ inside $D_m$, and $(\mathring{Q}_{m_1}, \mathring{R}_{m_2})_{m=1}^\aleph$ is the solution to the following system
\begin{equation}\notag
    \mathbb{A}_m \begin{pmatrix}
        \mathring{Q}_{m_1}\\
        \mathring{R}_{m_2}
    \end{pmatrix} - \sum_{j=1 \atop j\neq m}^\aleph  \mathbb{C}_{mj} \begin{pmatrix}
        \mathring{Q}_{j_1}\\
        \mathring{R}_{j_2}
    \end{pmatrix} = \begin{pmatrix}
        H^{Inc}(z_{m_1})\\
        E^{Inc}(z_{m_2})
    \end{pmatrix},
\end{equation}
where 
\begin{equation}\label{def-matrix-corollary}
{\small
\begin{aligned}
\mathbb{A}_m &:=\begin{pmatrix}
a^{h-3}\, \tfrac{\pm c_0}{i k \,\eta_0}\,\bigl({\bf P}_{0, 1}^{(1)}\bigr)^{-1} 
& i k\, \bigl(\nabla\Phi_k(z_{m_1}, z_{m_2})\times\bigr)\\
- k^2\, \bigl(\nabla\Phi_k(z_{m_1}, z_{m_2})\times\bigr) 
& a^{h-3}\, \tfrac{\pm d_0}{\eta_2}\, \bigl({\bf P}_{0, 2}^{(2)}\bigr)^{-1}
\end{pmatrix},\\[2pt]
\mathbb{C}_{mj} &:=\begin{pmatrix}
-i k^3 \, \Upsilon_k(z_{m_1}, z_{j_1}) 
& - i k\, \bigl(\nabla\Phi_k (z_{m_1}, z_{j_2})\times\bigr)\\
k^2\, \bigl(\nabla\Phi_k(z_{m_2}, z_{j_1})\times\bigr) 
& k^2\, \Upsilon_k(z_{m_2}, z_{j_2})
\end{pmatrix}.
\end{aligned}
}
\end{equation}
% \begin{eqnarray}\label{system-corollary}
%     &&\begin{pmatrix}
%         a^{h-3} \frac{\pm c_0}{i k \eta_0}({\bf P}_{0, 1}^{(1)})^{-1} & i k \nabla\Phi_k(z_{m_1}, z_{m_2})\times\\
%         - k^2 \nabla\Phi_k(z_{m_1}, z_{m_2})\times & a^{h-3} \frac{\pm d_0}{\eta_2} ({\bf P}_{0, 2}^{(2)})^{-1}
%     \end{pmatrix}\begin{pmatrix}
%         Q_{m_1}\\
%         R_{m_2}
%     \end{pmatrix}\notag\\
%     &-&\sum_{j=1 \atop j\neq m}^\aleph 
%     \begin{pmatrix}
%         -i k^3 \Upsilon_k(z_{m_0}, z_{j_0}) & - i k \nabla\Phi_k (z_{m_0}, z_{j_0})\times\\
%         k^2 \nabla\Phi_k(z_{m_0}, z_{j_0})\times & k^2 \Upsilon_k(z_{m_0}, z_{j_0})
%     \end{pmatrix}\begin{pmatrix}
%         Q_{j_1}\\
%         R_{j_2}
%     \end{pmatrix} = \begin{pmatrix}
%         H^{Inc}(z_{m_0})\\
%         E^{Inc}(z_{m_0})
%     \end{pmatrix}.
% \end{eqnarray}
\end{corollary}
% We refer to \textbf{Subsection SM-{\red{cite 1}}}
% % Appendix \ref{sec:appendix-coro} 
% for the detailed proof of Corollary \ref{corollary-main}.

%\begin{remark}
%    The matrix $\mathbb{A}_m$ in \eqref{def-matrix-corollary} implies the study on {\red{``bi-anisotropic"}} property of the medium.
%\end{remark}

%\begin{remark}
%    The coefficient matrix of the algebraic system \eqref{eq-al-D1-} and matrix $\mathbb{A}_m$ in \eqref{def-matrix-corollary} both demonstrate the ability to manipulate the electromagnetic wave by tuning the sign of $\pm c_0$ and $\pm d_0$ to generate single and eventually \red{double negativity} with proper incident frequencies.
%\end{remark}

%%%%%%%%%%%%%%%%%%%%%%%%%%%%%%

While Corollary~\ref{corollary-main} concentrates on the most dominant components $(Q_{m_1},R_{m_2})$ and leads
naturally to a bi-anisotropic effective model based on a $6\times 6$ polarizability matrix, the full
system \eqref{eq:FL-algebraic} contains the finer information on
$(Q_{m_2},R_{m_1})$ as well. These subdominant components can contribute to higher-order
corrections in homogenization (e.g.\ spatial dispersion, non-locality) and to the precise shape of
the band structure in periodic media. In particular:
\begin{itemize}
  \item keeping the full vector $\mathfrak{U}_m$ allows one to derive higher-order $\mathcal O(a^\alpha)$ corrections to the effective tensors $(\varepsilon_{\rm eff},\mu_{\rm eff},\xi,\zeta)$ which is introduced in \eqref{eq:bi-aniso-constitutive} and \eqref{eq:eff-tensors};
  \item the explicit structure of $\mathcal B$ makes it possible to explore systematically how
        choices of $(B_1,B_2)$, $(c_0,d_0)$ and the distribution of dimer centers influence the
        effective response.
\end{itemize}

\medskip

 Next, we provide a formal effective-medium interpretation of Corollary~\ref{corollary-main} and of the
algebraic system \eqref{eq:FL-algebraic}. The goal is to show how a cluster of hybrid dielectric--plasmonic
dimers may generate effective media with negative permittivity, negative permeability and
bi-anisotropic constitutive laws. A rigorous derivation of the corresponding homogenized model and
of the associated limiting constitutive parameters will be the subject of a subsequent work.

\subsection*{Dipole-level polarizability of a single hybrid dimer}

For each dimer $D_m = D_{m_1}\cup D_{m_2}$, let
\[
\mathcal{U}_m
:=
\begin{pmatrix}
\mathring Q_{m_1}\\[1mm] \mathring R_{m_2}
\end{pmatrix}\in\C^{6},
\qquad
\mathcal F_m^{\rm inc}
:=
\begin{pmatrix}
H^{\Inc}(z_{m_1})\\[1mm] E^{\Inc}(z_{m_2})
\end{pmatrix}\in\C^{6},
\]
where 
% $z_{m_0}$ is the intermediate point of the dimer and
$(\mathring Q_{m_1}, \mathring R_{m_2})$ are the dominant
degrees of freedom extracted in Corollary~\ref{corollary-main}. Up to the already controlled error term in
Corollary~\ref{corollary-main}, with notaiton \eqref{def-matrix-corollary}, the algebraic system can be written in the compact form
\begin{equation}\label{eq:discrete-system}
\mathbb A_m \mathcal U_m - \sum_{\substack{j=1\\ j\neq m}}^{\aleph} \mathbb C_{mj}\,\mathcal U_j
= \mathcal F_m^{\rm inc},
\qquad 1\leq m\leq\aleph.
\end{equation}
% with $A_m$ and $C_{mj}$ given by~(1.27). For convenience we recall
% \begin{equation}\label{eq:Am-explicit}
% A_m =
% \begin{pmatrix}
% a^{h-3}\dfrac{\pm c_0}{ik\eta_0}\big(P_{0,1}^{(1)}\big)^{-1}
% &
% ik\,\nabla\Phi_k(z_{m1},z_{m2})\times
% \\[2mm]
% -k^2\,\nabla\Phi_k(z_{m1},z_{m2})\times
% &
% a^{h-3}\dfrac{\pm d_0}{\eta_2}\big(P_{0,2}^{(2)}\big)^{-1}
% \end{pmatrix},
% \end{equation}
% and
% \begin{equation}\label{eq:Cmj-explicit}
% C_{mj} =
% \begin{pmatrix}
% -ik^3\,\Upsilon_k(z_{m0},z_{j0}) & -ik\,\nabla\Phi_k(z_{m0},z_{j0})\times\\[2mm]
% k^2\,\nabla\Phi_k(z_{m0},z_{j0})\times & k^2\,\Upsilon_k(z_{m0},z_{j0})
% \end{pmatrix}.
% \end{equation}

It is natural to interpret the interaction term in~\eqref{eq:discrete-system} as part of a
\emph{local field}. More precisely, we define the local field at intermediate point $z_{m_0}$ of the dimer by
\begin{equation}\label{eq:local-field}
\mathcal F_m^{\rm loc}
:=
\begin{pmatrix}
H^{\rm loc}(z_{m_0})\\[1mm] E^{\rm loc}(z_{m_0})
\end{pmatrix}
:=
\mathcal F_m^{\rm inc}
\;+\;
\sum_{\substack{j=1\\ j\neq m}}^{\aleph} \mathbb C_{mj}\,\mathcal U_j.
\end{equation}
Then~\eqref{eq:discrete-system} can be rewritten as
\begin{equation}\label{eq:single-dimer-law}
\mathbb A_m \mathcal U_m = \mathcal F_m^{\rm loc},
\qquad
\mathcal U_m = \mathcal P_m\,\mathcal F_m^{\rm loc},
\end{equation}
where
\begin{equation}\label{eq:single-dimer-polarizability}
\mathcal P_m := \mathbb A_m^{-1}
=
\begin{pmatrix}
\alpha^{(m)}_{HH} & \alpha^{(m)}_{HE}\\[1mm]
\alpha^{(m)}_{EH} & \alpha^{(m)}_{EE}
\end{pmatrix}
\in\C^{6\times 6}.
\end{equation}
The $3\times 3$ block matrices $\alpha^{(m)}_{HH},\alpha^{(m)}_{EE}$ describe, respectively, the 
magnetic and electric responses of the dimer, while the off-diagonal blocks
$\alpha^{(m)}_{HE},\alpha^{(m)}_{EH}$ encode magneto--electric coupling.

On the other hand, the far field expansion \eqref{far-field-corollary} can be written as
\begin{equation}\label{eq:farfield-dipoles}
E^\infty(\hat x,\theta,p)
=
\sum_{m=1}^{\aleph} e^{-ik\hat x\cdot z_{m_0}}
\big[
(I-\hat x\otimes\hat x)\, \mathring R_{m_2}
-ik\,\hat x\times \mathring Q_{m_1}
\big]
+ O\Big(a^{10-3h - \frac{23}{2}t_2}\Big).
\end{equation}
The factor $(I-\hat x\otimes\hat x) \mathring R_{m_2}$ is the standard electric-dipole contribution, while
$-ik\,\hat x\times \mathring Q_{m_1}$ is the magnetic-dipole contribution. Thus each dimer behaves, at the level
of~\eqref{eq:farfield-dipoles}, as a co-located electric and magnetic dipole with moments
$\mathring R_{m_2}$ and $\mathring Q_{m_1}$, linearly related to the local fields by~\eqref{eq:single-dimer-law}.

\subsection*{From discrete dimers to effective constitutive tensors}

We now give a \emph{formal} passage from the discrete dipole system~\eqref{eq:single-dimer-law} to a
macroscopic constitutive description. For simplicity, we assume that all dimers are identical and
identically oriented so that \(\mathcal P_m\equiv \mathcal  P := \begin{pmatrix}
    \alpha_{HH} & \alpha_{HE}\\[1mm]
\alpha_{EH} & \alpha_{EE}
\end{pmatrix}\) for all \(m\), and that the dimer
centers \(z_{m_0}\) form a sufficiently regular distribution in a reference region
\(\Omega\subset\R^3\).

Let \(V\subset\Omega\) be a small cube such that \(a\ll \diam(V)\ll 1\), and denote by
\begin{equation}\label{eq:density}
\rho_a := \frac{\#\{m:\ z_{m_0}\in V\}}{|V|}
\sim d_{\rm out}^{-3}
\end{equation}
the (local) number density of dimers, where \(d_{\rm out}\) is the typical inter-dimer distance.
We define the local averaged magnetic and electric polarizations by
\begin{equation}\label{eq:macro-M-P}
M(x) := \frac{1}{|V|}\sum_{z_{m_0}\in V} Q_{m_1},
\qquad
P(x) := \frac{1}{|V|}\sum_{z_{m_0}\in V} R_{m_2},
\qquad x\in V.
\end{equation}
Based on Corollary~\ref{corollary-main}, it is reasonable to approximate the local fields by smooth
macroscopic fields,
\[
H^{\rm loc}(z_{m_0})\approx H(x),\qquad E^{\rm loc}(z_{m_0})\approx E(x),\qquad x\in V.
\]
Using~\eqref{eq:single-dimer-law} and~\eqref{eq:density}--\eqref{eq:macro-M-P}, this leads to the
formal constitutive relation
\begin{equation}\label{eq:macro-chi}
\begin{pmatrix}
M(x)\\[1mm] P(x)
\end{pmatrix}
=
\rho_a\,\mathcal P
\begin{pmatrix}
H(x)\\[1mm] E(x)
\end{pmatrix}
=
\begin{pmatrix}
\chi_{HH} & \chi_{HE}\\[1mm]
\chi_{EH} & \chi_{EE}
\end{pmatrix}
\begin{pmatrix}
H(x)\\[1mm] E(x)
\end{pmatrix}.
\end{equation}
In the corresponding bi-anisotropic Maxwell model, one writes
\begin{equation}\label{eq:bi-aniso-constitutive}
D=\varepsilon_{\rm eff}E+\xi H,\qquad
B=\zeta E+\mu_{\rm eff}H,
\end{equation}
with effective tensors
\begin{equation}\label{eq:eff-tensors}
\varepsilon_{\rm eff}:=\varepsilon_0(I+\chi_{EE}),\quad
\mu_{\rm eff}:=\mu_0(I+\chi_{HH}),\quad
\xi:=\varepsilon_0\chi_{EH},\quad
\zeta:=\mu_0\chi_{HE}.
\end{equation}

The blocks of \(\mathcal P=\mathbb A_m^{-1}\) depend singularly on \(a\) through the resonant assumptions
and the geometric scalings. Combining the dominant polarizability scaling
\(|\alpha_{HH}|,|\alpha_{EE}|=\O(a^{3-h})\) and \(|\alpha_{HE}|,|\alpha_{EH}|=\O(a^{6-2h-2t_1})\) with
\(\rho_a\sim a^{-3t_2}\) yields the susceptibility scaling
\begin{equation}\label{eq:chi-scaling}
\chi_{HH},\chi_{EE}=\O(a^{3-h-3t_2}),\qquad
\chi_{HE},\chi_{EH}=\O(a^{6-2h-2t_1-3t_2}).
\end{equation}
In particular, the ``nontrivial'' diagonal regime \(3-h-3t_2=0\) leads to
\(\chi_{HE},\chi_{EH}=\O(a^{3t_2-2t_1})\), showing when magneto--electric coupling is order-one or small.

\medskip
This effective-medium interpretation is a consequence of the dipole-level Foldy--Lax description in
Corollary~\ref{corollary-main}; it is \emph{not} a homogenization theorem. Nevertheless, it indicates
that tuning the resonant parameters in \eqref{condition-on-k}, the anisotropy encoded in the polarization tensors
\(\mathbf P_{0,1}^{(1)}\), \(\mathbf P_{0,2}^{(2)}\), and the distribution scale \(d_{\rm out}\) can (formally) produce:
\begin{enumerate}[label=(\roman*)]
\item single- and double-negative windows (\(\Re\,\varepsilon_{\rm eff}<0\), \(\Re\,\mu_{\rm eff}<0\)),\cite{Veselago1968,Smith2000,Pendry2000}
\item bi-anisotropic (and potentially chiral) responses via \(\xi,\zeta\neq0\),\cite{Serdyukov2001,Lindell1994}
\item strongly anisotropic and possibly hyperbolic regimes,\cite{Poddubny2013}
\item near-zero behavior at frequencies where the relevant entries of \(\mathcal P\) cross through zero.\cite{AluEngheta2007,Maas2013}
\end{enumerate}
Explicit leading-order expressions supporting these implications
are recorded in Appendix~\ref{sec:appendix-susceptibilities}. A rigorous homogenization theory, in which the convergence of the discrete model to the effective
constitutive laws~\eqref{eq:bi-aniso-constitutive}--\eqref{eq:eff-tensors} is established together with
precise conditions for each regime, will be developed in a forthcoming work.
\bigskip

%%%%%%%%%%%%%%%%%%%%%%%%%%%%%%%

% \begin{remark}
%     For the particular case with only one pair of dimer, namely, $\aleph=1$, we refer to \cite{cao-ghandriche-sini-dimer}, where the inner distance $t_1$ fulfills 
%     \begin{equation}\label{condition-on-h-one dimer}
%         4-h-4t_1>0.
%     \end{equation}
% \end{remark}

% \blue{
% \begin{remark}
%     The condition $4-h-4t_1>0$ in \eqref{conditions-t-h} can be ensured by \cite[Theorem 1.2]{cao-ghandriche-sini-dimer} in the study for the electromagnetic waves generated by one hybrid dieletric-plasmonic dimer, where $t:= t_1$ therein stands for the inner distance inside one dimer.
% \end{remark} }

The rest of this paper is organized as follows. In Section \ref{sec-pre}, we introduce some preliminaries including the $\mathbb{L}^2$-decomposition, and the Lippmann-Schwinger equations of integral form for the solution to \eqref{model-m}. Based on the Lippmann-Schwinger system, then we present the associated a-priori estimates. In Section \ref{sec-la-system}, the precise form of the linear algebraic system whose solution formulates the Foldy-Lax approximation, is presented on the basis of the Lippmann-Schwinger equation. Section \ref{sec-proof-main} is devoted to prove our main theorem and the corollary for the asymptotic expansions of the Foldy-Lax approximation on the basis of Section \ref{sec-pre} and Section \ref{sec-la-system}. 
% Finally, we provide the rigorous proof of of the corresponding algebraic system in Section \ref{sec-proof-la}.
% Finally, in Section \ref{sec-proof-prior}, we provide the rigorous proof on deriving the a-priori estimates given in Section \ref{sec-pre}, and the detailedly proof of the corresponding algebraic system in Section \ref{sec-proof-la}.

\section{Some preliminaries and a-priori estimates}\label{sec-pre}

%\textbf{Spaces Utilised\\}

In this section, we present some necessary preliminaries and significant a-prior estimates. For the preliminaries, we cite some key points here for the completeness of the paper, please refer to \cite{cao-ghandriche-sini-highindex} for more details.

\subsection{Decomposition of $\mathbb{L}^2$.}

The following direct sum provide a useful decomposition of $\mathbb{L}^{2}$ (see \cite{Dautry-Lions}, page 314) 
\begin{equation}\label{L2-decomposition}
\mathbb{L}^{2} \, \equiv \, \mathbb{H}_{0}\left(\div=0 \right) \overset{\perp}{\oplus} \mathbb{H}_{0}\left(Curl=0 \right) \overset{\perp}{\oplus} \nabla \mathcal{H}armonic
\end{equation}   
where 
\begin{eqnarray*}
	\mathbb{H}_{0}\left(\div=0 \right) &:=& \left\lbrace E \in \left( \mathbb{L}^{2}(D)\right)^{3}, \, \div E = 0, \, \nu \cdot E = 0 \, \; \text{on} \;\, \partial D \right\rbrace ,\\
	\mathbb{H}_{0}\left(Curl =0 \right) &:=& \left\lbrace E \in \left( \mathbb{L}^{2}(D)\right)^{3}, \, Curl \, E = 0, \, \nu \times E = 0 \, \; \text{on} \;\, \partial D \right\rbrace,
\end{eqnarray*}
and 
\begin{equation*}
\nabla \mathcal{H}armonic := \left\lbrace E: \; E = \nabla \psi, \, \psi \in \mathbb{H}^{1}(D), \, \Delta\psi=0 \right\rbrace.
\end{equation*}
From the decomposition \eqref{L2-decomposition}, 
we define $\overset{1}{\mathbb{P}}, \overset{2}{\mathbb{P}}$ and $\overset{3}{\mathbb{P}}$ to be the natural projectors as follows
\begin{equation}\label{project}
	\overset{1}{\mathbb{P}} := \mathbb{L}^{2} \longrightarrow  \mathbb{H}_{0}\left(\div=0 \right), \;\;\;
	\overset{2}{\mathbb{P}} := \mathbb{L}^{2} \longrightarrow  \mathbb{H}_{0}\left(Curl = 0 \right) \;\; \text{and} \;\;
	\overset{3}{\mathbb{P}} := \mathbb{L}^{2} \longrightarrow  \nabla \mathcal{H}armonic.
\end{equation}

\begin{remark}
    Following the notations introduced in \textbf{Assumption \ref{ASass}}, we respectively denote 
    \begin{itemize}
        \item  $\left(\lambda_n^{(1)}(B_\ell); e_{n, B_\ell}^{(1)}\right)$ as the eigensystem of $N_{B_\ell}(\cdot)$ onto the first subspace $\mathbb{H}_0(\div=0)$;

        \item $\left(\lambda_n^{(2)}(B_\ell); e_{n, B_\ell}^{(2)}\right)$ as the eigensystem of $N_{B_\ell}(\cdot)$ onto the second subspace $\mathbb{H}_0(Curl=0)$;

        \item $\left(\lambda_n^{(3)}(B_\ell); e_{n, B_\ell}^{(3)}\right)$ as the eigensystem of $\nabla M_{B_\ell}(\cdot)$ onto the third subspace $\nabla \mathcal{H}armonic$.
    \end{itemize}
    Besides, the following vanishing integral character holds 
    \begin{equation}\label{VIP}
        \int_{B} e_{n, B}^{(\ell)}(x) \, dx \, = \, 0, \quad \text{for} \;\; n \in \mathbb{N} \quad \text{and} \quad \ell = 1, 2.
    \end{equation}
\end{remark}

\subsection{Lippmann-Schwinger integral formulation of the solutions.} 

Analogously to $(\ref{NM0})$, for any vector function $F$, we define the Newtonian potential operator $N^{k}_{D}(\cdot)$ and the Magnetization operator $\nabla M^{k}_{D}(\cdot)$ as follows:
\begin{equation}\label{N-M opera}
N^{k}_{D}(F)(x) \, := \, \int_{D} \Phi_{k}(x,y)F(y)\,dy \quad \text{and} \quad \nabla M^{k}_{D}(F)(x) \, := \, \underset{x}{\nabla}\int_{D}\underset{y}{\nabla}\Phi_{k}(x,y)\cdot F(y)\,dy,
\end{equation}
where $\Phi_{k}(\cdot,\cdot)$ is given by $(\ref{Helmholtz-kernel})$. The following Lippmann-Schwinger equation holds, constructed by the operators $N^{k}_{D}(\cdot)$ and $\nabla M^{k}_{D}(\cdot)$, whose solution solves $(\ref{model-m})$.

\begin{proposition}\label{prop-LS}
	The solution to the electromagnetic scattering problem \eqref{model-m} satisfies
	\begin{equation}\label{LS eq2}
		E^T(x) \, + \,  \nabla M^k(\eta \,E^T)(x) \, - \, k^2  \, N^{k}(\eta \, E^T)(x) \, =  \, E^{Inc}(x),\quad x\in D,
	\end{equation}
	where $\eta(\cdot)$ is defined by \eqref{def-eta}.
\end{proposition}
\begin{proof}
The proposition can be proved by utilizing the Stratton-Chu formula directly, see \cite[Theorem 6.1]{colton2019inverse} for more detailed discussions.
\end{proof}
 \begin{remark} Two remarks are in order. 
 \begin{enumerate}
     \item For simplification reasons, we use the notation $\nabla M(\cdot)$ (respectively, $N(\cdot)$, $E(\cdot)$) instead of $\nabla M^{0}(\cdot)$ (respectively, $N^{0}(\cdot)$, $E^{T}(\cdot)$) in the subsequent analyses.
     \item[]
     \item Under the $\mathbb{L}^2$-space decomposition introduced in \eqref{L2-decomposition}, by Green's formulas, we have
 	\begin{equation}\label{grad-M-1st-2nd}
 	\forall \; E \in \mathbb{H}_{0}\left(\div=0 \right) , \;\, \nabla M(E) = 0 \, \quad \text{and} \quad \forall \; E \in \mathbb{H}_{0}\left(Curl=0 \right) , \;\, \nabla M(E) = E.
 	\end{equation}
 	More properties for the Magnetization operator, such as self-adjointness, positivity, spectrum, boundedness ($\left\Vert \nabla M \right\Vert = 1$) , invariance of $\nabla \mathcal{H}armonic$, etc., can be found in \cite{friedman-pasciak,raevskii}.
 \end{enumerate}
 	
 \end{remark}

%%%%%%%%%%%%%%%%%%%%%%%%%%%%%%%%%%%%%%%%%%%%%

% \subsection{A-priori Estimates.}\label{subsec-eigen} Based on the decomposition \eqref{L2-decomposition}, we present here some necessary a-prior estimates derived from the Lippmann-Schwinger equation \eqref{LS eq2}, which play an important role in the proof of our main results. The proof of the propositions and lemmas in this subsection shall be given later in Section \ref{sec-proof-prior}.

\subsection{A-priori estimates.}

Suppose $E$ solves \eqref{model-m} with the integral formulation \eqref{LS eq2}. Then the projection of $E$ onto three subspaces $\mathbb{H}_0(\div=0)$, $\mathbb{H}_0(Curl=0)$ and $\nabla \mathcal{H}armonic$ can be respectively represented as $\overset{1}{\mathbb{P}}(E)$, $\overset{2}{\mathbb{P}}(E)$ and $\overset{3}{\mathbb{P}}(E)$.
Then we have the following estimates.

\vspace*{5px}

\begin{proposition}\label{lem-es-multi}
	Under Assumptions (\uppercase\expandafter{\romannumeral1}, \uppercase\expandafter{\romannumeral2}, \uppercase\expandafter{\romannumeral3}, \uppercase\expandafter{\romannumeral4}),
    consider the electromagnetic scattering problem \eqref{model-m} for the cluster of dimers $D=\cup_{m=1}^\aleph D_m=\cup_{m=1}^\aleph(D_{m_1}\cup D_{m_2})$, $m=1, 2, \cdots, \aleph$. Then for $t_1\in(0,1)$, $h\in\mathbb{R}^+$ and the constants $c_0, d_0$ satisfying  \begin{equation}\notag
\frac{9}{5}<h<\min\{2, 5-6t_1, 4-4t_1\}\quad\mbox{and}\quad \frac{k^4}{(4\pi)^4c_0^2 d_0^2}<1,
\end{equation}  
there holds the estimation 
	\begin{align}\label{max-P1P3}
\max_{1\le m\le \aleph}\left\lVert \overset{1}{\mathbb{P}}\bigl(\tilde{E}_{m_1}\bigr)\right\rVert_{\mathbb{L}^2(B_1)}
&\lesssim a^{4-2h} d_{\rm out}^{-\frac{9}{2}},\notag\\
\max_{1\le m\le \aleph}\left\lVert \overset{1}{\mathbb{P}}\bigl(\tilde{E}_{m_2}\bigr)\right\rVert_{\mathbb{L}^2(B_2)}
&\lesssim a^{7-2h}d_{\rm in}^{-3}d_{\rm out}^{-\frac{9}{2}}+ a^{7-2h} d_{\rm out}^{-\frac{17}{2}},\notag\\
\max_{1\le m\le \aleph}\left\lVert \overset{3}{\mathbb{P}}\bigl(\tilde{E}_{m_1}\bigr)\right\rVert_{\mathbb{L}^2(B_1)}
&\lesssim a^{8-2h}d_{\rm in}^{-3} d_{\rm out}^{-\frac{9}{2}},\notag\\
\max_{1\le m\le \aleph}\left\lVert \overset{3}{\mathbb{P}}\bigl(\tilde{E}_{m_2}\bigr)\right\rVert_{\mathbb{L}^2(B_2)}
&\lesssim a^{3-2h} d_{\rm out}^{-\frac{9}{2}}.
\end{align} 
    where $\overset{1}{\PP}(\tilde{E}_{m_\ell})$ and $\overset{3}{\PP}(\tilde{E}_{m_\ell})$, $\ell=1, 2$, are the wave fields after scaling $\overset{1}{\PP}({E}_{m_\ell})$ and $\overset{3}{\PP}({E}_{m_\ell})$ from $D_{m_\ell}$ to $B_\ell$ accordingly.
\end{proposition} 
\begin{proof}
    See {Section SM-1}.
    % Section \ref{sec-proof-prior}.
\end{proof}
% We postpone the proof of Proposition \ref{lem-es-multi} to Subsection \ref{sec-proof-prior} to clear the structure of the paper.

We introduce the vector functions $W_\ell$ and $V_\ell$, $\ell=1, 2$, from \cite{cao-ghandriche-sini-dimer} that

\begin{definition}\label{def-scattering coeff}
    	We define $W_{1}, W_{2}, V_{1}$ and $V_{2}$ to be solutions of,  
    	\begin{eqnarray}\label{AI1}
    		\left(I + \eta_1 \, \nabla M^{-k}_{D_{1}} \, - \, k^{2} \, \eta_1 \, N^{-k}_{D_{1}} \right)\left( W_1 \right)(x) &=& \mathcal{P}(x,z_1), \quad \; \, \; \quad x \in D_1, \label{DefW1}  \\
    		\left(I + \eta_2 \, \nabla M^{-k}_{D_{2}} \, - \, k^{2} \, \eta_2 \, N^{-k}_{D_{2}} \right)\left( W_2 \right)(x) &=& \overset{1}{\PP}\left(\mathcal{P}(x,z_2)\right), \quad x \in D_2,   \label{DefW2} \\
            \left(I + \eta_m \, \nabla M^{-k}_{D_{m}} \, - \, k^{2} \, \eta_m \, N^{-k}_{D_{m}} \right)\left( V_m \right)(x) &=& I_{3}, \quad \;\;  \qquad \qquad x \in D_m, \quad m = 1, 2, \label{def-Vm}
    	\end{eqnarray}
    	where, for $m=1,2$, the operators $\nabla M^{-k}_{D_{m}}(\cdot)$ and $N^{-k}_{D_{m}}(\cdot)$ are the adjoint operators to $\nabla M^k_{D_{m}}(\cdot)$ and $N^k_{D_{m}}(\cdot)$, introduced in \eqref{N-M opera}, and $\mathcal{P}(x, z)$ is the matrix expressed by
    	\begin{equation}\label{Def-matrix-p}
    	\mathcal{P}(x, z) \, := \, \left(\begin{array}{c}
    	(x-z)_1 \, I_{3}\\ 
    	(x-z)_2 \, I_{3}\\ 
    	(x-z)_3 \, I_{3}
    	\end{array} \right).
    	\end{equation} 
    \end{definition}
    \medskip 
    Then, there holds the following estimations w.r.t. $W_\ell, V_\ell$, $\ell=1, 2$, defined above.
\begin{proposition}\label{prop-scoeff}
        For $h \in (0,2)$, under assumption $(\ref{condition-on-k})$, the following estimations hold:
        \begin{enumerate}
            \item Regarding the scattering parameter $W_{1}$, defined by \eqref{AI1}, we have  \begin{equation}\label{*add0}
		\left\lVert \overset{1}{\mathbb{P}}\left(\tilde{W}_1\right)\right\rVert_{\mathbb{L}^2(B_1)}=\mathcal{O}(a^{1-h}) \quad \text{and} \quad
		\left\lVert \overset{j}{\mathbb{P}}\left(\tilde{W}_1\right)\right\rVert_{\mathbb{L}^2(B_1)}=\mathcal{O}(a^{3}), \quad \text{for} \quad j = 2, 3.	
		\end{equation}
            \item[]
            \item Regarding the scattering parameter $W_{2}$, defined by $(\ref{DefW2})$, we have  
        \begin{equation}\label{prop-es-W}
			\left\lVert \overset{1}{\mathbb{P}}\left(\tilde{W}_2\right)\right\rVert_{\mathbb{L}^2(B_2)}=\mathcal{O}(a),\quad 
			 \overset{2}{\mathbb{P}}\left(\tilde{W}_2\right) \, = \, 0 \quad \text{and} \quad
			\left\lVert\overset{3}{\mathbb{P}}\left(\tilde{W}_2\right)\right\rVert_{\mathbb{L}^2(B_2)}=\mathcal{O}(a^{5-h}).
		\end{equation}
        \item[] 
        \item Regarding the scattering parameter $V_{m}$, defined by \eqref{def-Vm}, for $m = 1, 2$, we have
        \begin{equation}\label{es-V12}
			\left\lVert \tilde{V}_1\right\rVert_{\mathbb{L}^2(B_1)}=\mathcal{O}\left(a^2\right)\quad\mbox{and}\quad\left\lVert \tilde{V}_2\right\rVert_{\mathbb{L}^2(B_2)}=\mathcal{O}\left(a^{-h}\right).
		\end{equation}
        \item[] 
        \end{enumerate}
	\end{proposition}
\begin{proof}
    The proof of {Proposition \ref{prop-scoeff}} follows the same argument to \cite[Proposition 2.7]{cao-ghandriche-sini-dimer} for the study for only one hybrid dieletric-plasmonic dimer case, and we refer the detailed proof therein.
\end{proof}

\begin{remark}
    From \cite[Proposition 2.5]{cao-ghandriche-sini-dimer} and \cite[Proposition 2.2]{cao-ghandriche-sini-highindex}, we can know that for $m=1, 2, \cdots, \aleph$, there holds
    \begin{equation}\notag
        \overset{2}{\PP}(E_{m_1})=\overset{2}{\PP}(E_{m_2})=0.
    \end{equation}
    Therefore, we can split $E_{m_\ell}$ as
    \begin{equation}\label{E-split}
        E_{m_\ell}=\overset{1}{\PP}\left(E_{m_\ell}\right)+ \overset{3}{\PP}\left(E_{m_\ell}\right)\quad \mbox{for}\quad m=1, 2, \cdots,\aleph\quad\mbox{and}\quad \ell=1, 2.
    \end{equation}
\end{remark}

\section{Linear algebraic system and its invertibility}\label{sec-la-system}

In this section, we present the precise expression of the linear algebraic system, whose solution contributes to formulate the far field asymptotic expansion of our main theorem, generated by the electromagnetic scattering by a cluster of hybric dieletric-plasnomic dimers.
% Suppose $m=1$ and $d:=d_{\rm in}$, for only one dimer scattering, based on \cite[Corollary 3.3]{cao-ghandriche-sini-dimer}, there holds the following expressions for the dominant terms associated with the corresponding linear algebraic system.

% \begin{proposition}\label{prop-coro-one dimer}
% 	For $t$ and $h$ fulfilling that
% 	\begin{equation*}
% 	7 \, - \, 2h \, - \, 7t \, > \, 0 \quad \text{and} \quad 3 \, - \,h \, - \, 2t \, > \, 0, 
% 	\end{equation*}
% 	then there holds the following expressions with respect to the dominant term $Q_{1}$ and $R_{2}$ as
% 	\begin{eqnarray}
% 	Q_{1}&=&\frac{i \eta_0 k}{\pm c_0}a^{3-h}{\bf P}_{0, 1}^{(1)}H^{inc}(z_0)+\O\left(k^2 a^{3-h}d\right)+\O\left( k^2 a^{13-4h}d^{-9}\right),\notag\\
% 	R_2&=&\frac{\eta_2}{\pm d_0}a^{3-h} {\bf P}_{0, 2}^{(2)}E^{inc}(z_0)+\O\left(k a^{3-h}d\right)+\O\left( a^{7-h}d^{-4}\right)+\O\left(k^2 a^{7-2h}d^{-3}\right).\notag
% 	\end{eqnarray}  
% 	where $z_0$ denotes the intermediate point between $z_1$ and $z_2$, and ${\bf P}_{0, 1}^{(1)}$, ${\bf P}_{0, 2}^{(2)}$ are given by \eqref{def-all-tensor}.
% \end{proposition}

% We refer to \cite[Section 6]{cao-ghandriche-sini-dimer} for the detailed proof of this proposition.

% For the cluster of dimers scattering $D_m$, $m=1, 2, \cdots, \aleph$, $\ell=1,2$, 
For $m=1, 2, \cdots, \aleph$, denote $\overset{1}{\mathbb{P}}(E_{m_\ell})$ and $\overset{3}{\mathbb{P}}(E_{m_\ell})$ as the projection of the total wave $E_{m_\ell}:=E|_{D_{m_\ell}}$ onto the subspaces $\mathbb{H}_0(\div=0)$ and $\nabla\mathcal{H}armonic$, respectively. Due to the relation $(\ref{Hdiv-curl})$,
%\begin{equation}\label{Hdiv-curl}
%\mathbb{H}_0(\div=0)=Curl(\mathbb{H}_0(Curl)\cap \mathbb{H}(\div=0)),
%\end{equation}
we can write
\begin{equation}\label{def-F_j}
\overset{1}{\mathbb{P}}(E_{m_\ell})=Curl (F_{m_\ell})\quad\mbox{with}\quad \nu\times F_{m_\ell}=0,\ \div(F_{m_\ell})=0.
\end{equation}
Set
\begin{equation}\label{Qj1}
\tilde{Q}_{m_\ell}:=\eta_{\ell}\int_{D_{m_{\ell}}}F_{m_\ell}(y)\,dy\quad\mbox{and}\quad \tilde{R}_{m_\ell}:=\eta_{\ell}\int_{D_{m_\ell}} \overset{3}{\mathbb{P}}\left(E_{m_\ell}\right)(y)\,dy,
\end{equation} 
where $\eta_m$ is defined in \eqref{def-eta}. Then we can derive the following linear algebraic system with respect to $\left(\tilde{\mathfrak{U}}_{m} \, := \, \left( \tilde{Q}_{m_1}, \tilde{R}_{m_1}, \tilde{Q}_{m_2}, \tilde{R}_{m_2} \right) \right)_{m=1}^{\aleph}$.

\begin{proposition}\label{prop-la}
	Under \textbf{Assumptions \ref{ASass}},
    %Assumptions (	\uppercase\expandafter{\romannumeral1}, 	\uppercase\expandafter{\romannumeral2}, 	\uppercase\expandafter{\romannumeral3}, 	\uppercase\expandafter{\romannumeral4}), 
    for $t_1\in (0, 1)$, $h\in \mathbb{R}_+$ and the constants $c_0, d_0$ such that 
	\begin{equation}\label{cond-prop-las}
		\frac{9}{5}<h<\min\{2, 5-8t_1\} \quad\mbox{and}\quad \frac{k^4}{(4\pi)^4 c_0^2 d_0^2}<1,
	\end{equation} 
	there holds the following linear algebraic system associated with the electromagnetic scattering problem \eqref{model-m} in the presence of a cluster of hybrid dieletric-plasmonic dimers $D \,:= \, \overset{\aleph}{ \underset{m=1}{\cup}} D_m$ as
    \begin{equation}\label{eq-al-D1}
            \begin{pmatrix}
                 \mathfrak{B_{1}} & - \, \Psi_{12} & \cdots & - \, \Psi_{1\aleph}  \\
                 - \, \Psi_{21} & \mathfrak{B_{2}} & \cdots & - \, \Psi_{2\aleph} \\
                 \vdots & \vdots & \ddots & \vdots \\
                 - \, \Psi_{\aleph1} & - \, \Psi_{\aleph2} & \cdots &  \mathfrak{B_{\aleph}} 
            \end{pmatrix} 
            \cdot
            \begin{pmatrix}
                \tilde{\mathfrak{U}}_{1} \\
                \tilde{\mathfrak{U}}_{2} \\
                \vdots \\
                \tilde{\mathfrak{U}}_{\aleph} 
            \end{pmatrix}
            = 
            \begin{pmatrix}
                \mathrm{S}_{1} \\
                \mathrm{S}_{2} \\
                \vdots \\
                \mathrm{S}_{\aleph}  
            \end{pmatrix} \, + \,
            \begin{pmatrix}
                \mathfrak{R}_{1} \\
                \mathfrak{R}_{2} \\
                \vdots \\
                \mathfrak{R}_{\aleph}  
            \end{pmatrix},
    %  \begin{pmatrix}
    % Error_1^{(1)} \\
    %  Error_1^{(2)} \\
    %  Error_2^{(1)} \\
    %  Error_2^{(2)}
    % \end{pmatrix},
        \end{equation}
        where $\left\{ \left\{\mathfrak{B_{m}} \right\}_{m=1}^{\aleph}; \left\{ \Psi_{mj}\right\}_{m,j=1}^{\aleph}; \left\{ \mathrm{S}_{m} \right\}_{m=1}^{\aleph} \right\}$ are given in $\textbf{Notations \ref{notbthm}},$
and for any $1\leq m \leq \aleph$, the error term is given by 
\begin{equation}\label{err-prop-la1}
\mathfrak{R}_{m}=
\begin{pmatrix}
    Error_{\tilde{Q}_{m_1}} \\
     Error_{\tilde{R}_{m_1}} \\
     Error_{\tilde{Q}_{m_2}} \\
     Error_{\tilde{R}_{m_2}}
\end{pmatrix} \, = \, \begin{pmatrix}
    \O\left(a^{10-3h} d_{\rm out}^{-\frac{17}{2}}\right) \\
    \O\left(a^{10-2h} d_{\rm in}^{-4} d_{\rm out}^{-\frac{9}{2}}\right) \\
    \O\left(a^{12-2h} d_{\rm in}^{-4} d_{\rm out}^{-\frac{9}{2}}\right) \\
    \O\left(a^{10-3h} d_{\rm out}^{-\frac{17}{2}}\right)
\end{pmatrix}.
% \quad\mbox{for any}\quad 1\leq m\leq \aleph.
\end{equation}
%\blue{The error term on the R.H.S. of \eqref{eq-al-D1} possesses
%	\begin{eqnarray}
%    Error_1^{(1)}&=&\mathcal{O}\left(a^{4-h}\right)+\O\left( a^{7-2h} d_{\rm out}^{-4}\right)+\O\left(a^{6-h}d^{-3}_{out}\right)\notag\\
%		Error_1^{(2)}&=&\mathcal{O}(a^4)+\O\left(a^{7-h} d_{\rm out}^{-4}\right)\notag\\
%        Error_2^{(1)}&=& \O\left(a^6\right) + \O\left(a^{9-h}d_{\rm in}^{-4}\right)\notag\\
%        Error_2^{(2)}&=& \O\left(a^{7-2h}d_{\rm out}^{-4}\right)+\O\left(a^{4-h}\right) + \O\left(a^{7-h}d_{\rm in}^{-4}\right)\notag\\
%        &+&  \O\left(a^{6}d_{\rm in}^{-3}\right)+\O\left(a^{6-h}d_{\rm out}^{-3}\right).
%	\end{eqnarray}  }
% 	Moreover, the linear algebraic system \eqref{eq-al-D1} is invertible under the condition
%    \begin{equation}\label{condi-invertibility-prop}
%     \frac{k^2}{2}\left( \frac{k^2 \eta_0}{c_0}|{\bf P}_{0, 1}^{(1)}|+\frac{\eta_2}{d_0}|{\bf P}_{0, 2}^{(2)}|\right)a^{3-h-3t_1}<1\, \mbox{ and }\,  \frac{k^2 \left( \frac{\max\{1, k^2\}\eta_0}{c_0}\left| {\bf P}_{0, 1}^{(1)} \right| + \frac{\eta_2}{d_0} \left| {\bf P}_{0, 2}^{(2)} \right| \right)}{1- \frac{k^2}{2} \left( \frac{k^2 \eta_0}{c_0}\left|{\bf P}_{0, 1}^{(1)} \right| +\frac{\eta_2}{d_0}\left| {\bf P}_{0, 2}^{(2)} \right| \right) a^{3-h-3t_1} } a^{3-h-3t_2} <1.
% \end{equation}
\end{proposition}
\begin{proof}
    See {Section SM-2}.
\end{proof}
\section{Proof of the main results}\label{sec-proof-main}

\subsection{Proof of Theorem \ref{main-1}}

With all the necessary knowledge presented in the previous sections, the proof {Theorem \ref{main-1}} shall be given as follows.
\begin{enumerate}
    \item Derivation of the scattered wave $E^{s}(\cdot,\theta, p, \omega)$. \\
    Repetition of the same arguments used in \cite[Section 4]{cao-ghandriche-sini-dimer}, we can derive the following approximation for the scattered field.
    \begin{equation*}
        E^{s}(x,\theta, p, \omega) \, = \, - \, k^{2} \, \sum_{m=1}^{\aleph} \sum_{\ell=1}^{2} \left[ \underset{y}{\nabla} \Phi_{k}(x,z_{m_{\ell}}) \times \tilde{Q}_{m_{\ell}} \, - \, \Upsilon_{k}(x,z_{m_{\ell}}) \cdot \tilde{R}_{m_{\ell}} \right] \, + \, Remainder, 
    \end{equation*}
    where $Remainder$ is the term given by 
    \begin{eqnarray*}
        Remainder \, &:=& \, - \, k^{2} \, \sum_{m=1}^{\aleph} \sum_{\ell=1}^{2} \eta_{\ell} \, \int_{D_{m_{\ell}}} \int_{0}^{1} \underset{y}{\nabla} \underset{y}{\nabla} \Phi_{k}(x,z_{m_{\ell}}+t(y-z_{m_{\ell}})) \cdot (y-z_{m_{\ell}}) \, dt \times F_{m_{\ell}}(y) \, dy \\
         &+&  k^{2} \, \sum_{m=1}^{\aleph} \sum_{\ell=1}^{2} \eta_{\ell} \, 
        \int_{D_{m_{\ell}}} \int_{0}^{1}  \underset{y}{\nabla} \Upsilon_{k}(x,z_{m_{\ell}}+t(y-z_{m_{\ell}})) \cdot \mathcal{P}(y, z_{m_{\ell}}) \, dt \cdot \overset{3}{\mathbb{P}}\left(E_{m_{\ell}}\right)(y) \, dy. 
    \end{eqnarray*}
    Next, we estimate the term $Remainder$. 
    \begin{equation*}
       \left\vert Remainder \right\vert \,  \lesssim  \, a^{4} \, \sum_{m=1}^{\aleph} \sum_{\ell=1}^{2} \left\vert \eta_{\ell} \right\vert \, \left[ \left\Vert \tilde{F}_{m_{\ell}} \right\Vert_{\mathbb{L}^{2}(B_{m_{\ell}})}  \, 
         +  \, 
        \left\Vert \overset{3}{\mathbb{P}}\left(\tilde{E}_{m_{\ell}}\right)\right\Vert_{\mathbb{L}^{2}(B_{m_{\ell}})}\right], 
    \end{equation*}
    which, by using $(\ref{d-a})$, $(\ref{def-eta})$, $(\ref{max-P1P3})$, and the formula
    \begin{equation}\label{def-F,A}
    \overset{1}{\PP}(E_{m_\ell}) \, = \, Curl\left( F_{m_\ell}\right)\quad\mbox{and}\quad \overset{1}{\PP}(W_1) \, = \, Curl\left(A_1\right) \quad \mbox{for}\quad m=1, 2, \cdots, \aleph\quad\mbox{and}\quad \ell=1, 2,
\end{equation}
    gives us 
    \begin{equation*}
        Remainder \, = \, \mathcal{O}\left( a^{7-2h} d_{\rm out}^{-\frac{15}{2}}\right).
    \end{equation*}
    Then,
    \begin{equation}\label{scattered field-foldy lax}
        E^{s}(x,\theta, p, \omega) \, = \, - \, k^{2} \, \sum_{m=1}^{\aleph} \sum_{\ell=1}^{2} \left[ \underset{y}{\nabla} \Phi_{k}(x,z_{m_{\ell}}) \times \tilde{Q}_{m_{\ell}} \, - \, \Upsilon_{k}(x,z_{m_{\ell}}) \cdot \tilde{R}_{m_{\ell}} \right] \, + \, \mathcal{O}\left( a^{7-2h-\frac{15}{2}t_2}\right), 
    \end{equation}
    where $\left(\tilde{\mathfrak{U}}_{m} \, := \, \left( \tilde{Q}_{m_1}, \tilde{R}_{m_1}, \tilde{Q}_{m_2}, \tilde{R}_{m_2} \right) \right)_{m=1}^{\aleph}$ is the solution of $(\ref{eq-al-D1})$.
    \item[]
    \item Derivation of the far field $E^{\infty}(\cdot,\theta, p, \omega)$. \\
    Similarly to the case of a single dimer, see \cite[Section 4]{cao-ghandriche-sini-dimer}, we can prove that the far field of the electromagnetic problem $(\ref{model-m})$ admits the following expansion,
	\begin{eqnarray}\label{far-expression-mid3}
	E^{\infty}(\hat{x},\theta, p) & = &  \left( I - \hat{x} \otimes \hat{x} \right) \cdot \sum_{m=1}^\aleph \sum_{i=1}^2 \eta_i  \,  e^{-i  k  \hat{x} \cdot z_{m_i}} \int_{D_{m_i}} \overset{3}{\mathbb{P}}\left(E_{m_i}\right)(y)  dy \notag\\
	&-& i  k  \left( I - \hat{x} \otimes \hat{x} \right) \cdot \sum_{m=1}^\aleph \sum_{i=1}^2 \eta_i \, e^{-i  k  \hat{x} \cdot z_{m_i}}\,\hat{x}\times\int_{D_{m_i}} F_{m_i}(y)  \,dy \, +\, Error_{0}, 
	\end{eqnarray}
	where
    \begin{align*}
Error_{0} \, :=\, &\left( I - \hat{x} \otimes \hat{x} \right) \cdot \sum_{m=1}^\aleph\, \sum_{i=1}^2 \, \eta_i 
\int_{D_{m_i}}\frac{k^2}{2}\bigl((y-z_{m_i})\cdot\hat{x}\bigr)^2
\int_{0}^1 (1-t)e^{- i k \hat{x}\cdot (z_{m_i}+t(y-z_{m_i}))}\,dt\; E_{m_i}(y)\,dy \\
&- i k\, \left( I - \hat{x} \otimes \hat{x} \right) \cdot \sum_{m=1}^\aleph \sum_{i=1}^2 \eta_i
\int_{D_{m_i}} e^{- i  k  \hat{x} \cdot z_{m_i}} \bigl((y-z_{m_i}) \cdot\hat{x}\bigr)\, \overset{3}{\mathbb{P}}\left(E_{m_i}\right)(y)\,dy.
\end{align*}
     Next, we estimate the above term. To do this, we take the absolute value on the both sides of $Error_0$ and use $(\ref{def-eta})$, to obtain
     \begin{eqnarray*}
        \left\vert Error_{0} \right\vert \, & \lesssim & \, \sum_{m=1}^\aleph \, \left\vert D_{m_1} \right\vert^{\frac{1}{2}} \,  \left\Vert E_{m_1}\right\Vert_{\mathbb{L}^2(D_{m_1})}\, + \, a^{2} \, \sum_{m=1}^\aleph \,  \left\vert D_{m_2} \right\vert^{\frac{1}{2}}\, \left\Vert E_{m_2}\right\Vert_{\mathbb{L}^2(D_{m_2})}\notag\\
        &+& a^{-1} \, \sum_{m=1}^{\aleph}  \left\vert D_{m_1} \right\vert^{\frac{1}{2}} \, \left\Vert \overset{3}{\PP}(E_{m_1})\right\Vert_{\mathbb{L}^2(D_{m_1})}\, +\, a \, \sum_{m=1}^{\aleph}  \left\vert D_{m_2} \right\vert^{\frac{1}{2}} \, \left\Vert \overset{3}{\PP}(E_{m_2})\right\Vert_{\mathbb{L}^2(D_{m_2})},
    \end{eqnarray*}
    which by scaling to $B_1, B_2$ gives us  
    \begin{eqnarray*}
       \left\vert Error_{0} \right\vert \, &\lesssim & \, a^{3} \, d_{\rm out}^{-3} \, \underset{1\leq m\leq \aleph}{\max} \left\lVert \tilde{E}_{m_1}\right\rVert_{\mathbb{L}^2(B_{1})}\, +\, a^{5} \, d_{\rm out}^{-3}\underset{1\leq m\leq\aleph}{\max}\left\lVert \tilde{E}_{m_2}\right\rVert_{\mathbb{L}^2(B_2)}   \\
        &+& \, a^{2} \, d_{\rm out}^{-3} \, \underset{1\leq m \leq \aleph}{\max}\left\lVert \overset{3}{\PP}\left(\tilde{E}_{m_1}\right)\right\rVert_{\mathbb{L}^2(B_1)}\, +\, a^{4}\, d_{\rm out}^{-3}\,\underset{1\leq m\leq \aleph}{\max}\left\lVert\overset{3}{\PP}\left(\tilde{E}_{m_2}\right)\right\rVert_{\mathbb{L}^2(B_2)}\\
        &\lesssim& \, a^{3} \, d_{\rm out}^{-3} \, \underset{1\leq m \leq \aleph}{\max}\left\lVert \overset{1}{\PP}\left(\tilde{E}_{m_1}\right)\right\rVert_{\mathbb{L}^2(B_1)}\, +\, a^{5}\, d_{\rm out}^{-3}\,\underset{1\leq m\leq \aleph}{\max}\left\lVert\overset{1}{\PP}\left(\tilde{E}_{m_2}\right)\right\rVert_{\mathbb{L}^2(B_2)}\\
        &+& \, a^{2} \, d_{\rm out}^{-3} \, \underset{1\leq m \leq \aleph}{\max}\left\lVert \overset{3}{\PP}\left(\tilde{E}_{m_1}\right)\right\rVert_{\mathbb{L}^2(B_1)}\, +\, a^{4}\, d_{\rm out}^{-3}\,\underset{1\leq m\leq \aleph}{\max}\left\lVert\overset{3}{\PP}\left(\tilde{E}_{m_2}\right)\right\rVert_{\mathbb{L}^2(B_2)}
     \end{eqnarray*}
     Thanks to {Proposition \ref{lem-es-multi}}, we deduce 
	\begin{equation}\label{FS-Eq-1202}
		Error_{0} \, = \, \mathcal{O}\left(a^{7-2h} \, d_{\rm out}^{-\frac{15}{2}} \right) \, \overset{(\ref{d-a})}{=} \, \mathcal{O}\left(a^{7-2h-\frac{15}{2}t_{2}} \right). 
	\end{equation} 
 By returning to $(\ref{far-expression-mid3})$ and using  the estimation $(\ref{FS-Eq-1202})$ we can derive 
 \begin{eqnarray}\label{FS-Equa1220}
 \nonumber
	E^{\infty}(\hat{x},\theta, p) & = &  \left( I - \hat{x} \otimes \hat{x} \right) \cdot \sum_{m=1}^\aleph \sum_{i=1}^2 \eta_i  \,  e^{- i  k  \hat{x} \cdot z_{m_i}} \int_{D_{m_i}} \overset{3}{\mathbb{P}}\left(E_{m_i}\right)(y)  dy \notag\\
	&-& i  k  \left( I - \hat{x} \otimes \hat{x} \right) \cdot \sum_{m=1}^\aleph \sum_{i=1}^2 \eta_i \, e^{- i  k  \hat{x} \cdot z_{m_i}}\,\hat{x}\times\int_{D_{m_i}} F_{m_i}(y)  \,dy \, +\, \mathcal{O}\left(a^{7-2h-\frac{15}{2}t_{2}} \right),
	\end{eqnarray}
    which with the notation \eqref{Qj1} further indicates that
    \begin{equation}\label{far-express-v1}
        E^{\infty}(\hat{x},\theta, p)=\sum_{m=1}^\aleph \sum_{i=1}^2 \left(e^{- i k \hat{x} \cdot z_{m_i}} (I-\hat{x}\otimes\hat{x})\tilde{R}_{m_i}- i k e^{- i k \hat{x}\cdot z_{m_i}} \hat{x}\times \tilde{Q}_{m_i}\right)+\O(a^{7-2h-\frac{15}{2}t_2}).
    \end{equation}
\item[]
    \item The gap between the perturbed solution and unperturbed solution related  to $(\ref{eq-al-D1})$. \\
    We recall that $\left(\tilde{\mathfrak{U}_{m}} \, := \, \left( \tilde{Q}_{m_1}, \tilde{R}_{m_1}, \tilde{Q}_{m_2}, \tilde{R}_{m_2} \right) \right)_{m=1}^{\aleph}$ is the solution to the following algebraic system  
    \begin{equation}\label{FA-Eq-perturbed}
            \begin{pmatrix}
                 \mathfrak{B_{1}} & - \, \Psi_{12} & \cdots & - \, \Psi_{1\aleph}  \\
                 - \, \Psi_{21} & \mathfrak{B_{2}} & \cdots & - \, \Psi_{2\aleph} \\
                 \vdots & \vdots & \ddots & \vdots \\
                 - \, \Psi_{\aleph1} & - \, \Psi_{\aleph2} & \cdots &  \mathfrak{B_{\aleph}} 
            \end{pmatrix} 
            \cdot
            \begin{pmatrix}
                \tilde{\mathfrak{U}}_{1} \\
                \tilde{\mathfrak{U}}_{2} \\
                \vdots \\
                \tilde{\mathfrak{U}}_{\aleph} 
            \end{pmatrix}
            = 
            \begin{pmatrix}
                \mathrm{S}_{1} \\
                \mathrm{S}_{2} \\
                \vdots \\
                \mathrm{S}_{\aleph}  
            \end{pmatrix} \, + \,
     \begin{pmatrix}
                \mathfrak{R}_{1} \\
                \mathfrak{R}_{2} \\
                \vdots \\
                \mathfrak{R}_{\aleph}  
            \end{pmatrix},
        \end{equation}
         and we denote by $\left({\mathfrak{U}}_{m} \, := \, \left( {Q}_{m_1}, {R}_{m_1}, {Q}_{m_2}, {R}_{m_2} \right) \right)_{m=1}^{\aleph}$ the solution to the following unperturbed algebraic system 
        \begin{equation}\label{FA-Eq-unperturbed}
            \begin{pmatrix}
                 \mathfrak{B_{1}} & - \, \Psi_{12} & \cdots & - \, \Psi_{1\aleph}  \\
                 - \, \Psi_{21} & \mathfrak{B_{2}} & \cdots & - \, \Psi_{2\aleph} \\
                 \vdots & \vdots & \ddots & \vdots \\
                 - \, \Psi_{\aleph1} & - \, \Psi_{\aleph2} & \cdots &  \mathfrak{B_{\aleph}} 
            \end{pmatrix} 
            \cdot
            \begin{pmatrix}
                {\mathfrak{U}}_{1} \\
                {\mathfrak{U}}_{2} \\
                \vdots \\
                {\mathfrak{U}}_{\aleph} 
            \end{pmatrix}
            = 
            \begin{pmatrix}
                \mathrm{S}_{1} \\
                \mathrm{S}_{2} \\
                \vdots \\
                \mathrm{S}_{\aleph}  
            \end{pmatrix},
        \end{equation}  
        where $\left\{ \left\{\mathfrak{B_{m}} \right\}_{m=1}^{\aleph}; \left\{ \Psi_{mj}\right\}_{m,j=1}^{\aleph}; \left\{ \mathrm{S}_{m} \right\}_{m=1}^{\aleph} \right\}$ are given in $\textbf{Notations \ref{notbthm}}$. Hence, by taking the difference $(\ref{FA-Eq-perturbed})$ and $(\ref{FA-Eq-unperturbed})$, we deduce that 
        \begin{equation}\label{FA-Eq-unperturbed-perturbed}
            \begin{pmatrix}
                 \mathfrak{B_{1}} & - \, \Psi_{12} & \cdots & - \, \Psi_{1\aleph}  \\
                 - \, \Psi_{21} & \mathfrak{B_{2}} & \cdots & - \, \Psi_{2\aleph} \\
                 \vdots & \vdots & \ddots & \vdots \\
                 - \, \Psi_{\aleph1} & - \, \Psi_{\aleph2} & \cdots &  \mathfrak{B_{\aleph}} 
            \end{pmatrix} 
            \cdot
            \begin{pmatrix}
                \left( \tilde{\mathfrak{U}}_{1} \, - \, {\mathfrak{U}}_{1} \right) \\
                \left( \tilde{\mathfrak{U}}_{2} \, - \,  {\mathfrak{U}}_{2} \right) \\
                \vdots \\
                \left( \tilde{\mathfrak{U}}_{\aleph} \, - \,  {\mathfrak{U}}_{\aleph} \right) 
            \end{pmatrix}
            = 
            \begin{pmatrix}
                \mathfrak{R}_{1} \\
                \mathfrak{R}_{2} \\
                \vdots \\
                \mathfrak{R}_{\aleph}  
            \end{pmatrix},
        \end{equation}
        which is invertible under the condition \eqref{condi-invertibility}.
        Now, component by component, by taking the Euclidean norm on both sides of the above equation, we obtain  
    \item[]
    \item The corresponding revised Foldy-Lax approximation.

    From \eqref{FA-Eq-unperturbed-perturbed}, we can directly obtain from the invertibility of the system that
    \begin{eqnarray}\label{far-error-diff}
        \left( \sum_{m=1}^\aleph |\mathfrak{U}_m-\tilde{\mathfrak{U}}_m|^2 \right)^{\frac{1}{2}} & \lesssim & \left(\sum_{m=1}^\aleph |\mathfrak{R}_m|^2\right)^\frac{1}{2} \lesssim \aleph^\frac{1}{2} \cdot |\mathfrak{R}_m|\notag\\
        &\lesssim& d_{\rm out}^{-\frac{3}{2}} \left( a^{10-3h-\frac{17}{2}t_2} + a^{10-2h-4t_1-\frac{9}{2}t_2} + a^{12-2h-4t_1-\frac{9}{2}t_2} \right)\notag\\
        &\lesssim& a^{10-3h-10t_2} + a^{10-2h-4t_1-6t_2},
    \end{eqnarray}
    by using \eqref{err-prop-la1} for the explicit expression of $\mathfrak{R}_m$.

    Now, we first recall the scattered wave expansion given by \eqref{scattered field-foldy lax}, which can be written equivalently as the matrix form
    \begin{equation}\label{scattered field-matrix form}
        E^s(x, \theta, p, \omega)= -k^2 \sum_{m=1}^\aleph \mathcal{G}_m\cdot \tilde{\mathfrak{U}}_m + \O(a^{7-2h-\frac{15}{2}t_2})\quad\mbox{with}\quad \mathcal{G}_m:=\begin{pmatrix}
            \underset{y}{\nabla}\Phi_k(x, z_{m_1})\times\\
            -\Upsilon_k(x, z_{m_1})\\
            \underset{y}{\nabla}\Phi_k(x, z_{m_2})\\
            -\Upsilon_k(x, z_{m_2})
        \end{pmatrix}^\mathsf{T}
    \end{equation}
%     the far-field expression \eqref{far-express-v1} by $\mathfrak{U}_m$ as 
%     \begin{equation}\label{far-express-vector}
%         E^\infty (\hat{x}, \theta, p, \omega)=\sum_{m=1}^\aleph \mathfrak{e}_m\cdot \tilde{\mathfrak{U}}_m + \O(a^{7-2h-\frac{15}{2}t_2})\quad\mbox{with}\quad \mathfrak{e}_m := \left(
%      \begin{array}{c}
% 	-ik e^{ik \hat{x}\cdot z_{m_1}} \hat{x}\times\\
% 	e^{i k \hat{x}\cdot z_{m_1}}(I-\hat{x}\otimes\hat{x}) \\
% 	-ik e^{ik \hat{x}\cdot z_{m_2}} \hat{x}\times\\
% 	e^{i k \hat{x}\cdot z_{m_2}}(I-\hat{x}\otimes\hat{x})
% \end{array}\right)^\mathsf{T}.
%     \end{equation}
    Then, combining with \eqref{far-error-diff} and the fact that $x$ is away from $D$, we can further obtain from \eqref{scattered field-matrix form} that
    \begin{equation}\label{far-mid}
        E^s ({x}, \theta, p, \omega)= - k^2 \sum_{m=1}^\aleph \mathcal{G}_m\cdot \mathfrak{U}_m - k^2 \sum_{m=1}^\aleph \mathcal{G}_m\cdot (\tilde{\mathfrak{U}}_m-\mathfrak{U}_m)+\O(a^{7-2h-\frac{15}{2}t_2}),
    \end{equation}
    where
    \begin{equation*}
        \left|\sum_{m=1}^\aleph \mathcal{G}_m\cdot (\tilde{\mathfrak{U}}_m-\mathfrak{U}_m)\right|\lesssim \left(\sum_{m=1}^\aleph |\mathcal{G}_m|^2\right)^\frac{1}{2}\cdot\left(\sum_{m=1}^\aleph |\tilde{\mathfrak{U}_m}-\mathfrak{U}_m|^2\right)^\frac{1}{2}\lesssim a^{10-3h-\frac{23}{2}t_2}+a^{10-2h-4t_1-\frac{15}{2}t_2}.
    \end{equation*}
    Therefore, we can derive the following Foldy-Lax approximation form with the dominant terms as
%     \begin{align*}\label{approximation-E-Foldy}
%     E^{\infty}(\hat{x}, \theta, p, \omega) \, 
%     % &=\sum_{m=1}^\aleph \sum_{i=1}^2 \left[e^{i k \hat{x}\cdot z_{m_i}}(I-\hat{x}\otimes\hat{x})R_{m_i}- i k e^{i k \hat{x}\cdot z_{m_i}} \hat{x}\times Q_{m_i}\right] + \O\left(a^{10-3h-\frac{23}{2}t_2}\right) \notag\\
%     &= \sum_{m=1}^\aleph \mathfrak{e}_m\cdot \mathfrak{U}_{m}+\O\left(a^{10-3h-\frac{23}{2}t_2}\right),
% \end{align*}
\begin{equation}\label{approximation-E-Foldy}
    E^{s}({x}, \theta, p, \omega) \, 
    = - k^2 \sum_{m=1}^\aleph \mathcal{G}_m\cdot \mathfrak{U}_{m}+\O\left(a^{10-3h-\frac{23}{2}t_2}\right),
\end{equation}
with $\mathcal{G}_m$ given by \eqref{scattered field-matrix form}. Similarly, by recalling the far field expression given by \eqref{far-express-v1} and using the estimates derived in \eqref{far-error-diff}, we can obtain
    \begin{align*}\label{approximation-E-Foldy}
    E^{\infty}(\hat{x}, \theta, p) \, 
    % &=\sum_{m=1}^\aleph \sum_{i=1}^2 \left[e^{i k \hat{x}\cdot z_{m_i}}(I-\hat{x}\otimes\hat{x})R_{m_i}- i k e^{i k \hat{x}\cdot z_{m_i}} \hat{x}\times Q_{m_i}\right] + \O\left(a^{10-3h-\frac{23}{2}t_2}\right) \notag\\
    &= \sum_{m=1}^\aleph \mathfrak{e}_m\cdot \mathfrak{U}_{m}+\O\left(a^{10-3h-\frac{23}{2}t_2}\right)\quad\mbox{with}\quad  \mathfrak{e}_m := \left(
     \begin{array}{c}
	-ik e^{- ik \hat{x}\cdot z_{m_1}} \hat{x}\times\\
	e^{- i k \hat{x}\cdot z_{m_1}}(I-\hat{x}\otimes\hat{x}) \\
	-ik e^{- ik \hat{x}\cdot z_{m_2}} \hat{x}\times\\
	e^{- i k \hat{x}\cdot z_{m_2}}(I-\hat{x}\otimes\hat{x})
\end{array}\right)^\mathsf{T},
\end{align*}
where $(\mathfrak{U}_m)_{m=1}^\aleph$ is the solution to the linear algebraic system \eqref{eq-al-D1-}. 
\end{enumerate}    
%%%%%%%%%%%%%%%%%%%%%%%%%%%%%%%%%%%%%%%%%%%%%%%%%%%%%%%%%%%%
%%%%%%%%%%%%%%%%%%%%%%%%%%%%%%%%%%%%%%%%%%%%%%%%%%%%%%%%%%%%

\subsection{Proof of Corollary \ref{corollary-main}}\label{sec:appendix-coro}
Recall the linear algebraic system \eqref{eq-al-D1-}. For any $1\leq m\leq \aleph$, there holds
\begin{equation}\label{system-coro-original}
    \mathfrak{B}_{m} \mathfrak{U}_m - \sum_{j=1 \atop j\neq m}^\aleph  \Psi_{mj}\mathfrak{U}_j  = \begin{pmatrix}
     \frac{i \, k \, \eta_{0}}{\pm \, c_{0}} \, a^{3-h} \, {\bf P}_{0, 1}^{(1)} \cdot H^{Inc}(z_{m_1}) \\
     a^{3} \, {\bf P}_{0, 1}^{(2)} \cdot E^{Inc}(z_{m_1}) \\
     i \, k \, a^{5} \, {\bf P}_{0, 2}^{(1)} \cdot H^{Inc}(z_{m_2}) \\
     \frac{\eta_{2}}{\pm \, d_{0}} \, a^{3-h} \, {\bf P}_{0, 2}^{(2)} \cdot E^{Inc}(z_{m_2})
    \end{pmatrix},
\end{equation}
with the R.H.S. reformulated as
\begin{equation}\notag
{\small
\begin{aligned}
\begin{pmatrix}
 \frac{i \, k \, \eta_{0}}{\pm \, c_{0}} \, a^{3-h} \, {\bf P}_{0, 1}^{(1)} \cdot H^{Inc}(z_{m_1}) \\
 a^{3} \, {\bf P}_{0, 1}^{(2)} \cdot E^{Inc}(z_{m_1}) \\
 i \, k \, a^{5} \, {\bf P}_{0, 2}^{(1)} \cdot H^{Inc}(z_{m_2}) \\
 \frac{\eta_{2}}{\pm \, d_{0}} \, a^{3-h} \, {\bf P}_{0, 2}^{(2)} \cdot E^{Inc}(z_{m_2})
\end{pmatrix}
&= {\bf T}_a\,\begin{pmatrix}
 H^{Inc}(z_{m_1})\\
 a^h E^{Inc}(z_{m_1})\\
 a^{2+h} H^{Inc}(z_{m_2})\\
 E^{Inc}(z_{m_2})
\end{pmatrix},\\[2pt]
{\bf T}_a &:= \begin{pmatrix}
 \frac{i k \eta_0}{\pm c_0} a^{3-h} {\bf P}_{0, 1}^{(1)} & 0 & 0 & 0\\
 0 & a^{3-h} {\bf P}_{0, 1}^{(2)} & 0 & 0\\
 0 & 0 & i k a^{3-h} {\bf P}_{0, 2}^{(1)} & 0\\
 0 & 0 & 0 & \frac{\eta_2}{\pm d_0}a^{3-h} {\bf P}_{0, 2}^{(2)}
\end{pmatrix}.
\end{aligned}
}
\end{equation}
By multiplying ${\bf T}_a^{-1}$ on the both sides of \eqref{system-coro-original}, we can get 
\begin{equation}\label{system-coro-reform1}
    \tilde{\mathfrak{B}}_m^a \mathfrak{U}_m - \sum_{j=1 \atop j\neq m}^\aleph \tilde{\Psi}_{mj}^a \mathfrak{U}_j = \begin{pmatrix}
        H^{Inc}(z_{m_1})\\
        a^h E^{Inc}(z_{m_1})\\
        a^{2+h} H^{Inc}(z_{m_2})\\
        E^{Inc}(z_{m_2})
    \end{pmatrix},
\end{equation}
where 
\begin{equation}\notag
    \tilde{\mathfrak{B}}_m^a:={\bf T}_a^{-1} \mathfrak{B}_m \quad\mbox{and}\quad \tilde{\Psi}_{mj}^a:= {\bf T}_a^{-1} \Psi_{mj}.
\end{equation}
Equivalently, \eqref{system-coro-reform1} indicates the following system of equations w.r.t. $(Q_{m_1}, R_{m_1}, Q_{m_2}, R_{m_2})$:
\begin{equation}\label{eq-system-1}
\begin{split}
&a^{h-3} \frac{\pm c_0}{i k \eta_0}\bigl({\bf P}_{0, 1}^{(1)}\bigr)^{-1} \Big(Q_{m_1} - \mathcal{B}_{13} Q_{m_2}- \mathcal{B}_{14} R_{m_2}\\
\qquad -& \sum_{j=1 \atop j\neq m}^\aleph \left(\mathcal{C}_{11} Q_{j_1} +\mathcal{C}_{12}R_{j_1} + \mathcal{C}_{13} Q_{j_2} + \mathcal{C}_{14}R_{j_2} \right)  \Big) = H^{Inc}(z_{m_1}),
\end{split}
\end{equation}
\begin{equation}\label{eq-system-2}
\begin{split}
&a^{h-3} \bigl({\bf P}_{0, 1}^{(2)}\bigr)^{-1} \Big(R_{m_1} - \mathcal{B}_{23} Q_{m_2}- \mathcal{B}_{24} R_{m_2}\\
\qquad -& \sum_{j=1 \atop j\neq m}^\aleph \left(\mathcal{C}_{21} Q_{j_1} +\mathcal{C}_{22}R_{j_1} + \mathcal{C}_{23} Q_{j_2} + \mathcal{C}_{24}R_{j_2} \right)  \Big) = a^h E^{Inc}(z_{m_1}),
\end{split}
\end{equation}
\begin{equation}\label{eq-system-3}
\begin{split}
&a^{h-3} \frac{i}{k} \bigl({\bf P}_{0, 2}^{(1)}\bigr)^{-1} \Big( \mathcal{B}_{31} Q_{m_1} + \mathcal{B}_{32} R_{m_1} - Q_{m_2}\\
\qquad +& \sum_{j=1 \atop j\neq m}^\aleph  \left(\mathcal{C}_{31} Q_{j_1} +\mathcal{C}_{32}R_{j_1} + \mathcal{C}_{33} Q_{j_2} + \mathcal{C}_{34}R_{j_2} \right)  \Big) = a^{h+2} H^{Inc}(z_{m_2}),
\end{split}
\end{equation}
\begin{equation}\label{eq-system-4}
\begin{split}
&a^{h-3} \frac{\mp d_0}{\eta_2} \bigl({\bf P}_{0, 2}^{(2)}\bigr)^{-1} \Big( \mathcal{B}_{41} Q_{m_1} + \mathcal{B}_{42} R_{m_1} - R_{m_2}\\
\qquad +& \sum_{j=1 \atop j\neq m}^\aleph  \left(\mathcal{C}_{41} Q_{j_1} +\mathcal{C}_{42}R_{j_1} + \mathcal{C}_{43} Q_{j_2} + \mathcal{C}_{44}R_{j_2} \right)  \Big) = E^{Inc}(z_{m_2}).
\end{split}
\end{equation}
By adding \eqref{eq-system-1} to \eqref{eq-system-3} and simplifying the corresponding expressions with \textbf{Notation \ref{notbthm}}, we can obtain that
\begin{eqnarray}\label{eq-sys1-modi}
    &&a^{h-3} \frac{\pm c_0}{i k \eta_0}\bigl({\bf P}_{0, 1}^{(1)}\bigr)^{-1} Q_{m_1} + \frac{1}{i k} a^{h-3} \bigl({\bf P}_{0, 2}^{(1)}\bigr)^{-1} Q_{m_2}\notag\\
    &+& i k^3 \Upsilon_k(z_{m_1}, z_{m_2})\left( \eta_2 a^{h+2} Q_{m_1}+Q_{m_2} \right) + i k \nabla\Phi_k(z_{m_1}, z_{m_2})\times \left(k^2 \eta_2 a^{h+2} R_{m_1} +R_{m_2} \right)\notag\\
    &+& i k \sum_{j=1 \atop j\neq m}^\aleph \Bigg[ k^2\left( \Upsilon_k(z_{m_1}, z_{j_1}) + \eta_2 a^{h+2} \Upsilon_k(z_{m_2}, z_{j_1})\right) Q_{j_1} + k^2 \left( \Upsilon_k(z_{m_1}, z_{j_2}) + \eta_2 a^{h+2} \Upsilon_k(z_{m_2}, z_{j_2}) \right) Q_{j_2}\notag\\
    &+& \left( {\nabla}\Phi_k(z_{m_1}, z_{j_1}) + \eta_2 a^{h+2} {\nabla}\Phi_k(z_{m_2}, z_{j_1})\right) \times R_{j_1} + \left( \nabla\Phi_k(z_{m_1}, z_{j_2}) + \eta_2 a^{h+2} \nabla\Phi_k(z_{m_2}, z_{j_2})\right)\times
    R_{j_2}  \Bigg]\notag\\
    &=& H^{Inc}(z_{m_1}) + a^{h+2} H^{Inc}(z_{m_2}),\notag\\
\end{eqnarray}
% It is obvious to see that
% \begin{eqnarray}\label{Error-coro-1}
%     |Error_{coro, 1}|&\lesssim& a^{t_1} \sum_{j=1 \atop j\neq m}^\aleph \left( d_{mj}^{-4} Q_{j_1}+d_{mj}^{-4} Q_{j_2} + d_{mj}^{-3} R_{j_1} +d_{mj}^{-3} R_{j_2} \right) + a^{t_1}| \nabla H^{Inc}(z_{m_0}+t(z_{m_1}-z_{m_0}))|\notag\\
%     &\lesssim& a^{t_1} \left( \sum_{j=1\atop j\neq m}^\aleph d_{mj}^{-8} \right)^\frac{1}{2} \sum_{\ell=1}^2 \left( \sum_{j=1 \atop j\neq m}^\aleph |Q_{j_\ell}|^2 \right)^\frac{1}{2} + a^{t_1} \left( \sum_{j=1\atop j\neq m}^\aleph d_{mj}^{-6} \right)^\frac{1}{2} \sum_{\ell=1}^2 \left( \sum_{j=1 \atop j\neq m}^\aleph |R_{j_\ell}|^2 \right)^\frac{1}{2}.
% \end{eqnarray}
% By virtue of Proposition SM-1.1, with the a-prior estimates derived in Proposition \ref{lem-es-multi}, we have 
% \begin{align}\label{es-QR12}
%      &\left(\sum_{m=1}^\aleph |Q_{m_1}|^2 \right)^\frac{1}{2}\lesssim a^{3-h} d_{\rm out}^{-\frac{3}{2}}, \qquad \; \; \left( \sum_{m=1}^\aleph |R_{m_1}|^2 \right)^\frac{1}{2}\lesssim a^{6-h} d_{\rm in}^{-3} d_{\rm out}^{-\frac{3}{2}},\notag\\
%     &\left(\sum_{m=1}^\aleph |Q_{m_2}|^2 \right)^\frac{1}{2}\lesssim a^{8-h} d_{\rm in}^{-3} d_{\rm out}^{-\frac{3}{2}}, \quad \left( \sum_{m=1}^\aleph |R_{m_2}|^2 \right)^\frac{1}{2}\lesssim a^{3-h} d_{\rm out}^{-\frac{3}{2}}.
% \end{align}
% Hence, taking \eqref{es-QR12} into \eqref{Error-coro-1},  we can derive the estimate by keeping the dominant term that
% \begin{equation}\label{es-Error-coro-1-final}
%     |Error_{coro, 1}|\lesssim a^{3-h+t_1-\frac{11}{2}t_2},
% \end{equation}
which further indicates by keeping the dominant terms that 
\begin{eqnarray}\label{eq-sys1-modi2}
    &&a^{h-3} \frac{\pm c_0}{i k \eta_0}\bigl({\bf P}_{0, 1}^{(1)}\bigr)^{-1} Q_{m_1} + i k \nabla\Phi_k(z_{m_1}, z_{m_2})\times R_{m_2} \notag\\
    &+& i k \sum_{j=1 \atop j\neq m}^\aleph \left( k^2 \Upsilon_k(z_{m_1}, z_{j_1}) Q_{j_1} + \nabla\Phi_k(z_{m_1}, z_{j_2})\times R_{j_2} \right) = H^{Inc}(z_{m_1}) + Error_{coro},
\end{eqnarray}
where 
\begin{eqnarray}
    Error_{coro}&:=&-  \left(\frac{1}{i k} a^{h-3} \bigl({\bf P}_{0, 2}^{(1)}\bigr)^{-1} + i k^3 \Upsilon_k(z_{m_1}, z_{m_2})\right) Q_{m_2} -  \sum_{j=1 \atop j\neq m}^\aleph i k^3 \Upsilon_k(z_{m_1}, z_{j_2}) Q_{j_2}  \notag\\
    &-& \sum_{j=1 \atop j\neq m}^\aleph i k \nabla\Phi_k(z_{m_1}, z_{j_1})\times R_{j_1} + a^{h+2} \Bigg[ H^{Inc}(z_{m_2}) - i k^3 \eta_2  \Big( \Upsilon_k(z_{m_1}, z_{m_2}) Q_{m_1} \notag\\
    &+& \nabla\Phi_k(z_{m_1}, z_{m_2})\times R_{m_1}\Big) - i k^3 \eta_2 \sum_{j=1\atop j\neq m}^\aleph \Big(\Upsilon_k(z_{m_2}, z_{j_1}) Q_{j_1} + \Upsilon_k(z_{m_2}, z_{j_2}) Q_{j_2}  \Big)\notag\\
    &-& i k \eta_2  \sum_{j=1\atop j\neq m}^\aleph \Big(\nabla\Phi_k(z_{m_2}, z_{j_1}) \times R_{j_1} + \nabla\Phi_k(z_{m_2}, z_{j_2}) \times R_{j_2} \Big)\Bigg].
\end{eqnarray}
% Now, we evaluate the R.H.S. of \eqref{eq-sys1-modi2}. By denoting 
% \begin{eqnarray*}
%     Error_{right}&:=& \frac{1}{i k} a^{h-3} \left[{\bf P}_{0, 2}^{(1)}\right]^{-1} Q_{m_2} + i k^3 \Upsilon_k(z_{m_1}, z_{m_2}) Q_{m_2} + i k^3 \eta_2 a^{2+h} \nabla\Phi_k(z_{m_1}, z_{m_2})\times R_{m_1}\notag\\
%     &+& \sum_{j=1 \atop j\neq m}^\aleph i k^3 \Upsilon_k(z_{m_0}, z_{j_0}) Q_{j_2}  + \sum_{j=1 \atop j\neq m}^\aleph i k \nabla\Phi_k(z_{m_0}, z_{j_0})\times R_{j_1},
% \end{eqnarray*}
% with direct calculation, it is straightforward to get from \eqref{es-QR12} that
% \begin{eqnarray*}
%     |Error_{right}| &\lesssim& \left( a^{h-3} \aleph^\frac{1}{2} + d_{\rm in}^{-3} \aleph^\frac{1}{2} + \left( \sum_{j=1 \atop j\neq m}^\aleph d_{mj}^{-6} \right)^\frac{1}{2} \right) \left(\sum_{m=1}^\aleph |Q_{m_2}|^2 \right)^\frac{1}{2} \notag\\
%     &+& \left( a^{2+h} d_{\rm in}^{-2} \aleph^\frac{1}{2} + \left( \sum_{j=1\atop j\neq m}^\aleph d_{mj}^{-4} \right)^\frac{1}{2} \right)\left(\sum_{m=1}^\aleph |R_{m_1}|^2 \right)^\frac{1}{2} \lesssim   a^{\min\{5-3t_1-3t_2, 8-h-6t_1-3t_2, 6-h-3t_1-\frac{7}{2}t_2\}}.
% \end{eqnarray*}
By virtue of Proposition SM-1.1, with the a-prior estimates derived in Proposition \ref{lem-es-multi}, we have 
\begin{align}\label{es-QR12}
     &\left(\sum_{m=1}^\aleph |Q_{m_1}|^2 \right)^\frac{1}{2}\lesssim a^{3-h} d_{\rm out}^{-\frac{3}{2}}, \qquad \; \; \left( \sum_{m=1}^\aleph |R_{m_1}|^2 \right)^\frac{1}{2}\lesssim a^{6-h} d_{\rm in}^{-3} d_{\rm out}^{-\frac{3}{2}},\notag\\
    &\left(\sum_{m=1}^\aleph |Q_{m_2}|^2 \right)^\frac{1}{2}\lesssim a^{8-h} d_{\rm in}^{-3} d_{\rm out}^{-\frac{3}{2}}, \quad \left( \sum_{m=1}^\aleph |R_{m_2}|^2 \right)^\frac{1}{2}\lesssim a^{3-h} d_{\rm out}^{-\frac{3}{2}}.
\end{align}
Combining with \eqref{es-QR12} and Proposition \ref{lem-es-multi}, utilizing the notation \eqref{Qj1}, we can obtain that
\begin{eqnarray*}
    |Error_{coro}| &\lesssim& a^{h+2} + \left(a^{h-3} + d_{\rm in}^{-3} \right)\max_m |Q_{m_2}| + \left( \sum_{j=1 \atop j\neq m}^\aleph d_{mj}^{-6} \right)^\frac{1}{2}\left(\sum_{m=1}^\aleph |Q_{m_2}|^2 \right)^\frac{1}{2} \notag\\
    &+& \left( \sum_{j=1 \atop j\neq m}^\aleph d_{mj}^{-4} \right)^\frac{1}{2}\left(\sum_{m=1}^\aleph |R_{m_1}|^2 \right)^\frac{1}{2} + a^{h+2} \left( d_{\rm in}^{-3}\max_{m} |Q_{m_1}| + d_{\rm in}^{-2} \max_m |R_{m_1}| \right)\notag\\
    &+& a^{h+2}  \left( \sum_{j=1 \atop j\neq m}^\aleph d_{mj}^{-6} \right)^\frac{1}{2} \left[\left(\sum_{m=1}^\aleph |Q_{m_1}|^2 \right)^\frac{1}{2} + \left(\sum_{m=1}^\aleph |Q_{m_2}|^2 \right)^\frac{1}{2} \right]\notag\\
    &+& a^{h+2}  \left( \sum_{j=1 \atop j\neq m}^\aleph d_{mj}^{-4} \right)^\frac{1}{2} \left[\left(\sum_{m=1}^\aleph |R_{m_1}|^2 \right)^\frac{1}{2} + \left(\sum_{m=1}^\aleph |R_{m_2}|^2 \right)^\frac{1}{2} \right]\notag\\
    &\lesssim& a^{h+2} + (a^{h-3} + d_{\rm in}^{-3}) a^4 \max_m \left\lVert\overset{1}{\PP}(\tilde{E}_{m_2})\right\rVert_{\LL^2(B_2)} + a^{6-h} d_{\rm in}^{-3} d_{\rm out}^{-\frac{7}{2}} \notag\\
    &+& a^{h+2} \left( a^2 d_{\rm in}^{-3} \left\lVert\overset{1}{\PP}(\tilde{E}_{m_1})\right\rVert_{\LL^2(B_1)} + a^3 d_{\rm in}^{-2} \max_m \left\lVert\overset{3}{\PP}(\tilde{E}_{m_1})\right\rVert_{\LL^2(B_1)}  \right) + a^5 d_{\rm out}^{-\frac{9}{2}}\notag\\
    &\lesssim& a^{6-h} d_{\rm in}^{-3} d_{\rm out}^{-\frac{7}{2}},
\end{eqnarray*}
by keeping the dominant term.
Hence, under the condition \eqref{conditions-t-h}, we can further simplify \eqref{eq-sys1-modi2} as
% \begin{eqnarray}\label{eq-system1-error final}
%     &&a^{h-3} \frac{\pm c_0}{i k \eta_0}\left[{\bf P}_{0, 1}^{(1)}\right]^{-1} Q_{m_1} + i k \nabla\Phi_k(z_{m_1}, z_{m_2})\times R_{m_2}\notag\\ 
%     &+& \sum_{j=1 \atop j\neq m}^\aleph i k \left( k^2 \Upsilon_k(z_{m_0}, z_{j_0}) Q_{j_1} + \nabla\Phi_k(z_{m_0}, z_{j_0})\times R_{j_2} \right)
%     = H^{Inc}(z_{m_0}) + \O(a^{4-h-\frac{11}{2}t_2}).
% \end{eqnarray}
\begin{eqnarray}\label{eq-system1-error final}
    &&a^{h-3} \frac{\pm c_0}{i k \eta_0}\bigl({\bf P}_{0, 1}^{(1)}\bigr)^{-1} Q_{m_1} + i k \nabla\Phi_k(z_{m_1}, z_{m_2})\times R_{m_2} \notag\\
    &+& i k \sum_{j=1 \atop j\neq m}^\aleph \left( k^2 \Upsilon_k(z_{m_1}, z_{j_1}) Q_{j_1} + \nabla\Phi_k(z_{m_1}, z_{j_2})\times R_{j_2} \right) = H^{Inc}(z_{m_1}) + \O(a^{6-h-3t_1-\frac{7}{2}t_2}).
\end{eqnarray}

Following the similar argument above, by adding \eqref{eq-system-2} to \eqref{eq-system-4}, with \textbf{Notation} \ref{notbthm}, we can obtain the expression related to the incident electric field under condition \eqref{conditions-t-h} that
\begin{eqnarray}\label{eq-system2-error-final}
    && - k^2 \nabla\Phi_k(z_{m_1}, z_{m_2}) \times Q_{m_1} +\frac{\pm d_0}{\eta_2} a^{h-3} \bigl({\bf P}_{0, 2}^{(2)}\bigr)^{-1} R_{m_2}\notag\\
    &-& \sum_{j=1 \atop j\neq m}^\aleph k^2 \left( \nabla\Phi_k(z_{m_2}, z_{j_1}) \times Q_{j_1} +\Upsilon_k(z_{m_2}, z_{j_2})R_{j_2} \right) = E^{Inc}(z_{m_2}) + \O(a^{6-h-3t_1-\frac{7}{2}t_2}).
\end{eqnarray}
From \eqref{eq-system1-error final} and \eqref{eq-system2-error-final}, suppose $(\mathring{Q}_{m_1}, \mathring{R}_{m_2})_{m=1}^\aleph$ is the solution to the following system of the matrix form 
\begin{eqnarray}\label{system-mathrix-error}
    &&\begin{pmatrix}
        a^{h-3} \frac{\pm c_0}{i k \eta_0} \bigl({\bf P}_{0, 1}^{-1} \bigr)^{-1} & i k \nabla\Phi_k(z_{m_1}, z_{m_2})\times\\
        - k^2 \nabla \Phi_k(z_{m_1}, z_{m_2})\times & a^{h-3} \frac{\pm d_0}{\eta_2} \bigl( {\bf P}_{0, 2}^{(2)} \bigr)^{-1}
    \end{pmatrix}\begin{pmatrix}
        \mathring{Q}_{m_1}\\
        \mathring{R}_{m_2}
    \end{pmatrix}\notag\\
    &-& \sum_{j=1 \atop j\neq m}^\aleph 
    \begin{pmatrix}
        - i k^3\Upsilon_k(z_{m_1}, z_{j_1}) & - i k \nabla\Phi_k(z_{m_1}, z_{j_2})\times\\
        k^2 \nabla\Phi_k(z_{m_2}, z_{j_1})\times & k^2 \Upsilon_k(z_{m_2}, z_{j_2})
    \end{pmatrix}\begin{pmatrix}
            \mathring{Q}_{j_1}\\
            \mathring{R}_{j_2}
        \end{pmatrix} = \begin{pmatrix}
        H^{Inc}(z_{m_1})\\
        E^{Inc}(z_{m_2})
    \end{pmatrix} .
    % + \begin{pmatrix}
    %     a^{4-h-\frac{11}{2}t_2}\\
    %      a^{4-h-\frac{11}{2}t_2}.
    % \end{pmatrix}
\end{eqnarray}

Recall the scattered field expansion with the form \eqref{scattered field-equivalent form} for $x\in \mathbb{R}^3\backslash D$. Taking Taylor expansion for the functions $\underset{y}{\nabla}\Phi_k(x, \cdot)$ and $\Upsilon_k(x, \cdot)$ near the intermediate point $z_{m_0}$, which is given by \eqref{def-intermediate point}, namely
\begin{equation*}
    \underset{y}{\nabla}\Phi_k(x, z_{m_\ell})=\underset{y}{\nabla}\Phi_k(x, z_{m_0}) + \O(a^{t_1})\quad\mbox{and}\quad \Upsilon_k(x, z_{m_\ell})=\Upsilon_k(x, z_{m_0})+\O(a^{t_1}),
\end{equation*}
% the far-field approximation \eqref{foldy-lax-final}, which can be equivalently written as
% \begin{equation}\label{far-field-original}
%     E^\infty(\hat{x}, \theta, p)= \sum_{m=1}^\aleph \sum_{i=1}^2 \left[ e^{- i k \hat{x}\cdot z_{m_i}} (I-\hat{x}\otimes \hat{x}) R_{m_i} - i k e^{- i k \hat{x}\cdot z_{m_i}} \hat{x}\times Q_{m_i} \right] +\O(a^{10-3h - \frac{23}{2}t_2}).
% \end{equation}
then by keeping the dominant terms, \eqref{scattered field-equivalent form} becomes
\begin{equation}\label{scattered-f1}
     E^s ({x}, \theta, p, \omega) = - k^2 \sum_{m=1}^\aleph \sum_{\ell=1}^2 \left[ \underset{y}{\nabla}\Phi_k(x, z_{m_0}) \times {Q}_{m_\ell} - \Upsilon_k(x, z_{m_0})\cdot {R}_{m_\ell} \right] + \O(a^{10-3h - \frac{23}{2}t_2}).
\end{equation}
Besides, due to the fact that
\begin{eqnarray}
    &&\left| k^2 \sum_{m=1}^\aleph \left( \underset{y}{\nabla}\Phi_k(x, z_{m_0})\times {Q}_{m_2} -\Upsilon_k(x, z_{m_0})\cdot {R}_{m_1} \right) \right|\notag\\
    &=& \left| k^2 \sum_{m=1}^\aleph \left( \underset{y}{\nabla}\Phi_k(x, z_{m_0})\times \eta_2 \int_{D_{m_2}} F_{m_2}(y)\,dy -\Upsilon_k(x, z_{m_0})\cdot \eta_1 \int_{D_{m_1}} \overset{3}{\PP}(E_{m_1})(y)\,dy \right) \right|\notag\notag\\
    &\lesssim& \sum_{m=1}^\aleph \left( a^4\left\lVert \overset{1}{\PP}(\tilde{E}_{m_2})\right\rVert_{\LL^2(B_2)} + a \left\lVert \overset{3}{\PP}(\tilde{E}_{m_1})\right\rVert_{\LL^2(B_1)} \right)\underset{\lesssim}{\eqref{max-P1P3}} \; a^{9-2h} d_{\rm in}^{-3} d_{\rm out}^{-\frac{15}{2}},
\end{eqnarray}
 we can further simplify \eqref{scattered-f1} as
 \begin{equation}\notag
     E^s ({x}, \theta, p, \omega) = - k^2 \sum_{m=1}^\aleph  \left[ \underset{y}{\nabla}\Phi_k(x, z_{m_0}) \times {Q}_{m_1} - \Upsilon_k(x, z_{m_0})\cdot {R}_{m_2} \right] + \O(a^{10-3h - \frac{23}{2}t_2}).
\end{equation}
With $(\mathring{Q}_{m_1}, \mathring{R}_{m_2})_{m=1}^\aleph$ given in \eqref{system-mathrix-error}, there holds
\begin{eqnarray}
     E^s ({x}, \theta, p, \omega) &=& - k^2 \sum_{m=1}^\aleph  \left[ \underset{y}{\nabla}\Phi_k(x, z_{m_0}) \times \mathring{Q}_{m_1} - \Upsilon_k(x, z_{m_0})\cdot \mathring{R}_{m_2} \right]\notag\\
     &-& k^2 \sum_{m=1}^\aleph  \left[ \underset{y}{\nabla}\Phi_k(x, z_{m_0}) \times \left( Q_{m_1} -\mathring{Q}_{m_1}\right) - \Upsilon_k(x, z_{m_0})\cdot \left(R_{m_2} - \mathring{R}_{m_2}\right) \right]+ \O(a^{10-3h - \frac{23}{2}t_2}),\notag
\end{eqnarray}
where 
\begin{eqnarray}
    && \left| k^2 \sum_{m=1}^\aleph  \left[ \underset{y}{\nabla}\Phi_k(x, z_{m_0}) \times \left( Q_{m_1} -\mathring{Q}_{m_1}\right) - \Upsilon_k(x, z_{m_0})\cdot \left(R_{m_2} - \mathring{R}_{m_2}\right) \right] \right|\notag\\
    &\lesssim& \aleph^\frac{1}{2} \left[ \left( \sum_{m=1}^\aleph |Q_{m_1}-\mathring{Q}_{m_1}|^2 \right)^\frac{1}{2} + \left(\sum_{m=1}^\aleph |R_{m_2} -\mathring{R}_{m_2} |^2\right)^\frac{1}{2} \right]\notag\\
    &\lesssim& \aleph \cdot a^{6-h-3t_1-\frac{7}{2}t_2} \lesssim a^{6-h-3t_1-\frac{13}{2}t_2}.\notag
\end{eqnarray}
Therefore, we can derive by keeping the dominant terms that 
\begin{equation}\label{scattered-coro-final}
      E^s ({x}, \theta, p, \omega) = - k^2 \sum_{m=1}^\aleph  \left[ \underset{y}{\nabla}\Phi_k(x, z_{m_0}) \times \mathring{Q}_{m_1} - \Upsilon_k(x, z_{m_0})\cdot \mathring{R}_{m_2} \right] + \O(a^{10-3h - \frac{23}{2}t_2}).
\end{equation}

% \begin{eqnarray}\label{far-field-taylor}
%     E^s (\hat{x}, \theta, p) &=& - k^2 \sum_{m=1}^\aleph \sum_{i=1}^2 \left[ \underset{y}{\nabla}\Phi_k(x, z_{m_0}) \times \tilde{Q}_{m_i} - \Upsilon_k(x, z_{m_0})\cdot \tilde{R}_{m_i} \right]\notag\\
%     &+& \sum_{m=1}^\aleph \sum_{i=1}^2 \frac{(- i k \hat{x} (z_{m_i}-z_{m_0}))^\ell}{\ell !}\left[ \underset{y}{\nabla}\Phi_k(x, z_{m_0}) \times \tilde{Q}_{m_i} - \Upsilon_k(x, z_{m_0})\cdot \tilde{R}_{m_i}  \right] +\O(a^{10-3h - \frac{23}{2}t_2}).
% \end{eqnarray}
%%%%%%%%%%%%%%%%%%%%%%%%%%%%%%%%%%%%%%%%%%%%%%%%%%%%%%%%far-field%%%%%%%%%%%%%%%%%%
For the far field approximation \eqref{foldy-lax-final}, which can be equivalently expressed as
\begin{equation}\label{far-field-original}
    E^\infty(\hat{x}, \theta, p)= \sum_{m=1}^\aleph \sum_{\ell=1}^2 \left[ e^{- i k \hat{x}\cdot z_{m_\ell}} (I-\hat{x}\otimes \hat{x}) R_{m_\ell} - i k e^{- i k \hat{x}\cdot z_{m_\ell}} \hat{x}\times Q_{m_\ell} \right] +\O(a^{10-3h - \frac{23}{2}t_2}),
\end{equation}
following the similar argument of deriving \eqref{scattered-coro-final}, by taking Taylor expansion at the intermediate point $z_{m_0}$ and keeping the dominant terms, with the notation $(\mathring{Q}_{m_1}, \mathring{R}_{m_2})_{m=1}^\aleph$, there holds
\begin{equation*}
    E^\infty (\hat{x}, \theta, p) = \sum_{m=1}^\aleph \left[ (I-\hat{x}\otimes\hat{x})\mathring{R}_{m_2} - i k \hat{x}\times \mathring{Q}_{m_1} \right] + \O(a^{10-3h - \frac{23}{2}t_2}).
\end{equation*}
The proof is complete.

\section{Appendix: Structure of the effective susceptibilities}
\label{sec:appendix-susceptibilities}

This appendix records the dominant leading-order expressions underlying the susceptibility scaling
\eqref{eq:chi-scaling} and the brief effective-medium discussion in the introduction. Throughout we
work under the assumptions of Theorem~\ref{main-1} and Corollary~\ref{corollary-main} and, for
simplicity, assume that all dimers are identical and identically oriented.

\subsection*{Single-dimer polarizability (dominant part)}

For each hybrid dimer \(D_m=D_{m_1}\cup D_{m_2}\), Corollary~\ref{corollary-main} yields the reduced
\(6\times6\) system
\[
\mathbb A_m
\binom{\mathring Q_{m_1}}{\mathring R_{m_2}}
=
\binom{H^{\rm loc}(z_{m_0})}{E^{\rm loc}(z_{m_0})},
\qquad
\binom{\mathring Q_{m1}}{\mathring R_{m2}}
=
\mathcal P_m
\binom{H^{\rm loc}(z_{m_0})}{E^{\rm loc}(z_{m_0})}
\quad\mbox{with}\quad \mathcal P_m:=\mathbb A_m^{-1},
\]
where \(z_{m_0}\) is the intermediate point and \(\mathbb A_m\) is given in \eqref{def-matrix-corollary}. Writing \(\mathbb A_m\) in \(3\times3\) blocks and setting
\[
G:=\nabla\Phi_k(z_{m_1},z_{m_2})\times,
\]
we have
\[
\mathbb A_{HH}=a^{h-3}\frac{\pm c_0}{ik\eta_0}\big(\mathbf P_{0,1}^{(1)}\big)^{-1},\quad
\mathbb A_{HE}=ik\,G,\quad
\mathbb A_{EH}=-k^2\,G,\quad
\mathbb A_{EE}=a^{h-3}\frac{\pm d_0}{\eta_2}\big(\mathbf P_{0,2}^{(2)}\big)^{-1}.
\]
Since \(|z_{m_1}-z_{m_2}|\sim d_{\rm in}=\alpha_0 a^{t_1}\), we have \(G=O(a^{-2t_1})\). Under the
parameter regime in Theorem~\ref{main-1} (in particular \(h<3-2t_1\)), the Schur-complement
corrections are lower order, and the standard block-inverse formula gives the dominant parts
\begin{equation}\label{eq:A31-alphaHH0-alphaEE0}
\alpha_{HH}^{(0)}
:=
\frac{ik\eta_0}{\pm c_0}\,a^{3-h}\, \mathbf P_{0,1}^{(1)},
\qquad
\alpha_{EE}^{(0)}
:=
\frac{\eta_2}{\pm d_0}\,a^{3-h}\, \mathbf P_{0,2}^{(2)},
\end{equation}
and
\begin{equation}\label{eq:A36-alphaHE0}
\alpha_{HE}^{(0)}
:=
k^2\,\frac{\eta_0\eta_2}{c_0d_0}\,a^{6-2h}\, \mathbf P_{0,1}^{(1)}\,G\, \mathbf P_{0,2}^{(2)},
\end{equation}
\begin{equation}\label{eq:A38-alphaEH0}
\alpha_{EH}^{(0)}
:=
ik^3\,\frac{\eta_0\eta_2}{c_0d_0}\,a^{6-2h}\,\mathbf P_{0,2}^{(2)}\,G\, \mathbf P_{0,1}^{(1)}.
\end{equation}
(Here and below we suppress harmless \(o(1)\) factors in operator norm as \(a\to0\).) In particular,
\[
\|\alpha_{HH}^{(0)}\|,\ \|\alpha_{EE}^{(0)}\|=O(a^{3-h}),
\qquad
\|\alpha_{HE}^{(0)}\|,\ \|\alpha_{EH}^{(0)}\|=O\big(a^{6-2h-2t_1}\big).
\]

\subsection*{From polarizabilities to macroscopic susceptibilities}

Let \(\rho_a\) denote the number density of dimers:
\begin{equation}\label{eq:A39-density}
\rho_a \sim d_{\rm out}^{-3}
= \beta_0^{-3}\,a^{-3t_2},
\end{equation}
where \(d_{\rm out}=\beta_0 a^{t_2}\) is the typical inter-dimer distance. The effective
susceptibilities are obtained by the standard (formal) identification
\begin{equation}\label{eq:A42-chi-def}
\chi_{HH}^{(0)} := \rho_a\,\alpha_{HH}^{(0)},\quad
\chi_{HE}^{(0)} := \rho_a\,\alpha_{HE}^{(0)},\quad
\chi_{EH}^{(0)} := \rho_a\,\alpha_{EH}^{(0)},\quad
\chi_{EE}^{(0)} := \rho_a\,\alpha_{EE}^{(0)}.
\end{equation}
Combining \eqref{eq:A39-density} with \eqref{eq:A31-alphaHH0-alphaEE0}, \eqref{eq:A36-alphaHE0} and
\eqref{eq:A38-alphaEH0}, it gives the dominant parts
\begin{equation}\label{eq:A43-chi-dominant}
\boxed{
\begin{aligned}
\chi_{HH}^{(0)}
&=
\beta_0^{-3}\,\frac{ik\eta_0}{\pm c_0}\,
a^{3-h-3t_2}\, \mathbf P_{0,1}^{(1)},\\[1mm]
\chi_{EE}^{(0)}
&=
\beta_0^{-3}\,\frac{\eta_2}{\pm d_0}\,
a^{3-h-3t_2}\,\mathbf P_{0,2}^{(2)},\\[1mm]
\chi_{HE}^{(0)}
&=
k^2\,\beta_0^{-3}\,\frac{\eta_0\eta_2}{c_0d_0}\,
a^{6-2h-3t_2}\, \mathbf P_{0,1}^{(1)}\,G\, \mathbf P_{0,2}^{(2)},\\[1mm]
\chi_{EH}^{(0)}
&=
ik^3\,\beta_0^{-3}\,\frac{\eta_0\eta_2}{c_0d_0}\,
a^{6-2h-3t_2}\,\mathbf P_{0,2}^{(2)}\,G\, \mathbf P_{0,1}^{(1)}.
\end{aligned}
}
\end{equation}
In terms of scaling exponents,
\begin{equation}\label{eq:A44-chi-scaling}
\|\chi_{HH}^{(0)}\|,\ \|\chi_{EE}^{(0)}\|
= O\big(a^{3-h-3t_2}\big),
\qquad
\|\chi_{HE}^{(0)}\|,\ \|\chi_{EH}^{(0)}\|
= O\big(a^{6-2h-2t_1-3t_2}\big).
\end{equation}

\subsection*{Remarks on regimes}

A natural ``non-trivial effective'' scaling for the diagonal response is
\[
3-h-3t_2=0,
\]
which keeps \(\chi_{HH}^{(0)}\) and \(\chi_{EE}^{(0)}\) of order one (away from the resonant
denominators). Under the same choice \(h=3-3t_2\), the off-diagonal exponent becomes
\(6-2h-2t_1-3t_2 = 3t_2-2t_1\). Hence magneto--electric coupling is:
(i) small if \(3t_2>2t_1\), (ii) order one if \(3t_2=2t_1\), and (iii) potentially dominant if
\(3t_2<2t_1\) (subject to the constraints of Theorem~\ref{main-1}).

Once \(\chi_{HH}^{(0)},\chi_{EE}^{(0)},\chi_{HE}^{(0)},\chi_{EH}^{(0)}\) are identified, the effective constitutive tensors
\(\varepsilon_{\rm eff},\mu_{\rm eff},\xi,\zeta\) follow from \eqref{eq:eff-tensors}. The signs and
anisotropy of \(\varepsilon_{\rm eff}\) and \(\mu_{\rm eff}\) are driven by the resonant parameters
in (\ref{condition-on-k}) and by the eigenstructure of \(\mathbf P_{0,1}^{(1)}\) and \(\mathbf P_{0,2}^{(2)}\); this yields the
single-/double-negative and hyperbolic scenarios discussed in the introduction, while
\(\chi_{HE}^{(0)},\chi_{EH}^{(0)}\) encode the bi-anisotropic coupling.


\begin{thebibliography}{26}

\bibitem{AluEngheta2007}
A.~Al\`u, M.~G.~Silveirinha, A.~Salandrino, and N.~Engheta,
\newblock Epsilon-near-zero metamaterials and electromagnetic sources: Tailoring the radiation phase pattern,
\newblock {\em Phys. Rev. B} {\bf 75} (2007), 155410.

\bibitem{AmmariKang_PMT}
H.~Ammari and H.~Kang,
\newblock {\em Polarization and Moment Tensors with Applications to Inverse Problems and Effective Medium Theory},
\newblock Applied Mathematical Sciences, vol.~162, Springer, New York, 2007.

\bibitem{AmmariEtAl_Maxwell}
H.~Ammari, M.~Ruiz, S.~Yu, and H.~Zhang,
\newblock Mathematical analysis of plasmonic resonances for nanoparticles: The full Maxwell equations,
\newblock {\em J. Differential Equations} {\bf 261} (2016), 3615--3669.

\bibitem{AmmariEtAl_SPR_ARMA}
H.~Ammari, P.~Millien, M.~Ruiz, and H.~Zhang,
\newblock Surface plasmon resonance of nanoparticles and applications,
\newblock {\em Arch. Ration. Mech. Anal.} {\bf 220} (2016), 109--153.

\bibitem{AmmariEtAl_DoubleNegAcoustic}
H.~Ammari, B.~Fitzpatrick, H.~Lee, E.~O.~Hiltunen, and S.~Yu,
\newblock Double-negative acoustic metamaterials,
\newblock {\em Quart. Appl. Math.} {\bf 77} (2019), 767--791.

\bibitem{AmmariWuYu_ChiralDoubleNeg}
H.~Ammari, W.~Wu, and S.~Yu,
\newblock Double-negative electromagnetic metamaterials due to chirality,
\newblock {\em Quart. Appl. Math.} {\bf 77} (2019), 105--130.

\bibitem{Ammai-Li-Zou-2023}
H. Ammari; B. Li; J. Zou,
\newblock{Mathematical analysis of electromagnetic scattering by dielectric nanoparticles with high refractive indices. Trans. Amer. Math. Soc. 376 (2023), no. 1, 39-90.}

\bibitem{amrouche-bernardi-dauge-girault}
C. Amrouche, C. Bernardi, M. Dauge and V. Girault,
\newblock{Vector potentials in three-dimensional non-smooth domains, Mathematical Methods in the Applied Sciences, 21 (9), 1998.}

\bibitem{BouchitteFelbacq_ArtificialMagnetism}
G.~Bouchitt\'e and D.~Felbacq,
\newblock Homogenization near resonances and artificial magnetism from dielectrics,
\newblock {\em C. R. Acad. Sci. Paris Ser. I} {\bf 339} (2004), 377--382.

\bibitem{BouchitteSchweizer_SplitRing}
G.~Bouchitt\'e and B.~Schweizer,
\newblock Homogenization of Maxwell's equations in a split ring geometry,
\newblock {\em Multiscale Model. Simul.} {\bf 8} (2010), 717--750.

\bibitem{cao-ghandriche-sini-dimer}
X. Cao, A. Ghandriche and M. Sini, \newblock{Electromagnetic Waves Generated by A Hybrid Dieletric-Plasmonic Dimer, SIAM J. Appl. Math., Vol  85, No. 5, pp. 1949--1975.}

\bibitem{cao-sini-effective}
X. Cao and M. Sini,
\newblock{The effective permittivity and permeability generated by a cluster of moderately contrasting nanoparticles. J. Differential Equations 367 (2023), 549-602.}

\bibitem{cao-ghandriche-sini-highindex}
X. Cao, A. Ghandriche and M. Sini,
\newblock{The electromagnetic waves generated by a cluster of nanoparticles with high refractive indices, Journal of the London Mathematical Society, 2023.}

\bibitem{ChenLipton_DoubleNeg}
Y.~Chen and R.~Lipton,
\newblock Resonance and double negative behavior in metamaterials,
\newblock {\em Arch. Ration. Mech. Anal.} {\bf 209} (2013), 835--868.

\bibitem{colton2019inverse}
D. Colton and R. Kress,
\newblock{Inverse acoustic and electromagnetic scattering theory, 93, 2019, Springer Nature.}

\bibitem{Dautry-Lions}
R. Dautray and J. L. Lions,
\newblock{Mathematical Analysis and Numerical Methods for Science and Technology Volume 3 Spectral Theory and Applications.}
\newblock{1st ed. Berlin, Heidelberg: Springer Berlin Heidelberg, 2000.}

\bibitem{friedman-pasciak}
M. J. Friedman and J. E. Pasciak,
\newblock{Spectral Properties for the Magnetization Integral Operator, Mathematics of Computation, 43 (168), 447--453, 1984.}

\bibitem{GS}
A. Ghandriche and M. Sini,
\newblock{Photo-acoustic inversion using plasmonic contrast agents: The full Maxwell model, J. Differential Equations, Vol. 341, 25 December 2022, pp 1--78.}

\bibitem{Lindell1994}
I.~V.~Lindell, A.~H.~Sihvola, S.~A.~Tretyakov, and A.~J.~Viitanen,
\newblock {\em Electromagnetic Waves in Chiral and Bi-Isotropic Media},
\newblock Artech House, Boston, 1994.

\bibitem{Maas2013}
R.~Maas, J.~Parsons, N.~Engheta, and A.~Polman,
\newblock Experimental realization of an epsilon-near-zero metamaterial at visible wavelengths,
\newblock {\em Nat. Photonics} {\bf 7} (2013), 907--912.

\bibitem{mitrea-mitrea-pipher}
D. Mitrea, M. Mitrea and J. Pipher,
\newblock{Vector potential theory on nonsmooth domains in $\mathbb{R}^{3}$ and applications to electromagnetic scattering. The Journal of Fourier Analysis and Applications 3, 131--192, 1997.}


\bibitem{nedelec}
J. ~C. N\'ed\'elec,
\newblock{Acoustic and electromagnetic equations
Appl. Math. Sci., 144
Springer-Verlag, New York, 2001, x+316 pp.}



\bibitem{Pendry2000}
J.~B.~Pendry,
\newblock Negative refraction makes a perfect lens,
\newblock {\em Phys. Rev. Lett.} {\bf 85} (2000), 3966--3969.

\bibitem{Poddubny2013}
A.~Poddubny, I.~Iorsh, P.~Belov, and Y.~Kivshar,
\newblock Hyperbolic metamaterials,
\newblock {\em Nat. Photonics} {\bf 7} (2013), 948--957.

\bibitem{raevskii}
V. Ya. Raevskii,
\newblock{Some properties of the operators of potential theory and their application to the investigation of the basic equation of electrostatics and magnetostatics,}
\newblock{Theoretical and Mathematical Physics, Volume 100, Number 3, pages 1040--1045, 1994.}

\bibitem{Serdyukov2001}
A.~Serdyukov, I.~Semchenko, S.~Tretyakov, and A.~Sihvola,
\newblock {\em Electromagnetics of Bi-anisotropic Materials: Theory and Applications},
\newblock Gordon and Breach Science Publishers, Amsterdam, 2001.

\bibitem{Smith2000}
D.~R.~Smith, W.~J.~Padilla, D.~C.~Vier, S.~C.~Nemat-Nasser, and S.~Schultz,
\newblock Composite medium with simultaneously negative permeability and permittivity,
\newblock {\em Phys. Rev. Lett.} {\bf 84} (2000), 4184--4187.

\bibitem{Veselago1968}
V.~G.~Veselago,
\newblock The electrodynamics of substances with simultaneously negative values of $\varepsilon$ and $\mu$,
\newblock {\em Soviet Phys. Usp.} {\bf 10} (1968), 509--514.

\end{thebibliography}
\end{document}

\sloppy

% \date{\today}

\allowdisplaybreaks

% \begin{abstract}
% We generate the electromagnetism by...
% \end{abstract}

% \subjclass[2010]{35R30, 35C20, 35Q60....}
% \keywords{...}

\maketitle

\tableofcontents

\section*{Supplementary material: purpose and organization}

This document collects technical material that supports the main paper \emph{``Electromagnetic Scattering by a Cluster of Hybrid Dielectric--Plasmonic Dimers''}. 
The goal is to keep the main paper focused on the core statements and their interpretation, while providing (i) detailed a-priori estimates for the Lippmann-Schwinger system and (ii) the longer algebraic/tensor manipulations used to obtain the simplified discrete systems and the key invertibility criteria.

\medskip

\begin{notation}
To make a cross reference to the main part document, we have added the prefix $MP-$ before the element that we want to make a reference to. For example, Proposition $MP-1.1$ refers to Proposition 1.1 in the main document, and (MP1.2) refers to equation (1.2) in the main document.
\end{notation}

\medskip

\noindent\textbf{How to read.}
\begin{itemize}
    \item \textbf{Proposition MP-3.4 (a-priori estimates):} the full proof starts at ``\textbf{Proof of Proposition MP-3.4}'' and provides the estimates used throughout.
    \item \textbf{Simplified linear systems and tensor reductions:} see the numbered item ``(5). Tensor analysis and the simplified linear algebraic systems.'' 
    \item \textbf{Invertibility / block-diagonal dominance arguments:} see the sections leading to the condition ensuring invertibility of the discrete system (MP1.25).
    \item \textbf{Additional technical proofs:} collected in the final \textbf{Appendix}.
\end{itemize}

\section{Proof of the a-priori estimates}\label{sec-proof-prior}

To justify the a-priori estimates in {Proposition MP-3.4}, we need the following preliminary result.

\begin{proposition}\label{prop-es-QR-domain}
    For any fixed $m$, $1\leq m\leq \aleph$, the linear algebraic system (MP4.4) admits
    \begin{equation}\label{eq-al-fix m}
        \mathfrak{B}_m \tilde{\mathfrak{U}}_m-\sum_{j=1\atop j\neq m}^\aleph \Psi_{mj} \tilde{\mathfrak{U}}_j = \mathrm{S}_m +\mathfrak{R}_m, \quad\mbox{with}\quad   \mathfrak{R}_m:=\begin{pmatrix}
    Error_{\tilde{Q}_{m_1}} \\
     Error_{\tilde{R}_{m_1}} \\
     Error_{\tilde{Q}_{m_2}} \\
     Error_{\tilde{R}_{m_2}}
\end{pmatrix}, 
    \end{equation}
    where $\left(\tilde{\mathfrak{U}}_{m} \, := \, \left( \tilde{Q}_{m_1}, \tilde{R}_{m_1}, \tilde{Q}_{m_2}, \tilde{R}_{m_2} \right) \right)_{m=1}^{\aleph}$ is the corresponding solution vector; $ \left\{\mathfrak{B_{m}} \right\}_{m=1}^{\aleph}$, $\left\{ \Psi_{mj}\right\}_{m,j=1}^{\aleph}$ and $\left\{ \mathrm{S}_{m} \right\}_{m=1}^{\aleph}$ are introduced in $\textbf{Notations MP-1.2}$, and $\mathfrak{R}_m$ signifies the remainder term given by \eqref{Equa-Error-Qm1-1229}.
    Then there holds the following estimation
    \begin{eqnarray}\label{es-Qm-main}
        \sum_{m=1}^\aleph |\tilde{Q}_{m_1}|^2 &\lesssim& a^{6-2h} d_{\rm out}^{-3}+ \left( a^{12-2h} d_{\rm out}^{-11} + a^{10} d_{\rm out}^{-9} \ln^2 d_{\rm out} +a^{22-2h} d_{\rm in}^{-14} d_{\rm out}^{-3} \right)\max_m \lVert \overset{1}{\PP}(\tilde{E}_{m_1})\rVert^2_{\LL^2(B_1)}\notag\\
        &+& \left( a^{16-2h} d_{\rm in}^{-8} d_{\rm out}^{-3} + a^{14} d_{\rm in}^{-6} d_{\rm out}^{-3} +a^{14} d_{\rm out}^{-9} \ln^2 d_{\rm out} \right)\max_m \lVert \overset{1}{\PP}(\tilde{E}_{m_2})\rVert^2_{\LL^2(B_2)}\notag\\
        &+& \Big( a^{4} d_{\rm out}^{-3} + a^{8} d_{\rm out}^{-9} + a^{10-2h} d_{\rm out}^{-9} \ln^2 d_{\rm out} + a^{14-2h} d_{\rm in}^{-6} d_{\rm out}^{-9} +  a^{16-4h} d_{\rm in}^{-8} d_{\rm out}^{-9}\notag\\
        &+& a^{14-2h} d_{\rm out}^{-15} \ln^2 d_{\rm out}\Big)\max_m \lVert \overset{3}{\PP}(\tilde{E}_{m_1})\rVert^2_{\LL^2(B_1)} + \Big( a^{14-2h} d_{\rm in}^{-6} d_{\rm out}^{-3} + a^{12} d_{\rm in}^{-4} d_{\rm out}^{-3} + a^{12} d_{\rm out}^{-9} \notag\\
        &+& a^{20-4h}  d_{\rm out}^{-17} +  a^{18-2h} d_{\rm out}^{-15} \ln^2 d_{\rm out} + a^{26-4h} d_{\rm in}^{-14} d_{\rm out}^{-9} \Big)\max_m \lVert \overset{3}{\PP}(\tilde{E}_{m_2})\rVert^2_{\LL^2(B_2)}, 
    \end{eqnarray}
        and
        \begin{eqnarray}\label{es-Rm-main}
        \sum_{m=1}^\aleph |\tilde{R}_{m_2}|^2 &\lesssim& a^{6-2h} d_{\rm out}^{-3}+ \Big( a^{10} d_{\rm in}^{-4} d_{\rm out}^{-3} + a^{10} d_{\rm out}^{-9} +a^{12-2h} d_{\rm in}^{-6} d_{\rm out}^{-3} + a^{12-2h} d_{\rm out}^{-9} \ln^2 d_{\rm out} \notag\\
        &+& a^{18-4h} d_{\rm out}^{-17} + a^{16-2h}d_{\rm out}^{-15}\ln^2 d_{\rm out} \Big)\max_m \lVert \overset{1}{\PP}(\tilde{E}_{m_1})\rVert^2_{\LL^2(B_1)}\notag\\
        &+& \Big( a^{14} d_{\rm out}^{-9} + a^{16-2h} d_{\rm out}^{-9} \ln^2 d_{\rm out} +a^{22-2h} d_{\rm in}^{-12} d_{\rm out}^{-3} + a^{22-4h} d_{\rm in}^{-8} d_{\rm out}^{-9} \notag\\
        &+& a^{20-2h} d_{\rm in}^{-6} d_{\rm out}^{-9} + a^{20-2h} d_{\rm out}^{-15}\ln^2 d_{\rm out} \Big)\max_m \lVert \overset{1}{\PP}(\tilde{E}_{m_2})\rVert^2_{\LL^2(B_2)}\notag\\
        &+& \Big( a^{8} d_{\rm in}^{-6} d_{\rm out}^{-3} + a^{8} d_{\rm out}^{-9} \ln^2 d_{\rm out} + a^{10-2h} d_{\rm out}^{-9} + a^{10-2h} d_{\rm in}^{-8} d_{\rm out}^{-3} \Big)\max_m \lVert \overset{3}{\PP}(\tilde{E}_{m_1})\rVert^2_{\LL^2(B_1)} \notag\\
        &+& \Big( a^{12} d_{\rm out}^{-9} \ln^2 d_{\rm out} + a^{14-2h} d_{\rm out}^{-11} + a^{20-4h} d_{\rm in}^{-6} d_{\rm out}^{-9} + a^{20-2h} d_{\rm in}^{-14}  d_{\rm out}^{-3} \Big)\max_m \lVert \overset{3}{\PP}(\tilde{E}_{m_2})\rVert^2_{\LL^2(B_2)}.
    \end{eqnarray}
\end{proposition}
See Appendix \ref{sec:Appendix-prop5.1} for the proof of this proposition.

\bigskip

\bigskip

\noindent \textbf{Proof of Proposition MP-3.4}.
For fixed $m$, the Lippmann-Schwinger equation (MP3.5) indicates that in the dimer $D_m=D_{m_1}\cup D_{m_2}$, there holds
		\begin{equation}\label{LS-Dm-sec5}
	 			\left(I+\eta \nabla M^k - k^2 \eta N^k\right)E_m(x)+\sum_{j=1\atop j\neq m}^\aleph \eta \nabla M^k(E_j)(x)-\sum_{j=1\atop j\neq m}^\aleph k^2 \eta N^k(E_j)(x)=E_m^{Inc}(x).
	 		\end{equation}
Denote
    \begin{equation}\label{denote-Tk}
        \TT_{k, \ell}=I \, +\, \eta_\ell\, \nabla M^k\, -\, k^2\, \eta_\ell\, N^k\quad\mbox{for}\quad \ell=1,2.
    \end{equation}
% \bigskip

\medskip

\noindent (1). \textit{Estimate for $\max_m \left\lVert \overset{1}{\PP}(\tilde{E}_{m_1})\right\rVert^2_{\LL^2(B_1)}$}.

\bigskip

In $D_{m_1}$, with the notation \eqref{denote-Tk}, the Lippmann-Schwinger equation \eqref{LS-Dm-sec5} becomes 
    \begin{eqnarray}\label{LS-Dm1}
        \TT_{k, 1}E_{m_1} (x) &+ &\eta_2 \left( \nabla M^k - k^2 N^k\right)(E_{m_2})(x)+ \sum_{j=1\atop j\neq m}^\aleph \eta_1 \left(\nabla M^k - k^2 N^k\right)(E_{j_1})(x)\notag\\
        & + & \sum_{j=1\atop j\neq m}^\aleph \eta_2 \left(\nabla M^k - k^2 N^k\right)(E_{j_2})(x)=E_{m_1}^{Inc}(x).
    \end{eqnarray}
%     By taking $\TT_{k, 1}$ on the both sides of \eqref{LS-Dm1}, with the notation \eqref{E-split}, we can obtain that
%     \begin{eqnarray}\label{LS-Dm1-P11-split}
%         \overset{1}{\PP}(E_{m_1})&+& \eta_2 \, \TT_{k,1}^{-1} \left((\nabla M^k - k^2 N^k)\overset{1}{\PP}(E_{m_2})\right)+ \sum_{j=1 \atop j\neq m}^\aleph \eta_1 \, \TT_{k,1}^{-1} \left((\nabla M^k - k^2 N^k)\overset{1}{\PP}(E_{j_1})\right)\notag\\
%         &+& \sum_{j=1 \atop j\neq m}^\aleph \eta_2 \, \TT_{k,1}^{-1} \left((\nabla M^k - k^2 N^k)\overset{1}{\PP}(E_{j_2})\right)\notag\\
%         &=& \TT_{k, 1}^{-1} (E_{m_1}^{Inc}) - \overset{3}{\PP}(E_{m_1}) -\eta_2 \TT_{k,1}^{-1} \left((\nabla M^k - k^2 N^k)\overset{3}{\PP}(E_{m_2})\right)\notag\\
%         &-& \sum_{j=1 \atop j\neq m}^\aleph \eta_1 \, \TT_{k,1}^{-1} \left((\nabla M^k - k^2 N^k)\overset{3}{\PP}(E_{j_1})\right) - \sum_{j=1 \atop j\neq m}^\aleph \eta_2 \, \TT_{k,1}^{-1} \left((\nabla M^k - k^2 N^k)\overset{3}{\PP}(E_{j_2})\right).
%     \end{eqnarray}
% Since 
%%%%%%%%%%%%%%%%%%%%%%%%%%%%%%%%%%%%%%%%%%%%%%%%%%%%%%%%%%

By using the split (MP3.15) and the fact that    
\begin{equation}\label{M|div=0}
    \nabla M^k \left(\overset{1}{\PP}(E_{m_\ell})\right)=0\quad\mbox{for any}\quad m=1, 2, \cdots, \aleph\quad\mbox{and}\quad \ell=1, 2,
\end{equation}
we shall reformulate \eqref{LS-Dm1} as
\begin{eqnarray}\label{LS-Dm1-P11-Phi}
    \TT_{k, 1} \left(\overset{1}{\PP}(E_{m_1})\right)(x) &+& \TT_{k, 1} \left(\overset{3}{\PP}(E_{m_1})\right)(x) - k^2 \eta_2 N^k\left(\overset{1}{\PP}(E_{m_2})\right)(x) + \eta_2 \left(\nabla M^k-k^2 N^k\right)\left(\overset{3}{\PP}(E_{m_2})\right)(x)\notag\\
    &-& k^2 \eta_1\sum_{j=1 \atop j\neq m}^\aleph  N^k\left(\overset{1}{\PP}(E_{j_1})\right)(x) + \eta_1 \sum_{j=1 \atop j\neq m}^\aleph (\nabla M^k -k^2 N^k)\left(\overset{3}{\PP}(E_{j_1})\right)(x)\notag\\
    &-& k^2 \eta_2 \sum_{j=1 \atop j\neq m}^\aleph N^k\left(\overset{1}{\PP}(E_{j_2})\right)(x) + \eta_2 \sum_{j=1 \atop j\neq m}^\aleph (\nabla M^k- k^2 N^k)\left(\overset{3}{\PP}(E_{j_2})\right)(x)= E_{m_1}^{Inc}(x).
\end{eqnarray}
By taking the $\mathbb{L}^2$-inner product w.r.t. $e_{n}^{(1)}(D_{m_1})$ on the both sides of \eqref{LS-Dm1-P11-Phi}, where 
\begin{equation*}
    e_n^{(1)}(x) := e_n^{(1)}\left(\frac{x-z_{m_1}}{a}\right)a^{-\frac{3}{2}},\quad x\in D_{m_1},
\end{equation*}
due to the fact that for any vector function $F$, there holds
\begin{equation}\label{M-project1=0}
    \left\langle\nabla M^k(F), e_n^{(1)}\right\rangle=\left\langle F, \nabla M^{-k}(e_n^{(1)}) \right\rangle=0,
\end{equation}
we can derive from \eqref{LS-Dm1-P11-Phi} that 
\begin{eqnarray}\label{LS-Dm1-P11-inner-first}
    &&\left\langle\TT_{k, 1}\left(\overset{1}{\PP}(E_{m_1})\right), e_n^{(1)} \right\rangle + \left\langle \TT_{k, 1}\left(\overset{3}{\PP}(E_{m_1})\right), e_n^{{(1)}} \right\rangle- k^2 \eta_2 \left\langle N^k\left(\overset{1}{\PP}(E_{m_2})\right), e_n^{(1)} \right\rangle\notag\\
    &-& k^2 \eta_2 \left\langle N^k\left(\overset{3}{\PP}(E_{m_2})\right), e_n^{(1)}\right\rangle - k^2 \eta_1 \sum_{j=1 \atop j\neq m}^\aleph \left\langle N^k\left(\overset{1}{\PP}(E_{j_1})\right), e_n^{(1)} \right\rangle - k^2 \eta_1 \sum_{j=1 \atop j\neq m}^\aleph \left\langle N^k\left(\overset{3}{\PP}(E_{j_1})\right), e_n^{(1)} \right\rangle\notag\\
    &-& k^2 \eta_2 \sum_{j=1 \atop j\neq m}^\aleph \left\langle N^k\left(\overset{1}{\PP}(E_{j_2})\right), e_n^{(1)} \right\rangle - k^2 \eta_2 \sum_{j=1 \atop j\neq m}^\aleph \left\langle N^k\left(\overset{3}{\PP}(E_{j_2})\right), e_n^{(1)}\right\rangle = \left\langle E_{m_1}^{Inc}, e_n^{(1)} \right\rangle.\notag\\
\end{eqnarray}
It is direct to see that the second term on the L.H.S. of \eqref{LS-Dm1-P11-inner-first} fulfills
\begin{equation*}
    \left\langle \TT_{k, 1}\left(\overset{3}{\PP}(E_{m_1})\right), e_n^{(1)} \right\rangle = \left\langle \overset{3}{\PP}(E_{m_1}), \TT_{-k, 1}(e_n^{(1)})  \right\rangle, 
\end{equation*}
where 
\begin{eqnarray*}
    \TT_{-k, 1}&:=&(I+\eta_1 \nabla M^{-k}-k^2 \eta_1  N^{-k})(e_n^{(1)})\notag\\
    &=& (I- k^2 \eta_1 N)(e_n^{(1)}) - k^2 \eta_1 (N^{-k}-N)(e_n^{(1)})\notag\\
    &=&(1-k^2 \eta_1 \lambda_n^{(1)})(e_n^{(1)})-\underbrace{\frac{k^2 \eta_1}{4\pi}\sum_{\ell\geq 1}\frac{(-ik)^{\ell+1}}{(\ell+1)!}\int_{D_{m_1}}\lVert\cdot-y\rVert^\ell e_n^{(1)}(y)\,dy}_{=:r_{k, 1}(x)},
\end{eqnarray*}
by using \cite[Remark 2.2]{CGS-one dimer}. Then 
\begin{equation}\label{LHS-2-simply}
    \left\langle \TT_{k, 1}\left(\overset{3}{\PP}(E_{m_1})\right), e_n^{(1)} \right\rangle = (1-k^2 \eta_1 \lambda_n^{(1)})\left\langle \overset{3}{\PP}(E_{m_1}), e_n^{(1)} \right\rangle- \left\langle\overset{3}{\PP}(E_{m_1}), r_{k, 1}\right\rangle=- \left\langle\overset{3}{\PP}(E_{m_1}), r_{k, 1}\right\rangle,
\end{equation}
which indicates that \eqref{LS-Dm1-P11-inner-first} becomes
\begin{eqnarray}\label{LS-Dm1-P11-inner-first-simple}
    &&\left\langle\TT_{k, 1}\left(\overset{1}{\PP}(E_{m_1})\right), e_n^{(1)} \right\rangle - \left\langle \overset{3}{\PP}(E_{m_1}), r_{k,1}\right\rangle - k^2 \eta_2 \left\langle N^k\left(\overset{1}{\PP}(E_{m_2})\right), e_n^{(1)} \right\rangle\notag\\
    &-& k^2 \eta_2 \left\langle N^k\left(\overset{3}{\PP}(E_{m_2})\right), e_n^{(1)}\right\rangle - k^2 \eta_1 \sum_{j=1 \atop j\neq m}^\aleph \left\langle N^k\left(\overset{1}{\PP}(E_{j_1})\right), e_n^{(1)} \right\rangle - k^2 \eta_1 \sum_{j=1 \atop j\neq m}^\aleph \left\langle N^k\left(\overset{3}{\PP}(E_{j_1})\right), e_n^{(1)} \right\rangle\notag\\
    &-& k^2 \eta_2 \sum_{j=1 \atop j\neq m}^\aleph \left\langle N^k\left(\overset{1}{\PP}(E_{j_2})\right), e_n^{(1)} \right\rangle - k^2 \eta_2 \sum_{j=1 \atop j\neq m}^\aleph \left\langle N^k\left(\overset{3}{\PP}(E_{j_2})\right), e_n^{(1)}\right\rangle = \left\langle E_{m_1}^{Inc}, e_n^{(1)} \right\rangle.\notag\\
\end{eqnarray}

Taking the Taylor expansion for $\Phi_k(x, y)$ in $N^k\left(\overset{1}{\PP}(E_{m_2})\right)$ and $N^k\left(\overset{3}{\PP}(E_{m_2})\right)$, we respectively have
\begin{eqnarray*}
    N^k\left(\overset{1}{\PP}(E_{m_2})\right)&=& \int_{D_{m_2}}\Phi_k(x, y)\overset{1}{\PP}(E_{m_2})(y)\,dy\notag\\
    &=& \underset{y}{\nabla}\Phi_k(z_{m_1}, z_{m_2})\int_{D_{m_2}}\mathcal{P}(y, z_{m_2})\overset{1}{\PP}(E_{m_2})(y)\,dy\notag\\
    &+&\underset{x}{\nabla}\underset{y}{\nabla}\Phi_k(z_{m_1}, z_{m_2})(x-z_{m_1})\int_{D_{m_2}}\mathcal{P}(y, z_{m_2})\overset{1}{\PP}(E_{m_2})(y)\,dy\notag\\
    &+& \int_{D_{m_2}} \int_0^1 (1-t)(y-z_{m_2})^\perp\underset{y}{\nabla}\underset{y}{\nabla}\Phi_k(z_{m_1}, z_{m_2}+t(y-z_{m_2}))(y-z_{m_2})\,dt\overset{1}{\PP}(E_{m_2})(y)\,dy\notag\\
    &+& \int_{D_{m_2}} \int_0^1 (x-z_{m_1})^\perp \underset{x}{\nabla} (\underset{x}{\nabla}\underset{y}{\nabla}\Phi_k)(z_{m_1}+t(x-z_{m_1}), z_{m_2})(x-z_{m_1})\,dt\mathcal{P}(y, z_{m_2})\overset{1}{\PP}(E_{m_2})(y)\,dy\notag\\
    &+&\int_{D_{m_2}}\int_0^1 \int_0^1 (1-t) (y-z_{m_2})^\perp \underset{x}{\nabla}(\underset{y}{\nabla}\underset{y}{\nabla}\Phi_k)(z_{m_1}+s(x-z_{m_1}), z_{m_2}+t(y-z_{m_2}))\notag\\
    &&\qquad\qquad\qquad\qquad\qquad \cdot (y-z_{m_2})\,dt\cdot\mathcal{P}(x,z_{m_1})\,ds\overset{1}{\PP}(E_{m_2})\,dy,
\end{eqnarray*}
and 
\begin{eqnarray*}
    N^k\left(\overset{3}{\PP}(E_{m_2})\right)&=&\int_{D_{m_2}}\Phi_k(x, y)\overset{3}{\PP}(E_{m_2})(y)\,dy\notag\\
    &=&\Phi_k(z_{m_1}, z_{m_2})\int_{D_{m_2}}\overset{3}{\PP}(E_{m_2})\,dy+\underset{x}{\nabla}\Phi_{k}(z_{m_1}, z_{m_2})\mathcal{P}(x, z_{m_1})\int_{D_{m_2}}\overset{3}{\PP}(E_{m_2})(y)\,dy\notag\\
    &+&\underset{y}{\nabla}\Phi_k(z_{m_1}, z_{m_2})\int_{D_{m_2}}\mathcal{P}(y, z_{m_2})\overset{3}{\PP}(E_{m_2})(y)\,dy\notag\\
    &+& \int_{D_{m_2}}\int_0^1 (x-z_{m_1})^\perp \underset{x}{\nabla}\underset{x}{\nabla}\Phi_k(z_{m_1}+t(x-z_{m_1}), z_{m_2})\cdot (x-z_{m_1})\,dt\cdot\overset{3}{\PP}(E_{m_2})(y)\,dy\notag\\
    &+& \int_{D_{m_2}} \int_0^1 \underset{x}{\nabla}\underset{y}{\nabla} \Phi_k(z_{m_1}+t(x-z_{m_1}), z_{m_2})\mathcal{P}(y, z_{m_2})\cdot(x-z_{m_1})\,dt\cdot\overset{3}{\PP}(E_{m_2})(y)\,dy\notag\\
    &+&\int_{D_{m_2}}\int_0^1(1-t) (y-z_{m_2})^\perp \underset{y}{\nabla}\underset{y}{\nabla}\Phi_k(z_{m_1}, z_{m_2}+t(y-z_{m_2}))\cdot(y-z_{m_2})\,dt\cdot\overset{3}{\PP}(E_{m_2})(y)\,dy\notag\\
    &+& \int_{D_{m_2}}\int_0^1 \int_0^1 (1-t) (y-z_{m_2})^\perp \underset{x}{\nabla}(\underset{y}{\nabla}\underset{y}{\nabla}\Phi_k)(z_{m_1}+s(x-z_{m_1}), z_{m_2}+t(y-z_{m_2}))\notag\\
    &&\qquad\qquad\qquad\qquad\qquad \cdot (y-z_{m_2})\,dt\cdot\mathcal{P}(x, z_{m_1})\,ds\cdot\overset{3}{\PP}(E_{m_2})(y)\,dy,
\end{eqnarray*}
which directly yields that
\begin{eqnarray}\label{inner p-P1Em2}
    \left\langle N^k\left(\overset{1}{\PP}(E_{m_2})\right), e_n^{(1)} \right\rangle &=& \left\langle \underset{x}{\nabla}\underset{y}{\nabla}\Phi_k(z_{m_1}, z_{m_2})(x-z_{m_1})\int_{D_{m_2}}\mathcal{P}(y, z_{m_2})\overset{1}{\PP}(E_{m_2})(y)\,dy, e_n^{(1)} \right\rangle\notag\\
    &+& \left\langle Err_1\left(\overset{1}{\PP}(E_{m_2})\right), e_n^{(1)} \right\rangle,
\end{eqnarray}
with 
\begin{eqnarray}\label{Error_1-P1m2}
    Err_1\left(\overset{1}{\PP}(E_{m_2})\right)&=& \int_{D_{m_2}} \int_0^1 (x-z_{m_1})^\perp \underset{x}{\nabla} (\underset{x}{\nabla}\underset{y}{\nabla}\Phi_k)(z_{m_1}+t(x-z_{m_1}), z_{m_2})(x-z_{m_1})\,dt\mathcal{P}(y, z_{m_2})\overset{1}{\PP}(E_{m_2})(y)\,dy\notag\\
    &+&\int_{D_{m_2}}\int_0^1 \int_0^1 (1-t) (y-z_{m_2})^\perp \underset{x}{\nabla}(\underset{y}{\nabla}\underset{y}{\nabla}\Phi_k)(z_{m_1}+s(x-z_{m_1}), z_{m_2}+t(y-z_{m_2}))\notag\\
    &&\qquad\qquad\qquad\qquad\qquad \cdot (y-z_{m_2})\,dt\cdot\mathcal{P}(x,z_{m_1})\,ds\overset{1}{\PP}(E_{m_2})\,dy,
\end{eqnarray}
and 
\begin{eqnarray}\label{inner p-P3Em2}
    \left\langle N^k\left(\overset{3}{\PP}(E_{m_2})\right), e_n^{(1)} \right\rangle &=& \left\langle \underset{x}{\nabla}\Phi_k(z_{m_1}, z_{m_2})\mathcal{P}(x, z_{m_1})\int_{D_{m_2}}\overset{3}{\PP}(E_{m_2})(y)\,dy, e_n^{(1)} \right\rangle\notag\\
    &+& \left\langle Err_1\left(\overset{3}{\PP}(E_{m_2})\right), e_n^{(1)} \right\rangle,
\end{eqnarray}
with
\begin{eqnarray}\label{Error_1-P3m2}
    Err_1\left(\overset{3}{\PP}(E_{m_2})\right)&=& \int_{D_{m_2}}\int_0^1 (x-z_{m_1})^\perp \underset{x}{\nabla}\underset{x}{\nabla}\Phi_k(z_{m_1}+t(x-z_{m_1}), z_{m_2})\cdot (x-z_{m_1})\,dt\cdot\overset{3}{\PP}(E_{m_2})(y)\,dy\notag\\
    &+& \int_{D_{m_2}} \int_0^1 \underset{x}{\nabla}\underset{y}{\nabla} \Phi_k(z_{m_1}+t(x-z_{m_1}), z_{m_2})\mathcal{P}(y, z_{m_2})\cdot(x-z_{m_1})\,dt\cdot\overset{3}{\PP}(E_{m_2})(y)\,dy\notag\\
    &+& \int_{D_{m_2}}\int_0^1 \int_0^1 (1-t) (y-z_{m_2})^\perp \underset{x}{\nabla}(\underset{y}{\nabla}\underset{y}{\nabla}\Phi_k)(z_{m_1}+s(x-z_{m_1}), z_{m_2}+t(y-z_{m_2}))\notag\\
    &&\qquad\qquad\qquad\qquad\qquad \cdot (y-z_{m_2})\,dt\cdot\mathcal{P}(x, z_{m_1})\,ds\cdot\overset{3}{\PP}(E_{m_2})(y)\,dy.
\end{eqnarray}
In a similar manner, we can get that
\begin{eqnarray}\label{inner p-P1j1j2}
    \left\langle N^k\left(\overset{1}{\PP}(E_{j_1})\right), e_n^{(1)}\right\rangle  &=&  \left\langle \underset{x}{\nabla}\underset{y}{\nabla}\Phi_k(z_{m_1}, z_{j_1})(x-z_{m_1})\int_{D_{j_1}}\mathcal{P}(y, z_{j_1})\overset{1}{\PP}(E_{j_1})(y)\,dy, e_n^{(1)} \right\rangle\notag\\
    &+& \left\langle Err_1\left(\overset{1}{\PP}(E_{j_1})\right), e_n^{(1)} \right\rangle,\notag\\
     \left\langle N^k\left(\overset{3}{\PP}(E_{j_1})\right), e_n^{(1)} \right\rangle &=& \left\langle \underset{x}{\nabla}\Phi_k(z_{m_1}, z_{j_1})\mathcal{P}(x, z_{m_1})\int_{D_{j_1}}\overset{3}{\PP}(E_{j_1})(y)\,dy, e_n^{(1)} \right\rangle\notag\\
    &+& \left\langle Err_1\left(\overset{3}{\PP}(E_{j_1})\right), e_n^{(1)} \right\rangle,\notag\\
    \left\langle N^k\left(\overset{1}{\PP}(E_{j_2})\right), e_n^{(1)}\right\rangle  &=&  \left\langle \underset{x}{\nabla}\underset{y}{\nabla}\Phi_k(z_{m_1}, z_{j_2})(x-z_{m_1})\int_{D_{j_2}}\mathcal{P}(y, z_{j_2})\overset{1}{\PP}(E_{j_2})(y)\,dy, e_n^{(1)} \right\rangle\notag\\
    &+& \left\langle Err_1\left(\overset{1}{\PP}(E_{j_2})\right), e_n^{(1)} \right\rangle,\notag\\
     \left\langle N^k\left(\overset{3}{\PP}(E_{j_2})\right), e_n^{(1)} \right\rangle &=& \left\langle \underset{x}{\nabla}\Phi_k(z_{m_1}, z_{j_2})\mathcal{P}(x, z_{m_1})\int_{D_{j_2}}\overset{3}{\PP}(E_{j_2})(y)\,dy, e_n^{(1)} \right\rangle\notag\\
    &+& \left\langle Err_1\left(\overset{3}{\PP}(E_{j_2})\right), e_n^{(1)} \right\rangle,
\end{eqnarray}
where
\begin{eqnarray}\label{Err_1-P1j1j2}
     Err_1\left(\overset{1}{\PP}(E_{j_1})\right)&=& \int_{D_{j_1}} \int_0^1 (x-z_{m_1})^\perp \underset{x}{\nabla} (\underset{x}{\nabla}\underset{y}{\nabla}\Phi_k)(z_{m_1}+t(x-z_{m_1}), z_{j_1})(x-z_{m_1})\,dt\mathcal{P}(y, z_{j_1})\overset{1}{\PP}(E_{j_1})(y)\,dy\notag\\
    &+&\int_{D_{j_1}}\int_0^1 \int_0^1 (1-t) (y-z_{j_1})^\perp \underset{x}{\nabla}(\underset{y}{\nabla}\underset{y}{\nabla}\Phi_k)(z_{m_1}+s(x-z_{m_1}), z_{j_1}+t(y-z_{j_1}))\notag\\
    &&\qquad\qquad\qquad\qquad\qquad \cdot (y-z_{j_1})\,dt\cdot\mathcal{P}(x,z_{m_1})\,ds\cdot\overset{1}{\PP}(E_{j_1})\,dy,\notag\\
    Err_1\left(\overset{3}{\PP}(E_{j_1})\right)&=& \int_{D_{j_1}}\int_0^1 (x-z_{m_1})^\perp \underset{x}{\nabla}\underset{x}{\nabla}\Phi_k(z_{m_1}+t(x-z_{m_1}), z_{j_1})\cdot (x-z_{m_1})\,dt\cdot\overset{3}{\PP}(E_{j_1})(y)\,dy\notag\\
    &+& \int_{D_{j_1}} \int_0^1 \underset{x}{\nabla}\underset{y}{\nabla} \Phi_k(z_{m_1}+t(x-z_{m_1}), z_{j_1})\mathcal{P}(y, z_{j_1})\cdot(x-z_{m_1})\,dt\cdot\overset{3}{\PP}(E_{j_1})(y)\,dy\notag\\
    &+& \int_{D_{j_1}}\int_0^1 \int_0^1 (1-t) (y-z_{j_1})^\perp \underset{x}{\nabla}(\underset{y}{\nabla}\underset{y}{\nabla}\Phi_k)(z_{m_1}+s(x-z_{m_1}), z_{j_1}+t(y-z_{j_1}))\notag\\
    &&\qquad\qquad\qquad\qquad\qquad \cdot (y-z_{j_1})\,dt\cdot\mathcal{P}(x, z_{m_1})\,ds\cdot\overset{3}{\PP}(E_{j_1})(y)\,dy,\notag\\
     Err_1\left(\overset{1}{\PP}(E_{j_2})\right)&=& \int_{D_{j_2}} \int_0^1 (x-z_{m_1})^\perp \underset{x}{\nabla} (\underset{x}{\nabla}\underset{y}{\nabla}\Phi_k)(z_{m_1}+t(x-z_{m_1}), z_{j_2})(x-z_{m_1})\,dt\mathcal{P}(y, z_{j_2})\overset{1}{\PP}(E_{j_2})(y)\,dy\notag\\
    &+&\int_{D_{j_2}}\int_0^1 \int_0^1 (1-t) (y-z_{j_2})^\perp \underset{x}{\nabla}(\underset{y}{\nabla}\underset{y}{\nabla}\Phi_k)(z_{m_1}+s(x-z_{m_1}), z_{j_2}+t(y-z_{j_2}))\notag\\
    &&\qquad\qquad\qquad\qquad\qquad \cdot (y-z_{j_2})\,dt\cdot\mathcal{P}(x,z_{m_1})\,ds\cdot\overset{1}{\PP}(E_{j_2})\,dy,\notag\\
    Err_1\left(\overset{3}{\PP}(E_{j_2})\right)&=& \int_{D_{j_2}}\int_0^1 (x-z_{m_1})^\perp \underset{x}{\nabla}\underset{x}{\nabla}\Phi_k(z_{m_1}+t(x-z_{m_1}), z_{j_2})\cdot (x-z_{m_1})\,dt\cdot\overset{3}{\PP}(E_{j_2})(y)\,dy\notag\\
    &+& \int_{D_{j_2}} \int_0^1 \underset{x}{\nabla}\underset{y}{\nabla} \Phi_k(z_{m_1}+t(x-z_{m_1}), z_{j_2})\mathcal{P}(y, z_{j_2})\cdot(x-z_{m_1})\,dt\cdot\overset{3}{\PP}(E_{j_2})(y)\,dy\notag\\
    &+& \int_{D_{j_2}}\int_0^1 \int_0^1 (1-t) (y-z_{j_2})^\perp \underset{x}{\nabla}(\underset{y}{\nabla}\underset{y}{\nabla}\Phi_k)(z_{m_1}+s(x-z_{m_1}), z_{j_2}+t(y-z_{j_2}))\notag\\
    &&\qquad\qquad\qquad\qquad\qquad \cdot (y-z_{j_2})\,dt\cdot\mathcal{P}(x, z_{m_1})\,ds\cdot\overset{3}{\PP}(E_{j_2})(y)\,dy.
\end{eqnarray}
By taking \eqref{inner p-P1Em2}, \eqref{inner p-P3Em2} and \eqref{inner p-P1j1j2} into \eqref{LS-Dm1-P11-inner-first-simple}, it gives us 
\begin{eqnarray}\label{LS-Dm1-P11-inner-after taylor}
    &&\left\langle \TT_{k, 1}\left(\overset{1}{\PP}(E_{m_1})\right), e_n^{(1)}\right\rangle - \left\langle \overset{3}{\PP}(E_{m_1}) , r_{k, 1}\right\rangle\notag\\
    &-& k^2 \eta_2 \left\langle \underset{x}{\nabla}\underset{y}{\nabla}\Phi_k(z_{m_1}, z_{m_2})(x-z_{m_1})\int_{D_{m_2}}\mathcal{P}(y, z_{m_2})\overset{1}{\PP}(E_{m_2})(y)\,dy, e_n^{(1)} \right\rangle\notag\\
    &-& k^2 \eta_2 \left\langle \underset{x}{\nabla}\Phi_k(z_{m_1}, z_{m_2})\mathcal{P}(x, z_{m_1})\int_{D_{m_2}}\overset{3}{\PP}(E_{m_2})(y)\,dy, e_n^{(1)} \right\rangle\notag\\
    &-& k^2 \eta_1 \sum_{j=1 \atop j\neq m}^\aleph \left\langle \underset{x}{\nabla}\underset{y}{\nabla}\Phi_k(z_{m_1}, z_{j_1})(x-z_{m_1})\int_{D_{j_1}}\mathcal{P}(y, z_{j_1})\overset{1}{\PP}(E_{j_1})(y)\,dy, e_n^{(1)} \right\rangle\notag\\
    &-& k^2 \eta_1 \sum_{j=1 \atop j\neq m}^\aleph \left\langle \underset{x}{\nabla}\Phi_k(z_{m_1}, z_{j_1})\mathcal{P}(x, z_{m_1})\int_{D_{j_1}}\overset{3}{\PP}(E_{j_1})(y)\,dy, e_n^{(1)} \right\rangle\notag\\
    &-& k^2 \eta_2 \sum_{j=1 \atop j\neq m}^\aleph \left\langle \underset{x}{\nabla}\underset{y}{\nabla}\Phi_k(z_{m_1}, z_{j_2})(x-z_{m_1})\int_{D_{j_2}}\mathcal{P}(y, z_{j_2})\overset{1}{\PP}(E_{j_2})(y)\,dy, e_n^{(1)} \right\rangle\notag\\
    &-& k^2 \eta_2 \sum_{j=1 \atop j\neq m}^\aleph \left\langle \underset{x}{\nabla}\Phi_k(z_{m_1}, z_{j_2})\mathcal{P}(x, z_{m_1})\int_{D_{j_2}}\overset{3}{\PP}(E_{j_2})(y)\,dy, e_n^{(1)}\right\rangle\notag\\
    &=& \left\langle E_{m_1}^{Inc}, e_n^{(1)} \right\rangle + k^2 \eta_2 \left\langle Err_1\left( \overset{1}{\PP}(E_{m_2}) \right), e_n^{(1)} \right\rangle + k^2 \eta_2 \left\langle Err_1\left( \overset{3}{\PP}(E_{m_2}) \right), e_n^{(1)} \right\rangle\notag\\
    &+& k^2 \eta_1 \sum_{j=1\atop j\neq m}^\aleph \left\langle Err_1\left( \overset{1}{\PP}(E_{j_1}) \right), e_n^{(1)} \right\rangle + k^2 \eta_1 \sum_{j=1\atop j\neq m}^\aleph \left\langle Err_1\left( \overset{3}{\PP}(E_{j_1}) \right), e_n^{(1)} \right\rangle\notag\\
    &+& k^2 \eta_2 \sum_{j=1\atop j\neq m}^\aleph \left\langle Err_1\left( \overset{1}{\PP}(E_{j_2}) \right), e_n^{(1)} \right\rangle + k^2 \eta_2 \sum_{j=1\atop j\neq m}^\aleph \left\langle Err_1\left( \overset{3}{\PP}(E_{j_2}) \right), e_n^{(1)} \right\rangle
\end{eqnarray}
After scaling \eqref{LS-Dm1-P11-inner-after taylor} to $B_1$ and $B_2$, there holds
\begin{eqnarray}\label{LS-Dm1-P11-after scaling}
    \left\langle \TT_{ka, 1}\left(\overset{1}{\PP}(\tilde{E}_{m_1})\right), e_n^{(1)}\right\rangle &=& \left\langle \overset{3}{\PP}(\tilde{E}_{m_1}) , \tilde{r}_{k, 1}\right\rangle + \left\langle \tilde{E}_{m_1}^{Inc}, e_n^{(1)} \right\rangle\notag\\
    &+& k^2 \eta_2 a^5 \left\langle \underset{x}{\nabla}\underset{y}{\nabla}\Phi_k(z_{m_1}, z_{m_2})\cdot x\cdot\int_{B_2}\mathcal{P}(y, 0)\overset{1}{\PP}(\tilde{E}_{m_2})(y)\,dy, e_n^{(1)} \right\rangle\notag\\
    &+& k^2 \eta_2 a^4 \left\langle \underset{x}{\nabla}\Phi_k(z_{m_1}, z_{m_2})\mathcal{P}(x, 0)\int_{B_2}\overset{3}{\PP}(\tilde{E}_{m_2})(y)\,dy, e_n^{(1)} \right\rangle\notag\\
    &+& k^2 \eta_1 a^5 \sum_{j=1 \atop j\neq m}^\aleph \left\langle \underset{x}{\nabla}\underset{y}{\nabla}\Phi_k(z_{m_1}, z_{j_1})\cdot x\cdot\int_{B_1}\mathcal{P}(y, 0)\overset{1}{\PP}(\tilde{E}_{j_1})(y)\,dy, e_n^{(1)} \right\rangle\notag\\
    &+& k^2 \eta_1 a^4 \sum_{j=1 \atop j\neq m}^\aleph \left\langle \underset{x}{\nabla}\Phi_k(z_{m_1}, z_{j_1})\mathcal{P}(x, 0)\int_{B_1}\overset{3}{\PP}(\tilde{E}_{j_1})(y)\,dy, e_n^{(1)} \right\rangle\notag\\
    &+& k^2 \eta_2 a^5 \sum_{j=1 \atop j\neq m}^\aleph \left\langle \underset{x}{\nabla}\underset{y}{\nabla}\Phi_k(z_{m_1}, z_{j_2})\cdot x\cdot\int_{B_2}\mathcal{P}(y, 0)\overset{1}{\PP}(\tilde{E}_{j_2})(y)\,dy, e_n^{(1)} \right\rangle\notag\\
    &+& k^2 \eta_2 a^4 \sum_{j=1 \atop j\neq m}^\aleph \left\langle \underset{x}{\nabla}\Phi_k(z_{m_1}, z_{j_2})\mathcal{P}(x, 0)\int_{B_2}\overset{3}{\PP}(\tilde{E}_{j_2})(y)\,dy, e_n^{(1)}\right\rangle\notag\\
    &+& k^2 \eta_2 \left\langle \widetilde{Err}_1\left( \overset{1}{\PP}(E_{m_2}) \right), e_n^{(1)} \right\rangle + k^2 \eta_2 \left\langle \widetilde{Err}_1\left( \overset{3}{\PP}(E_{m_2}) \right), e_n^{(1)} \right\rangle\notag\\
    &+& k^2 \eta_1 \sum_{j=1\atop j\neq m}^\aleph \left\langle \widetilde{Err}_1\left( \overset{1}{\PP}(E_{j_1}) \right), e_n^{(1)} \right\rangle + k^2 \eta_1 \sum_{j=1\atop j\neq m}^\aleph \left\langle \widetilde{Err}_1\left( \overset{3}{\PP}(E_{j_1}) \right), e_n^{(1)} \right\rangle\notag\\
    &+& k^2 \eta_2 \sum_{j=1\atop j\neq m}^\aleph \left\langle \widetilde{Err}_1\left( \overset{1}{\PP}(E_{j_2}) \right), e_n^{(1)} \right\rangle + k^2 \eta_2 \sum_{j=1\atop j\neq m}^\aleph \left\langle \widetilde{Err}_1\left( \overset{3}{\PP}(E_{j_2}) \right), e_n^{(1)} \right\rangle.
\end{eqnarray}
Due to the fact that 
\begin{equation}\label{adjoint-T_{ka,1}}
    \left\langle \TT_{ka, 1}\left(\overset{1}{\PP}(\tilde{E}_{m_1})\right), e_n^{(1)} \right\rangle= \left\langle \overset{1}{\PP}(\tilde{E}_{m_1}),  \TT_{-ka, 1}(e_n^{(1)}) \right\rangle,
\end{equation}
where by (MP1.16)
\begin{equation}\label{def-T_{-ka,1}}
    \TT_{-ka, 1}(e_{n}^{(1)})=(1-k^2 \eta_1 a^2 \lambda_n^{(1)})e_n^{(1)}-\frac{1\mp c_0 a^h}{4\pi \lambda_{n_0}^{(1)}}\sum_{\ell\geq 1}\frac{(-ika)^{\ell+1}}{(\ell+1)!}\int_{B_1}\lVert\cdot -y\rVert^\ell e_n^{(1)}(y)\,dy,
\end{equation}
together with the expression \eqref{adjoint-T_{ka,1}}, we can derive from \eqref{LS-Dm1-P11-after scaling} by taking the square modulus on the both sides and summing w.r.t $n$ that
\begin{eqnarray}\label{eq-P11-norm}
    \left\lVert\overset{1}{\PP}(\tilde{E}_{m_1})\right\rVert^2_{\LL^2({B_1})}&\lesssim& \frac{1}{|1-k^2\eta_1 a^2 \lambda_n^{(1)}|^2}\Bigg[ a^{10} \left\lVert \underset{x}{\nabla}\underset{y}{\nabla}\Phi_k(z_{m_1}, z_{m_2})\cdot x\cdot \int_{B_2}\mathcal{P}(y, 0)\overset{1}{\PP}(\tilde{E}_{m_2})(y)\,dy \right\rVert_{\LL^2(B_1)}^2\notag\\
    &+& a^8 \left\lVert \underset{x}{\nabla}\Phi_k(z_{m_1}, z_{m_2})\mathcal{P}(x, 0)\int_{B_2} \overset{3}{\PP}(\tilde{E}_{m_2})(y)\,dy \right\rVert_{\LL^2(B_1)}^2 \notag\\
    &+& a^6\sum_{j=1 \atop j\neq m}^\aleph  \left\lVert \underset{x}{\nabla}\underset{y}{\nabla}\Phi_k(z_{m_1}, z_{j_1})\cdot x\cdot \int_{B_1}\mathcal{P}(y, 0)\overset{1}{\PP}(\tilde{E}_{j_1})(y)\,dy \right\rVert_{\LL^2(B_1)}^2\notag\\
    &+& a^6 \sum_n \sum_{i=1 \atop i\neq m}^\aleph\left| \left\langle\underset{x}{\nabla}\underset{y}{\nabla}\Phi_k(z_{m_1}, z_{i_1})\cdot x\cdot \int_{B_1}\mathcal{P}(y, 0)\overset{1}{\PP}(\tilde{E}_{i_1})(y)\,dy, e_n^{(1)}\right\rangle \right|\notag\\
    &&\qquad\qquad\qquad\qquad\cdot \sum_{\ell>i\atop \ell\neq m}^\aleph \left|\left\langle \underset{x}{\nabla}\underset{y}{\nabla}\Phi_k(z_{m_1}, z_{\ell_1})\cdot x\cdot \int_{B_1}\mathcal{P}(y, 0)\overset{1}{\PP}(\tilde{E}_{\ell_1})(y)\,dy, e_n^{(1)} \right\rangle\right|\notag\\
    &+& a^4 \sum_{j=1 \atop j\neq m}^\aleph \left\lVert \underset{x}{\nabla}\Phi_k(z_{m_1}, z_{j_1})\mathcal{P}(x, 0)\int_{B_1} \overset{3}{\PP}(\tilde{E}_{j_1})(y)\,dy \right\rVert_{\LL^2(B_1)}^2\notag\\
    &+& a^4 \sum_n \sum_{i=1 \atop i\neq m}^\aleph \left|\left\langle \underset{x}{\nabla}\Phi_k(z_{m_1}, z_{i_1})\mathcal{P}(x, 0)\int_{B_1} \overset{3}{\PP}(\tilde{E}_{i_1})(y)\,dy, e_n^{(1)}\right\rangle \right|\notag\\
    &&\qquad\qquad\qquad\qquad\cdot \sum_{\ell>i\atop \ell\neq m}^\aleph \left|\left\langle \underset{x}{\nabla}\Phi_k(z_{m_1}, z_{\ell_1})\mathcal{P}(x, 0)\int_{B_1} \overset{3}{\PP}(\tilde{E}_{\ell_1})(y)\,dy, e_n^{(1)} \right\rangle\right|\notag\\
    &+& a^{10}\sum_{j=1 \atop j\neq m}^\aleph  \left\lVert \underset{x}{\nabla}\underset{y}{\nabla}\Phi_k(z_{m_1}, z_{j_2})\cdot x\cdot \int_{B_2}\mathcal{P}(y, 0)\overset{1}{\PP}(\tilde{E}_{j_2})(y)\,dy \right\rVert_{\LL^2(B_1)}^2\notag\\
    &+& a^{10} \sum_n \sum_{i=1 \atop i\neq m}^\aleph\left| \left\langle\underset{x}{\nabla}\underset{y}{\nabla}\Phi_k(z_{m_1}, z_{i_2})\cdot x\cdot \int_{B_2}\mathcal{P}(y, 0)\overset{1}{\PP}(\tilde{E}_{i_2})(y)\,dy, e_n^{(1)}\right\rangle \right|\notag\\
    &&\qquad\qquad\qquad\qquad\cdot \sum_{\ell>i\atop \ell\neq m}^\aleph \left|\left\langle \underset{x}{\nabla}\underset{y}{\nabla}\Phi_k(z_{m_1}, z_{\ell_2})\cdot x\cdot \int_{B_2}\mathcal{P}(y, 0)\overset{1}{\PP}(\tilde{E}_{\ell_2})(y)\,dy, e_n^{(1)} \right\rangle\right|\notag\\
     &+& a^8 \sum_{j=1 \atop j\neq m}^\aleph \left\lVert \underset{x}{\nabla}\Phi_k(z_{m_1}, z_{j_2})\mathcal{P}(x, 0)\int_{B_2} \overset{3}{\PP}(\tilde{E}_{j_2})(y)\,dy \right\rVert_{\LL^2(B_1)}^2\notag\\
    &+& a^8 \sum_n \sum_{i=1 \atop i\neq m}^\aleph \left|\left\langle \underset{x}{\nabla}\Phi_k(z_{m_1}, z_{i_2})\mathcal{P}(x, 0)\int_{B_2} \overset{3}{\PP}(\tilde{E}_{i_2})(y)\,dy, e_n^{(1)}\right\rangle \right|\notag\\
    &&\qquad\qquad\qquad\qquad\cdot \sum_{\ell>i\atop \ell\neq m}^\aleph \left|\left\langle \underset{x}{\nabla}\Phi_k(z_{m_1}, z_{\ell_2})\mathcal{P}(x, 0)\int_{B_2} \overset{3}{\PP}(\tilde{E}_{\ell_2})(y)\,dy, e_n^{(1)} \right\rangle\right|\notag\\
    &+& \sum_n\left| \left\langle \overset{1}{\PP}(\tilde{E}_{m_1}), \frac{1\mp c_0 a^h}{4\pi \lambda_{n_0}^{(1)}}\sum_{\ell\geq 1}\frac{(-ika)^{\ell+1}}{(\ell+1)!}\int_{B_1}\left\lVert\cdot - y\right\rVert^\ell e_n^{(1)}(y)\,dy \right\rangle \right|^2\notag\\
    &+& \sum_n \left|\left\langle \overset{3}{\PP}(\tilde{E}_{m_1}), \tilde{r}_{k, 1} \right\rangle\right|^2 + \sum_n \left|\left\langle \tilde{E}_{m_1}^{Inc}, e_n^{(1)} \right\rangle\right|^2 + \sum_n \left|\left\langle \widetilde{Err}_1\left(\overset{1}{\PP}(E_{m_2})\right), e_n^{(1)} \right\rangle\right|^2\notag\\
    &+& \sum_n\left| \left\langle \widetilde{Err}_1\left(\overset{3}{\PP}(E_{m_2})\right), e_n^{(1)} \right\rangle \right|^2 + a^{-4} \sum_n \sum_{j=1 \atop j\neq m}^\aleph \left| \left\langle \widetilde{Err}_1\left(\overset{1}{\PP}(E_{j_1})\right), e_n^{(1)} \right\rangle \right|^2\notag\\
    &+& a^{-4} \sum_n \sum_{i=1 \atop i\neq m}^\aleph \left| \left\langle \widetilde{Err}_1\left(\overset{1}{\PP}(E_{i_1})\right), e_n^{(1)} \right\rangle \right|\cdot \sum_{\ell>i \atop \ell\neq m}^\aleph \left| \left\langle \widetilde{Err}_1\left(\overset{1}{\PP}(E_{\ell_1})\right), e_n^{(1)} \right\rangle \right|\notag\\
    &+& a^{-4} \sum_n \sum_{j=1 \atop j\neq m}^\aleph \left| \left\langle \widetilde{Err}_1\left(\overset{3}{\PP}(E_{j_1})\right), e_n^{(1)} \right\rangle \right|^2\notag\\
    &+&a^{-4} \sum_n \sum_{i=1 \atop i\neq m}^\aleph \left| \left\langle \widetilde{Err}_1\left(\overset{3}{\PP}(E_{i_1})\right), e_n^{(1)} \right\rangle \right|\cdot \sum_{\ell>i \atop \ell\neq m}^\aleph \left| \left\langle \widetilde{Err}_1\left(\overset{3}{\PP}(E_{\ell_1})\right), e_n^{(1)} \right\rangle \right|\notag\\
    &+& \sum_n \sum_{j=1 \atop j\neq m}^\aleph \left| \left\langle \widetilde{Err}_1\left(\overset{1}{\PP}(E_{j_2})\right), e_n^{(1)} \right\rangle \right|^2 + \sum_n \sum_{j=1 \atop j\neq m}^\aleph \left| \left\langle \widetilde{Err}_1\left(\overset{3}{\PP}(E_{j_2})\right), e_n^{(1)} \right\rangle \right|^2\notag\\
    &+& \sum_n \sum_{i=1 \atop i\neq m}^\aleph \left| \left\langle \widetilde{Err}_1\left(\overset{1}{\PP}(E_{i_2})\right), e_n^{(1)} \right\rangle \right|\cdot \sum_{\ell>i \atop \ell\neq m}^\aleph \left| \left\langle \widetilde{Err}_1\left(\overset{1}{\PP}(E_{\ell_2})\right), e_n^{(1)} \right\rangle \right|\notag\\
    &+&\sum_n \sum_{i=1 \atop i\neq m}^\aleph \left| \left\langle \widetilde{Err}_1\left(\overset{3}{\PP}(E_{i_2})\right), e_n^{(1)} \right\rangle \right|\cdot \sum_{\ell>i \atop \ell\neq m}^\aleph \left| \left\langle \widetilde{Err}_1\left(\overset{3}{\PP}(E_{\ell_2})\right), e_n^{(1)} \right\rangle \right|\Bigg].
\end{eqnarray}
Now, we are in a position to evaluate the R.H.S. of \eqref{eq-P11-norm} respectively as follows.

\bigskip

\noindent (a). Estimate for the terms with respect to $\underset{x}{\nabla}\Phi_k$ and $\underset{x}{\nabla}\underset{y}{\nabla}\Phi_k$:

\bigskip

Without loss of generality, we evaluate the first four terms on the R.H.S. of \eqref{eq-P11-norm}, the other terms can be evaluated similarly. It is clear to see that
\begin{eqnarray*}
    &&a^{10} \left\lVert \underset{x}{\nabla} \underset{y}{\nabla}\Phi_k(z_{m_1}, z_{m_2})\cdot x \cdot \int_{B_2}\mathcal{P}(y, 0)\overset{1}{\PP}(\tilde{E}_{m_2})(y)\,dy \right\rVert_{\LL^2(B_1)}^2\lesssim a^{10} d_{\rm in}^{-6}\left\lVert \overset{1}{\PP}(\tilde{E}_{m_2})\right\rVert_{\LL^2(B_2)}^2,\notag\\
    &&a^8 \left\lVert \underset{x}{\nabla}\Phi_k(z_{m_1}, z_{m_2})\mathcal{P}(x, 0)\int_{B_2}\overset{3}{\PP}(\tilde{E}_{m_2})(u)\,dy \right\rVert_{\LL^2(B_1)}^2\lesssim a^8 d_{\rm in}^{-4} \left\lVert \overset{3}{\PP}(\tilde{E}_{m_2})\right\rVert_{\LL^2(B_2)}^2,
\end{eqnarray*}
and 
\begin{eqnarray}\label{es-grad_Phi-term}
    &&a^{6} \sum_{j=1 \atop j\neq m}^\aleph \left\lVert \underset{x}{\nabla} \underset{y}{\nabla}\Phi_k(z_{m_1}, z_{j_1})\cdot x \cdot \int_{B_1}\mathcal{P}(y, 0)\overset{1}{\PP}(\tilde{E}_{j_1})(y)\,dy \right\rVert_{\LL^2(B_1)}^2\notag\\
    &\lesssim& a^6 \sum_{j=1 \atop j\neq m}^\aleph d_{mj}^{-6}\left\lVert \overset{1}{\PP}(\tilde{E}_{j_1})\right\rVert_{\LL^2(B_1)}^2 \lesssim a^{6} d_{\rm out}^{-6}\max_j \left\lVert \overset{1}{\PP}(\tilde{E}_{j_1})\right\rVert_{\LL^2(B_1)}^2.
\end{eqnarray}
Besides, for the crossing term
\begin{eqnarray}\label{grad_Phi-term-crossing}
     I^*&:=& a^{6} \sum_n \sum_{i=1 \atop i\neq m}^\aleph \left|\left\langle \underset{x}{\nabla} \underset{y}{\nabla}\Phi_k(z_{m_1}, z_{i_1})\cdot x \cdot \int_{B_1}\mathcal{P}(y, 0)\overset{1}{\PP}(\tilde{E}_{i_1})(y)\,dy, e_n^{(1)} \right\rangle\right|\notag\\
    &&\qquad\qquad\qquad\qquad \cdot \; \sum_{\ell>i \atop \ell\neq m}^\aleph \left|\left\langle \underset{x}{\nabla} \underset{y}{\nabla}\Phi_k(z_{m_1}, z_{\ell_1})\cdot x \cdot \int_{B_1}\mathcal{P}(y, 0)\overset{1}{\PP}(\tilde{E}_{\ell_1})(y)\,dy, e_n^{(1)} \right\rangle\right|,
\end{eqnarray}
since
\begin{eqnarray*}
    \int_{B_1}\mathcal{P}(y, 0)\overset{1}{\PP}(\tilde{E}_{i_1})(y)\,dy&=& a^{-4} \int_{D_{i_1}} \mathcal{P}(y, z_{i_1}) \overset{1}{\PP}(E_{i_1})(y)\,dy = a^{-4} \int_{D_{i_1}} Curl \mathcal{P}(y, z_{i_1}) F_{i_1}(y)\,dy\notag\\
    &=&a^{-4} \mathcal{M}\cdot \int_{D_{i_1}} F_{i_1}(y)\,dy = \frac{a^{-2}}{\eta_0} \mathcal{M}\cdot (\eta_0 a^{-2}) \int_{D_{i_1}} F_{i_1}(y)\,dy \notag\\
    &=& \frac{a^{-2}}{\eta_0} \mathcal{M} \cdot \eta_1 \int_{D_{i_1}} F_{i_1}(y)\,dy  = \frac{a^{-2}}{\eta_0} \mathcal{M}\cdot \tilde{Q}_{i_1},
\end{eqnarray*}
where we use the relation (MP4.1) and 
\begin{equation*}
    \mathcal{M} := Curl \mathcal{P}(y, z_{i_1})= \left(	\begin{array}{ccccccccc}
		0 & 0 & 0 & 0 & 0 & 1 & 0 & 1 & 0 \\
		0 & 0 & 1 & 0 & 0 & 0 & -1 & 0 & 0 \\
		0 & -1 & 0 & 1 & 0 & 0 & 0 & 0 & 0
	\end{array}\right)^\mathsf{T}.
\end{equation*}
Then with the help of Proposition \ref{prop-es-QR-domain}, there holds
\begin{eqnarray}
    |I^*|&\lesssim& a^{6} \sum_{i=1 \atop i\neq m}^\aleph \left\lVert \underset{x}{\nabla} \underset{y}{\nabla}\Phi_k(z_{m_1}, z_{i_1})\cdot x \cdot \int_{B_1}\mathcal{P}(y, 0)\overset{1}{\PP}(\tilde{E}_{i_1})(y)\,dy\right\rVert_{\LL^2(B_1)}\notag\\
    &&\qquad\qquad\qquad\qquad \cdot \; \sum_{\ell>i \atop \ell\neq m}^\aleph \left\lVert \underset{x}{\nabla} \underset{y}{\nabla}\Phi_k(z_{m_1}, z_{\ell_1})\cdot x \cdot \int_{B_1}\mathcal{P}(y, 0)\overset{1}{\PP}(\tilde{E}_{\ell_1})(y)\,dy\right\rVert_{\LL^2(B_1)}\notag\\
    &\lesssim& a^6 \left( \sum_{i=1 \atop i\neq m}^\aleph d_{mi}^{-3}\cdot a^{-2} |\tilde{Q}_{i_1}| \right) \left( \sum_{\ell>i \atop \ell\neq m} d_{m\ell}^{-3}\cdot a^{-2} |\tilde{Q}_{\ell_1}| \right)\lesssim a^2 \left( \sum_{i=1\atop i\neq m}^\aleph d_{mi}^{-3} |\tilde{Q}_{i_1}|  \right)^2\lesssim a^2 d_{\rm out}^{-6} \sum_{i=1 \atop i\neq m}^\aleph |\tilde{Q}_{i_1}|^2\notag\\
    &\underset{\lesssim}{\eqref{es-Qm-main}}& a^{8-2h} d_{\rm out}^{-9} + \left( a^{14-2h} d_{\rm out}^{-17} + a^{12} d_{\rm out}^{-15} \ln^2 d_{\rm out} + a^{24-2h} d_{\rm in}^{-14} d_{\rm out}^{-9} \right)\max_m \left\lVert \overset{1}{\PP}(\tilde{E}_{m_1}) \right\rVert_{\LL^2(B_1)}^2\notag\\
    &+&\left( a^{18-2h} d_{\rm in}^{-8} d_{\rm out}^{-9} + a^{16} d_{\rm in}^{-6} d_{\rm out}^{-9} + a^{16} d_{\rm out}^{-15} \ln^2 d_{\rm out} \right)\max_m \left\lVert \overset{1}{\PP}(\tilde{E}_{m_2}) \right\rVert_{\LL^2(B_2)}^2 + \Big( a^6 d_{\rm out}^{-9} + a^{10} d_{\rm out}^{-15}\notag\\
    &+& a^{12-2h} d_{\rm out}^{-15} \ln^2 d_{\rm out} + a^{16-2h} d_{\rm in}^{-6} d_{\rm out}^{-15} + a^{18-4h} d_{\rm in}^{-8} d_{\rm out}^{-15} + a^{16-2h} d_{\rm out}^{-21} \ln^2 d_{\rm out} \Big)\max_m \left\lVert \overset{3}{\PP}(\tilde{E}_{m_1}) \right\rVert_{\LL^2(B_1)}\notag\\
    &+& \Big(a^{16-2h} d_{\rm in}^{-6} d_{\rm out}^{-9} + a^{14} d_{\rm in}^{-4} d_{\rm out}^{-9} + a^{14} d_{\rm out}^{-15} + a^{22-4h} d_{\rm out}^{-23} + a^{20-2h} d_{\rm out}^{-21} \ln^2 d_{\rm out}\notag\\
    &+& a^{28-4h} d_{\rm in}^{-14} d_{\rm out}^{-15} \Big)\max_m \left\lVert \overset{3}{\PP}(\tilde{E}_{m_2}) \right\rVert_{\LL^2(B_2)},\notag\\
\end{eqnarray}
by using Cauchy-Schwartz inequality and the counting lemma \cite[Lemma 3.4]{CS-JDE}.

\bigskip

\noindent (b) Estimate for 
\begin{eqnarray}\label{es-P11-remind-r}
    &&\sum_n\left| \left\langle \overset{1}{\PP}(\tilde{E}_{m_1}), \frac{1\mp c_0 a^h}{4\pi \lambda_{n_0}^{(1)}}\sum_{\ell\geq 1}\frac{(-ika)^{\ell+1}}{(\ell+1)!}\int_{B_1}\left\lVert\cdot - y\right\rVert^\ell e_n^{(1)}(y)\,dy \right\rangle \right|^2\notag\\
    &\lesssim& a^4 \sum_n\left|  \left\langle \overset{1}{\PP}(\tilde{E}_{m_1}), \sum_{\ell\geq 1}\int_{B_1}\frac{\left\lVert\cdot - y\right\rVert^\ell}{(\ell+1)!} e_n^{(1)}(y)\,dy \right\rangle \right|^2 \lesssim a^4 \sum_n \sum_{\ell\geq 1} \left|  \left\langle\int_{B_1} \overset{1}{\PP}(\tilde{E}_{m_1})(x)\frac{\left\lVert\cdot - x\right\rVert^\ell}{(\ell+1)!}\,dx, e_n^{(1)} \right\rangle \right|^2\notag\\
    &\lesssim& a^4 \sum_{\ell\geq 1}\left\lVert \int_{B_1} \overset{1}{\PP}(\tilde{E}_{m_1})(x)\frac{\left\lVert\cdot - x\right\rVert^\ell}{(\ell+1)!}\,dx \right\rVert_{\LL^2(B_1)}^2\lesssim a^4 \left\lVert\overset{1}{\PP}(\tilde{E}_{m_1})\right\rVert_{\LL^2(B_1)}^2.\notag\\
\end{eqnarray}
Similarly, we get
\begin{eqnarray}\label{es-P31-remind-r}
    \sum_n \left|\left\langle \overset{3}{\PP}(\tilde{E}_{m_1}), \tilde{r}_{k, 1} \right\rangle\right|^2&=&\sum_n  \left| \left\langle \overset{3}{\PP}(\tilde{E}_{m_1}), \frac{k^2 a^2 \eta_1}{4\pi} \sum_{\ell\geq 1}\frac{(-ika)^{\ell+1}}{(\ell+1)!}\int_{B_1}\left\lVert\cdot - y\right\rVert^\ell e_n^{(1)}(y)\,dy \right\rangle \right|^2\notag\\
    &\lesssim& a^4 \sum_{\ell\geq 1}\left\lVert \int_{B_1} \overset{3}{\PP}(\tilde{E}_{m_1})(x)\frac{\left\lVert\cdot - x\right\rVert^\ell}{(\ell+1)!}\,dx \right\rVert_{\LL^2(B_1)}^2\lesssim a^4 \left\lVert\overset{3}{\PP}(\tilde{E}_{m_1})\right\rVert_{\LL^2(B_1)}^2.
\end{eqnarray}

\bigskip

\noindent (c) Estimate for $\sum_n\left|\left\langle \tilde{E}_{m_1}^{Inc}, e_n^{(1)} \right\rangle\right|^2$:

\begin{equation}\label{es-E^inc_m1-P11}
    \sum_n\left|\left\langle \tilde{E}_{m_1}^{Inc}, e_n^{(1)} \right\rangle\right|^2\lesssim a^2 \sum_n \left|\left\langle \tilde{H}_{m_1}^{Inc}, \phi_n \right\rangle\right|^2 \lesssim a^2,
\end{equation}
with $Curl \phi_n=e_n^{(1)}$.

\bigskip

\noindent (d) Estimate for $\sum_n \left| \left\langle \widetilde{Err_1}\left(\overset{1}{\PP}(E_{m_2})\right), e_n^{(1)} \right\rangle \right|^2$:

By the definition for ${Err_1}\left(\overset{1}{\PP}(E_{m_2})\right)$ in \eqref{Error_1-P1m2}, we have 
\begin{eqnarray}\label{es-Err_1-P1m2}
    && \sum_n \left| \left\langle \widetilde{Err_1}\left(\overset{1}{\PP}(E_{m_2})\right), e_n^{(1)} \right\rangle \right|^2\notag\\
    &\lesssim& a^{12} \left\lVert \int_{B_2}\int_0^1 x^\perp \underset{x}{\nabla} (\underset{x}{\nabla}\underset{y}{\nabla}\Phi_k)(z_{m_1}+tax, z_{m_2})\cdot x\cdot \mathcal{P}(y, 0)\cdot\overset{1}{\PP}(\tilde{E}_{m_2})(y)\,dy \right\rVert_{\LL^2(B_1)}^2\notag\\
    &+& a^{12}\left\lVert \int_{B_2}\int_0^1\int_0^1 (1-t)y^\perp \underset{x}{\nabla}(\underset{y}{\nabla
    }\underset{y}{\nabla}\Phi_k)(z_{m_1}+sax, z_{m_2}+tay)\cdot y\,dt\cdot\mathcal{P}(x, 0)\,ds\cdot\overset{1}{\PP}(\tilde{E}_{m_2})(y)\,dy\right\rVert_{\LL^2(B_1)}^2\notag\\
    &\lesssim& a^{12} d_{\rm in}^{-8} \left\lVert\overset{1}{\PP}(\tilde{E}_{m_2})\right\rVert_{\LL^2(B_2)}^2.\notag\\
\end{eqnarray}

% \bigskip

\noindent (e) Estimate for $\sum_n \left| \left\langle \widetilde{Err_1}\left(\overset{3}{\PP}(E_{m_2})\right), e_n^{(1)} \right\rangle \right|^2$:

From \eqref{Error_1-P3m2}, we can know that
\begin{eqnarray}\label{es-Err_1-P3m2}
    &&\sum_n \left| \left\langle \widetilde{Err_1}\left(\overset{3}{\PP}(E_{m_2})\right), e_n^{(1)} \right\rangle \right|^2\notag\\
    &\lesssim& a^{10} \left\lVert \int_{B_2}\int_0^1 x^\perp \underset{x}{\nabla}\underset{x}{\nabla}\Phi_k(z_{m_1}+tax, z_{m_2})\cdot x\,dt\cdot\overset{3}{\PP}(\tilde{E}_{m_2})(y)\,dy \right\rVert_{\LL^2(B_1)}^2\notag\\
    &+& a^{10} \left\lVert \int_{B_2}\int_0^1 \underset{x}{\nabla}\underset{y}{\nabla}\Phi_k(z_{m_1}+tax, z_{m_2})\cdot\mathcal{P}(y, 0)\cdot x\,dt\cdot \overset{3}{\PP}(\tilde{E}_{m_2})(y)\,dy\right\rVert_{\LL^2(B_1)}^2\notag\\
    &+& a^{12} \left\lVert \int_{B_2}\int_0^1\int_0^1(1-t) y^\perp\underset{x}{\nabla}(\underset{y}{\nabla}\underset{y}{\nabla}\Phi_k)(z_{m_1}+sax, z_{m_2}+tay)\cdot y\,dt\cdot \mathcal{P}(x, 0)\,ds\cdot\overset{3}{\PP}(\tilde{E}_{m_2})(y)\,dy \right\rVert_{\LL^2(B_1)}^2\notag\\
    &\lesssim& a^{10} d_{\rm in}^{-6} \left\lVert \overset{3}{\PP}(\tilde{E}_{m_2})\right\rVert_{\LL^2(B_2)}^2+a^{12} d_{\rm in}^{-8} \left\lVert \overset{3}{\PP}(\tilde{E}_{m_2})\right\rVert_{\LL^2(B_2)}^2.\notag\\
\end{eqnarray}

% \bigskip

\noindent (f) Estimate for $a^{-4}\sum_n\sum_{j=1 \atop j\neq m}^\aleph \left| \left\langle \widetilde{Err_1}\left(\overset{1}{\PP}(E_{j_1})\right), e_n^{(1)} \right\rangle \right|^2$:

\begin{eqnarray}\label{es-Err_1-P1j1}
    && a^{-4}\sum_n\sum_{j=1 \atop j\neq m}^\aleph \left| \left\langle \widetilde{Err_1}\left(\overset{1}{\PP}(E_{j_1})\right), e_n^{(1)} \right\rangle \right|^2\notag\\
     &\lesssim& a^{8} \sum_{j=1\atop j\neq m}^\aleph \left\lVert \int_{B_1}\int_0^1 x^\perp \underset{x}{\nabla} (\underset{x}{\nabla}\underset{y}{\nabla}\Phi_k)(z_{m_1}+tax, z_{j_1})\cdot x\cdot \mathcal{P}(y, 0)\cdot\overset{1}{\PP}(\tilde{E}_{j_1})(y)\,dy \right\rVert_{\LL^2(B_1)}^2\notag\\
    &+& a^{8} \sum_{j=1\atop j\neq m}^\aleph \left\lVert \int_{B_1}\int_0^1\int_0^1 (1-t)y^\perp \underset{x}{\nabla}(\underset{y}{\nabla
    }\underset{y}{\nabla}\Phi_k)(z_{m_1}+sax, z_{j_1}+tay)\cdot y\,dt\cdot\mathcal{P}(x, 0)\,ds\cdot\overset{1}{\PP}(\tilde{E}_{j_1})(y)\,dy\right\rVert_{\LL^2(B_1)}^2\notag\\
    &\lesssim& a^{8}\sum_{j=1\atop j\neq m}^\aleph d_{mj}^{-8} \left\lVert\overset{1}{\PP}(\tilde{E}_{j_1})\right\rVert_{\LL^2(B_2)}^2 \lesssim a^8 d_{\rm out}^{-8} \max_j \left\lVert\overset{1}{\PP}(\tilde{E}_{j_1})\right\rVert_{\LL^2(B_1)}^2.\notag\\
\end{eqnarray}

% \bigskip

\noindent (g) Estimate for 
\begin{eqnarray}\label{es-Err_1-crossing-P11}
    && a^{-4}\sum_n\sum_{i=1 \atop i\neq m}^\aleph \left| \left\langle \widetilde{Err_1}\left(\overset{1}{\PP}(E_{i_1})\right), e_n^{(1)} \right\rangle \right|\cdot\sum_{\ell>i \atop \ell\neq m}^\aleph \left| \left\langle \widetilde{Err_1}\left(\overset{1}{\PP}(E_{\ell_1})\right), e_n^{(1)} \right\rangle \right| \notag\\
    &\lesssim& a^{-4} \sum_{i=1\atop i\neq m}^\aleph \left\lVert \widetilde{Err}_1 \left(\overset{1}{\PP}(E_{i_1})\right)\right\rVert_{\LL^2(B_1)}\cdot \sum_{\ell>i \atop \ell\neq m}^\aleph \left\lVert \widetilde{Err}_1 \left(\overset{1}{\PP}(E_{\ell_1})\right)\right\rVert_{\LL^2(B_1)}\notag\\
    &\underset{\lesssim}{\eqref{Err_1-P1j1j2}}& a^{-4} \left( \sum_{i=1\atop i\neq m}^\aleph a^6 d_{mi}^{-4} \left\lVert \overset{1}{\PP}(\tilde{E}_{i_1})\right\rVert_{\LL^2(B_1)} \right)\cdot \left( \sum_{\ell> i \atop \ell \neq m}^\aleph a^6 d_{m\ell}^{-4} \left\lVert \overset{1}{\PP}(\tilde{E}_{\ell_1})\right\rVert_{\LL^2(B_1)} \right)\notag\\
    &\lesssim& a^8 \left( \sum_{i=1 \atop i\neq m}^\aleph d_{mi}^{-4} \left\lVert \overset{1}{\PP}(\tilde{E}_{i_1})\right\rVert_{\LL^2(B_1)}  \right)^2 \lesssim a^8 d_{\rm out}^{-8} \max_j \left\lVert \overset{1}{\PP}(\tilde{E}_{j_1})\right\rVert_{\LL^2(B_1)}^2.
\end{eqnarray}

% \bigskip

\noindent (h) Estimate for $a^{-4}\sum_n\sum_{j=1 \atop j\neq m}^\aleph \left| \left\langle \widetilde{Err_1}\left(\overset{3}{\PP}(E_{j_1})\right), e_n^{(1)} \right\rangle \right|^2$:

\begin{eqnarray}\label{es-Err_1-P3j1}
    && a^{-4}\sum_n\sum_{j=1 \atop j\neq m}^\aleph \left| \left\langle \widetilde{Err_1}\left(\overset{3}{\PP}(E_{j_1})\right), e_n^{(1)} \right\rangle \right|^2\notag\\
    &\lesssim& a^{6} \sum_{j=1 \atop j\neq m}^\aleph \left\lVert \int_{B_1}\int_0^1 x^\perp \underset{x}{\nabla}\underset{x}{\nabla}\Phi_k(z_{m_1}+tax, z_{j_1})\cdot x\,dt\cdot\overset{3}{\PP}(\tilde{E}_{j_1})(y)\,dy \right\rVert_{\LL^2(B_1)}^2\notag\\
    &+& a^{6}\sum_{j=1 \atop j\neq m}^\aleph \left\lVert \int_{B_1}\int_0^1 \underset{x}{\nabla}\underset{y}{\nabla}\Phi_k(z_{m_1}+tax, z_{j_1})\cdot\mathcal{P}(y, 0)\cdot x\,dt\cdot \overset{3}{\PP}(\tilde{E}_{j_1})(y)\,dy\right\rVert_{\LL^2(B_1)}^2\notag\\
    &+& a^{8}\sum_{j=1 \atop j\neq m}^\aleph \left\lVert \int_{B_1}\int_0^1\int_0^1(1-t) y^\perp\underset{x}{\nabla}(\underset{y}{\nabla}\underset{y}{\nabla}\Phi_k)(z_{m_1}+sax, z_{j_1}+tay)\cdot y\,dt\cdot \mathcal{P}(x, 0)\,ds\cdot\overset{3}{\PP}(\tilde{E}_{j_1})(y)\,dy \right\rVert_{\LL^2(B_1)}^2\notag\\
    &\lesssim& a^{6}\sum_{j=1 \atop j\neq m}^\aleph d_{mj}^{-6} \left\lVert \overset{3}{\PP}(\tilde{E}_{j_1})\right\rVert_{\LL^2(B_1)}^2+a^{8}\sum_{j=1 \atop j\neq m}^\aleph d_{mj}^{-8} \left\lVert \overset{3}{\PP}(\tilde{E}_{j_1})\right\rVert_{\LL^2(B_1)}^2 \lesssim a^6 d_{\rm out}^{-6} \max_j\left\lVert \overset{3}{\PP}(\tilde{E}_{j_1})\right\rVert_{\LL^2(B_1)}^2,\notag\\
\end{eqnarray}
by keeping the dominant term.

\medskip

\noindent (i) Estimate for 
\begin{eqnarray}\label{es-Err_1-crossing-P31}
     && a^{-4}\sum_n\sum_{i=1 \atop i\neq m}^\aleph \left| \left\langle \widetilde{Err_1}\left(\overset{3}{\PP}(E_{i_1})\right), e_n^{(1)} \right\rangle \right|\cdot\sum_{\ell>i \atop \ell\neq m}^\aleph \left| \left\langle \widetilde{Err_1}\left(\overset{3}{\PP}(E_{\ell_1})\right), e_n^{(1)} \right\rangle \right| \notag\\
    &\lesssim& a^{-4} \sum_{i=1\atop i\neq m}^\aleph \left\lVert \widetilde{Err}_1 \left(\overset{3}{\PP}(E_{i_1})\right)\right\rVert_{\LL^2(B_1)}\cdot \sum_{\ell>i \atop \ell\neq m}^\aleph \left\lVert \widetilde{Err}_1 \left(\overset{3}{\PP}(E_{\ell_1})\right)\right\rVert_{\LL^2(B_1)}\notag\\
    &\lesssim& a^{-4} \left( a^5 \sum_{i=1 \atop i\neq m}^\aleph d_{mi}^{-3} \left\lVert \overset{3}{\PP}(\tilde{E}_{i_1})\right\rVert_{\LL^2(B_1)} + a^6 \sum_{i=1 \atop i\neq m}^\aleph d_{mi}^{-4} \left\lVert \overset{3}{\PP}(\tilde{E}_{i_1})\right\rVert_{\LL^2(B_1)}\right)\notag\\
    &&\qquad\qquad\qquad\qquad\cdot \left( a^5 \sum_{\ell>i \atop \ell\neq m}^\aleph d_{m\ell}^{-3} \left\lVert \overset{3}{\PP}(\tilde{E}_{\ell_1})\right\rVert_{\LL^2(B_1)} + a^6 \sum_{\ell>i \atop \ell\neq m}^\aleph d_{m\ell}^{-4} \left\lVert \overset{3}{\PP}(\tilde{E}_{\ell_1})\right\rVert_{\LL^2(B_1)}\right)\notag\\
    &\lesssim& a^{-4} \left( a^5 d_{\rm out}^{-3} |\ln (d_{\rm out})| \max_j \left\lVert \overset{3}{\PP}(\tilde{E}_{j_1})\right\rVert_{\LL^2(B_1)} + a^6 d_{\rm out}^{-4}\max_j \left\lVert \overset{3}{\PP}(\tilde{E}_{j_1})\right\rVert_{\LL^2(B_1)}\right)^2\notag\\
    &\lesssim& a^6 d_{\rm out}^{-6} \ln^2 (d_{\rm out}) \max_j \left\lVert \overset{3}{\PP}(\tilde{E}_{j_1})\right\rVert_{\LL^2(B_1)}^2 + a^8 d_{\rm out}^{-8} \max_j \left\lVert \overset{3}{\PP}(\tilde{E}_{j_1})\right\rVert_{\LL^2(B_1)}^2.
\end{eqnarray}

% \bigskip

% \medskip

\noindent (j) Estimate for $\sum_n\sum_{j=1 \atop j\neq m}^\aleph \left| \left\langle \widetilde{Err_1}\left(\overset{1}{\PP}(E_{j_2})\right), e_n^{(1)} \right\rangle \right|^2$:

\begin{eqnarray}\label{es-Err_1-P1j2}
    && \sum_n\sum_{j=1 \atop j\neq m}^\aleph \left| \left\langle \widetilde{Err_1}\left(\overset{1}{\PP}(E_{j_2})\right), e_n^{(1)} \right\rangle \right|^2\notag\\
     &\underset{\lesssim}{\eqref{Err_1-P1j1j2}}& a^{12} \sum_{j=1\atop j\neq m}^\aleph \left\lVert \int_{B_2}\int_0^1 x^\perp \underset{x}{\nabla} (\underset{x}{\nabla}\underset{y}{\nabla}\Phi_k)(z_{m_1}+tax, z_{j_2})\cdot x\,dt\cdot \mathcal{P}(y, 0)\cdot\overset{1}{\PP}(\tilde{E}_{j_2})(y)\,dy \right\rVert_{\LL^2(B_1)}^2\notag\\
    &+& a^{12} \sum_{j=1\atop j\neq m}^\aleph \left\lVert \int_{B_2}\int_0^1\int_0^1 (1-t)y^\perp \underset{x}{\nabla}(\underset{y}{\nabla
    }\underset{y}{\nabla}\Phi_k)(z_{m_1}+sax, z_{j_2}+tay)\cdot y\,dt\cdot\mathcal{P}(x, 0)\,ds\cdot\overset{1}{\PP}(\tilde{E}_{j_2})(y)\,dy\right\rVert_{\LL^2(B_1)}^2\notag\\
    &\lesssim& a^{12}\sum_{j=1\atop j\neq m}^\aleph d_{mj}^{-8} \left\lVert\overset{1}{\PP}(\tilde{E}_{j_2})\right\rVert_{\LL^2(B_2)}^2 \lesssim a^{12} d_{\rm out}^{-8} \max_j \left\lVert\overset{1}{\PP}(\tilde{E}_{j_2})\right\rVert_{\LL^2(B_2)}^2.\notag\\
\end{eqnarray}

% \bigskip
% \medskip

\noindent (k) Estimate for 
\begin{eqnarray}\label{es-Err_1-crossing-P12}
     && \sum_n\sum_{i=1 \atop i\neq m}^\aleph \left| \left\langle \widetilde{Err_1}\left(\overset{1}{\PP}(E_{i_2})\right), e_n^{(1)} \right\rangle \right|\cdot\sum_{\ell>i \atop \ell\neq m}^\aleph \left| \left\langle \widetilde{Err_1}\left(\overset{1}{\PP}(E_{\ell_2})\right), e_n^{(1)} \right\rangle \right| \notag\\
    &\lesssim& \sum_{i=1\atop i\neq m}^\aleph \left\lVert \widetilde{Err}_1 \left(\overset{1}{\PP}(E_{i_2})\right)\right\rVert_{\LL^2(B_1)}\cdot \sum_{\ell>i \atop \ell\neq m}^\aleph \left\lVert \widetilde{Err}_1 \left(\overset{1}{\PP}(E_{\ell_2})\right)\right\rVert_{\LL^2(B_1)}\notag\\
    &\lesssim& \left( a^6 \sum_{i=1 \atop i\neq m}^\aleph d_{mi}^{-4} \left\lVert \overset{1}{\PP}(\tilde{E}_{i_2})\right\rVert_{\LL^2(B_2)} \right)\cdot \left( a^6 \sum_{\ell>i \atop \ell\neq m}^\aleph d_{m\ell}^{-4} \left\lVert \overset{1}{\PP}(\tilde{E}_{\ell_2})\right\rVert_{\LL^2(B_2)}\right)\lesssim a^{12} d_{\rm out}^{-8} \max_j \left\lVert \overset{1}{\PP}(\tilde{E}_{j_2})\right\rVert_{\LL^2(B_2)}^2.\notag\\
\end{eqnarray}

\bigskip

\noindent (l) Estimate for $\sum_n\sum_{j=1 \atop j\neq m}^\aleph \left| \left\langle \widetilde{Err_1}\left(\overset{3}{\PP}(E_{j_2})\right), e_n^{(1)} \right\rangle \right|^2$:

\begin{eqnarray}\label{es-Err_1-P3j2}
    && \sum_n\sum_{j=1 \atop j\neq m}^\aleph \left| \left\langle \widetilde{Err_1}\left(\overset{3}{\PP}(E_{j_2})\right), e_n^{(1)} \right\rangle \right|^2\notag\\
    &\lesssim& a^{10} \sum_{j=1 \atop j\neq m}^\aleph \left\lVert \int_{B_2}\int_0^1 x^\perp \underset{x}{\nabla}\underset{x}{\nabla}\Phi_k(z_{m_1}+tax, z_{j_2})\cdot x\,dt\cdot\overset{3}{\PP}(\tilde{E}_{j_2})(y)\,dy \right\rVert_{\LL^2(B_1)}^2\notag\\
    &+& a^{10}\sum_{j=1 \atop j\neq m}^\aleph \left\lVert \int_{B_2}\int_0^1 \underset{x}{\nabla}\underset{y}{\nabla}\Phi_k(z_{m_1}+tax, z_{j_2})\cdot\mathcal{P}(y, 0)\cdot x\,dt\cdot \overset{3}{\PP}(\tilde{E}_{j_2})(y)\,dy\right\rVert_{\LL^2(B_1)}^2\notag\\
    &+& a^{12}\sum_{j=1 \atop j\neq m}^\aleph \left\lVert \int_{B_2}\int_0^1\int_0^1(1-t) y^\perp\underset{x}{\nabla}(\underset{y}{\nabla}\underset{y}{\nabla}\Phi_k)(z_{m_1}+sax, z_{j_2}+tay)\cdot y\,dt\cdot \mathcal{P}(x, 0)\,ds\cdot\overset{3}{\PP}(\tilde{E}_{j_2})(y)\,dy \right\rVert_{\LL^2(B_1)}^2\notag\\
    &\lesssim& a^{10}\sum_{j=1 \atop j\neq m}^\aleph d_{mj}^{-6} \left\lVert \overset{3}{\PP}(\tilde{E}_{j_2})\right\rVert_{\LL^2(B_2)}^2+a^{12}\sum_{j=1 \atop j\neq m}^\aleph d_{mj}^{-8} \left\lVert \overset{3}{\PP}(\tilde{E}_{j_2})\right\rVert_{\LL^2(B_2)}^2 \lesssim a^{10} d_{\rm out}^{-6} \max_j\left\lVert \overset{3}{\PP}(\tilde{E}_{j_2})\right\rVert_{\LL^2(B_2)}^2,\notag\\
\end{eqnarray}
by keeping the dominant term.

\medskip

\noindent (m) Estimate for 
\begin{eqnarray}\label{es-Err_1-crossing-P32}
     && \sum_n\sum_{i=1 \atop i\neq m}^\aleph \left| \left\langle \widetilde{Err_1}\left(\overset{3}{\PP}(E_{i_2})\right), e_n^{(1)} \right\rangle \right|\cdot\sum_{\ell>i \atop \ell\neq m}^\aleph \left| \left\langle \widetilde{Err_1}\left(\overset{3}{\PP}(E_{\ell_2})\right), e_n^{(1)} \right\rangle \right| \notag\\
    &\lesssim& \sum_{i=1\atop i\neq m}^\aleph \left\lVert \widetilde{Err}_1 \left(\overset{3}{\PP}(E_{i_2})\right)\right\rVert_{\LL^2(B_1)}\cdot \sum_{\ell>i \atop \ell\neq m}^\aleph \left\lVert \widetilde{Err}_1 \left(\overset{3}{\PP}(E_{\ell_2})\right)\right\rVert_{\LL^2(B_1)}\notag\\
    &\lesssim& \left( a^5 \sum_{i=1 \atop i\neq m}^\aleph d_{mi}^{-3} \left\lVert \overset{3}{\PP}(\tilde{E}_{i_2})\right\rVert_{\LL^2(B_2)} + a^6 \sum_{i=1 \atop i\neq m}^\aleph d_{mi}^{-4} \left\lVert \overset{3}{\PP}(\tilde{E}_{i_2})\right\rVert_{\LL^2(B_2)}\right)\notag\\
    &&\qquad\qquad\qquad\qquad\cdot \left( a^5 \sum_{\ell>i \atop \ell\neq m}^\aleph d_{m\ell}^{-3} \left\lVert \overset{3}{\PP}(\tilde{E}_{\ell_2})\right\rVert_{\LL^2(B_2)} + a^6 \sum_{\ell>i \atop \ell\neq m}^\aleph d_{m\ell}^{-4} \left\lVert \overset{3}{\PP}(\tilde{E}_{\ell_2})\right\rVert_{\LL^2(B_2)}\right)\notag\\
    &\lesssim& \left( a^5 d_{\rm out}^{-3} |\ln (d_{\rm out})| \max_j \left\lVert \overset{3}{\PP}(\tilde{E}_{j_2})\right\rVert_{\LL^2(B_2)} + a^6 d_{\rm out}^{-4}\max_j \left\lVert \overset{3}{\PP}(\tilde{E}_{j_2})\right\rVert_{\LL^2(B_2)}\right)^2\notag\\
    &\lesssim& a^{10} d_{\rm out}^{-6} \ln^2 d_{\rm out} \max_j \left\lVert \overset{3}{\PP}(\tilde{E}_{j_2})\right\rVert_{\LL^2(B_2)}^2 + a^{12} d_{\rm out}^{-8} \max_j \left\lVert \overset{3}{\PP}(\tilde{E}_{j_2})\right\rVert_{\LL^2(B_2)}^2.
\end{eqnarray}

Following the similar argument to \eqref{es-grad_Phi-term}, combining with the estimates \eqref{es-P11-remind-r}--\eqref{es-Err_1-crossing-P32} above, we shall generate from \eqref{eq-P11-norm} by taking $\max_m(\cdot)$ on the both sides that
\begin{eqnarray}\label{es-norm-P1m1}
    \max_{m} \left\lVert \overset{1}{\PP}(\tilde{E}_{m_1})\right\rVert_{\LL^2(B_1)}^2&\lesssim& a^{8-4h} d_{\rm out}^{-9} + a^{2-2h} \notag\\
    &+& \left( a^{10-2h} d_{\rm in}^{-6} + a^{10-2h} d_{\rm out}^{-6} \ln^2 d_{\rm out} + a^{18-4h} d_{\rm in}^{-8} d_{\rm out}^{-9} \right) \max_m \left\lVert\overset{1}{\PP}(\tilde{E}_{m_2})\right\rVert_{\LL^2(B_2)}^2\notag\\
    &+& \left( a^{4-2h} d_{\rm out}^{-6} + a^{16-4h} d_{\rm in}^{-6} d_{\rm out}^{-15}  + a^{18-6h} d_{\rm in}^{-8} d_{\rm out}^{-15} \right) \max_m \left\lVert\overset{3}{\PP}(\tilde{E}_{m_1})\right\rVert_{\LL^2(B_1)}^2\notag\\
    &+& \left( a^{8-2h} d_{\rm in}^{-4} + a^{8-2h} d_{\rm out}^{-6} \right) \max_m \left\lVert\overset{3}{\PP}(\tilde{E}_{m_2})\right\rVert_{\LL^2(B_2)}^2
\end{eqnarray}
by keeping the dominant terms under the condition that
\begin{equation}\label{cond-on-t2}
    \frac{9}{5}< h < \min\{ 2, 15-14t_1 \}.
\end{equation}

\bigskip

\noindent (2). \textit{Estimate for $\max_m \left\lVert \overset{3}{\PP}(\tilde{E}_{m_1}) \right\rVert^2_{\LL^2(B_1)}$}.

\bigskip

Recall the Lippmann-Schwinger equation \eqref{LS-Dm1}. With (MP3.15) and \eqref{M|div=0}, there holds
\begin{eqnarray*}
    && \left(I+\eta_1 \nabla M^k\right)\overset{3}{\PP}(E_{m_1})\notag\\
    &=& E_{m_1}^{Inc} - \overset{1}{\PP}(E_{m_1}) + k^2 \eta_1 N^k \left(\overset{1}{\PP}(E_{m_1})\right) + k^2 \eta_1 N^k \left(\overset{3}{\PP}(E_{m_1})\right) + k^2 \eta_2 \int_{D_{m_2}}\Phi_k(x, y)\overset{1}{\PP}(E_{m_2})(y)\,dy\notag\\
    &+& k^2 \eta_1 \sum_{j=1\atop j\neq m}^\aleph \int_{D_{j_1}}\Phi_k(x, y)\overset{1}{\PP}(E_{j_1})(y)\,dy + k^2 \eta_2 \sum_{j=1\atop j\neq m}^\aleph \int_{D_{j_2}}\Phi_k(x, y)\overset{1}{\PP}(E_{j_2})(y)\,dy + k^2 \eta_2 \int_{D_{m_2}} \Upsilon_k(x, y)\overset{3}{\PP}(E_{m_2})(y)\,dy\notag\\
    &+& k^2 \eta_1 \sum_{j=1 \atop j\neq m}^\aleph \int_{D_{j_1}} \Upsilon_k(x, y)\overset{3}{\PP}(E_{j_1})(y)\,dy + k^2 \eta_2 \sum_{j=1 \atop j\neq m}^\aleph \int_{D_{j_2}} \Upsilon_k(x, y)\overset{3}{\PP}(E_{j_2})(y)\,dy,
\end{eqnarray*}
where $\Upsilon_k(x, y)$ is the dyadic Green function given by (MP1.17). By taking the Taylor expansion of $\Phi_k(x, y)$ and $\Upsilon_k(x, y)$ at $(z_{m_1}, z_{m_2})$ and $(z_{m_1}, z_{j_i}), i=1,2$, respectively, we can get that
\begin{eqnarray}\label{Ls-Dm1-P31-taylor}
    && \left(I+\eta_1 \nabla M \right)\overset{3}{\PP}(E_{m_1})\notag\\
    &=& E_{m_1}^{Inc} - \overset{1}{\PP}(E_{m_1}) + k^2 \eta_1 N^k \left(\overset{1}{\PP}(E_{m_1})\right) + k^2 \eta_1 N^k \left(\overset{3}{\PP}(E_{m_1})\right) - \eta_1 \left( \nabla M^k -\nabla M \right)\overset{3}{\PP}(E_{m_1})\notag\\
    &+& k^2 \eta_2 \underset{y}{\nabla} (\Phi_kI)(z_{m_1}, z_{m_2}) \int_{D_{m_2}}\mathcal{P}(y, z_{m_2})\cdot\overset{1}{\PP}(E_{m_2})(y)\,dy + k^2 \eta_1 \sum_{j=1\atop j\neq m}^\aleph \underset{y}{\nabla}(\Phi_k I) (z_{m_1}, z_{j_1}) \int_{D_{j_1}}\mathcal{P}(y, z_{j_1}) \overset{1}{\PP}(E_{j_1})(y)\,dy\notag\\
    &+& k^2 \eta_2 \sum_{j=1\atop j\neq m}^\aleph \underset{y}{\nabla}(\Phi_k I)(z_{m_1}, z_{j_2}) \int_{D_{j_2}} \mathcal{P}(y, z_{j_2}) \overset{1}{\PP}(E_{j_2})(y)\,dy + k^2 \eta_2 \Upsilon_k(z_{m_1}, z_{m_2}) \int_{D_{m_2}} \overset{3}{\PP}(E_{m_2})(y)\,dy\notag\\
    &+& k^2 \eta_1 \sum_{j=1 \atop j\neq m}^\aleph \Upsilon_k(z_{m_1}, z_{j_1})\int_{D_{j_1}} \overset{3}{\PP}(E_{j_1})(y)\,dy + k^2 \eta_2 \sum_{j=1 \atop j\neq m}^\aleph\Upsilon_k(z_{m_1}, z_{j_2}) \int_{D_{j_2}} \overset{3}{\PP}(E_{j_2})(y)\,dy\notag\\
    &+& Err_2\left( \overset{1}{\PP}(E_{m_2}) \right) + Err_2\left( \overset{3}{\PP}(E_{m_2}) \right) + \sum_{j=1 \atop j\neq m}^\aleph Err_2\left( \overset{1}{\PP}(E_{j}) \right) + \sum_{j=1 \atop j\neq m}^\aleph Err_2\left( \overset{3}{\PP}(E_{j}) \right),\notag\\
\end{eqnarray}
where
\begin{eqnarray}\label{def-Err_2-P13mj}
    &&Err_2\left(\overset{1}{\PP}(E_{m_2})\right)\notag\\
     &=& \eta_2 k^2  \int_{D_{m_2}} \int_0^1 (1-t) (y-z_{m_2})^\perp \underset{y}{Hess}\Phi_k(z_{m_1}, z_{m_2}+t(y-z_{m_2}))\cdot(y-z_{m_2})\,dt\cdot\overset{1}{\PP}(E_{m_2})(y)\,dy\notag\\
     &-& \eta_2 k^2 \int_{D_{m_2}} \int_0^1  \underset{y}{Hess}\Phi_k(z_{m_1}+t(x-z_{m_1}), z_{m_2})\cdot(y-z_{m_2}) \cdot(x-z_{m_1})\,dt\cdot\overset{1}{\PP}(E_{m_2})(y)\,dy\notag\\
     &+& \eta_2 k^2 \int_{D_{m_2}} \int_0^1 \int_0^1 (1-t) (y-z_{m_2})^\perp \underset{x}{\nabla}(\underset{y}{Hess}\Phi_k)(z_{m_1}+t(x-z_{m_1}), z_{m_2}+s(y-z_{m_2}))\notag\\
     &&\qquad\qquad\qquad\qquad\cdot(y-z_{m_2})\,ds\cdot\mathcal{P}(x, z_{m_1})\,dt\cdot\overset{1}{\PP}(E_{m_2})(y)\,dy\notag\\
     &&Err_2\left(\overset{1}{\PP}(E_{j})\right)\notag\\
     &=& \eta_1 k^2 \int_{D_{j_1}} \int_0^1 (1-t) (y-z_{j_1})^\perp \underset{y}{Hess}\Phi_k(z_{m_1}, z_{j_1}+t(y-z_{j_1}))\cdot(y-z_{j_1})\,dt\cdot\overset{1}{\PP}(E_{j_1})(y)\,dy\notag\\
     &-& \eta_1 k^2 \int_{D_{j_1}} \int_0^1  \underset{y}{Hess}\Phi_k(z_{m_1}+t(x-z_{m_1}), z_{j_1})\cdot(y-z_{j_1})\cdot(x-z_{m_1})\,dt\cdot\overset{1}{\PP}(E_{j_1})(y)\,dy\notag\\
     &+& \eta_1 k^2 \int_{D_{j_1}} \int_0^1 \int_0^1 (1-t) (y-z_{j_1})^\perp \underset{x}{\nabla}(\underset{y}{Hess}\Phi_k)(z_{m_1}+t(x-z_{m_1}), z_{j_1}+s(y-z_{j_1}))\notag\\
     &&\qquad\qquad\qquad\qquad\cdot(y-z_{j_1})\,ds\cdot\mathcal{P}(x, z_{m_1})\,dt\cdot\overset{1}{\PP}(E_{j_1})(y)\,dy\notag\\
     &+& \eta_2 k^2 \int_{D_{j_2}} \int_0^1 (1-t) (y-z_{j_2})^\perp \underset{y}{Hess}\Phi_k(z_{m_1}, z_{j_2}+t(y-z_{j_2}))\cdot(y-z_{j_2})\,dt\cdot\overset{1}{\PP}(E_{j_2})(y)\,dy\notag\\
     &-& \eta_2 k^2 \int_{D_{j_2}} \int_0^1  \underset{y}{Hess}\Phi_k(z_{m_1}+t(x-z_{m_1}), z_{j_2})\cdot(y-z_{j_2})\cdot(x-z_{m_1})\,dt\cdot\overset{1}{\PP}(E_{j_2})(y)\,dy\notag\\
     &+& \eta_2 k^2 \int_{D_{j_2}} \int_0^1 \int_0^1 (1-t) (y-z_{j_2})^\perp \underset{x}{\nabla}(\underset{y}{Hess}\Phi_k)(z_{m_1}+t(x-z_{m_1}), z_{j_2}+s(y-z_{j_2}))\notag\\
     &&\qquad\qquad\qquad\qquad\cdot(y-z_{j_2})\,ds\cdot\mathcal{P}(x, z_{m_1})\,dt\cdot\overset{1}{\PP}(E_{j_2})(y)\,dy\notag\\
      &&Err_2\left(\overset{3}{\PP}(E_{m_2})\right)\notag\\
     &=& \eta_2 k^2  \int_{D_{m_2}} \int_0^1  \underset{x}{\nabla}\Upsilon_k(z_{m_1}+t(x-z_{m_1}), z_{m_2})\cdot\mathcal{P}(x, z_{m_1})\,dt\cdot\overset{3}{\PP}(E_{m_2})(y)\,dy\notag\\
     &+ & \eta_2 k^2 \int_{D_{m_2}} \int_0^1  \underset{y}{\nabla}\Upsilon_k(z_{m_1}, z_{m_2}+t(y-z_{m_2}))\cdot\mathcal{P}(y, z_{m_2})\,dt\cdot\overset{3}{\PP}(E_{m_2})(y)\,dy\notag\\
      &&Err_2\left(\overset{3}{\PP}(E_{j})\right)\notag\\
     &=& \eta_1 k^2  \int_{D_{j_1}} \int_0^1  \underset{x}{\nabla}\Upsilon_k(z_{m_1}+t(x-z_{m_1}), z_{j_1})\cdot\mathcal{P}(x, z_{m_1})\,dt\cdot\overset{3}{\PP}(E_{j_1})(y)\,dy\notag\\
     &+ & \eta_1 k^2 \int_{D_{j_1}} \int_0^1  \underset{y}{\nabla}\Upsilon_k(z_{m_1}, z_{j_1}+t(y-z_{j_1}))\cdot\mathcal{P}(y, z_{j_1})\,dt\cdot\overset{3}{\PP}(E_{j_1})(y)\,dy\notag\\
    &+& \eta_2 k^2 \int_{D_{j_2}} \int_0^1  \underset{x}{\nabla}\Upsilon_k(z_{m_1}+t(x-z_{m_1}), z_{j_2})\cdot\mathcal{P}(x, z_{m_1})\,dt\cdot\overset{3}{\PP}(E_{j_2})(y)\,dy\notag\\
     &+ & \eta_2 k^2 \int_{D_{j_2}} \int_0^1  \underset{y}{\nabla}\Upsilon_k(z_{m_1}, z_{j_2}+t(y-z_{j_2}))\cdot\mathcal{P}(y, z_{j_2})\,dt\cdot\overset{3}{\PP}(E_{j_2})(y)\,dy.\notag\\
\end{eqnarray}

Scaling \eqref{Ls-Dm1-P31-taylor} to $B_1$ and $B_2$, and taking the inner product w.r.t $e_n^{(3)}$, we can know that
\begin{eqnarray}\label{LS-Dm1-P31-inner}
    && \left\langle \overset{3}{\PP}(\tilde{E}_{m_1}), e_n^{(3)} \right\rangle\notag\\
    &=&\frac{1}{1+\eta_1 \lambda_n^{(3)}} \Bigg[\left\langle \tilde{E}_{m_1}^{Inc}, e_n^{(3)} \right\rangle + k^2 \eta_1 a^2 \left\langle N^{ka} \left(\overset{1}{\PP}(\tilde{E}_{m_1})\right), e_n^{(3)} \right\rangle + k^2 \eta_1 a^2\left\langle N^{ka} \left(\overset{3}{\PP}(\tilde{E}_{m_1})\right), e_n^{(3)}\right\rangle\notag\\
    &-& \eta_1 \left\langle (\nabla M^{ka} - \nabla M)\left(\overset{3}{\PP}(\tilde{E}_{m_1})\right), e_n^{(3)} \right\rangle + k^2 \eta_2 a^4 \left\langle \underset{y}{\nabla} (\Phi_kI)(z_{m_1}, z_{m_2}) \int_{B_{2}}\mathcal{P}(y, 0)\cdot\overset{1}{\PP}(\tilde{E}_{m_2})(y)\,dy, e_n^{(3)} \right\rangle \notag\\
    &+& k^2 \eta_1 a^4 \sum_{j=1 \atop j\neq m}^\aleph \left\langle  \underset{y}{\nabla}(\Phi_k I) (z_{m_1}, z_{j_1}) \int_{B_{1}}\mathcal{P}(y, 0) \overset{1}{\PP}(\tilde{E}_{j_1})(y)\,dy, e_n^{(3)}\right\rangle\notag\\
    &+ &k^2 \eta_2 a^4 \sum_{j=1 \atop j\neq m}^\aleph \left\langle  \underset{y}{\nabla}(\Phi_k I) (z_{m_1}, z_{j_2}) \int_{B_{2}}\mathcal{P}(y, 0) \overset{1}{\PP}(\tilde{E}_{j_2})(y)\,dy, e_n^{(3)}\right\rangle + k^2 \eta_2 a^3 \left\langle \Upsilon_k(z_{m_1}, z_{m_2})\int_{B_2} \overset{3}{\PP}(\tilde{E}_{m_2})(y)\,dy, e_n^{(3)}\right\rangle \notag\\
    &+& k^2 \eta_1 a^3 \sum_{j=1 \atop j\neq m}^\aleph  \left\langle \Upsilon_k(z_{m_1}, z_{j_1})\int_{B_1} \overset{3}{\PP}(\tilde{E}_{j_1})(y)\,dy, e_n^{(3)} \right\rangle  + k^2 \eta_2 a^3 \sum_{j=1 \atop j\neq m}^\aleph  \left\langle \Upsilon_k(z_{m_1}, z_{j_2})\int_{B_2} \overset{3}{\PP}(\tilde{E}_{j_2})(y)\,dy, e_n^{(3)} \right\rangle \notag\\
    &+& \left\langle \widetilde{Err}_2 \left(\overset{1}{\PP}(E_{m_2})\right), e_n^{(3)} \right\rangle + \left\langle \widetilde{Err}_2 \left(\overset{3}{\PP}(E_{m_2})\right), e_n^{(3)} \right\rangle + \sum_{j=1 \atop j\neq m}^\aleph \left\langle \widetilde{Err}_2 \left(\overset{1}{\PP}(E_{j})\right) + \widetilde{Err}_2 \left(\overset{3}{\PP}(E_{j})\right), e_n^{(3)} \right\rangle \Bigg].\notag\\
\end{eqnarray}

Observing from \cite[Equations (2.7) and (2.8)]{CGS}, we have the expansions
\begin{eqnarray}\label{expan-M^k-M}
    && \eta_1 \left\langle (\nabla M^{ka}- \nabla M)\left( \overset{3}{\PP}(\tilde{E}_{m_1})\right), e_n^{(3)} \right\rangle\notag\\
    &=& \frac{1\mp c_0 a^h}{2\lambda_{n_0}^{(1)}}\left\langle N\left(\overset{3}{\PP}(\tilde{E}_{m_1})\right), e_n^{(3)} \right\rangle + \eta_1 \frac{i (ka)^3}{12\pi} \left\langle \int_{B_1} \overset{3}{\PP}(\tilde{E}_{m_1})(y)\,dy, e_n^{(3)} \right\rangle\notag\\ \nonumber
    &-& \frac{1\mp c_0 a^h}{2\lambda_{n_0}^{(1)}}\left\langle \int_{B_1}\Phi_0(x, y)\frac{A(x, y)\overset{3}{\PP}(\tilde{E}_{m_1})(y)}{\lVert\cdot- y \rVert^2}\,dy, e_n^{(3)} \right\rangle \\ &-& \frac{\eta_1}{4\pi} \sum_{\ell\geq 3}(ika)^{\ell+1} \left\langle \int_{B_1} \frac{Hess(\lVert\cdot-y\rVert^\ell)}{(\ell+1)!}\overset{3}{\PP}(\tilde{E}_{m_1})(y)\,dy, e_n^{(3)} \right\rangle,
\end{eqnarray}
where $A(x, y)=(x-y)\otimes (x-y)$, and
\begin{eqnarray}\label{expan-N^ka}
    k^2 \eta_1 a^2 \left\langle N^{ka} \left( \overset{1}{\PP}(\tilde{E}_{m_1}) \right), e_n^{(3)} \right\rangle & = &  k^2 \eta_1 a^2 \left\langle N \left( \overset{1}{\PP}(\tilde{E}_{m_1}) \right), e_n^{(3)} \right\rangle +  k^2 \eta_1 a^2 \left\langle (N^{ka}-N) \left( \overset{1}{\PP}(\tilde{E}_{m_1}) \right), e_n^{(3)} \right\rangle \notag\\
    &=& \frac{1\mp c_0 a^h}{4\pi \lambda_{n_0}^{(1)}}\sum_{\ell\geq 1}(i ka)^{\ell+1}\left\langle \int_{B_1}\frac{\lVert\cdot - y \rVert^\ell}{(\ell+1)!}  \overset{1}{\PP}(\tilde{E}_{m_1}), e_n^{(3)} \right\rangle.
\end{eqnarray}

With the help of \eqref{expan-M^k-M} and \eqref{expan-N^ka}, we can derive from \eqref{LS-Dm1-P31-inner} by taking the square modulus and summing w.r.t. $n$ that
\begin{eqnarray}\label{LS-Dm1-P31-norm-formulate}
    && \left\lVert \overset{3}{\PP}(\tilde{E}_{m_1})\right\rVert_{\LL^2(B_1)}^2\notag\\
    &\lesssim& \sum_n\frac{1}{|1+\eta_1 \lambda_n^{(3)}|^2}\Bigg[ \left|\left\langle \tilde{E}_{m_1}^{Inc}, e_n^{(3)} \right\rangle\right|^2
    + a^4 \left|\sum_{\ell\geq 1}\left\langle \int_{B_1}\frac{\lVert\cdot - y \rVert^\ell}{(\ell+1)!}  \overset{1}{\PP}(\tilde{E}_{m_1}), e_n^{(3)} \right\rangle\right|^2 + \left| \left\langle N\left(\overset{3}{\PP}(\tilde{E}_{m_1})\right), e_n^{(3)} \right\rangle  \right|^2\notag\\
    &+& a^4 \left|\sum_{\ell\geq 1}\left\langle \int_{B_1}\frac{\lVert\cdot - y \rVert^\ell}{(\ell+1)!}  \overset{3}{\PP}(\tilde{E}_{m_1}), e_n^{(3)} \right\rangle\right|^2 + a^2 \left| \left\langle \int_{B_1}\overset{3}{\PP}(\tilde{E}_{m_1})(y)\,dy, e_n^{(3)} \right\rangle \right|^2\notag\\
    &+& \left| \left\langle \int_{B_1}\Phi_0(x, y)\frac{A(x, y)\overset{3}{\PP}(\tilde{E}_{m_1})(y)}{\lVert\cdot- y \rVert^2}\,dy, e_n^{(3)} \right\rangle \right|^2 + a^4 \left| \sum_{\ell\geq 3}\left\langle \int_{B_1} \frac{Hess(\lVert\cdot-y\rVert^\ell)}{(\ell+1)!}\overset{3}{\PP}(\tilde{E}_{m_1})(y)\,dy, e_n^{(3)} \right\rangle \right|^2\notag\\
    &+& a^8 \left| \left\langle \underset{y}{\nabla} (\Phi_kI)(z_{m_1}, z_{m_2}) \int_{B_{2}}\mathcal{P}(y, 0)\cdot\overset{1}{\PP}(\tilde{E}_{m_2})(y)\,dy, e_n^{(3)} \right\rangle \right|^2\notag\\
    &+& a^4  \sum_{j=1 \atop j\neq m}^\aleph \left| \left\langle  \underset{y}{\nabla}(\Phi_k I) (z_{m_1}, z_{j_1}) \int_{B_{1}}\mathcal{P}(y, 0) \overset{1}{\PP}(\tilde{E}_{j_1})(y)\,dy, e_n^{(3)}\right\rangle \right|^2 \notag\\
    &+& a^4 \sum_{i=1 \atop i\neq m}^\aleph \left| \left\langle  \underset{y}{\nabla}(\Phi_k I) (z_{m_1}, z_{i_1}) \int_{B_{1}}\mathcal{P}(y, 0) \overset{1}{\PP}(\tilde{E}_{i_1})(y)\,dy, e_n^{(3)}\right\rangle \right|\notag\\
    && \qquad\qquad\qquad\qquad \cdot \,  \sum_{\ell>i \atop \ell\neq m}^\aleph \left| \left\langle  \underset{y}{\nabla}(\Phi_k I) (z_{m_1}, z_{\ell_1}) \int_{B_{1}}\mathcal{P}(y, 0) \overset{1}{\PP}(\tilde{E}_{\ell_1})(y)\,dy, e_n^{(3)}\right\rangle \right|\notag\\
    &+&  a^8 \sum_{j=1 \atop j\neq m}^\aleph \left| \left\langle  \underset{y}{\nabla}(\Phi_k I) (z_{m_1}, z_{j_2}) \int_{B_{2}}\mathcal{P}(y, 0) \overset{1}{\PP}(\tilde{E}_{j_2})(y)\,dy, e_n^{(3)}\right\rangle \right|^2\notag\\
    &+&  a^8 \sum_{i=1 \atop i\neq m}^\aleph \left| \left\langle  \underset{y}{\nabla}(\Phi_k I) (z_{m_1}, z_{i_2}) \int_{B_{2}}\mathcal{P}(y, 0) \overset{1}{\PP}(\tilde{E}_{i_2})(y)\,dy, e_n^{(3)}\right\rangle \right|\notag\\
    &&\qquad\qquad\qquad\qquad \cdot \,  a^8 \sum_{\ell>i \atop \ell\neq m}^\aleph \left| \left\langle  \underset{y}{\nabla}(\Phi_k I) (z_{m_1}, z_{\ell_2}) \int_{B_{2}}\mathcal{P}(y, 0) \overset{1}{\PP}(\tilde{E}_{\ell_2})(y)\,dy, e_n^{(3)}\right\rangle \right|\notag\\
    &+& a^6 \left| \left\langle \Upsilon_k(z_{m_1}, z_{m_2})\int_{B_2} \overset{3}{\PP}(\tilde{E}_{m_2})(y)\,dy, e_n^{(3)}\right\rangle \right|^2 + a^2 \sum_{j=1 \atop j\neq m}^\aleph \left| \left\langle \Upsilon_k(z_{m_1}, z_{j_1})\int_{B_1} \overset{3}{\PP}(\tilde{E}_{j_1})(y)\,dy, e_n^{(3)} \right\rangle \right|^2 \notag\\
    &+& a^2 \sum_{\ell=1 \atop \ell\neq m}^\aleph \left| \left\langle \Upsilon_k(z_{m_1}, z_{\ell_1})\int_{B_1} \overset{3}{\PP}(\tilde{E}_{\ell_1})(y)\,dy, e_n^{(3)} \right\rangle \right|\cdot \sum_{i>\ell\atop i\neq m}^\aleph \left| \left\langle \Upsilon_k(z_{m_1}, z_{i_1})\int_{B_1} \overset{3}{\PP}(\tilde{E}_{i_1})(y)\,dy, e_n^{(3)} \right\rangle \right| \notag\\
    &+& a^6 \sum_{j=1 \atop j\neq m}^\aleph \left|  \left\langle \Upsilon_k(z_{m_1}, z_{j_2})\int_{B_2} \overset{3}{\PP}(\tilde{E}_{j_2})(y)\,dy, e_n^{(3)} \right\rangle \right|^2 \notag\\
    &+& a^6 \sum_{\ell=1 \atop \ell\neq m}^\aleph \left|  \left\langle \Upsilon_k(z_{m_1}, z_{\ell_2})\int_{B_2} \overset{3}{\PP}(\tilde{E}_{\ell_2})(y)\,dy, e_n^{(3)} \right\rangle \right| \cdot \sum_{i>\ell \atop i\neq m}^\aleph \left|  \left\langle \Upsilon_k(z_{m_1}, z_{i_2})\int_{B_2} \overset{3}{\PP}(\tilde{E}_{i_2})(y)\,dy, e_n^{(3)} \right\rangle \right|\notag\\
    &+& \left|\left\langle \widetilde{Err}_2 \left(\overset{1}{\PP}(E_{m_2})\right), e_n^{(3)} \right\rangle\right|^2 +  \left|\left\langle \widetilde{Err}_2 \left(\overset{3}{\PP}(E_{m_2})\right), e_n^{(3)} \right\rangle\right|^2 + \sum_{j=1 \atop j\neq m}^\aleph  \left|\left\langle \widetilde{Err}_2 \left(\overset{1}{\PP}(E_{j})\right), e_n^{(3)} \right\rangle\right|^2\notag\\
    &+& \sum_{\ell=1 \atop \ell\neq m}^\aleph  \left|\left\langle \widetilde{Err}_2 \left(\overset{1}{\PP}(E_{\ell})\right), e_n^{(3)} \right\rangle\right|\cdot \sum_{i>\ell \atop i\neq m}^\aleph  \left|\left\langle \widetilde{Err}_2 \left(\overset{1}{\PP}(E_{i})\right), e_n^{(3)} \right\rangle\right| + \sum_{j=1 \atop j\neq m}^\aleph  \left|\left\langle \widetilde{Err}_2 \left(\overset{3}{\PP}(E_{j})\right), e_n^{(3)} \right\rangle\right|^2 \notag\\
    &+& \sum_{\ell=1 \atop \ell\neq m}^\aleph  \left|\left\langle \widetilde{Err}_2 \left(\overset{3}{\PP}(E_{\ell})\right), e_n^{(3)} \right\rangle\right|\cdot \sum_{i>\ell \atop i\neq m}^\aleph  \left|\left\langle \widetilde{Err}_2 \left(\overset{3}{\PP}(E_{i})\right), e_n^{(3)} \right\rangle\right| \Bigg].\notag\\
\end{eqnarray}

We are now in a position to estimate the R.H.S. of \eqref{LS-Dm1-P31-norm-formulate} term by term.

\bigskip 

\noindent (a). Estimate for :
\begin{eqnarray}\label{es-P31-RHS-1}
    && \sum_{n} \frac{1}{|1+\eta_1 \lambda_n^{(3)}|^2} a^8 \left| \left\langle \underset{y}{\nabla} (\Phi_kI)(z_{m_1}, z_{m_2}) \int_{B_{2}}\mathcal{P}(y, 0)\cdot\overset{1}{\PP}(\tilde{E}_{m_2})(y)\,dy, e_n^{(3)} \right\rangle \right|^2\notag\\
    &\lesssim& a^{12} \left\lVert \underset{y}{\nabla} (\Phi_kI)(z_{m_1}, z_{m_2}) \int_{B_{2}}\mathcal{P}(y, 0)\cdot\overset{1}{\PP}(\tilde{E}_{m_2})(y)\,dy \right\rVert_{\LL^2(B_1)}^2 \, \lesssim \, a^{12} \, d_{\rm in}^{-4} \left\lVert\overset{1}{\PP}(\tilde{E}_{m_2})\right\rVert^2_{\LL^2(B_2)}.
\end{eqnarray}

\noindent (b). Estimate for :
\begin{eqnarray}\label{es-P31-RHS-2}
    && \sum_{n} \frac{1}{|1+\eta_1 \lambda_n^{(3)}|^2} a^4  \sum_{j=1 \atop j\neq m}^\aleph \left| \left\langle \underset{y}{\nabla} (\Phi_kI)(z_{m_1}, z_{j_1}) \int_{B_{1}}\mathcal{P}(y, 0)\cdot\overset{1}{\PP}(\tilde{E}_{j_1})(y)\,dy, e_n^{(3)} \right\rangle \right|^2\notag\\
    &\lesssim& a^{8} \sum_{j=1 \atop j\neq m}^\aleph \left\lVert \underset{y}{\nabla} (\Phi_kI)(z_{m_1}, z_{j_1}) \int_{B_{1}}\mathcal{P}(y, 0)\cdot\overset{1}{\PP}(\tilde{E}_{j_1})(y)\,dy \right\rVert_{\LL^2(B_1)}^2 \notag\\
    &\lesssim&  a^{8} \sum_{j=1 \atop j\neq m}^\aleph d_{mj}^{-4} \left\lVert\overset{1}{\PP}(\tilde{E}_{j_1})\right\rVert^2_{\LL^2(B_1)}\, \lesssim \, a^{8} d_{\rm out}^{-4} \max_j\left\lVert\overset{1}{\PP}(\tilde{E}_{j_1})\right\rVert^2_{\LL^2(B_1)}. 
\end{eqnarray}

\noindent (c) Estimate for 
\begin{eqnarray}\label{es-P31-RHS-2-crossing}
    && \sum_{n} \frac{1}{|1+\eta_1 \lambda_n^{(3)}|^2} a^4  \sum_{i=1 \atop i\neq m}^\aleph \left| \left\langle \underset{y}{\nabla} (\Phi_kI)(z_{m_1}, z_{i_1}) \int_{B_{1}}\mathcal{P}(y, 0)\cdot\overset{1}{\PP}(\tilde{E}_{i_1})(y)\,dy, e_n^{(3)} \right\rangle \right|\notag\\
    &&\qquad\qquad\qquad\qquad\cdot\, \sum_{\ell>i \atop \ell\neq m}^\aleph \left| \left\langle \underset{y}{\nabla} (\Phi_kI)(z_{m_1}, z_{\ell_1}) \int_{B_{1}}\mathcal{P}(y, 0)\cdot\overset{1}{\PP}(\tilde{E}_{\ell_1})(y)\,dy, e_n^{(3)} \right\rangle \right|\notag\\
    &\lesssim& a^{8} \left( \sum_{i=1 \atop i\neq m}^\aleph d_{mi}^{-2}\left\lVert\overset{1}{\PP}(\tilde{E}_{i_1})\right\rVert_{\LL^2(B_1)} \right)\cdot \left( \sum_{\ell>i \atop \ell\neq m}^\aleph d_{m\ell}^{-2}\left\lVert\overset{1}{\PP}(\tilde{E}_{\ell_1})\right\rVert_{\LL^2(B_1)} \right)\notag\\
    &\lesssim&  a^{8} \left( d_{\rm out}^{-3}\max_j \left\lVert\overset{1}{\PP}(\tilde{E}_{j_1})\right\rVert_{\LL^2(B_1)}\right)^2 \lesssim a^{8} d_{\rm out}^{-6} \max_j\left\lVert\overset{1}{\PP}(\tilde{E}_{j_1})\right\rVert^2_{\LL^2(B_1)}. 
\end{eqnarray} 

\noindent (d). Estimate for:
\begin{eqnarray}\label{es-P31-RHS-3}
    && \sum_{n} \frac{1}{|1+\eta_1 \lambda_n^{(3)}|^2} a^8 \sum_{j=1 \atop j\neq m}^\aleph \left| \left\langle \underset{y}{\nabla} (\Phi_kI)(z_{m_1}, z_{j_2}) \int_{B_{2}}\mathcal{P}(y, 0)\cdot\overset{1}{\PP}(\tilde{E}_{j_2})(y)\,dy, e_n^{(3)} \right\rangle \right|^2\notag\\
    &\lesssim& a^{12} \sum_{j=1 \atop j\neq m}^\aleph \left\lVert \underset{y}{\nabla} (\Phi_kI)(z_{m_1}, z_{j_2}) \int_{B_{2}}\mathcal{P}(y, 0)\cdot\overset{1}{\PP}(\tilde{E}_{j_2})(y)\,dy \right\rVert_{\LL^2(B_1)}^2 \notag\\
    &\lesssim&  a^{12} \sum_{j=1 \atop j\neq m}^\aleph d_{mj}^{-4} \left\lVert\overset{1}{\PP}(\tilde{E}_{j_2})\right\rVert^2_{\LL^2(B_2)}\, \lesssim \, a^{12} d_{\rm out}^{-4} \max_j\left\lVert\overset{1}{\PP}(\tilde{E}_{j_2})\right\rVert^2_{\LL^2(B_2)}. 
\end{eqnarray}

\noindent (e) Estimate for 
\begin{eqnarray}\label{es-P31-RHS-3-crossing}
    && \sum_{n} \frac{1}{|1+\eta_1 \lambda_n^{(3)}|^2} a^8 \sum_{i=1 \atop i\neq m}^\aleph \left| \left\langle \underset{y}{\nabla} (\Phi_kI)(z_{m_1}, z_{i_2}) \int_{B_{2}}\mathcal{P}(y, 0)\cdot\overset{1}{\PP}(\tilde{E}_{i_2})(y)\,dy, e_n^{(3)} \right\rangle \right|\notag\\
    &&\qquad\qquad\qquad\qquad\cdot\, \sum_{\ell>i \atop \ell\neq m}^\aleph \left| \left\langle \underset{y}{\nabla} (\Phi_kI)(z_{m_1}, z_{\ell_2}) \int_{B_{2}}\mathcal{P}(y, 0)\cdot\overset{1}{\PP}(\tilde{E}_{\ell_2})(y)\,dy, e_n^{(3)} \right\rangle \right|\notag\\
    &\lesssim& a^{12} \left( \sum_{i=1 \atop i\neq m}^\aleph d_{mi}^{-2}\left\lVert\overset{1}{\PP}(\tilde{E}_{i_2})\right\rVert_{\LL^2(B_2)} \right)\cdot \left( \sum_{\ell>i \atop \ell\neq m}^\aleph d_{m\ell}^{-2}\left\lVert\overset{1}{\PP}(\tilde{E}_{\ell_2})\right\rVert_{\LL^2(B_2)} \right)\notag\\
    &\lesssim&  a^{12} \left( d_{\rm out}^{-3}\max_j \left\lVert\overset{1}{\PP}(\tilde{E}_{j_2})\right\rVert_{\LL^2(B_2)}\right)^2 \lesssim a^{12} d_{\rm out}^{-6} \max_j\left\lVert\overset{1}{\PP}(\tilde{E}_{j_2})\right\rVert^2_{\LL^2(B_2)}. 
\end{eqnarray} 

\noindent (f). Estimate for:
\begin{eqnarray}\label{es-P31-RHS-Upsilon1}
    &&  \sum_n \frac{1}{|1+\eta_1 \lambda_n^{(3)}|^2} a^6 \left| \left\langle \Upsilon_k(z_{m_1}, z_{m_2}) \int_{B_2}\overset{3}{\PP}(\tilde{E}_{m_2})(y)\,dy, e_n^{(3)} \right\rangle\right|^2\notag\\
    &\lesssim & a^{10} \left\lVert \Upsilon_k(z_{m_1}, z_{m_2}) \int_{B_2}\overset{3}{\PP}(\tilde{E}_{m_2})(y)\,dy \right\rVert^2_{\LL^2(B_1)}\, \lesssim \, a^{10} \, d_{\rm in}^{-6} \left\lVert \overset{3}{\PP}(\tilde{E}_{m_2})\right\rVert^2_{\LL^2(B_2)}.
\end{eqnarray}

\noindent (g). Estimate for:
\begin{eqnarray}\label{es-P31-RHS-Upsilon2}
    &&  \sum_n \frac{1}{|1+\eta_1 \lambda_n^{(3)}|^2} a^2 \sum_{j=1 \atop j\neq m}^\aleph \left| \left\langle \Upsilon_k(z_{m_1}, z_{j_1}) \int_{B_1}\overset{3}{\PP}(\tilde{E}_{j_1})(y)\,dy, e_n^{(3)} \right\rangle\right|^2\notag\\
    &\lesssim & a^{6} \sum_{j=1 \atop j\neq m}^\aleph \left\lVert \Upsilon_k(z_{m_1}, z_{j_1}) \int_{B_1}\overset{3}{\PP}(\tilde{E}_{j_1})(y)\,dy \right\rVert^2_{\LL^2(B_1)}\, \lesssim \, a^6\sum_{j=1 \atop j\neq m}^\aleph d_{mj}^{-6} \left\lVert \overset{3}{\PP}(\tilde{E}_{j_1})\right\rVert_{\LL^2(B_1)}^2 \notag\\
    &\lesssim & a^{6} \, d_{\rm out}^{-6} \max_j\left\lVert \overset{3}{\PP}(\tilde{E}_{j_1})\right\rVert^2_{\LL^2(B_1)}.
\end{eqnarray}

\noindent (h). Estimate for:
\begin{eqnarray}\label{es-P31-RHS-Upsilon3}
    &&  \sum_n \frac{1}{|1+\eta_1 \lambda_n^{(3)}|^2} a^2 \sum_{\ell=1 \atop \ell\neq m}^\aleph \left| \left\langle \Upsilon_k(z_{m_1}, z_{\ell_1}) \int_{B_1}\overset{3}{\PP}(\tilde{E}_{\ell_1})(y)\,dy, e_n^{(3)} \right\rangle\right|\notag\\
    &&\qquad\qquad\qquad\qquad\cdot \sum_{i>\ell \atop i\neq m}^\aleph \left| \left\langle \Upsilon_k(z_{m_1}, z_{i_1}) \int_{B_1}\overset{3}{\PP}(\tilde{E}_{i_1})(y)\,dy, e_n^{(3)} \right\rangle\right|\notag\\
    &\lesssim & a^{6} \left(\sum_{\ell=1 \atop \ell\neq m}^\aleph d_{m\ell}^{-3} \left\lVert \overset{3}{\PP}(\tilde{E}_{\ell_1})\right\rVert_{\LL^2(B_1)}\right) \cdot \left(\sum_{i>\ell \atop i\neq m}^\aleph d_{mi}^{-3} \left\lVert \overset{3}{\PP}(\tilde{E}_{i_1})\right\rVert_{\LL^2(B_1)}\right)\notag\\
    &\lesssim& a^6 \left( d_{\rm out}^{-3} |\ln (d_{\rm out})|\max_j \left\lVert\overset{3}{\PP}(\tilde{E}_{j_1})\right\rVert_{\LL^2(B_1)} \right)^2\, \lesssim\, a^6 d_{\rm out}^{-6} \ln^2 d_{\rm out} \max_j\left\lVert\overset{3}{\PP}(\tilde{E}_{j_1})\right\rVert_{\LL^2(B_1)}^2.
\end{eqnarray}

\noindent (i). Similar to \eqref{es-P31-RHS-Upsilon2} and \eqref{es-P31-RHS-Upsilon3}, there holds the estimate for:
\begin{equation}\label{es-P31-RHS-Upsilon4}
      \sum_n \frac{1}{|1+\eta_1 \lambda_n^{(3)}|^2} a^6 \sum_{j=1 \atop j\neq m}^\aleph \left| \left\langle \Upsilon_k(z_{m_1}, z_{j_2}) \int_{B_2}\overset{3}{\PP}(\tilde{E}_{j_2})(y)\,dy, e_n^{(3)} \right\rangle\right|^2 \lesssim  a^{10} \, d_{\rm out}^{-6} \max_j\left\lVert \overset{3}{\PP}(\tilde{E}_{j_2})\right\rVert^2_{\LL^2(B_2)},
\end{equation}
and 
\begin{eqnarray}\label{es-P31-RHS-Upsilon5}
     &&  \sum_n \frac{1}{|1+\eta_1 \lambda_n^{(3)}|^2} a^6 \sum_{\ell=1 \atop \ell\neq m}^\aleph \left| \left\langle \Upsilon_k(z_{m_1}, z_{\ell_2}) \int_{B_2}\overset{3}{\PP}(\tilde{E}_{\ell_2})(y)\,dy, e_n^{(3)} \right\rangle\right|\notag\\
    &&\qquad\qquad\qquad\qquad\cdot \sum_{i>\ell \atop i\neq m}^\aleph \left| \left\langle \Upsilon_k(z_{m_1}, z_{i_2}) \int_{B_2}\overset{3}{\PP}(\tilde{E}_{i_2})(y)\,dy, e_n^{(3)} \right\rangle\right|\notag\\
    &\lesssim& a^{10} d_{\rm out}^{-6} \ln^2 d_{\rm out} \max_j \left\lVert \overset{3}{\PP}(\tilde{E}_{j_2})\right\rVert^2_{\LL^2(B_2)}.
\end{eqnarray}

\noindent (j) Estimate for $\sum_n\frac{1}{|1+\eta_1 \lambda_n^{(3)}|^2}\left| \left\langle \widetilde{Err}_2\left( \overset{1}{\PP}(E_{m_2})\right), e_n^{(3)} \right\rangle \right|^2$:

From the definition of $ \widetilde{Err}_2\left( \overset{1}{\PP}(E_{m_2})\right)$ given in \eqref{def-Err_2-P13mj}, we can know that
\begin{eqnarray}\label{es-P13-RHS-Err_2-P1m2}
    && \sum_n\frac{1}{|1+\eta_1 \lambda_n^{(3)}|^2}\left| \left\langle \widetilde{Err}_2\left( \overset{1}{\PP}(E_{m_2})\right), e_n^{(3)} \right\rangle \right|^2 \notag\\
    &\lesssim& a^{14} \left\lVert \int_{B_{2}} \int_0^1 (1-t) y^\perp \underset{y}{Hess}\Phi_k(z_{m_1}, z_{m_2}+tay)\cdot y\,dt\cdot\overset{1}{\PP}(\tilde{E}_{m_2})(y)\,dy \right\rVert^2_{\LL^2(B_1)}\notag\\
    &+& a^{14} \left\lVert \int_{B_{2}} \int_0^1  \underset{y}{Hess}\Phi_k(z_{m_1}+tax, z_{m_2})\cdot y \cdot x\,dt\cdot\overset{1}{\PP}(\tilde{E}_{m_2})(y)\,dy \right\rVert^2_{\LL^2(B_1)}\notag\\
    &+& a^{16} \left\lVert \int_{B_{2}} \int_0^1 \int_0^1 (1-t) y^\perp \underset{x}{\nabla}(\underset{y}{Hess}\Phi_k)(z_{m_1}+tax, z_{m_2}+say)\cdot y\,ds\cdot\mathcal{P}(x, 0)\,dt\cdot\overset{1}{\PP}(\tilde{E}_{m_2})(y)\,dy \right\rVert^2_{\LL^2(B_2)}\notag\\
    &\lesssim& a^{14} \, d_{\rm in}^{-6} \, \left\lVert \overset{1}{\PP}(\tilde{E}_{m_2})\right\rVert^2_{\LL^2(B_2)}.\notag\\
\end{eqnarray}

\noindent (k) Estimate for $\sum_n\frac{1}{|1+\eta_1 \lambda_n^{(3)}|^2}\left| \left\langle \widetilde{Err}_2\left( \overset{3}{\PP}(E_{m_2})\right), e_n^{(3)} \right\rangle \right|^2$:

From the definition of $\widetilde{Err}_2\left( \overset{3}{\PP}(E_{m_2})\right)$ in \eqref{def-Err_2-P13mj}, we directly have
\begin{eqnarray}\label{es-P13-RHS-Err_2-P3m2}
    &&\sum_n\frac{1}{|1+\eta_1 \lambda_n^{(3)}|^2}\left| \left\langle \widetilde{Err}_2\left( \overset{3}{\PP}(E_{m_2})\right), e_n^{(3)} \right\rangle \right|^2\notag\\
    &\lesssim& a^{12} \left\lVert  \int_{B_{2}} \int_0^1  \underset{x}{\nabla}\Upsilon_k(z_{m_1}+tax, z_{m_2})\cdot\mathcal{P}(x, 0)\,dt\cdot\overset{3}{\PP}(\tilde{E}_{m_2})(y)\,dy\ \right\rVert_{\LL^2(B_1)}^2\notag\\
    &+& a^{12} \left\lVert \int_{B_{2}} \int_0^1  \underset{y}{\nabla}\Upsilon_k(z_{m_1}, z_{m_2}+tay)\cdot\mathcal{P}(y, 0)\,dt\cdot\overset{3}{\PP}(\tilde{E}_{m_2})(y)\,dy \right\rVert^2_{\LL^2(B_1)}\, \lesssim\, a^{12}\, d_{\rm in}^{-8} \left\lVert \overset{3}{\PP}(\tilde{E}_{m_2})\right\rVert^2_{\LL^2(B_2)}.
\end{eqnarray}

\noindent (l) Estimate for $\sum_n\frac{1}{|1+\eta_1 \lambda_n^{(3)}|^2}\sum_{j=1 \atop j\neq m}^\aleph \left| \left\langle \widetilde{Err}_2\left( \overset{1}{\PP}(E_{j_1})\right), e_n^{(3)} \right\rangle \right|^2$:

With the definition of $\widetilde{Err}_2\left( \overset{1}{\PP}(E_{j_1})\right)$ in \eqref{def-Err_2-P13mj}, we obtain that
\begin{eqnarray}\label{es-P13-RHS-Err_2-P1j1}
    &&\sum_n\frac{1}{|1+\eta_1 \lambda_n^{(3)}|^2}\sum_{j=1 \atop j\neq m}^\aleph \left| \left\langle \widetilde{Err}_2\left( \overset{1}{\PP}(E_{j})\right), e_n^{(3)} \right\rangle \right|^2\notag\\
    &\lesssim& a^4\cdot a^6 \sum_{j=1\atop j\neq m}^\aleph \left\lVert \int_{B_{1}} \int_0^1 (1-t) y^\perp \underset{y}{Hess}\Phi_k(z_{m_1}, z_{j_1}+tay)\cdot y\,dt\cdot\overset{1}{\PP}(\tilde{E}_{j_1})(y)\,dy\right\rVert^2_{\LL^2(B_1)}\notag\\
    &+& a^4\cdot a^6 \sum_{j=1\atop j\neq m}^\aleph \left\lVert \int_{B_{1}} \int_0^1  \underset{y}{Hess}\Phi_k(z_{m_1}+tax, z_{j_1})\cdot y\cdot x\,dt\cdot\overset{1}{\PP}(\tilde{E}_{j_1})(y)\,dy \right\rVert_{\LL^2(B_1)}^2\notag\\
    &+& a^4 \cdot a^{8}\sum_{j=1\atop j\neq m}^\aleph \left\lVert \int_{B_{1}} \int_0^1 \int_0^1 (1-t) y^\perp \underset{x}{\nabla}(\underset{y}{Hess}\Phi_k)(z_{m_1}+tax, z_{j_1}+say)\cdot y\,ds\cdot\mathcal{P}(x, 0)\,dt\cdot\overset{1}{\PP}(\tilde{E}_{j_1})(y)\,dy \right\rVert_{\LL^2(B_1)}^2\notag\\
    &+& a^4\cdot a^{10} \sum_{j=1\atop j\neq m}^\aleph \left\lVert \int_{B_{2}} \int_0^1 (1-t) y^\perp \underset{y}{Hess}\Phi_k(z_{m_1}, z_{j_2}+tay)\cdot y\,dt\cdot\overset{1}{\PP}(\tilde{E}_{j_2})(y)\,dy\right\rVert^2_{\LL^2(B_1)}\notag\\
    &+& a^4\cdot a^{10} \sum_{j=1\atop j\neq m}^\aleph \left\lVert \int_{B_{2}} \int_0^1  \underset{y}{Hess}\Phi_k(z_{m_1}+tax, z_{j_2})\cdot y\cdot x\,dt\cdot\overset{1}{\PP}(\tilde{E}_{j_2})(y)\,dy \right\rVert_{\LL^2(B_1)}^2\notag\\
    &+& a^4 \cdot a^{12} \sum_{j=1\atop j\neq m}^\aleph \left\lVert \int_{B_{2}} \int_0^1 \int_0^1 (1-t) y^\perp \underset{x}{\nabla}(\underset{y}{Hess}\Phi_k)(z_{m_1}+tax, z_{j_2}+say)\cdot y\,ds\cdot\mathcal{P}(x, 0)\,dt\cdot\overset{1}{\PP}(\tilde{E}_{j_2})(y)\,dy \right\rVert_{\LL^2(B_1)}^2\notag\\
    &\lesssim& a^{10} \sum_{j=1 \atop j\neq m}^\aleph d_{mj}^{-6} \left\lVert \overset{1}{\PP}(\tilde{E}_{j_1})\right\rVert_{\LL^2(B_1)}^2 \, + \, a^{14} \sum_{j=1 \atop j\neq m}^\aleph d_{mj}^{-6} \left\lVert \overset{1}{\PP}(\tilde{E}_{j_2})\right\rVert_{\LL^2(B_2)}^2\notag\\
    &\lesssim& a^{10} \, d_{\rm out}^{-6} \max_j \left\lVert \overset{1}{\PP}(\tilde{E}_{j_1})\right\rVert_{\LL^2(B_1)}^2 \, + \, a^{14} d_{\rm out}^{-6} \max_j \left\lVert \overset{1}{\PP}(\tilde{E}_{j_2})\right\rVert_{\LL^2(B_2)}^2.\notag\\
\end{eqnarray}

\noindent (m) Estimate for $\sum_n\frac{1}{|1+\eta_1 \lambda_n^{(3)}|^2}\sum_{\ell=1 \atop \ell\neq m}^\aleph \left| \left\langle \widetilde{Err}_2\left( \overset{1}{\PP}(E_{\ell})\right), e_n^{(3)} \right\rangle \right|\cdot \sum_{i>\ell\atop i\neq m}^\aleph \left| \left\langle \widetilde{Err}_2\left( \overset{1}{\PP}(E_{i})\right), e_n^{(3)} \right\rangle \right|$:

By \eqref{def-Err_2-P13mj}, we shall write 
\begin{equation*}
    \widetilde{Err}_2\left( \overset{1}{\PP}(E_{j})\right)=\widetilde{Err}_2\left( \overset{1}{\PP}(E_{j_1})\right)+ \widetilde{Err}_2\left( \overset{1}{\PP}(E_{j_1})\right),
\end{equation*}
where
\begin{eqnarray*}
    &&\widetilde{Err}_2\left( \overset{1}{\PP}(E_{j_1})\right)\notag\\
    &:=& \eta_1 k^2 \int_{B_{1}} \int_0^1 (1-t) y^\perp \underset{y}{Hess}\Phi_k(z_{m_1}, z_{j_1}+tay)\cdot y\,dt\cdot\overset{1}{\PP}(\tilde{E}_{j_1})(y)\,dy\notag\\
     &-& \eta_1 k^2 \int_{B_{1}} \int_0^1  \underset{y}{Hess}\Phi_k(z_{m_1}+tax, z_{j_1})\cdot y\cdot x\,dt\cdot\overset{1}{\PP}(\tilde{E}_{j_1})(y)\,dy\notag\\
     &+& \eta_1 k^2 \int_{B_{1}} \int_0^1 \int_0^1 (1-t) y^\perp \underset{x}{\nabla}(\underset{y}{Hess}\Phi_k)(z_{m_1}+tax, z_{j_1}+say)\cdot y\,ds\cdot\mathcal{P}(x, 0)\,dt\cdot\overset{1}{\PP}(\tilde{E}_{j_1})(y)\,dy\notag\\
     &&\widetilde{Err}_2\left( \overset{1}{\PP}(E_{j_1})\right)\notag\\
     &:=& \eta_2 k^2 \int_{B_{2}} \int_0^1 (1-t) y^\perp \underset{y}{Hess}\Phi_k(z_{m_1}, z_{j_2}+tay)\cdot y\,dt\cdot\overset{1}{\PP}(\tilde{E}_{j_2})(y)\,dy\notag\\
     &-& \eta_2 k^2 \int_{B_{2}} \int_0^1  \underset{y}{Hess}\Phi_k(z_{m_1}+tax, z_{j_2})\cdot y\cdot x\,dt\cdot\overset{1}{\PP}(\tilde{E}_{j_2})(y)\,dy\notag\\
     &+& \eta_2 k^2 \int_{B_{2}} \int_0^1 \int_0^1 (1-t) y^\perp \underset{x}{\nabla}(\underset{y}{Hess}\Phi_k)(z_{m_1}+tax, z_{j_2}+say)\cdot y\,ds\cdot\mathcal{P}(x, 0)\,dt\cdot\overset{1}{\PP}(\tilde{E}_{j_2})(y)\,dy.
\end{eqnarray*}
Hence, there holds
\begin{eqnarray*}
    &&\sum_n\frac{1}{|1+\eta_1 \lambda_n^{(3)}|^2}\sum_{\ell=1 \atop \ell\neq m}^\aleph \left| \left\langle \widetilde{Err}_2\left( \overset{1}{\PP}(E_{\ell})\right), e_n^{(3)} \right\rangle \right|\cdot \sum_{i>\ell\atop i\neq m}^\aleph \left| \left\langle \widetilde{Err}_2\left( \overset{1}{\PP}(E_{i})\right), e_n^{(3)} \right\rangle \right|\notag\\
    &=& \sum_n\frac{1}{|1+\eta_1 \lambda_n^{(3)}|^2}\sum_{\ell=1 \atop \ell\neq m}^\aleph \left| \left\langle \widetilde{Err}_2\left( \overset{1}{\PP}(E_{\ell_1})\right), e_n^{(3)} \right\rangle + \left\langle \widetilde{Err}_2\left( \overset{1}{\PP}(E_{\ell_2})\right), e_n^{(3)} \right\rangle \right| \notag\\
    &&\qquad\qquad\qquad\qquad\cdot \sum_{i>\ell\atop i\neq m}^\aleph \left| \left\langle \widetilde{Err}_2\left( \overset{1}{\PP}(E_{i_1})\right), e_n^{(3)} \right\rangle + \left\langle \widetilde{Err}_2\left( \overset{1}{\PP}(E_{i_2})\right), e_n^{(3)} \right\rangle \right|\notag\\
    &\lesssim& a^4 \sum_{\ell=1\atop \ell\neq m}^\aleph \left( \left\lVert \widetilde{Err}_2\left(\overset{1}{\PP}(E_{\ell_1})\right)\right\rVert_{\LL^2(B_1)} + \left\lVert \widetilde{Err}_2\left(\overset{1}{\PP}(E_{\ell_2})\right)\right\rVert_{\LL^2(B_1)} \right)\notag\\
    &&\qquad\qquad\qquad\qquad\cdot \sum_{i>\ell \atop i\neq m}^\aleph \left( \left\lVert \widetilde{Err}_2\left(\overset{1}{\PP}(E_{i_1})\right)\right\rVert_{\LL^2(B_1)} + \left\lVert \widetilde{Err}_2\left(\overset{1}{\PP}(E_{i_2})\right)\right\rVert_{\LL^2(B_1)} \right),
\end{eqnarray*}
where
\begin{eqnarray*}
    &&\left\lVert \widetilde{Err}_2\left(\overset{1}{\PP}(E_{j_1})\right)\right\rVert_{\LL^2(B_1)} = \left( \sum_n \left| \left\langle \widetilde{Err}_2\left(\overset{1}{\PP}(E_{j_1})\right), e_n^{(3)} \right\rangle \right|^2 \right)^{\frac{1}{2}}\notag\\
    &\lesssim& a^3 \left\lVert \int_{B_{1}} \int_0^1 (1-t) y^\perp \underset{y}{Hess}\Phi_k(z_{m_1}, z_{j_1}+tay)\cdot y\,dt\cdot\overset{1}{\PP}(\tilde{E}_{j_1})(y)\,dy \right\rVert_{\LL^2(B_1)}\notag\\
    &+& a^3 \left\lVert \int_{B_{1}} \int_0^1  \underset{y}{Hess}\Phi_k(z_{m_1}+tax, z_{j_1})\cdot y\cdot x\,dt\cdot\overset{1}{\PP}(\tilde{E}_{j_1})(y)\,dy \right\rVert_{\LL^2(B_1)}\notag\\
    &+& a^4 \left\lVert \int_{B_{1}} \int_0^1 \int_0^1 (1-t) y^\perp \underset{x}{\nabla}(\underset{y}{Hess}\Phi_k)(z_{m_1}+tax, z_{j_1}+say)\cdot y\,ds\cdot\mathcal{P}(x, 0)\,dt\cdot\overset{1}{\PP}(\tilde{E}_{j_1})(y)\,dy \right\rVert_{\LL^2(B_1)}\notag\\
    &\lesssim& a^3 \, d_{mj}^{-3} \, \left\lVert \overset{1}{\PP}(\tilde{E}_{j_1}) \right\rVert_{\LL^2(B_1)},
\end{eqnarray*}
and 
\begin{eqnarray*}
     &&\left\lVert \widetilde{Err}_2\left(\overset{1}{\PP}(E_{j_2})\right)\right\rVert_{\LL^2(B_1)} = \left( \sum_n \left| \left\langle \widetilde{Err}_2\left(\overset{1}{\PP}(E_{j_2})\right), e_n^{(3)} \right\rangle \right|^2 \right)^{\frac{1}{2}}\notag\\
     &\lesssim& a^5 \left\lVert \int_{B_{2}} \int_0^1 (1-t) y^\perp \underset{y}{Hess}\Phi_k(z_{m_1}, z_{j_2}+tay)\cdot y\,dt\cdot\overset{1}{\PP}(\tilde{E}_{j_2})(y)\,dy \right\rVert_{\LL^2(B_1)}\notag\\
     &+& a^5 \left\lVert \int_{B_{2}} \int_0^1  \underset{y}{Hess}\Phi_k(z_{m_1}+tax, z_{j_2})\cdot y\cdot x\,dt\cdot\overset{1}{\PP}(\tilde{E}_{j_2})(y)\,dy \right\rVert_{\LL^2(B_1)}\notag\\
     &+& a^6 \left\lVert  \int_{B_{2}} \int_0^1 \int_0^1 (1-t) y^\perp \underset{x}{\nabla}(\underset{y}{Hess}\Phi_k)(z_{m_1}+tax, z_{j_2}+say)\cdot y\,ds\cdot\mathcal{P}(x, 0)\,dt\cdot\overset{1}{\PP}(\tilde{E}_{j_2})(y)\,dy \right\rVert_{\LL^2(B_1)}\notag\\
     &\lesssim& a^5 \, d_{mj}^{-3} \, \left\lVert \overset{1}{\PP}(\tilde{E}_{j_2}) \right\rVert_{\LL^2(B_2)}.
\end{eqnarray*}
Therefore, we can obtain that
\begin{eqnarray}\label{es-P13-RHS-P1j-crossing}
    &&\sum_n\frac{1}{|1+\eta_1 \lambda_n^{(3)}|^2}\sum_{\ell=1 \atop \ell\neq m}^\aleph \left| \left\langle \widetilde{Err}_2\left( \overset{1}{\PP}(E_{\ell})\right), e_n^{(3)} \right\rangle \right|\cdot \sum_{i>\ell\atop i\neq m}^\aleph \left| \left\langle \widetilde{Err}_2\left( \overset{1}{\PP}(E_{i})\right), e_n^{(3)} \right\rangle \right|\notag\\
    &\lesssim& a^4 \sum_{\ell=1 \atop\ell\neq m}^\aleph \left(a^3 d_{m\ell}^{-3} \left\lVert \overset{1}{\PP}(\tilde{E}_{\ell_1}) \right\rVert_{\LL^2(B_1)} + a^5  d_{m\ell}^{-3}  \left\lVert \overset{1}{\PP}(\tilde{E}_{\ell_2}) \right\rVert_{\LL^2(B_2)} \right)\notag\\
    &&\qquad\qquad\qquad\cdot \sum_{i>\ell \atop i\neq m}^\aleph \left(a^3 d_{mi}^{-3} \left\lVert \overset{1}{\PP}(\tilde{E}_{i_1}) \right\rVert_{\LL^2(B_1)} + a^5  d_{mi}^{-3}  \left\lVert \overset{1}{\PP}(\tilde{E}_{i_2}) \right\rVert_{\LL^2(B_2)} \right)\notag\\
    &\lesssim& a^4 \left( a^3\, d_{\rm out}^{-3} |\ln (d_{\rm out})| \max_j \left\lVert \overset{1}{\PP}(\tilde{E}_{j_1}) \right\rVert_{\LL^2(B_1)} + a^5\, d_{\rm out}^{-3} |\ln (d_{\rm out})| \max_j \left\lVert \overset{1}{\PP}(\tilde{E}_{j_2}) \right\rVert_{\LL^2(B_2)} \right)^2\notag\\
    &\lesssim& a^{10}\, d_{\rm out}^{-6}\, \ln^2 d_{\rm out} \max_j \left\lVert \overset{1}{\PP}(\tilde{E}_{j_1}) \right\rVert_{\LL^2(B_1)}^2 + a^{14} \, d_{\rm out}^{-6}\, \ln^2 d_{\rm out} \max_j \left\lVert \overset{1}{\PP}(\tilde{E}_{j_2}) \right\rVert_{\LL^2(B_2)}^2.
\end{eqnarray}

\noindent (n) Estimate for $\sum_n\frac{1}{|1+\eta_1 \lambda_n^{(3)}|^2}\sum_{j=1 \atop j\neq m}^\aleph \left| \left\langle \widetilde{Err}_2\left( \overset{3}{\PP}(E_{j_1})\right), e_n^{(3)} \right\rangle \right|^2$:

From $\widetilde{Err}_2\left( \overset{3}{\PP}(E_{j_1})\right)$ in \eqref{def-Err_2-P13mj}, we obtain that
\begin{eqnarray}\label{es-P13-RHS-Err_2-P3j1}
    &&\sum_n\frac{1}{|1+\eta_1 \lambda_n^{(3)}|^2}\sum_{j=1 \atop j\neq m}^\aleph \left| \left\langle \widetilde{Err}_2\left( \overset{3}{\PP}(E_{j})\right), e_n^{(3)} \right\rangle \right|^2\notag\\
    &\lesssim& a^4 \cdot a^4 \sum_{j=1 \atop j\neq m}^\aleph \left\lVert \int_{B_{1}} \int_0^1  \underset{x}{\nabla}\Upsilon_k(z_{m_1}+tax, z_{j_1})\cdot\mathcal{P}(x, 0)\,dt\cdot\overset{3}{\PP}(\tilde{E}_{j_1})(y)\,dy \right\rVert^2_{\LL^2(B_1)}\notag\\
    &+& a^4 \cdot a^4 \sum_{j=1 \atop j\neq m}^\aleph \left\lVert \int_{B_{1}} \int_0^1  \underset{y}{\nabla}\Upsilon_k(z_{m_1}, z_{j_1}+tay)\cdot\mathcal{P}(y, 0)\,dt\cdot\overset{3}{\PP}(\tilde{E}_{j_1})(y)\,dy \right\rVert_{\LL^2(B_1)}\notag\\
    &+& a^4 \cdot a^8 \sum_{j=1 \atop j\neq m}^\aleph \left\lVert \int_{B_{2}} \int_0^1  \underset{x}{\nabla}\Upsilon_k(z_{m_1}+tax, z_{j_2})\cdot\mathcal{P}(x, 0)\,dt\cdot\overset{3}{\PP}(\tilde{E}_{j_2})(y)\,dy \right\rVert^2_{\LL^2(B_1)}\notag\\
    &+& a^4 \cdot a^8 \sum_{j=1 \atop j\neq m}^\aleph \left\lVert \int_{B_{2}} \int_0^1  \underset{y}{\nabla}\Upsilon_k(z_{m_1}, z_{j_2}+tay)\cdot\mathcal{P}(y, 0)\,dt\cdot\overset{3}{\PP}(\tilde{E}_{j_2})(y)\,dy \right\rVert^2_{\LL^2(B_1)}\notag\\
    &\lesssim& a^8 \sum_{j=1 \atop j\neq m}^\aleph d_{mj}^{-8}  \left\lVert \overset{3}{\PP}(\tilde{E}_{j_1}) \right\rVert_{\LL^2(B_1)}^2 \, +\, a^{12} \sum_{j=1 \atop j\neq m}^\aleph d_{mj}^{-8}  \left\lVert \overset{3}{\PP}(\tilde{E}_{j_2}) \right\rVert_{\LL^2(B_2)}^2\notag\\
    &\lesssim& a^8 d_{\rm out}^{-8} \max_{j} \left\lVert \overset{3}{\PP}(\tilde{E}_{j_1}) \right\rVert_{\LL^2(B_1)}^2 \, +\, a^{12} d_{\rm out}^{-8} \max_{j} \left\lVert \overset{3}{\PP}(\tilde{E}_{j_2})\right\rVert_{\LL^2(B_2)}^2
\end{eqnarray}

\noindent (o) Estimate for $\sum_n\frac{1}{|1+\eta_1 \lambda_n^{(3)}|^2}\sum_{\ell=1 \atop \ell\neq m}^\aleph \left| \left\langle \widetilde{Err}_2\left( \overset{3}{\PP}(E_{\ell})\right), e_n^{(3)} \right\rangle \right|\cdot \sum_{i>\ell\atop i\neq m}^\aleph \left| \left\langle \widetilde{Err}_2\left( \overset{3}{\PP}(E_{i})\right), e_n^{(3)} \right\rangle \right|$:

Based on \eqref{def-Err_2-P13mj}, rewrite 
\begin{equation*}
    \widetilde{Err}_2\left( \overset{3}{\PP}(E_{j})\right)=\widetilde{Err}_2\left( \overset{3}{\PP}(E_{j_1})\right)+ \widetilde{Err}_2\left( \overset{3}{\PP}(E_{j_2})\right),
\end{equation*}
where
\begin{eqnarray*}
    \widetilde{Err}_2\left( \overset{3}{\PP}(E_{j_1})\right)&:=&\eta_1 k^2 a^4 \int_{B_{1}} \int_0^1  \underset{x}{\nabla}\Upsilon_k(z_{m_1}+tax, z_{j_1})\cdot\mathcal{P}(x, 0)\,dt\cdot\overset{3}{\PP}(\tilde{E}_{j_1})(y)\,dy\notag\\
     &+ & \eta_1 k^2 a^4 \int_{B_{1}} \int_0^1  \underset{y}{\nabla}\Upsilon_k(z_{m_1}, z_{j_1}+tay)\cdot\mathcal{P}(y, 0)\,dt\cdot\overset{3}{\PP}(\tilde{E}_{j_1})(y)\,dy,\notag\\
     \widetilde{Err}_2\left( \overset{3}{\PP}(E_{j_2})\right)&:=& \eta_2 k^2 a^4 \int_{B_{2}} \int_0^1  \underset{x}{\nabla}\Upsilon_k(z_{m_1}+tax, z_{j_2})\cdot\mathcal{P}(x, 0)\,dt\cdot\overset{3}{\PP}(\tilde{E}_{j_2})(y)\,dy\notag\\
     &+ & \eta_2 k^2 a^4 \int_{B_{2}} \int_0^1  \underset{y}{\nabla}\Upsilon_k(z_{m_1}, z_{j_2}+tay)\cdot\mathcal{P}(y, 0)\,dt\cdot\overset{3}{\PP}(\tilde{E}_{j_2})(y)\,dy.
\end{eqnarray*}
Hence, there holds
\begin{eqnarray*}
    &&\sum_n\frac{1}{|1+\eta_1 \lambda_n^{(3)}|^2}\sum_{\ell=1 \atop \ell\neq m}^\aleph \left| \left\langle \widetilde{Err}_2\left( \overset{3}{\PP}(E_{\ell})\right), e_n^{(3)} \right\rangle \right|\cdot \sum_{i>\ell\atop i\neq m}^\aleph \left| \left\langle \widetilde{Err}_2\left( \overset{3}{\PP}(E_{i})\right), e_n^{(3)} \right\rangle \right|\notag\\
    &=& \sum_n\frac{1}{|1+\eta_1 \lambda_n^{(3)}|^2}\sum_{\ell=1 \atop \ell\neq m}^\aleph \left| \left\langle \widetilde{Err}_2\left( \overset{3}{\PP}(E_{\ell_1})\right), e_n^{(3)} \right\rangle + \left\langle \widetilde{Err}_2\left( \overset{3}{\PP}(E_{\ell_2})\right), e_n^{(3)} \right\rangle \right| \notag\\
    &&\qquad\qquad\qquad\qquad\cdot \sum_{i>\ell\atop i\neq m}^\aleph \left| \left\langle \widetilde{Err}_2\left( \overset{3}{\PP}(E_{i_1})\right), e_n^{(3)} \right\rangle + \left\langle \widetilde{Err}_2\left( \overset{3}{\PP}(E_{i_2})\right), e_n^{(3)} \right\rangle \right|\notag\\
    &\lesssim& a^4 \sum_{\ell=1\atop \ell\neq m}^\aleph \left( \left\lVert \widetilde{Err}_2\left(\overset{3}{\PP}(E_{\ell_1})\right)\right\rVert_{\LL^2(B_1)} + \left\lVert \widetilde{Err}_2\left(\overset{3}{\PP}(E_{\ell_2})\right)\right\rVert_{\LL^2(B_1)} \right)\notag\\
    &&\qquad\qquad\qquad\qquad\cdot \sum_{i>\ell \atop i\neq m}^\aleph \left( \left\lVert \widetilde{Err}_2\left(\overset{3}{\PP}(E_{i_1})\right)\right\rVert_{\LL^2(B_1)} + \left\lVert \widetilde{Err}_2\left(\overset{3}{\PP}(E_{i_2})\right)\right\rVert_{\LL^2(B_1)} \right),
\end{eqnarray*}
where
\begin{eqnarray*}
    &&\left\lVert \widetilde{Err}_2\left(\overset{3}{\PP}(E_{j_1})\right)\right\rVert_{\LL^2(B_1)} = \left( \sum_n \left| \left\langle \widetilde{Err}_2\left(\overset{3}{\PP}(E_{j_1})\right), e_n^{(3)} \right\rangle \right|^2 \right)^{\frac{1}{2}}\notag\\
    &\lesssim& a^2 \left\lVert \int_{B_{1}} \int_0^1  \underset{x}{\nabla}\Upsilon_k(z_{m_1}+tax, z_{j_1})\cdot\mathcal{P}(x, 0)\,dt\cdot\overset{3}{\PP}(\tilde{E}_{j_1})(y)\,dy \right\rVert_{\LL^2(B_1)}\notag\\
    &+& a^2 \left\lVert \int_{B_{1}} \int_0^1  \underset{y}{\nabla}\Upsilon_k(z_{m_1}, z_{j_1}+tay)\cdot\mathcal{P}(y, 0)\,dt\cdot\overset{3}{\PP}(\tilde{E}_{j_1})(y)\,dy \right\rVert_{\LL^2(B_1)}\lesssim a^2 \, d_{mj}^{-4} \left\lVert \overset{3}{\PP}(\tilde{E}_{j_1})\right\rVert_{\LL^2(B_1)}
\end{eqnarray*}
and 
\begin{eqnarray*}
     &&\left\lVert \widetilde{Err}_2\left(\overset{3}{\PP}(E_{j_2})\right)\right\rVert_{\LL^2(B_1)} = \left( \sum_n \left| \left\langle \widetilde{Err}_2\left(\overset{3}{\PP}(E_{j_2})\right), e_n^{(3)} \right\rangle \right|^2 \right)^{\frac{1}{2}}\notag\\
     &\lesssim& a^4\left\lVert \int_{B_{2}} \int_0^1  \underset{x}{\nabla}\Upsilon_k(z_{m_1}+tax, z_{j_2})\cdot\mathcal{P}(x, 0)\,dt\cdot\overset{3}{\PP}(\tilde{E}_{j_2})(y)\,dy \right\rVert_{\LL^2(B_1)}\notag\\
     &+& a^4 \left\lVert \int_{B_{2}} \int_0^1  \underset{y}{\nabla}\Upsilon_k(z_{m_1}, z_{j_2}+tay)\cdot\mathcal{P}(y, 0)\,dt\cdot\overset{3}{\PP}(\tilde{E}_{j_2})(y)\,d \right\rVert_{\LL^2(B_1)}\lesssim a^4 \, d_{mj}^{-4} \left\lVert \overset{3}{\PP}(\tilde{E}_{j_2})\right\rVert_{\LL^2(B_2)}
\end{eqnarray*}
Therefore, we can obtain that
\begin{eqnarray}\label{es-P13-RHS-P3j-crossing}
    &&\sum_n\frac{1}{|1+\eta_1 \lambda_n^{(3)}|^2}\sum_{\ell=1 \atop \ell\neq m}^\aleph \left| \left\langle \widetilde{Err}_2\left( \overset{3}{\PP}(E_{\ell})\right), e_n^{(3)} \right\rangle \right|\cdot \sum_{i>\ell\atop i\neq m}^\aleph \left| \left\langle \widetilde{Err}_2\left( \overset{3}{\PP}(E_{i})\right), e_n^{(3)} \right\rangle \right|\notag\\
    &\lesssim& a^4 \sum_{\ell=1 \atop\ell\neq m}^\aleph \left(a^2 d_{m\ell}^{-4} \left\lVert \overset{3}{\PP}(\tilde{E}_{\ell_1}) \right\rVert_{\LL^2(B_1)} + a^4  d_{m\ell}^{-4}  \left\lVert \overset{3}{\PP}(\tilde{E}_{\ell_2}) \right\rVert_{\LL^2(B_2)} \right)\notag\\
    &&\qquad\qquad\qquad\cdot \sum_{i>\ell \atop i\neq m}^\aleph \left(a^2 d_{mi}^{-4} \left\lVert \overset{3}{\PP}(\tilde{E}_{i_1}) \right\rVert_{\LL^2(B_1)} + a^4  d_{mi}^{-4}  \left\lVert \overset{3}{\PP}(\tilde{E}_{i_2}) \right\rVert_{\LL^2(B_2)} \right)\notag\\
    &\lesssim& a^4 \left( a^2\, d_{\rm out}^{-4} \max_j \left\lVert \overset{3}{\PP}(\tilde{E}_{j_1}) \right\rVert_{\LL^2(B_1)} + a^4\, d_{\rm out}^{-4}  \max_j \left\lVert \overset{3}{\PP}(\tilde{E}_{j_2}) \right\rVert_{\LL^2(B_2)} \right)^2\notag\\
    &\lesssim& a^{8}\, d_{\rm out}^{-8}\,  \max_j \left\lVert \overset{3}{\PP}(\tilde{E}_{j_1}) \right\rVert_{\LL^2(B_1)}^2 + a^{12} \, d_{\rm out}^{-8}\,  \max_j \left\lVert \overset{3}{\PP}(\tilde{E}_{j_2}) \right\rVert_{\LL^2(B_2)}^2.
\end{eqnarray}

Combining with the estimates \eqref{es-P31-RHS-1}--\eqref{es-P13-RHS-P3j-crossing}, we can derive from \eqref{LS-Dm1-P31-norm-formulate} that
\begin{eqnarray*}
    &&\left\lVert\overset{3}{\PP}(\tilde{E}_{m_1})\right\rVert_{\LL^2(B_1)}^2 \notag\\
    &\lesssim& a^4 \left\lVert\tilde{E}_{m_1}^{Inc}\right\rVert^2_{\LL^2(B_1)} + a^8 \left\lVert\overset{1}{\PP}(\tilde{E}_{m_1})\right\rVert_{\LL^2(B_1)}^2 + 
 a^4 \left\lVert\overset{3}{\PP}(\tilde{E}_{m_1})\right\rVert_{\LL^2(B_1)}^2 + a^{12}\, d_{\rm in}^{-4} \left\lVert\overset{1}{\PP}(\tilde{E}_{m_2})\right\rVert_{\LL^2(B_2)}^2\notag\\
 &+& a^{8} \, d_{\rm out}^{-4} \max_j \left\lVert\overset{1}{\PP}(\tilde{E}_{j_1})\right\rVert_{\LL^2(B_1)}^2 + a^{12} \, d_{\rm out}^{-4} \max_j \left\lVert\overset{1}{\PP}(\tilde{E}_{j_2})\right\rVert_{\LL^2(B_2)}^2 + a^{10} \, d_{\rm in}^{-6} \left\lVert\overset{3}{\PP}(\tilde{E}_{m_2})\right\rVert_{\LL^2(B_2)}^2\notag\\
 &+& a^{6} \, d_{\rm out}^{-6} \max_j \left\lVert\overset{3}{\PP}(\tilde{E}_{j_1})\right\rVert_{\LL^2(B_1)}^2 + a^{10} \, d_{\rm out}^{-6} \max_j \left\lVert\overset{3}{\PP}(\tilde{E}_{j_2})\right\rVert_{\LL^2(B_2)}^2 + a^{14} \, d_{\rm in}^{-6} \left\lVert\overset{1}{\PP}(\tilde{E}_{m_2})\right\rVert_{\LL^2(B_2)}^2\notag\\
 &+& a^{12} \, d_{\rm in}^{-8} \max_j \left\lVert\overset{3}{\PP}(\tilde{E}_{m_2})\right\rVert_{\LL^2(B_2)}^2 + a^{10} \, d_{\rm out}^{-6} \max_j \left\lVert\overset{1}{\PP}(\tilde{E}_{j_1})\right\rVert_{\LL^2(B_1)}^2 + a^{14} \, d_{\rm out}^{-6} \max_j \left\lVert\overset{1}{\PP}(\tilde{E}_{j_2})\right\rVert_{\LL^2(B_2)}^2\notag\\
 &+&  a^{8} \, d_{\rm out}^{-8} \max_j \left\lVert\overset{3}{\PP}(\tilde{E}_{j_1})\right\rVert_{\LL^2(B_1)}^2 + a^{12} \, d_{\rm out}^{-8} \max_j \left\lVert\overset{3}{\PP}(\tilde{E}_{j_2})\right\rVert_{\LL^2(B_2)}^2\notag\\
 &+& a^{8} \, d_{\rm out}^{-6} \max_j \left\lVert\overset{1}{\PP}(\tilde{E}_{j_1})\right\rVert_{\LL^2(B_1)}^2 + a^{12} \, d_{\rm out}^{-6} \max_j \left\lVert\overset{1}{\PP}(\tilde{E}_{j_2})\right\rVert_{\LL^2(B_2)}^2
\end{eqnarray*}
which by taking $\max_m(\cdot)$ and $\max_j(\cdot)$ on the both sides gives 
\begin{eqnarray}\label{es-norm-P31}
      \left\lVert\overset{3}{\PP}(\tilde{E}_{m_1})\right\rVert_{\LL^2(B_1)}^2 &\lesssim& a^4 + a^{8} d_{\rm out}^{-6} \max_m \left\lVert\overset{1}{\PP}(\tilde{E}_{m_1})\right\rVert_{\LL^2(B_1)}^2 \notag\\
      &+& a^{12} \left( d_{\rm in}^{-4} +  d_{\rm out}^{-6} \right) \max_m \left\lVert\overset{1}{\PP}(\tilde{E}_{m_2})\right\rVert_{\LL^2(B_2)}^2 + a^{10}\, d_{\rm in}^{-6} \max_m \left\lVert\overset{3}{\PP}(\tilde{E}_{m_2})\right\rVert_{\LL^2(B_2)}^2.
\end{eqnarray}

\medskip

\noindent (3) \textit{Estimate for $\max_m\left\lVert \overset{1}{\PP}(\tilde{E}_{m_2})\right\rVert^2_{\LL^2(B_2)}$}.

\medskip

We start from the Lippmann-Schwinger equation on $D_{m_2}$ and project onto $\mathbb{H}_0(\rm{div}=0)$.
Using the decomposition of the electric field and the moment notation (cf. \textbf{Notation MP-1.2}), one arrives at the projected identity

\begin{eqnarray}\label{LS-Dm2-P12-Phi}
    \TT_{k, 2} \left(\overset{1}{\PP}(E_{m_2})\right)(x) &+& \TT_{k, 2} \left(\overset{3}{\PP}(E_{m_2})\right)(x) - k^2 \eta_1 N^k\left(\overset{1}{\PP}(E_{m_1})\right)(x) + \eta_1 \left(\nabla M^k-k^2 N^k\right)\left(\overset{3}{\PP}(E_{m_1})\right)(x)\notag\\
    &-& k^2 \eta_1\sum_{j=1 \atop j\neq m}^\aleph  N^k\left(\overset{1}{\PP}(E_{j_1})\right)(x) + \eta_1 \sum_{j=1 \atop j\neq m}^\aleph (\nabla M^k -k^2 N^k)\left(\overset{3}{\PP}(E_{j_1})\right)(x)\notag\\
    &-& k^2 \eta_2 \sum_{j=1 \atop j\neq m}^\aleph N^k\left(\overset{1}{\PP}(E_{j_2})\right)(x) + \eta_2 \sum_{j=1 \atop j\neq m}^\aleph (\nabla M^k- k^2 N^k)\left(\overset{3}{\PP}(E_{j_2})\right)(x)= E_{m_2}^{Inc}(x).
\end{eqnarray}

The interaction terms are then expanded around the centers by Taylor's formula, and the resulting Taylor remainders are estimated using the separation assumptions and the bounds on $\nabla \Phi_k$, $\operatorname{Hess} \Phi_k$ and $\nabla\Upsilon_k$.
After rescaling to $B_1,B_2$ and applying Cauchy--Schwarz in the modal index, we obtain the norm inequality

\begin{eqnarray}\label{eq-P12-norm-simple}
    \left\lVert \overset{1}{\PP}(\tilde{E}_{m_2}) \right\rVert_{\LL^2(B_2)}^2 &\lesssim& a^2 + a^{8} \left\lVert\overset{1}{\PP}(\tilde{E}_{m_2})\right\rVert_{\LL^2(B_2)}^2 + a^{8} \left\lVert\overset{3}{\PP}(\tilde{E}_{m_2})\right\rVert_{\LL^2(B_2)}^2 + a^{6} d_{\rm in}^{-6} \left\lVert\overset{1}{\PP}(\tilde{E}_{m_1})\right\rVert_{\LL^2(B_1)}^2\notag\\
    &+& a^{4} d_{\rm in}^{-4} \left\lVert\overset{3}{\PP}(\tilde{E}_{m_1})\right\rVert_{\LL^2(B_1)}^2 + a^{6} d_{\rm out}^{-6} \max_j \left\lVert\overset{1}{\PP}(\tilde{E}_{j_1})\right\rVert_{\LL^2(B_1)}^2 + a^{4} d_{\rm out}^{-4} \max_j \left\lVert\overset{3}{\PP}(\tilde{E}_{j_1})\right\rVert_{\LL^2(B_1)}^2\notag\\
    &+& a^{6} d_{\rm out}^{-6} \ln^2 d_{\rm out} \max_j \left\lVert\overset{1}{\PP}(\tilde{E}_{j_1})\right\rVert_{\LL^2(B_1)}^2 + a^{4} d_{\rm out}^{-6} \max_j \left\lVert\overset{3}{\PP}(\tilde{E}_{j_1})\right\rVert_{\LL^2(B_1)}^2 \notag\\
    &+& a^{10} d_{\rm out}^{-6} \max_j \left\lVert\overset{1}{\PP}(\tilde{E}_{j_2})\right\rVert_{\LL^2(B_2)}^2 + a^{10} d_{\rm out}^{-6} \ln^2 d_{\rm out} \max_j \left\lVert\overset{1}{\PP}(\tilde{E}_{j_2})\right\rVert_{\LL^2(B_2)}^2\notag\\
    &+& a^{8} d_{\rm out}^{-4} \max_j \left\lVert\overset{3}{\PP}(\tilde{E}_{j_2})\right\rVert_{\LL^2(B_2)}^2 + a^{8} d_{\rm out}^{-6} \max_j  \left\lVert\overset{3}{\PP}(\tilde{E}_{j_2})\right\rVert_{\LL^2(B_2)}^2 + a^{8} d_{\rm in}^{-8} \left\lVert\overset{1}{\PP}(\tilde{E}_{m_1})\right\rVert_{\LL^2(B_1)}^2\notag\\
    &+& a^{6} d_{\rm in}^{-6} \left\lVert\overset{3}{\PP}(\tilde{E}_{m_1})\right\rVert_{\LL^2(B_1)}^2+ a^{8} d_{\rm out}^{-8}\max_j \left\lVert\overset{1}{\PP}(\tilde{E}_{j_1})\right\rVert_{\LL^2(B_1)}^2 + a^{6} d_{\rm out}^{-6} \max_j \left\lVert\overset{3}{\PP}(\tilde{E}_{j_1})\right\rVert_{\LL^2(B_1)}^2\notag\\
    &+&  a^{6} d_{\rm out}^{-6} \ln^2 d_{\rm out} \max_j \left\lVert\overset{3}{\PP}(\tilde{E}_{j_1})\right\rVert_{\LL^2(B_1)}^2 + a^{8} d_{\rm out}^{-8}  \max_j \left\lVert\overset{3}{\PP}(\tilde{E}_{j_1})\right\rVert_{\LL^2(B_1)}^2\notag\\
    &+& a^{12} d_{\rm out}^{-8} \max_j \left\lVert\overset{1}{\PP}(\tilde{E}_{j_2})\right\rVert_{\LL^2(B_2)}^2+ a^{10} d_{\rm out}^{-6} \max_j \left\lVert\overset{3}{\PP}(\tilde{E}_{j_2})\right\rVert_{\LL^2(B_2)}^2\notag\\
    &+& a^{10} d_{\rm out}^{-6} \ln^2 d_{\rm out} \max_j \left\lVert\overset{3}{\PP}(\tilde{E}_{j_2})\right\rVert_{\LL^2(B_2)}^2+ a^{12} d_{\rm out}^{-8} \max_j \left\lVert\overset{3}{\PP}(\tilde{E}_{j_2})\right\rVert_{\LL^2(B_2)}^2.
\end{eqnarray}

Taking $\max_j$ and $\max_m$ on both sides yields

\begin{eqnarray}\label{es-norm-P1m2}
    \max_{m} \left\lVert \overset{1}{\PP}(\tilde{E}_{m_2})\right\rVert_{\LL^2(B_2)}^2&\lesssim& a^{2} + a^{6} d_{\rm in}^{-6} \max_m \left\lVert\overset{1}{\PP}(\tilde{E}_{m_1})\right\rVert_{\LL^2(B_1)}^2 + a^{4}\left(d_{\rm in}^{-4} + d_{\rm out}^{-6}\right) \max_m \left\lVert\overset{3}{\PP}(\tilde{E}_{m_1})\right\rVert_{\LL^2(B_1)}^2\notag\\
    &+& a^{8} d_{\rm out}^{-6}\max_m\left\lVert\overset{3}{\PP}(\tilde{E}_{m_2})\right\rVert_{\LL^2(B_2)}^2.
\end{eqnarray}

\medskip

\noindent (4). \textit{Estimate for $\max_m \left\lVert \overset{3}{\PP}(\tilde{E}_{m_2}) \right\rVert^2_{\LL^2(B_2)}$}.

\medskip

Projecting the Lippmann-Schwinger equation on $D_{m_2}$ onto $\nabla \mathcal{H}armonic$ and applying $\TT_{k,2}^{-1}$ leads, after Taylor expansion around the centers and collecting the Taylor remainders as $Err_4(\cdot)$, to the following identity
\begin{eqnarray}\label{Ls-Dm2-P32-taylor}
    && \left(I+\eta_2 \nabla M \right)\overset{3}{\PP}(E_{m_2})\notag\\
    &=& E_{m_2}^{Inc} - \overset{1}{\PP}(E_{m_2}) + k^2 \eta_2 N^k \left(\overset{1}{\PP}(E_{m_2})\right) + k^2 \eta_1 N^k \left(\overset{3}{\PP}(E_{m_2})\right) - \eta_2 \left( \nabla M^k -\nabla M \right)\overset{3}{\PP}(E_{m_2})\notag\\
    &+& k^2 \eta_1 \underset{y}{\nabla} (\Phi_kI)(z_{m_2}, z_{m_1}) \int_{D_{m_1}}\mathcal{P}(y, z_{m_1})\cdot\overset{1}{\PP}(E_{m_1})(y)\,dy + k^2 \eta_1 \sum_{j=1\atop j\neq m}^\aleph \underset{y}{\nabla}(\Phi_k I) (z_{m_2}, z_{j_1}) \int_{D_{j_1}}\mathcal{P}(y, z_{j_1}) \overset{1}{\PP}(E_{j_1})(y)\,dy\notag\\
    &+& k^2 \eta_2 \sum_{j=1\atop j\neq m}^\aleph \underset{y}{\nabla}(\Phi_k I)(z_{m_2}, z_{j_2}) \int_{D_{j_2}} \mathcal{P}(y, z_{j_2}) \overset{1}{\PP}(E_{j_2})(y)\,dy + k^2 \eta_1 \Upsilon_k(z_{m_2}, z_{m_1}) \int_{D_{m_1}} \overset{3}{\PP}(E_{m_1})(y)\,dy\notag\\
    &+& k^2 \eta_1 \sum_{j=1 \atop j\neq m}^\aleph \Upsilon_k(z_{m_2}, z_{j_1})\int_{D_{j_1}} \overset{3}{\PP}(E_{j_1})(y)\,dy + k^2 \eta_2 \sum_{j=1 \atop j\neq m}^\aleph\Upsilon_k(z_{m_2}, z_{j_2}) \int_{D_{j_2}} \overset{3}{\PP}(E_{j_2})(y)\,dy\notag\\ 
    &+& Err_4\left( \overset{1}{\PP}(E_{m_1}) \right) + Err_4\left( \overset{3}{\PP}(E_{m_1}) \right) + \sum_{j=1 \atop j\neq m}^\aleph Err_4\left( \overset{1}{\PP}(E_{j}) \right) + \sum_{j=1 \atop j\neq m}^\aleph Err_4\left( \overset{3}{\PP}(E_{j}) \right),\notag\\
\end{eqnarray}

The Taylor remainder $Err_4(\cdot)$ is controlled by the same kernel derivative bounds and the separation/scaling assumptions; hence it contributes higher-order terms that can be absorbed.
Proceeding as above (rescaling to $B_1,B_2$, taking $\LL^2$-inner products against the modes, and summing in the modal index) gives the estimate
\begin{eqnarray}\label{es-norm-P32}
      \left\lVert\overset{3}{\PP}(\tilde{E}_{m_2})\right\rVert_{\LL^2(B_2)}^2 &\lesssim& a^{-2h} + a^{6-4h} d_{\rm out}^{-9} + a^{4-2h} \left( d_{\rm in}^{-4} + d_{\rm out}^{-6} \right) \max_m \left\lVert\overset{1}{\PP}(\tilde{E}_{m_1})\right\rVert_{\LL^2(B_1)}^2 \notag\\
      &+&  \left(a^{8-2h}\, d_{\rm out}^{-6} + a^{22-6h} d_{\rm in}^{-8} d_{\rm out}^{-15} ) \max_m \left\lVert\overset{1}{\PP}(\tilde{E}_{m_2}\right)\right\rVert_{\LL^2(B_2)}^2 \notag\\
      &+& \left(a^{2-2h}\, d_{\rm in}^{-6} + a^{2-2h} d_{\rm out}^{-6} \ln^2 d_{\rm out} + a^{10-4h} d_{\rm in}^{-8} d_{\rm out}^{-9} \right) \max_m \left\lVert\overset{3}{\PP}(\tilde{E}_{m_1})\right\rVert_{\LL^2(B_1)}^2,
\end{eqnarray}
provided that
\begin{equation}\label{condition-es-P32}
   \frac{9}{5} < h < \min\{ 2, 5-6t_1 \}.
\end{equation}

Therefore, combining with the estimates \eqref{es-norm-P1m1}, \eqref{es-norm-P31}, \eqref{es-norm-P1m2} and \eqref{es-norm-P32} derived above, under the sufficient condition that
\begin{equation}\label{condition-a_priori es}
    \frac{9}{5} < h <\min\{ 2, 5-6t_1, 4-4t_1\}\quad \mbox{and}\quad \frac{k^4}{(4\pi)^4 c_0^2 d_0^2}<1,
\end{equation} 
we can obtain 
\begin{align*}
   & \max_m \left\lVert\overset{1}{\PP}(\tilde{E}_{m_1}) \right\rVert_{\LL^2(B_1)}^2\leq a^{4-2h} d_{\rm out}^{-\frac{9}{2}} && \max_m \left\lVert\overset{1}{\PP}(\tilde{E}_{m_2}) \right\rVert_{\LL^2(B_2)}^2\leq a^{7-2h} d_{\rm in}^{-3} d_{\rm out}^{-\frac{9}{2}} + a^{7-2h} d_{\rm out}^{-\frac{17}{2}} \\
   & \max_m \left\lVert\overset{3}{\PP}(\tilde{E}_{m_1}) \right\rVert_{\LL^2(B_1)}^2\leq a^{8-2h} d_{\rm in}^{-3} d_{\rm out}^{-\frac{9}{2}} && \max_m \left\lVert\overset{3}{\PP}(\tilde{E}_{m_2}) \right\rVert_{\LL^2(B_2)}^2\leq  a^{3-2h} d_{\rm out}^{-\frac{9}{2}},
\end{align*} 
which justifies (MP3.7) in Proposition MP-3.4. Note that the sufficient condition \eqref{condition-a_priori es} meets \eqref{cond-on-t2} and \eqref{condition-es-P32}.

\section{Derivation of the linear algebraic system and corresponding invertibility}\label{sec-proof-la}

\subsection{Derivation of the discrete system (summary)}\label{sec-derivation-summary}

The full derivation of the discrete system (MP1.24) (equivalently, (MP4.4)) from the Lippmann-Schwinger formulation is presented in the \emph{full} supplementary material, and the same Taylor/moment expansions are also used inside the proof of Proposition MP-3.4 in Section~\ref{sec-proof-prior} above.
In this (reduced) supplement we only record the structural output needed for the invertibility analysis:

\begin{itemize}
    \item Testing the Lippmann-Schwinger equation on each inclusion against the projector kernels and using the moment definitions leads to the unknown vector
    \(
        \mathfrak{U}_m=(Q_{m_1},R_{m_1},Q_{m_2},R_{m_2})^\mathsf{T}.
    \)
    \item The resulting algebraic system has the block form
    \(
        \mathfrak{B}_m \mathfrak{U}_m-\sum_{j\neq m}\Psi_{mj}\mathfrak{U}_j=\mathrm{S}_m,
    \)
    where the blocks $\mathfrak{B}_m$, $\Psi_{mj}$ and the source term $\mathrm{S}_m$ are those of \textbf{Notation MP-1.2}.
    \item The remainders produced by the Taylor expansions are controlled under the separation assumptions on $d_{\rm in},d_{\rm out}$ and the scaling constraints on $(t_1,t_2,h)$; these bounds are exactly the ones used in the proof of Proposition MP-3.4 and in the full supplement.
\end{itemize}

We now focus on the key point needed by the main paper: an explicit quantitative condition ensuring the invertibility of (MP1.24).

\subsection{Invertibility of the algebraic system (MP1.24)}\label{Inv-A-S-119}

This section is devoted to justify the invertibility of the linear algebraic system (MP1.24). First, we rewrite (MP1.24) as
\begin{equation}\label{system-matrix form}
    \begin{pmatrix}
        \mathfrak{B}_1 & 0 & \cdots & 0\\
        0 & \mathfrak{B}_2 & \cdots & 0\\
        \vdots & \vdots & \ddots & \vdots\\
        0 & 0 & \cdots & \mathfrak{B}_\aleph
    \end{pmatrix}\begin{pmatrix}
        \mathfrak{U}_1\\
        \mathfrak{U}_2\\
        \vdots\\
        \mathfrak{U}_\aleph
    \end{pmatrix}-\begin{pmatrix}
                 0 & - \, \Psi_{12} & \cdots & - \, \Psi_{1\aleph}  \\
                 - \, \Psi_{21} & 0 & \cdots & - \, \Psi_{2\aleph} \\
                 \vdots & \vdots & \ddots & \vdots \\
                 - \, \Psi_{\aleph1} & - \, \Psi_{\aleph2} & \cdots &  0 
            \end{pmatrix}\begin{pmatrix}
        \mathfrak{U}_1\\
        \mathfrak{U}_2\\
        \vdots\\
        \mathfrak{U}_\aleph
    \end{pmatrix} = \begin{pmatrix}
                \mathrm{S}_{1} \\
                \mathrm{S}_{2} \\
                \vdots \\
                \mathrm{S}_{\aleph}  
            \end{pmatrix}.
\end{equation}
With \textbf{Notation MP-1.2}, for each $1\leq m\leq \aleph$, we can reformulate $\mathfrak{B}_m$ on the L.H.S. of \eqref{system-matrix form} in the way 
\begin{equation}\label{split-B_m}
    \mathfrak{B}_m=\underbrace{\begin{pmatrix}
        I_3 & 0 & -\mathcal{B}_{13} & 0\\
        0 & I_3 & 0 & 0\\
        0 & 0 & I_3 & 0\\
        0 & -\mathcal{B}_{42} & 0 & I_3
    \end{pmatrix}}_{\mathfrak{B}_m^d} + \underbrace{\begin{pmatrix}
        0 & 0 & 0 & -\mathcal{B}_{14}\\
        0 & 0 & -\mathcal{B}_{23} & -\mathcal{B}_{24}\\
        -\mathcal{B}_{31} & -\mathcal{B}_{32} & 0 & 0\\
        -\mathcal{B}_{41} & 0 & 0 & 0
    \end{pmatrix}}_{\mathfrak{B}_m^r},
\end{equation}
where we can obtain directly from (MP1.19) that
\begin{equation*}
    |\mathfrak{B}_m^r| \lesssim \max\{ a^{3-h} d_{\rm in}^{-2},  a^3 d_{\rm in}^{-3}  \}\cdot \left|\begin{pmatrix}
        0 & 0 & 0 & I_3\\
        0 & 0 & I_3 & I_3\\
        I_3 & I_3 & 0 & 0\\
        I_3 & 0 & 0 & 0
    \end{pmatrix}\right| \lesssim \max\{ a^{3-h} d_{\rm in}^{-2},  a^3 d_{\rm in}^{-3}\},
\end{equation*}
performing as the remainder term for $h$ and $t_1$ satisfy (MP1.21). Then, \eqref{split-B_m} can be simplified with only the dominant part as 
\begin{equation}\label{B_m-express-domin}
    \mathfrak{B}_m = \mathfrak{B}_m^d + \O(\max\{ a^{3-h} d_{\rm in}^{-2},  a^3 d_{\rm in}^{-3}\}).
\end{equation}
Since 
\begin{eqnarray}\label{B_m-domin-quadratic}
      (\mathfrak{B}_m^d \mathfrak{U}_m)^\mathsf{T}\cdot \mathfrak{U}_m &=& \left[ \begin{pmatrix}
        I_3 & 0 & -\mathcal{B}_{13} & 0\\
        0 & I_3 & 0 & 0\\
        0 & 0 & I_3 & 0\\
        0 & -\mathcal{B}_{42} & 0 & I_3
    \end{pmatrix} \begin{pmatrix}
        Q_{m_1}\\
        R_{m_1}\\
        Q_{m_2}\\
        R_{m_2}
    \end{pmatrix} \right]^\mathsf{T}\cdot \begin{pmatrix}
         Q_{m_1}\\
        R_{m_1}\\
        Q_{m_2}\\
        R_{m_2}
    \end{pmatrix}=\begin{pmatrix}
        Q_{m_1}-\mathcal{B}_{31}Q_{m_2}\\
        R_{m_1}\\
        Q_{m_2}\\
        -\mathcal{B}_{42}R_{m_1}+R_{m_2}
    \end{pmatrix}^\mathsf{T}\cdot\begin{pmatrix}
         Q_{m_1}\\
        R_{m_1}\\
        Q_{m_2}\\
        R_{m_2}
    \end{pmatrix}\notag\\
    &=& (Q_{m_1}-\mathcal{B}_{13}Q_{m_2})Q_{m_1} + |R_{m_1}|^2 + | Q_{m_2} |^2 + (-\mathcal{B}_{42} R_{m_1}+ R_{m_2})R_{m_2}\notag\\
    &=& | Q_{m_1} |^2 + | R_{m_1} |^2 + |Q_{m_2} |^2 + |R_{m_2} |^2 -\mathcal{B}_{13} Q_{m_2} Q_{m_1} - \mathcal{B}_{42} R_{m_1} R_{m_2},
\end{eqnarray}
with 
\begin{equation}\label{upper-B13 term}
    | \mathcal{B}_{13} Q_{m_2} Q_{m_1} |\lesssim a^{3-h} d_{\rm in}^{-3} \frac{k^4 \eta_0}{c_0} \left|{\bf P}_{0, 1}^{(1)} \right|\cdot |Q_{m_2} Q_{m_1} |\lesssim a^{3-h} d_{\rm in}^{-3} \frac{k^4 \eta_0}{2 c_0} \left|{\bf P}_{0, 1}^{(1)} \right| \left( |Q_{m_2}|^2 +|Q_{m_1}|^2 \right)
\end{equation}
and 
\begin{equation}\label{upper-B42 term}
     | \mathcal{B}_{42} R_{m_1} R_{m_2} |\lesssim a^{3-h} d_{\rm in}^{-3} \frac{k^2 \eta_2}{d_0} \left|{\bf P}_{0, 2}^{(2)} \right|\cdot |R_{m_2} R_{m_1} |\lesssim a^{3-h} d_{\rm in}^{-3} \frac{k^2 \eta_2}{2 d_0} \left|{\bf P}_{0, 2}^{(2)} \right| \left( |R_{m_2}|^2 +|R_{m_1}|^2 \right).
\end{equation}
Denote
\begin{equation}\label{def-alpha-beta}
    \alpha:= \frac{k^4 \eta_0}{2 c_0} a^{3-h} d_{\rm in}^{-3}  \left|{\bf P}_{0, 1}^{(1)} \right|\quad\mbox{and}\quad \beta:= \frac{k^2 \eta_2}{2 d_0}  a^{3-h} d_{\rm in}^{-3} \left|{\bf P}_{0, 2}^{(2)} \right|,
\end{equation}
then we can obtain from \eqref{B_m-domin-quadratic}, \eqref{upper-B13 term} and \eqref{upper-B42 term} that
\begin{eqnarray}\label{es-invert-lower bound}
    \left| (\mathfrak{B}_m\mathfrak{U}_m)^\mathsf{T}\cdot\mathfrak{U}_m \right|&\gtrsim& |Q_{m_1} |^2 + |R_{m_1}|^2 + |Q_{m_2}|^2 + |R_{m_2}|^2 - \alpha (|Q_{m_1}|^2 + |Q_{m_2}|^2) -\beta (|R_{m_1}|^2 + |R_{m_2}|^2) \notag\\
    &\gtrsim& (1-\max\{\alpha, \beta\}) \left(|Q_{m_1} |^2 + |R_{m_1}|^2 + |Q_{m_2}|^2 + |R_{m_2}|^2\right) \notag\\
    &\gtrsim& (1-\max\{\alpha, \beta\})|\mathfrak{U}_m|^2,
\end{eqnarray}
under the condition that
\begin{equation}\label{condition-invert-lower bound}
    \frac{k^2}{2} a^{3-h} d_{\rm in}^{-3} \left( \frac{k^2 \eta_0}{c_0}\left|{\bf P}_{0, 1}^{(1)} \right| +\frac{\eta_2}{d_0}\left| {\bf P}_{0, 2}^{(2)} \right| \right)<1.
\end{equation}
Together with \eqref{B_m-express-domin} and \eqref{es-invert-lower bound}, we can derive that
\begin{equation}\label{eq-invert-lower bound-final}
    |\mathfrak{B}_m| \simeq |\mathfrak{B}_m^d| \gtrsim 1-\max\{\alpha, \beta\}= 1- \frac{k^2}{2} a^{3-h} d_{\rm in}^{-3} \left( \frac{k^2 \eta_0}{c_0}\left|{\bf P}_{0, 1}^{(1)} \right| +\frac{\eta_2}{d_0}\left| {\bf P}_{0, 2}^{(2)} \right| \right).
\end{equation}

Next, we evaluate the second term on the L.H.S. of \eqref{system-matrix form}. For each entry $\Psi_{mj}$, $1\leq m\neq j\leq \aleph$, with (MP1.20), we can rewrite $\Psi_{mj}$ as
\begin{equation}\notag
    \Psi_{mj}={\bf P}_{cons}\cdot \hat{\Psi}_{mj},
\end{equation}
where 
\begin{equation}\label{def-P-constant}
{\bf P}_{cons}:=
    \begin{pmatrix}
        \frac{k^2 \eta_0}{\pm c_0} {\bf P}_{0, 1}^{(1)} & 0 & 0 & 0\\
        0 & k^2 {\bf P}_{0, 1}^{(2)} & 0 & 0 \\
        0 & 0 & k^2 \eta_2 {\bf P}_{0, 2}^{(1)} & 0\\
        0 & 0 & 0 & \frac{k^2 \eta_2}{\pm d_0}{\bf P}_{0, 2}^{(2)}
    \end{pmatrix}
\end{equation}
and 
\begin{equation}\notag
\hat{\Psi}_{mj}:=
    \begin{pmatrix}
        k^2 a^{3-h} \Upsilon_k(z_{m_1}, z_{j_1}) & a^{3-h} \nabla\Phi_k (z_{m_1}, z_{j_1})\times & k^2 a^{3-h} \Upsilon_k(z_{m_1}, z_{j_2}) & a^{3-h} \nabla\Phi_k(z_{m_1}, z_{j_2})\times\\
        a^3 \nabla\Phi_k(z_{m_1}, z_{j_1})\times & a^3 \Upsilon_k(z_{m_1}, z_{j_1}) & a^3 \nabla\Phi_k(z_{m_1}, z_{j_2})\times & a^3 \Upsilon_k(z_{m_1}, z_{j_2})\\
        k^2 a^5 \Upsilon_k(z_{m_2}, z_{j_1}) & a^5 \nabla\Phi_k(z_{m_2}, z_{j_1})\times & k^2 a^5\Upsilon_k(z_{m_2}, z_{j_2}) & a^5 \nabla\Phi_k(z_{m_2}, z_{j_2})\times\\
        a^{3-h} \nabla\Phi_k(z_{m_2}, z_{j_1})\times & a^{3-h} \Upsilon_k(z_{m_2}, z_{j_1}) & a^{3-h} \nabla\Phi_k(z_{m_2}, z_{j_2})\times & a^{3-h} \Upsilon_k(z_{m_2}, z_{j_2})
    \end{pmatrix}.
\end{equation}
Then \eqref{system-matrix form} becomes
\begin{equation}
    \begin{pmatrix}\label{system-split-2nd}
        \mathfrak{B}_1\mathfrak{U}_1\\
         \mathfrak{B}_2\mathfrak{U}_2\\
         \vdots\\
          \mathfrak{B}_\aleph\mathfrak{U}_\aleph
    \end{pmatrix} - \begin{pmatrix}
        0 & {\bf P}_{cons}\hat{\Psi}_{12} & \cdots & {\bf P}_{cons}\hat{\Psi}_{1\aleph}\\
        {\bf P}_{cons}\hat{\Psi}_{21} & 0 & \cdots & {\bf P}_{cons}\hat{\Psi}_{2\aleph}\\
        \vdots & \vdots & \ddots & \vdots \\
        {\bf P}_{cons}\hat{\Psi}_{\aleph1} & {\bf P}_{cons}\hat{\Psi}_{\aleph2} & \cdots & 0
    \end{pmatrix}\begin{pmatrix}
        \mathfrak{U}_1 \\
        \mathfrak{U}_2\\
        \vdots\\
        \mathfrak{U}_\aleph
    \end{pmatrix}=\begin{pmatrix}
                \mathrm{S}_{1} \\
                \mathrm{S}_{2} \\
                \vdots \\
                \mathrm{S}_{\aleph}  
            \end{pmatrix}.
\end{equation}
Taking the $\LL^2$-inner product of \eqref{system-split-2nd} w.r.t. $(\mathfrak{U}_1, \cdots, \mathfrak{U}_\aleph)^\mathsf{T}$, it gives
\begin{equation}\label{system-inner product}
    \sum_{m=1}^\aleph \left\langle \mathfrak{B}_m\mathfrak{U}_m, \mathfrak{U}_m \right\rangle - \sum_{m=1}^\aleph \left\langle \sum_{j=1\atop j\neq m}^\aleph {\bf P}_{cons}\hat{\Psi}_{mj}\mathfrak{U}_j, \mathfrak{U}_m \right\rangle =\sum_{m=1}^\aleph \left\langle \mathrm{S}_m, \mathfrak{U}_m \right\rangle.
\end{equation}
To extract the dominant terms of $\hat{\Psi}_{mj}$, we reformulate
\begin{equation}\notag
    \hat{\Psi}_{mj}=a^{3-h} \left( \hat{\Psi}_{mj, 0} +  \hat{\Psi}_{mj, k} + \hat{\Psi}_{mj, r_1} \right) + a^3 \hat{\Psi}_{mj, r_2},
\end{equation}
with
\begin{equation}\label{def-Phi-mj-0}
    \hat{\Psi}_{mj, 0}:=\begin{pmatrix}
        k^2 \Upsilon_0(z_{m_1}, z_{j_1}) & 0 & k^2 \Upsilon_0(z_{m_1}, z_{j_2}) & 0\\
        0 & 0 & 0 & 0\\
         0 & 0 & 0 & 0\\
         0 & \Upsilon_0(z_{m_2}, z_{j_1}) & 0 & \Upsilon_0(z_{m_2}, z_{j_2})
    \end{pmatrix},
\end{equation}
\begin{equation}\label{def-Phi-mj-k}
    \hat{\Psi}_{mj, k}:=\begin{pmatrix}
        k^2 (\Upsilon_k-\Upsilon_0)(z_{m_1}, z_{j_1}) & 0 & k^2 (\Upsilon_k-\Upsilon_0)(z_{m_1}, z_{j_2}) & 0\\
        0 & 0 & 0 & 0\\
         0 & 0 & 0 & 0\\
         0 & (\Upsilon_k-\Upsilon_0)(z_{m_2}, z_{j_1}) & 0 & (\Upsilon_k-\Upsilon_0)(z_{m_2}, z_{j_2})
    \end{pmatrix},
\end{equation}
\begin{equation}\label{def-Phi-mj-r1}
    \hat{\Psi}_{mj, r_1}:=\begin{pmatrix}
        0 & \nabla\Phi_k(z_{m_1}, z_{j_1})\times & 0 & \nabla\Phi_k(z_{m_1}, z_{j_2})\times \\
        0 & 0&0&0\\
        0&0&0&0\\
        \nabla\Phi_k(z_{m_2}, z_{j_1})\times & 0 & \nabla\Phi_k(z_{m_2}, z_{j_2})\times & 0
    \end{pmatrix},
\end{equation}
and
\begin{equation}\label{def-psi-mj-r2}
    \hat{\Psi}_{mj, r_2}:=\begin{pmatrix}
        0&0&0&0\\
        \nabla\Phi_k(z_{m_1}, z_{j_1})\times & \Upsilon_k(z_{m_1}, z_{j_1}) & \nabla\Phi_k(z_{m_1}, z_{j_2})\times & \Upsilon_k(z_{m_1}, z_{j_2})\\
        k^2 a^2 \Upsilon_k(z_{m_2}, z_{j_1}) & a^2 \nabla\Phi_k(z_{m_2}, z_{j_1})\times & k^2 a^2 \Upsilon_k(z_{m_2}, z_{j_2}) & a^2 \nabla\Phi_k(z_{m_2}, z_{j_2})\times\\
        0&0&0&0
    \end{pmatrix}.
\end{equation}
Then \eqref{system-inner product} equivalently indicates
\begin{eqnarray}\label{eq-system-split0kr12}
    &&\sum_{m=1}^\aleph\left\langle \mathfrak{B}_m\mathfrak{U}_m, \mathfrak{U}_m \right\rangle - a^{3-h} \Bigg( \sum_{m=1}^\aleph\left\langle \sum_{j=1 \atop j\neq m}^\aleph {\bf P}_{cons}\hat{\Psi}_{mj, 0}\mathfrak{U}_j, \mathfrak{U}_m \right\rangle + \sum_{m=1}^\aleph\left\langle \sum_{j=1 \atop j\neq m}^\aleph {\bf P}_{cons}\hat{\Psi}_{mj, k}\mathfrak{U}_j, \mathfrak{U}_m \right\rangle\notag\\
    &+& \sum_{m=1}^\aleph\left\langle \sum_{j=1 \atop j\neq m}^\aleph {\bf P}_{cons}\hat{\Psi}_{mj, r_1}\mathfrak{U}_j, \mathfrak{U}_m \right\rangle \Bigg) = \sum_{m=1}^\aleph\left\langle \mathrm{S}_m, \mathfrak{U}_m \right\rangle + a^3 \sum_{m=1}^\aleph\left\langle \sum_{j=1 \atop j\neq m}^\aleph {\bf P}_{cons}\hat{\Psi}_{mj, r_2}\mathfrak{U}_j, \mathfrak{U}_m \right\rangle.\notag\\
\end{eqnarray}
We first evaluate the R.H.S. of \eqref{eq-system-split0kr12}. Denote
\begin{equation}\label{def-invert-I1-0}
    \mathcal{I}_1:=a^3 \sum_{m=1}^\aleph\left\langle \sum_{j=1 \atop j\neq m}^\aleph {\bf P}_{cons}\hat{\Psi}_{mj, r_2}\mathfrak{U}_j, \mathfrak{U}_m \right\rangle.
\end{equation}
From \eqref{def-psi-mj-r2}, we can know that by taking the intermediate point $z_{m_0}$ in $D_m$ (resp. $z_{j_0}$ in $D_j$), there holds
\begin{equation}\label{es-psi-mj-r2}
    \left| \hat{\Psi}_{mj, r_2} \right|\lesssim \left| \Upsilon_k(z_{m_0}, z_{j_0}) \mathbb{I}_a \right|\quad\mbox{with}\quad \mathbb{I}_a:=\begin{pmatrix}
        0&0&0&0\\
        I_3 & I_3& I_3 & I_3\\
        k^2 a^2 I_3 & a^2 I_3 & k^2 a^2 I_3  & a^2 I_3\\
        0&0&0&0
    \end{pmatrix},
\end{equation}
which further implies
\begin{eqnarray}
    |\mathcal{I}_1|&\lesssim& a^3 \left| \sum_{m=1}^\aleph\left\langle \sum_{j=1\atop j\neq m}^\aleph {\bf P}_{cons} \Upsilon_k(z_{m_0}, z_{j_0})\mathbb{I}_a\mathfrak{U}_j, \mathfrak{U}_m \right\rangle \right|\notag\\
    &\lesssim&  \underbrace{ a^3 \left| \sum_{m=1}^\aleph\left\langle \sum_{j=1\atop j\neq m}^\aleph {\bf P}_{cons} \mathbb{I}_a \Upsilon_0(z_{m_0}, z_{j_0}) \mathfrak{U}_j, \mathfrak{U}_m \right\rangle \right|}_{\mathcal{I}_{11}} + \underbrace{ a^3 \left| \sum_{m=1}^\aleph\left\langle \sum_{j=1\atop j\neq m}^\aleph {\bf P}_{cons} (\Upsilon_k-\Upsilon_0)(z_{m_0}, z_{j_0})\mathbb{I}_a\mathfrak{U}_j, \mathfrak{U}_m \right\rangle \right|}_{\mathcal{I}_{12}}.\notag
\end{eqnarray}
For $\mathcal{I}_{11}$, we can know from (MP1.17) that since $\Upsilon_0(\cdot, \cdot)$ is harmonic, by the Mean Value Theorem, there holds
\begin{equation}\label{mvt-upsilon0}
    \Upsilon_0(z_{m_0}, z_{j_0})=\frac{1}{|B_{m_0}|}\frac{1}{|B_{j_0}|}\int_{B_{m_0}} \int_{B_{j_0}} \Upsilon_0(x, y)\,dy\,dx,
\end{equation}
where $B_{m_0}$ and $B_{j_0}$ are the balls with radius $d_{\rm out}/2$, centered respectively at $z_{m_0}$ and $z_{j_0}$. Denote 
\begin{equation}\label{notation-integral form1-invert}
    \Omega=\underset{m=1}{\overset{\aleph}{\cup}}B_{m_0},\quad \mathbb{P}_{a}(x)=\begin{cases}
        {\bf P}_{cons}\mathbb{I}_a &\mbox{in} \; B_{m_0}\\
        0 & \mbox{Otherwise}
    \end{cases},
    \quad
    \Gamma(x)=\begin{cases}
        \mathfrak{U}_m & \mbox{in} \; B_{m_0}\\
        0& \mbox{Otherwise}
    \end{cases},
\end{equation}
then combining with (MP1.17) and (MP3.4), we have
\begin{eqnarray*}
    |\mathcal{I}_{11}|&\lesssim& a^{3}d_{\rm out}^{-6} \left(\left|\int_{\Omega}\int_\Omega \left\langle\mathbb{P}_a \Upsilon_0(x, y) \Gamma(y), \Gamma(x)\right\rangle\,dy\,dx\right|+ \left| \sum_{m=1}^\aleph \int_{B_{m_0}}\int_{B_{m_0}} \left\langle {\bf P}_{cons} \mathbb{I}_a \Upsilon_0(x, y)\mathfrak{U}_m, \mathfrak{U}_m \right\rangle \,dy\,dx\right|\right)\notag\\
    &\lesssim& a^{3}d_{\rm out}^{-6} \left( |\mathbb{P}_a|\cdot\left| \int_\Omega \left\langle\nabla M(\Gamma)(x), \Gamma(x)\right\rangle\,dx \right| + | {\bf P}_{cons} \mathbb{I}_a| d_{\rm out}^{3} \left| \sum_{m=1}^\aleph \left\langle \mathfrak{U}_m, \mathfrak{U}_m \right\rangle \right| \right),
\end{eqnarray*}
where we use the fact that $ \lVert \nabla M\rVert_{\mathcal{L}(\mathbb{L}^2(\Omega); \mathbb{L}^2(\Omega))}=1$ and 
\begin{equation}\label{eq-JLMS-d^3}
    \int_{B_{m_1}} \int_{B_{m_1}} \Upsilon_0(x, y)\,dy\,dx=-\frac{\pi d^3}{18} I_3.
\end{equation}
see \cite[eq. (6.7)]{CGS}. Besides, from \eqref{def-P-constant} and \eqref{es-psi-mj-r2}, it gives
\begin{equation}\notag
    {\bf P}_{cons}\mathbb{I}_a=\begin{pmatrix}
        0&0&0&0\\
        k^2 {\bf P}_{0, 1}^{(2)} & k^2 {\bf P}_{0, 1}^{(2)} & k^2 {\bf P}_{0, 1}^{(2)} & k^2 {\bf P}_{0, 1}^{(2)}\\
        k^4 \eta_2 a^2 {\bf P}_{0, 2}^{(1)} & k^2 \eta_2 a^2 {\bf P}_{0, 2}^{(1)} & k^4 \eta_2 a^2 {\bf P}_{0, 2}^{(1)} & k^2 \eta_2 a^2 {\bf P}_{0, 2}^{(1)}\\
        0&0&0&0
\end{pmatrix}.
\end{equation}
Hence by keeping the dominant terms, we can further obtain that
\begin{eqnarray}\label{es-invert-I_11}
    |\mathcal{I}_{11}| &\lesssim& a^3 d_{\rm out}^{-6} \left( |\mathbb{P}_a| \left\lVert\nabla M\right\rVert_{\mathcal{L}(\LL^2(\Omega); \LL^2(\Omega))} \left\lVert \Gamma \right\rVert_{\LL^2(\Omega)}^2 + |{\bf P}_{cons}\mathbb{I}_a| \left|\left\langle \Gamma, \Gamma \right\rangle\right| \right)\notag\\
    &\lesssim& |{\bf P}_{cons}\mathbb{I}_a|a^3 d_{\rm out}^{-6} \left\lVert \Gamma\right\rVert^2_{\LL^2(\Omega)} \lesssim k^2 \left|{\bf P}_{0, 1}^{(2)}\right| a^3 d_{\rm out}^{-6} \left\lVert \Gamma\right\rVert^2_{\LL^2(\Omega)},
\end{eqnarray}
where we use
\begin{equation}\label{transform-dis-conti}
d_{\rm out}^3 \sum_{m=1}^\aleph \left\langle \mathfrak{U}_m, \mathfrak{U}_m \right\rangle \lesssim  |B_{m_0}| \sum_{m=1}^\aleph \left\langle \mathfrak{U}_m, \mathfrak{U}_m \right\rangle \lesssim  \sum_{m=1}^\aleph \int_{B_{m_0}} \left\langle \mathfrak{U}_m, \mathfrak{U}_m \right\rangle \,dx =  \left\langle \Gamma(x), \Gamma(x)\right\rangle.
\end{equation}
%%%%%%%%%%%%%%%%%%%%%%%%%%%%%%%%%%%%%%%%%%%%%%%%%%

For $\mathcal{I}_{12}$, since 
\begin{equation}\label{equiv-diff-Upsilon}
    (\Upsilon_k-\Upsilon_0)(x, y)=\frac{1}{k^2} \underset{x}{\nabla} \underset{x}{\nabla}\Phi_k -\underset{x}{\nabla} \underset{x}{\nabla} \Phi_0 +\Phi_k(x, y) I_3 \sim \Phi_0(x, y) I_3,
\end{equation}
and $\Phi_0(x, y)$ is harmonic, similar to \eqref{mvt-upsilon0}, we can get that
\begin{eqnarray}\label{es-invert-I-12}
    |\mathcal{I}_{12}| &\lesssim& a^3 \left| \sum_{m=1}^\aleph\left\langle \sum_{j=1\atop j\neq m}^\aleph {\bf P}_{cons} \Phi_0(z_{m_0}, z_{j_0})\mathbb{I}_a\mathfrak{U}_j, \mathfrak{U}_m \right\rangle \right|\notag\\
     &\lesssim& a^{3} d_{\rm out}^{-6} \left| \sum_{m=1}^\aleph \sum_{j=1\atop j\neq m}^\aleph  \int_{B_{m_0}}\int_{B_{j_0}} \left\langle{\bf P}_{cons} \mathbb{I}_a  \Phi_0(x, y) \mathfrak{U}_j, \mathfrak{U}_m \right\rangle \,dy\,dx \right|\notag\\
    &\lesssim& a^{3} d_{\rm out}^{-6} | {\bf P}_{cons}\mathbb{I}_a |\cdot \left| \int_{\Omega} \left\langle N\left( \Gamma\right)(x), \Gamma(x) \right\rangle \,dx \right| + a^{3} d_{\rm out}^{-1} | {\bf P}_{cons}\mathbb{I}_a | \sum_{m=1}^\aleph \left| \left\langle \mathfrak{U}_m, \mathfrak{U}_m \right\rangle \right|\notag\\
    &\lesssim& k^2 \left|{\bf P}_{0, 1}^{(2)}\right|  a^{3} d_{\rm out}^{-6}   \left(\lVert N\rVert_{\mathcal{L}(\mathbb{L}^2(\Omega); \LL^2(\Omega))} +  d_{\rm out}^2\right) \left\lVert \Gamma \right\rVert^2_{\LL^2(\Omega)} \notag\\
    &\lesssim&  k^2 \left|{\bf P}_{0, 1}^{(2)}\right| a^3 d_{\rm out}^{-6} \left\lVert \Gamma\right\rVert^2_{\LL^2(\Omega)},
\end{eqnarray}
by keeping the dominant terms, where we use the fact that
\begin{equation*}
    \int_{B_{m_0}}\int_{B_{m_0}}\Phi_0(x, y)\,dy\,dx =\frac{\pi d_{\rm out}^5}{60},
\end{equation*}
see \cite[eq. (6.11)]{CGS}. Then together with \eqref{es-invert-I_11} and \eqref{es-invert-I-12}, we can obtain that
\begin{equation}\notag
    |\mathcal{I}_1|\lesssim k^2 \left|{\bf P}_{0, 1}^{(2)}\right| a^3 d_{\rm out}^{-6} \left\lVert \Gamma\right\rVert^2_{\LL^2(\Omega)},
\end{equation}
which indicates that \eqref{eq-system-split0kr12} can be simplified as  
\begin{eqnarray}\label{eq-system-split0kr12- with error}
    &&\sum_{m=1}^\aleph\left\langle \mathfrak{B}_m\mathfrak{U}_m, \mathfrak{U}_m \right\rangle - a^{3-h} \Bigg( \sum_{m=1}^\aleph\left\langle \sum_{j=1 \atop j\neq m}^\aleph {\bf P}_{cons}\hat{\Psi}_{mj, 0}\mathfrak{U}_j, \mathfrak{U}_m \right\rangle + \sum_{m=1}^\aleph\left\langle \sum_{j=1 \atop j\neq m}^\aleph {\bf P}_{cons}\hat{\Psi}_{mj, k}\mathfrak{U}_j, \mathfrak{U}_m \right\rangle\notag\\
    &+& \sum_{m=1}^\aleph\left\langle \sum_{j=1 \atop j\neq m}^\aleph {\bf P}_{cons}\hat{\Psi}_{mj, r_1}\mathfrak{U}_j, \mathfrak{U}_m \right\rangle \Bigg) = \sum_{m=1}^\aleph\left\langle \mathrm{S}_m, \mathfrak{U}_m \right\rangle + \O\left( k^2 \left|{\bf P}_{0, 1}^{(2)}\right| a^3 d_{\rm out}^{-6} \left\lVert \Gamma\right\rVert^2_{\LL^2(\Omega)} \right).\notag\\
\end{eqnarray}

Now, we investigate the L.H.S. of \eqref{eq-system-split0kr12- with error} respectively as follows.

\medskip

\noindent (i) Estimate for 
\begin{equation*}
    \mathcal{I}_2:= a^{3-h} \sum_{m=1}^\aleph\left\langle \sum_{j=1 \atop j\neq m}^\aleph {\bf P}_{cons}\hat{\Psi}_{mj, 0}\mathfrak{U}_j, \mathfrak{U}_m \right\rangle.
\end{equation*}
Considering $\hat{\Psi}_{mj, 0}$ in \eqref{def-Phi-mj-0} with substitution by the intermediate points $(z_{m_0}, z_{j_0})$, we have
\begin{equation}\label{def-I-k,1}
    \hat{\Psi}_{mj, 0} \sim \Upsilon_0(z_{m_0}, z_{j_0}) \mathbb{I}_{k}\quad\mbox{with}\quad \mathbb{I}_{k}:=\begin{pmatrix}
        k^2 I_3 & 0 & k^2 I_3 & 0\\
        0&0&0&0\\
        0&0&0&0\\
        0 & I_3 & 0 & I_3
    \end{pmatrix}.
\end{equation}
Then from \eqref{eq-JLMS-d^3} and the notation \eqref{notation-integral form1-invert}, together with \eqref{eq-JLMS-d^3}, we can know by following the similar argument to \eqref{es-invert-I_11} that
\begin{eqnarray}\label{es-invert-I2}
    |\mathcal{I}_2| &\lesssim& a^{3-h}\left|\sum_{m=1}^\aleph\left\langle \sum_{j=1 \atop j\neq m}^\aleph {\bf P}_{cons} \mathbb{I}_{k} \Upsilon_0(z_{m_0}, z_{j_0}) \mathfrak{U}_j, \mathfrak{U}_m \right\rangle\right|\notag\\
    &\lesssim& a^{3-h} d_{\rm out}^{-6} \left| {\bf P}_{cons}\mathbb{I}_{k} \right| \left( \left| \int_\Omega\int_\Omega \left\langle\Upsilon_0(x, y) \Gamma(y), \Gamma(x)\right\rangle \,dx\,dy\right| + \sum_{m=1}^\aleph \left|\int_{B_{m_0}}\int_{B_{m_0}}  \left\langle \Upsilon_0(x, y)\mathfrak{U}_m, \mathfrak{U}_m \right\rangle \,dx\,dy\right|\right)\notag\\
    &\lesssim& k^2 \left( \frac{k^2 \eta_0}{c_0} \left| {\bf P}_{0, 1}^{(1)} \right| + \frac{\eta_2}{d_0} {\bf P}_{0, 2}^{(2)} \right) a^{3-h} d_{\rm out}^{-6} \left\lVert \Gamma \right\rVert^2_{\LL^2(\Omega)}.\notag\\
\end{eqnarray}

\medskip

\noindent (ii) Estimate for 
\begin{equation*}
    \mathcal{I}_3:=a^{3-h} \sum_{m=1}^\aleph\left\langle \sum_{j=1 \atop j\neq m}^\aleph {\bf P}_{cons}\hat{\Psi}_{mj, k}\mathfrak{U}_j, \mathfrak{U}_m \right\rangle.
\end{equation*}
In \eqref{def-Phi-mj-k}, by using the intermediate points $(z_{m_0}, z_{j_0})$ and \eqref{equiv-diff-Upsilon}, we have
\begin{equation*}
    \hat{\Psi}_{mj, k}\sim \Phi_0(z_{m_0}, z_{j_0})\mathbb{I}_{k},
\end{equation*}
where $\mathbb{I}_{k}$ is given in \eqref{def-I-k,1}. Hence, following the similar argument to \eqref{es-invert-I-12}, it yields
\begin{eqnarray}\label{es-invert-I3}
    |\mathcal{I}_3|&\lesssim& a^{3-h}\left|\sum_{m=1}^\aleph\left\langle \sum_{j=1 \atop j\neq m}^\aleph {\bf P}_{cons} \mathbb{I}_{k}\Phi_0(z_{m_0}, z_{j_0}) \mathfrak{U}_j, \mathfrak{U}_m \right\rangle\right|\notag\\
    &\lesssim& k^2 \left( \frac{k^2 \eta_0}{c_0} \left| {\bf P}_{0, 1}^{(1)} \right| + \frac{\eta_2}{d_0} {\bf P}_{0, 2}^{(2)} \right) a^{3-h} d_{\rm out}^{-6} \left\lVert \Gamma \right\rVert^2_{\LL^2(\Omega)}.
\end{eqnarray}

\medskip

\noindent (iii) Estimate for 
\begin{equation*}
    \mathcal{I}_4:=a^{3-h} \sum_{m=1}^\aleph\left\langle \sum_{j=1 \atop j\neq m}^\aleph {\bf P}_{cons}\hat{\Psi}_{mj, r_1}\mathfrak{U}_j, \mathfrak{U}_m \right\rangle.
\end{equation*}
By using the intermediate points $(z_{m_0}, z_{j_0})$, we can see directly from \eqref{def-Phi-mj-r1} that
\begin{equation}\label{split-nabla Phi_k}
    \hat{\Psi}_{mj, r_1}\sim \left(\nabla\Phi_k(z_{m_0}, z_{j_0})\times\right) \mathbb{I}_{0}\quad\mbox{with}\quad \mathbb{I}_0:=\begin{pmatrix}
        0 & I_3 & 0 & I_3\\
        0&0&0&0\\
        0&0&0&0\\
        I_3 & 0 & I_3 & 0
    \end{pmatrix}.
\end{equation}
Since 
\begin{equation}\label{split-nabla-Phi-k}
    \nabla\Phi_k(z_{m_0}, z_{j_0})=\nabla\Phi_0(z_{m_0}, z_{j_0})+(\nabla\Phi_k-\nabla\Phi_0)(z_{m_0}, z_{j_0}),
\end{equation}
with $\nabla\Phi_0(\cdot, \cdot)$ being harmonic and 
\begin{equation}\notag
    (\nabla\Phi_k -\nabla\Phi_0) (z_{m_0}, z_{j_0}) = \frac{ik}{4\pi} \int_0^1 i k t\cdot e^{i k t|z_{m_0}-z_{j_0} |} \,dt \cdot\frac{z_{m_0}-z_{j_0}}{|z_{m_0}-z_{j_0}|},
\end{equation} 
we can obtain that
\begin{eqnarray}\label{es-invert-I4}
    |\mathcal{I}_4|&\lesssim& a^{3-h} \left| \sum_{m=1}^\aleph\left\langle \sum_{j=1 \atop j\neq m}^\aleph {\bf P}_{cons} \mathbb{I}_0 \nabla\Phi_k(z_{m_0}, z_{j_0})\times \mathfrak{U}_j, \mathfrak{U}_m \right\rangle \right|\notag\\
    % &\lesssim& a^{3-h} \left| \sum_{m=1}^\aleph\left\langle \sum_{j=1 \atop j\neq m}^\aleph {\bf P}_{cons} \mathbb{I}_0 \nabla\Phi_0(z_{m_0}, z_{j_0})\times \mathfrak{U}_j, \mathfrak{U}_m \right\rangle \right| \notag\\
    % &+& a^{3-h} \left| \sum_{m=1}^\aleph\left\langle \sum_{j=1 \atop j\neq m}^\aleph {\bf P}_{cons} \mathbb{I}_0 (\nabla\Phi_k-\nabla\Phi_0)(z_{m_0}, z_{j_0})\times \mathfrak{U}_j, \mathfrak{U}_m \right\rangle \right| \notag\\
    &\lesssim& a^{3-h} d_{\rm out}^{-6} |{\bf P}_{cons}\mathbb{I}_0| \Bigg(\left| \int_\Omega\int_{\Omega} \left\langle \nabla\Phi_0(x, y)\times \Gamma(y), \Gamma(x) \right\rangle \,dy\,dx \right|\notag\\
    &+& \left| \sum_{m=1}^\aleph \int_{B_{m_0}} \int_{B_{m_0}} \left\langle \nabla\Phi_0(x, y)\times\mathfrak{U}_m, \mathfrak{U}_m \right\rangle\,dx\,dy \right| + d_{\rm out}^3 \sum_{m=1}^\aleph\left\langle \mathfrak{U}_m, \mathfrak{U}_m \right\rangle  \Bigg)\notag\\
    &\lesssim& k^2 \left( \frac{\eta_0}{c_0}\left| {\bf P}_{0, 1}^{(1)} \right| + \frac{\eta_2}{d_0} \left| {\bf P}_{0, 2}^{(2)} \right| \right) a^{3-h} d_{\rm out}^{-6} \left( \left\lVert \nabla N \right\rVert_{\mathcal{L}(\LL^2(\Omega); \LL^2(\Omega))} + 1 \right) \left\lVert \Gamma\right\rVert^2_{\LL^2(\Omega)}\notag\\
    &\lesssim& k^2 \left( \frac{\eta_0}{c_0}\left| {\bf P}_{0, 1}^{(1)} \right| + \frac{\eta_2}{d_0} \left| {\bf P}_{0, 2}^{(2)} \right| \right) a^{3-h} d_{\rm out}^{-6} \left\lVert \Gamma\right\rVert^2_{\LL^2(\Omega)}
\end{eqnarray}
due to the fact that 
\begin{equation*}
    \left|\int_{B_{m_0}} \int_{B_{m_0}}\nabla \Phi_0(x, y) \,dx\,dy \right| \lesssim d_{\rm out}^3, 
\end{equation*}
see \cite[eq. (1.13)]{GS}.

Moreover, for the first term on the L.H.S. of \eqref{eq-system-split0kr12- with error}, denote
\begin{equation*}
    \mathbb{B}:=\begin{pmatrix}
        \mathfrak{B}_1 & 0 & \cdots & 0\\
        0 & \mathfrak{B}_2 & \cdots & 0\\
        \vdots & \vdots & \ddots & \vdots\\
        0 & 0 & \cdots & \mathfrak{B}_\aleph
    \end{pmatrix}.
\end{equation*}
By virtue of \eqref{es-invert-lower bound} and \eqref{eq-invert-lower bound-final}, we can derive that
\begin{equation}\label{invert-lower bound}
    \left| \sum_{m=1}^\aleph \left\langle \mathfrak{B}_m \mathfrak{U}_m, \mathfrak{U}_m \right\rangle \right| \underset{\gtrsim}{\eqref{transform-dis-conti}} d_{\rm out}^{-3} \left| \left\langle \mathbb{B} \Gamma(x), \Gamma(x) \right\rangle \right| \underset{\gtrsim}{\eqref{eq-invert-lower bound-final}}d_{\rm out}^{-3} \left(1- \max\{\alpha, \beta\} \right)\left\lVert \Gamma\right\rVert^2_{\LL^2(\Omega)}.
\end{equation}
Combining with \eqref{es-invert-I2}, \eqref{es-invert-I3},  \eqref{es-invert-I4} and \eqref{invert-lower bound}, we can derive that it is sufficient to let
\begin{equation*}
    k^2 \left( \frac{\max\{1, k^2\}\eta_0}{c_0}\left| {\bf P}_{0, 1}^{(1)} \right| + \frac{\eta_2}{d_0} \left| {\bf P}_{0, 2}^{(2)} \right| \right) a^{3-h} d_{\rm out}^{-6} \left\lVert \Gamma\right\rVert^2_{\LL^2(\Omega)} < d_{\rm out}^{-3} \left(1- \max\{\alpha, \beta\} \right)\left\lVert \Gamma\right\rVert^2_{\LL^2(\Omega)},
    % d_{\rm out}^{-3} \left[1- \frac{k^2}{2} a^{3-h} d_{\rm in}^{-3} \left( \frac{k^2 \eta_0}{c_0}\left|{\bf P}_{0, 1}^{(1)} \right| +\frac{\eta_2}{d_0}\left| {\bf P}_{0, 2}^{(2)} \right| \right) \right]\left\lVert \Gamma\right\rVert^2_{\LL^2(\Omega)}
\end{equation*}
such that the algebraic system (MP1.24) is invertible. Therefore, we obtain under condition \eqref{condition-invert-lower bound} that
\begin{equation*}
    \frac{k^2 \left( \frac{\max\{1, k^2\}\eta_0}{c_0}\left| {\bf P}_{0, 1}^{(1)} \right| + \frac{\eta_2}{d_0} \left| {\bf P}_{0, 2}^{(2)} \right| \right)}{1- \frac{k^2}{2} a^{3-h} d_{\rm in}^{-3} \left( \frac{k^2 \eta_0}{c_0}\left|{\bf P}_{0, 1}^{(1)} \right| +\frac{\eta_2}{d_0}\left| {\bf P}_{0, 2}^{(2)} \right| \right)} a^{3-h} d_{\rm out}^{-3} <1.
\end{equation*}

\section{Appendix}

\subsection{Proof of Proposition \ref{prop-es-QR-domain}}\label{sec:Appendix-prop5.1}

For any fixed $m$, $1\leq m\leq \aleph$, the algebraic system (MP4.4) indicates \eqref{eq-al-fix m}, which can be written more explicitly with the notation presented in \textbf{Notation MP-1.2} as 
\begin{equation}\label{a-system-fix m-1}
    \mathfrak{B}_m \tilde{\mathfrak{U}}_m-\sum_{j=1\atop j\neq m}^\aleph \Psi_{mj} \tilde{\mathfrak{U}}_j = \begin{pmatrix}
     \frac{i \, k \, \eta_{0}}{\pm \, c_{0}} \, a^{3-h} \, {\bf P}_{0, 1}^{(1)} \cdot H^{Inc}(z_{m_1}) \\
     a^{3} \, {\bf P}_{0, 1}^{(2)} \cdot E^{Inc}(z_{m_1}) \\
     i \, k \, a^{5} \, {\bf P}_{0, 2}^{(1)} \cdot H^{Inc}(z_{m_2}) \\
     \frac{\eta_{2}}{\pm \, d_{0}} \, a^{3-h} \, {\bf P}_{0, 2}^{(2)} \cdot E^{Inc}(z_{m_2})
    \end{pmatrix} + \begin{pmatrix}
    Error_{\tilde{Q}_{m_1}} \\
     Error_{\tilde{R}_{m_1}} \\
     Error_{\tilde{Q}_{m_2}} \\
     Error_{\tilde{R}_{m_2}}
\end{pmatrix},
\end{equation}
% \begin{eqnarray}\label{a-system-fix m-1}
%      \begin{pmatrix}
%     I_{3} & 0 & - \, \mathcal{B}_{13} & - \, \mathcal{B}_{14} \\
%     0 & I_{3} & - \, \mathcal{B}_{23}  & -\, \mathcal{B}_{24} \\
%     - \, \mathcal{B}_{31} & - \, \mathcal{B}_{32} & I_{3} & 0 \\  
%     - \, \mathcal{B}_{41} & - \, \mathcal{B}_{42} & 0 & I_{3}
%     \end{pmatrix}_{m}  \begin{pmatrix}
%         \tilde{Q}_{m_1}\\
%         \tilde{R}_{m_1}\\
%         \tilde{Q}_{m_2}\\
%         \tilde{R}_{m_2}
%     \end{pmatrix}&-&\sum_{j=1\atop j\neq m}^\aleph \begin{pmatrix}
%         \mathcal{C}_{11} & \mathcal{C}_{12} & \mathcal{C}_{13} & \mathcal{C}_{14}\\
%         \mathcal{C}_{21} & \mathcal{C}_{22} & \mathcal{C}_{23} & \mathcal{C}_{24} \\
%         \mathcal{C}_{31} & \mathcal{C}_{32} & \mathcal{C}_{33} & \mathcal{C}_{34} \\
%         \mathcal{C}_{41} & \mathcal{C}_{42} & \mathcal{C}_{43} & \mathcal{C}_{44} \\
%     \end{pmatrix}_{mj} \begin{pmatrix}
%         \tilde{Q}_{j_1}\\
%         \tilde{R}_{j_1}\\
%         \tilde{Q}_{j_2}\\
%         \tilde{R}_{j_2}
%     \end{pmatrix}\notag\\
%     &=& \begin{pmatrix}
%      \frac{i \, k \, \eta_{0}}{\pm \, c_{0}} \, a^{3-h} \, {\bf P}_{0, 1}^{(1)} \cdot H^{Inc}(z_{m_1}) \\
%      a^{3} \, {\bf P}_{0, 1}^{(2)} \cdot E^{Inc}(z_{m_1}) \\
%      i \, k \, a^{5} \, {\bf P}_{0, 2}^{(1)} \cdot H^{Inc}(z_{m_2}) \\
%      \frac{\eta_{2}}{\pm \, d_{0}} \, a^{3-h} \, {\bf P}_{0, 2}^{(2)} \cdot E^{Inc}(z_{m_2})
%     \end{pmatrix} + \begin{pmatrix}
%     Error_{Q_{m_1}} \\
%      Error_{R_{m_1}} \\
%      Error_{Q_{m_2}} \\
%      Error_{R_{m_2}}
% \end{pmatrix},
% \end{eqnarray}
with $Error_{\tilde{Q}_{m_1}}$, $Error_{\tilde{R}_{m_1}}$, $Error_{\tilde{Q}_{m_2}}$ and $Error_{\tilde{R}_{m_2}}$ respectively given by 
% \eqref{es-norm-P1m1}, \eqref{es-norm-P31}, \eqref{es-norm-P1m2} and \eqref{es-norm-P32}.
 \begin{eqnarray}\label{Equa-Error-Qm1-1229}
    \nonumber
     Error_{\tilde{Q}_{m_1}}  \, & = & \O\Bigg( a^3\, +\, a^{4-h} \, + \, \left( a^{6-h} \, d_{out}^{-4} + a^{5} \, d_{out}^{-3} \, \left\vert \ln(d_{out}) \right\vert \right) \max_{m} \left\lVert \overset{1}{\PP}(\tilde{E}_{m_1})\right\rVert_{\LL^2(B_1)} \\ \nonumber &+& \, \left( a^{8-h} \, d_{in}^{-4} \, + \, a^{7} \, (d_{in}^{-3}+ d_{out}^{-3} |\ln (d_{out})|)  \right) \, \max_{m} \left\lVert \overset{1}{\PP}(\tilde{E}_{m_2})\right\rVert_{\LL^2(B_2)} \notag\\
    &+& \left( a^{2} \, + \, a^{4} \, d_{out}^{-3}\, + \, a^{5-h}\, d_{out}^{-3} |\ln (d_{out})|  \right) \, \max_{m} \left\lVert \overset{3}{\PP}(\tilde{E}_{m_1})\right\rVert_{\LL^2(B_1)}\notag \\ &+& \, \left( a^{7-h} \, d_{in}^{-3} + a^{6} \, d_{in}^{-2} \, + a^{6} d_{out}^{-3} \right) \, \max_{m} \left\lVert \overset{3}{\PP}(\tilde{E}_{m_2})\right\rVert_{\LL^2(B_2)}\Bigg),\notag\\
    Error_{\tilde{R}_{m_1}}  &=& \O\Bigg( a^{4} \, %+ \,  a^8 \, d_{in}^{-3} \left\lVert \overset{3}{\PP}(\tilde{E}_{m_2})\right\rVert_{\LL^2(B_2)} \, 
    + \, a^5 \, d_{out}^{-4}\, \max_j\left\lVert\overset{3}{\PP}(\tilde{E}_{j_1})\right\rVert_{\LL^2(B_1)} \, + \, a^{7} \, d_{in}^{-4} \max_j\left\lVert\overset{3}{\PP}(\tilde{E}_{j_2})\right\rVert_{\LL^2(B_2)}\notag\\
    &+& %a^9 \, d_{in}^{-2} \left\lVert\overset{1}{\PP}(\tilde{E}_{m_2})\right\rVert_{\LL^2(B_2)} \,+ 
    \, a^{6} \, d_{out}^{-3} \, \left\vert \ln\left(d_{out}\right) \right\vert \, \max_j \left\lVert \overset{1}{\PP}(\tilde{E}_{j_1})\right\rVert_{\LL^2(B_1)} \, +\, a^{8} \, (d_{in}^{-3}+ d_{out}^{-3} |\ln (d_{out})|) \max_j \left\lVert \overset{1}{\PP}(\tilde{E}_{j_2})\right\rVert_{\LL^2(B_2)}\Bigg),\notag\\
     Error_{\tilde{Q}_{m_2}} &=& \O\Bigg(a^{6} 
    + \, \left(a^{8} \, d_{in}^{-4} + a^{9-h} d_{in}^{-3}\right) \max_j\left\lVert \overset{1}{\PP}(\tilde{E}_{j_1})\right\rVert_{\LL^2(B_1)} \, + \, a^{10} \, d_{out}^{-4} \max_j\left\lVert \overset{1}{\PP}(\tilde{E}_{j_2})\right\rVert_{\LL^2(B_2)}\notag\\
    &+& a^{7} \, \left(d_{in}^{-3} + d_{out}^{-3} |\ln (d_{out})|\right) \max_j\left\lVert \overset{3}{\PP}(\tilde{E}_{j_1})\right\rVert_{\LL^2(B_1)} \, + \, a^{9} \, d_{out}^{-3} \, \left\vert \ln\left( d_{out} \right) \right\vert \, \max_j\left\lVert \overset{3}{\PP}(\tilde{E}_{j_2})\right\rVert_{\LL^2(B_2)}\Bigg),\notag\\
     Error_{\tilde{R}_{m_2}} & = & \O\Bigg(a^{3} \, + \, a^{4-h} \, + \,  %a^4 \, d_{in}^{-3} \left\lVert \overset{3}{\PP}(\tilde{E}_{m_1})\right\rVert_{\LL^2(B_1)} \, +
    \,\left( a^4 \, \left(d_{in}^{-3} + d_{out}^{-3} |\ln (d_{out})| \right) + a^{5-h} \, d_{in}^{-4} \right) \max_j\left\lVert \overset{3}{\PP}(\tilde{E}_{j_1})\right\rVert_{\LL^2(B_1)} \notag\\ &+& \, \left(a^{7-h} \, d_{out}^{-4} \, + \, a^{6} \, d_{out}^{-3} \, \left\vert \ln(d_{out}) \right\vert \right) \max_j\left\lVert \overset{3}{\PP}(\tilde{E}_{j_2})\right\rVert_{\LL^2(B_2)}\notag\\
    &+&\left(a^{5} \, \left(d_{in}^{-2} + d_{out}^{-3}\right) \, + \, a^{6-h} \left(d_{in}^{-3} +d_{out}^{-3} |\ln (d_{out})| \right) \right) \max_j\left\lVert \overset{1}{\PP}(\tilde{E}_{j_1})\right\rVert_{\LL^2(B_1)} \notag\\
    &+& \left( a^{7} \, d_{out}^{-3} \, +\, a^{8-h}\, d_{out}^{-3} |\ln (d_{out})| \right) \max_j\left\lVert \overset{1}{\PP}(\tilde{E}_{j_2})\right\rVert_{\LL^2(B_2)}\Bigg).
\end{eqnarray}

Since the first term on the R.H.S. of \eqref{a-system-fix m-1} is equivalent to 
\begin{equation*}
    \begin{pmatrix}
     \frac{i \, k \, \eta_{0}}{\pm \, c_{0}} \, a^{3-h} \, {\bf P}_{0, 1}^{(1)} \cdot H^{Inc}(z_{m_1}) \\
     a^{3} \, {\bf P}_{0, 1}^{(2)} \cdot E^{Inc}(z_{m_1}) \\
     i \, k \, a^{5} \, {\bf P}_{0, 2}^{(1)} \cdot H^{Inc}(z_{m_2}) \\
     \frac{\eta_{2}}{\pm \, d_{0}} \, a^{3-h} \, {\bf P}_{0, 2}^{(2)} \cdot E^{Inc}(z_{m_2})
    \end{pmatrix} = \underbrace{\begin{pmatrix}
    \frac{i \, k \, \eta_{0}}{\pm \, c_{0}} \, a^{3-h} \, {\bf P}_{0, 1}^{(1)} & 0 & 0 & 0 \\
    0 & a^{3} \, {\bf P}_{0, 1}^{(2)}& 0  &0 \\
    0 & 0 & i \, k \, a^{5} \, {\bf P}_{0, 2}^{(1)} & 0 \\  
    0 & 0 & 0 & \frac{\eta_{2}}{\pm \, d_{0}} \, a^{3-h} \, {\bf P}_{0, 2}^{(2)}
    \end{pmatrix}}_{{\bf P}_a}  \begin{pmatrix}
      H^{Inc}(z_{m_1}) \\
      E^{Inc}(z_{m_1}) \\
      H^{Inc}(z_{m_2}) \\
      E^{Inc}(z_{m_2})
    \end{pmatrix}, 
\end{equation*}
by multiplying ${\bf P}_a$ on the both sides of \eqref{a-system-fix m-1}, we have
\begin{equation}\label{a-sys-fix m-matrix}
    \mathfrak{B}_m^a \cdot \tilde{\mathfrak{U}}^a_m -\sum_{j=1 \atop j\neq m}^\aleph \Psi_{mj}^a \tilde{\mathfrak{U}}_j^a =\begin{pmatrix}
      H^{Inc}(z_{m_1}) \\
      E^{Inc}(z_{m_1}) \\
      H^{Inc}(z_{m_2}) \\
      E^{Inc}(z_{m_2})
    \end{pmatrix}+ \begin{pmatrix}
    Error_{Q_{m_1}}^a \\
     Error_{R_{m_1}}^a \\
     Error_{Q_{m_2}}^a \\
     Error_{R_{m_2}}^a
\end{pmatrix}=:\begin{pmatrix}
   L_{H, m_1} \\
     L_{E, m_1} \\
    L_{H, m_2} \\
     L_{E, m_2}
\end{pmatrix},
\end{equation}
where 
\begin{equation}\label{notation-ABU}
    \mathfrak{B}_m^a:={\bf P}_a^{-1} \mathfrak{B}_m {\bf P}_a, \quad \Psi_{mj}^a:= {\bf P}_a^{-1} \Psi_{mj} {\bf P}_a, \quad \tilde{\mathfrak{U}}_m^a:={\bf P}_a^{-1} \tilde{\mathfrak{U}}_m = \begin{pmatrix}
   \tilde{Q}_{m_1}^a\\
     \tilde{R}_{m_1}^a \\
   \tilde{Q}_{m_2}^a \\
     \tilde{R}_{m_2}^a
\end{pmatrix}
\end{equation}
and
\begin{equation}\notag
   \begin{pmatrix}
    Error_{Q_{m_1}}^a \\
     Error_{R_{m_1}}^a \\
     Error_{Q_{m_2}}^a \\
     Error_{R_{m_2}}^a
\end{pmatrix}= {\bf P}_a^{-1} \begin{pmatrix}
    Error_{\tilde{Q}_{m_1}} \\
     Error_{\tilde{R}_{m_1}} \\
     Error_{\tilde{Q}_{m_2}} \\
     Error_{\tilde{R}_{m_2}}
\end{pmatrix}.
\end{equation}

From the matrix form \eqref{a-sys-fix m-matrix}, we can extract the precise $\ell_2$-estimates for $(\tilde{Q}_{m_1}^a, \tilde{R}_{m_1}^a, \tilde{Q}_{m_2}^a, \tilde{R}_{m_2}^a)_{m=1}^\aleph$ and eventually for $(\tilde{Q}_{m_1}, \tilde{R}_{m_1}, \tilde{Q}_{m_2}, \tilde{R}_{m_2})_{m=1}^\aleph$ respectively as follows:

\medskip

\noindent (i) \emph{Estimate for $\lVert \tilde{Q}_{m_1}\rVert_{\ell_2}$}:

\medskip

By injecting the accurate forms of $\mathfrak{B}_m$ and $\Psi_{mj}$ given by \textbf{Notation MP-1.2}, we can get from \eqref{a-sys-fix m-matrix} the following quadratic form
\begin{eqnarray}\label{Qm1-quadratic}
    &&\sum_{m=1}^\aleph \tilde{Q}_{m_1}^a\cdot \tilde{Q}_{m_1}^a - k^4 a^5 \sum_{m=1}^\aleph \Upsilon_k(z_{m_1}, z_{m_2}){\bf P}_{0, 2}^{(1)} \tilde{Q}_{m_2}^a\cdot \tilde{Q}_{m_1}^a\notag\\
    &-&\frac{k \eta_2}{\pm i d_0} a^{3-h} \sum_{m=1}^\aleph \nabla\Phi_k(z_{m_1}, z_{m_2})\times{\bf P}_{0, 2}^{(2)}\tilde{R}_{m_2}^a\cdot\tilde{Q}_{m_1}^a\notag\\
    &-& \sum_{m=1}^\aleph \sum_{j=1\atop j\neq m}^\aleph \frac{k^4 \eta_0}{\pm c_0}a^{3-h} \Upsilon_k(z_{m_1}, z_{j_1}){\bf P}_{0, 1}^{(1)}\tilde{Q}_{j_1}^a\cdot\tilde{Q}_{m_1}^a-\sum_{m=1}^\aleph \sum_{j=1\atop j\neq m}^\aleph \frac{k}{i}a^3 \nabla\Phi_k(z_{m_1}, z_{j_1})\times{\bf P}_{0, 1}^{(2)}\tilde{R}_{j_1}^a\cdot\tilde{Q}_{m_1}^a\notag\\
    &-& \sum_{m=1}^\aleph \sum_{j=1\atop j\neq m}^\aleph k^4 a^5  \Upsilon_k(z_{m_1}, z_{j_2}){\bf P}_{0, 2}^{(1)}\tilde{Q}_{j_2}^a\cdot\tilde{Q}_{m_1}^a-\sum_{m=1}^\aleph \sum_{j=1\atop j\neq m}^\aleph \frac{ k \eta_2}{\pm i d_0} a^{3-h} \nabla\Phi_k(z_{m_1}, z_{j_2})\times{\bf P}_{0, 2}^{(2)}\tilde{R}_{j_2}^a\cdot\tilde{Q}_{m_1}^a\notag\\
    & =& \sum_{m=1}^\aleph L_{h, m_1}\cdot\tilde{Q}_{m_1}^a.
\end{eqnarray}
Now, we evaluate the L.H.S. of \eqref{Qm1-quadratic} term by term.

% \medskip 

\noindent (a) Estimate for 
\begin{eqnarray}\label{es-L1}
     L_1&:=& k^4 a^5 \sum_{m=1}^\aleph \Upsilon_k(z_{m_1}, z_{m_2}){\bf P}_{0, 2}^{(1)} \tilde{Q}_{m_2}^a\cdot \tilde{Q}_{m_1}^a\notag\\
     |L_1|&\lesssim& k^4 a^5 d_{\rm in}^{-3} |{\bf P}_{0, 2}^{(1)}| \sum_{m=1}^\aleph  |\tilde{Q}_{m_2}^a|\cdot |\tilde{Q}_{m_1}^a|\lesssim k^4 a^5 d_{\rm in}^{-3} |{\bf P}_{0, 2}^{(1)}| \left(\sum_{m=1}^\aleph  |\tilde{Q}_{m_2}^a|^2\right)^\frac{1}{2}\left(\sum_{m=1}^\aleph|\tilde{Q}_{m_1}^a|^2\right)^\frac{1}{2}.
\end{eqnarray}

\noindent (b) Estimate for 
\begin{eqnarray}\label{es-L2}
    L_2&:=&  \frac{k \eta_2}{\pm i d_0} a^{3-h} \sum_{m=1}^\aleph \nabla\Phi_k(z_{m_1}, z_{m_2})\times{\bf P}_{0, 2}^{(2)}\tilde{R}_{m_2}^a\cdot\tilde{Q}_{m_1}^a\notag\\
    |L_2|&\lesssim& \frac{k \eta_2}{d_0} a^{3-h} d_{\rm in}^{-2} |{\bf P}_{0, 2}^{(2)}| \sum_{m=1}^\aleph |\tilde{R}_{m_2}^a|\cdot|\tilde{Q}_{m_1}^a|\lesssim \frac{k \eta_2}{d_0} a^{3-h} d_{\rm in}^{-2} |{\bf P}_{0, 2}^{(2)}| \left(\sum_{m=1}^\aleph |\tilde{R}_{m_2}^a|^2\right)^\frac{1}{2}\left(\sum_{m=1}^\aleph|\tilde{Q}_{m_1}^a|^2\right)^\frac{1}{2}.\notag\\
\end{eqnarray}

\noindent (c) Estimate for 
\begin{equation*}
   L_3:=  \sum_{m=1}^\aleph \sum_{j=1\atop j\neq m}^\aleph \frac{k^4 \eta_0}{\pm c_0}a^{3-h} \Upsilon_k(z_{m_1}, z_{j_1}){\bf P}_{0, 1}^{(1)}\tilde{Q}_{j_1}^a\cdot\tilde{Q}_{m_1}^a.
\end{equation*}
   % Write 
   % \begin{equation}\label{split-Upsilon}
   %     \Upsilon_k(z_{m_1}, z_{j_1})=\Upsilon_0(z_{m_1}, z_{j_1})+(\Upsilon_k-\Upsilon_0)(z_{m_1}, z_{j_1}).
   % \end{equation}
% Since $\Upsilon_0(\cdot, \cdot)$ is harmonic, 
% by the Mean Value Theorem, we have
% \begin{equation}\label{mvt-upsilon0}
%     \Upsilon_0(z_{m_1}, z_{j_1})=\frac{1}{|B_{m_1}|}\frac{1}{|B_{j_1}|}\int_{B_{m_1}} \int_{B_{j_1}} \Upsilon_0(x, y)\,dy\,dx,
% \end{equation}
By using \eqref{mvt-upsilon0}, similar to $\mathcal{I}_1$, we can write
\begin{eqnarray}\label{def-L3}
    L_3&=& \underbrace{\frac{k^4 \eta_0}{\pm c_0} a^{3-h} \sum_{m=1}^\aleph \sum_{j=1\atop j\neq m}^\aleph \frac{1}{|B_{m_1}|}\frac{1}{|B_{j_1}|}\int_{B_{m_1}} \int_{B_{j_1}} \Upsilon_0(x, y)\,dy\,dx\cdot {\bf P}_{0, 1}^{(1)}\tilde{Q}_{j_1}^a\cdot \tilde{Q}_{m_1}^a}_{L_3^0}\notag\\
    &+& \underbrace{\frac{k^4 \eta_0}{\pm c_0}a^{3-h}\sum_{m=1}^\aleph \sum_{j=1\atop j\neq m}^\aleph  (\Upsilon_k-\Upsilon_0)(z_{m_1}, z_{j_1}){\bf P}_{0, 1}^{(1)}\tilde{Q}_{j_1}^a\cdot\tilde{Q}_{m_1}^a}_{L_3^r},
\end{eqnarray}
where $B_{m_1}, B_{j_1}$ are two balls with radius $\frac{d_{\rm out}}{2}$, centered respectively at $z_{m_1}, z_{j_1}$. Denote 
\begin{equation}\label{notation-integral form1}
    \Omega=\underset{m=1}{\overset{\aleph}{\cup}}B_{m_1},\quad \Lambda_{1}(x)=\begin{cases}
        \tilde{Q}_{m_1}^a &\mbox{in} \; B_{m_1}\\
        0 & \mbox{Otherwise}
    \end{cases},
    \quad
    \mathbb{P}_1^{(1)}=\begin{cases}
        {\bf P}_{0, 1}^{(1)} & \mbox{in} \; B_{m_1}\\
        0& \mbox{Otherwise}
    \end{cases},
\end{equation}
then 
\begin{eqnarray*}
    L_3^0 &=& \frac{36}{\pi}\frac{k^4 \eta_0}{\pm c_0}a^{3-h}d_{\rm out}^{-6} \int_{\Omega}\int_\Omega \Upsilon_0(x, y) \mathbb{P}_{1}^{(1)}\cdot\Lambda_1(y) \cdot\Lambda_1(x)\,dy\,dx\\
    &-&\frac{36}{\pi}\frac{k^4 \eta_0}{\pm c_0}a^{3-h}d_{\rm out}^{-6} \sum_{m=1}^\aleph \int_{B_{m_1}}\int_{B_{m_1}} \Upsilon_0(x, y) \,dy\,dx\cdot {\bf P}_{0, 1}^{(1)}\tilde{Q}_{m_1}^a\cdot \tilde{Q}_{m_1}^a,
\end{eqnarray*}
which implies by using \eqref{eq-JLMS-d^3} that
% \begin{equation*}
%      |L_3^0| \lesssim \frac{k^4 \eta_0}{c_0}a^{3-h} d_{\rm out}^{-6} \cdot  \left|\int_{\Omega}\int_\Omega \Upsilon_0(x, y) \cdot \mathbb{P}_{1}^{(1)}\Lambda_1(y) \cdot\Lambda_1(x)\,dy\,dx\right| + \frac{k^4 \eta_0}{\pm c_0}a^{3-h}d_{\rm out}^{-6} \cdot d_{\rm out}^3 \cdot |{\bf P}_{0, 1}^{(1)}| \cdot \sum_{m=1}^\aleph |\tilde{Q}_{m_1}^a|^2,
% \end{equation*}
% by using the fact that 
% \begin{equation*}
%     \int_{B_{m_1}} \int_{B_{m_1}} \Upsilon_0(x, y)\,dy\,dx=-\frac{\pi d^3}{18} I_3.
% \end{equation*}
% see \cite[eq. (6.7)]{CGS}. Then
\begin{eqnarray}\label{es-L30-final}
% |L_3^0| &\lesssim& \frac{k^4 \eta_0}{c_0}a^{3-h} d_{\rm out}^{-6} \cdot  \left|\int_{\Omega}\int_\Omega \Upsilon_0(x, y) \cdot \mathbb{P}_{1}^{(1)}\Lambda_1(y) \cdot\Lambda_1(x)\,dy\,dx\right| + \frac{k^4 \eta_0}{\pm c_0}a^{3-h}d_{\rm out}^{-6} \cdot d_{\rm out}^3 \cdot |{\bf P}_{0, 1}^{(1)}| \cdot \sum_{m=1}^\aleph |\tilde{Q}_{m_1}^a|^2\notag\\
   |L_3^0|   &\lesssim& \frac{k^4 \eta_0}{c_0}a^{3-h} d_{\rm out}^{-6} \cdot  \left|\int_{\Omega} \nabla M( \mathbb{P}_{1}^{(1)} \Lambda_1)(x) \cdot\Lambda_1(x)\,dx\right| + \frac{k^4 \eta_0}{\pm c_0}a^{3-h}d_{\rm out}^{-3} \cdot |{\bf P}_{0, 1}^{(1)}| \cdot \sum_{m=1}^\aleph |\tilde{Q}_{m_1}^a|^2\notag\\
      &\lesssim& \frac{k^4 \eta_0}{c_0}a^{3-h} d_{\rm out}^{-6} \lVert\nabla M\rVert_{\mathcal{L}(\mathbb{L}^2(\Omega); \mathbb{L}^2(\Omega))} \cdot |\mathbb{P}_{1}^{(1)}|\cdot \lVert \Lambda_1\rVert_{\mathbb{L}^2(\Omega)}^2 + \frac{k^4 \eta_0}{\pm c_0}a^{3-h}d_{\rm out}^{-3} \cdot |{\bf P}_{0, 1}^{(1)}| \cdot \sum_{m=1}^\aleph |\tilde{Q}_{m_1}^a|^2\notag\\
      &\lesssim&  \frac{k^4 \eta_0}{ c_0} |{\bf P}_{0, 1}^{(1)}| a^{3-h}d_{\rm out}^{-3}  \cdot \sum_{m=1}^\aleph |\tilde{Q}_{m_1}^a|^2,
\end{eqnarray}
due to the fact that $ \lVert \nabla M\rVert_{\mathcal{L}(\mathbb{L}^2(\Omega); \mathbb{L}^2(\Omega))}=1$ and
\begin{equation*}
   \lVert \Lambda_1\rVert_{\mathbb{L}^2(\Omega)}=\int_\Omega \Lambda_1(x)\cdot\Lambda_1(x)\,dx=\sum_{m}\int_{B_{m_1}}\tilde{Q}_{m_1}^a\cdot\tilde{Q}_{m_1}^a\,dx=|B_{m_1}|\sum_m \tilde{Q}_{m_1}^a\cdot\tilde{Q}_{m_1}^a.
\end{equation*}

For $L_3^r$, with the help of \eqref{equiv-diff-Upsilon},
% since 
% \begin{equation}\notag
%     (\Upsilon_k-\Upsilon_0)(x, y)=\frac{1}{k^2} \underset{x}{\nabla} \underset{x}{\nabla} -\underset{x}{\nabla} \underset{x}{\nabla} \Phi_0 +\Phi_k(x, y) I_3 \sim \Phi_0(x, y) I_3,
% \end{equation}
similar to \eqref{es-invert-I-12}, we can get that
\begin{eqnarray}\label{es-L3r-final}
    % |L_3^r| &=& \frac{k^4 \eta_0}{c_0}a^{3-h}\left| \sum_{m=1}^\aleph \sum_{j=1\atop j\neq m}^\aleph  (\Upsilon_k-\Upsilon_0)(z_{m_1}, z_{j_1}){\bf P}_{0, 1}^{(1)}\tilde{Q}_{j_1}^a\cdot\tilde{Q}_{m_1}^a\right|\notag\\
    |L_3^r|&\lesssim& \frac{k^4 \eta_0}{c_0}a^{3-h} d_{\rm out}^{-6} \left| \sum_{m=1}^\aleph \sum_{j=1\atop j\neq m}^\aleph  \int_{B_{m_1}}\int_{B_{j_1}} \Phi_0(x, y)\,dy\,dx {\bf P}_{0, 1}^{(1)}\tilde{Q}_{j_1}^a\cdot\tilde{Q}_{m_1}^a\right|\notag\\
    &\lesssim& \frac{k^4 \eta_0}{c_0}a^{3-h} d_{\rm out}^{-6} \left| \int_{\Omega} N\left( \mathbb{P}_1^{(1)} \Lambda_1\right)(x)\cdot\Lambda_1(x)\,dx \right| + \frac{k^4 \eta_0}{c_0}a^{3-h} d_{\rm out}^{-1} |{\bf P}_{0, 1}^{(1)}| \sum_{m=1}^\aleph |\tilde{Q}_{m_1}^a|^2\notag\\
    &\lesssim& \frac{k^4 \eta_0}{c_0}a^{3-h} d_{\rm out}^{-6} \lVert N\rVert_{\mathcal{L}(\mathbb{L}^2(\Omega); \LL^2(\Omega))}|{\bf P}_{0, 1}^{(1)}| \cdot d_{\rm out}^3 \cdot \sum_{m=1}^\aleph |\tilde{Q}_{m_1}^a|^2\notag\\
    &\lesssim& \frac{k^4 \eta_0}{c_0}  |{\bf P}_{0, 1}^{(1)}| a^{3-h} d_{\rm out}^{-3} \sum_{m=1}^\aleph |\tilde{Q}_{m_1}^a|^2,
\end{eqnarray}
% by keeping the dominant terms, where we use the fact that
% \begin{equation*}
%     \int_{B_{m_1}}\int_{B_{m_1}}\Phi_0(x, y)\,dy\,dx =\frac{\pi d_{\rm out}^5}{60},
% \end{equation*}
% see {\color{red}{check the reference}}. 
Together with \eqref{es-L30-final}, we can obtain that
\begin{equation}\label{es-L3-final}
    |L_3|\lesssim |L_3^0|+|L_3^r|\lesssim \frac{k^4 \eta_0}{c_0} |{\bf P}_{0, 1}^{(1)}| a^{3-h} d_{\rm out}^{-3} \sum_{m=1}^\aleph |\tilde{Q}_{m_1}^a|^2.
\end{equation}

\medskip

\noindent (d) Estimate for 
\begin{equation*}
    L_4:= \sum_{m=1}^\aleph \sum_{j=1\atop j\neq m}^\aleph \frac{k}{i} a^3 \nabla\Phi_k(z_{m_1}, z_{j_1})\times {\bf P}_{0, 1}^{(1)} \tilde{R}_{j_1}^a\cdot \tilde{Q}_{m_1}^a.
\end{equation*}
With the expression \eqref{split-nabla-Phi-k} 
% Since 
% \begin{equation*}
%     \nabla \Phi_k(z_{m_1}, z_{j_1})=\nabla \Phi_0(z_{m_1}, z_{j_1})+(\nabla \Phi_k -\nabla\Phi_0)(z_{m_1}, z_{j_1})
% \end{equation*}
for $\nabla\Phi_0$ being harmonic, by using the notation \eqref{notation-integral form1}, there holds
\begin{eqnarray*}
    L_4 & = & \underbrace{\frac{k}{i} a^3 \sum_{m=1}^\aleph \sum_{j=1\atop j\neq m}^\aleph \frac{1}{|B_{m_1}|} \frac{1}{|B_{j_1}|}\int_{B_{m_1}}\int_{B_{j_1}} \nabla\Phi_0(x, y) \,dy\,dx \times {\bf P}_{0,1}^{(1)} \tilde{R}_{j_1}^a \cdot\tilde{Q}_{m_1}^a}_{L_4^0}\notag\\
    &+& \underbrace{\frac{k}{i} a^3 \sum_{m=1}^\aleph \sum_{j=1\atop j\neq m}^\aleph (\nabla\Phi_k-\nabla\Phi_0)(z_{m_1}, z_{j_1}) \times {\bf P}_{0,1}^{(1)} \tilde{R}_{j_1}^a \cdot\tilde{Q}_{m_1}^a}_{L_4^r}.
\end{eqnarray*}
Denote
\begin{equation}\notag
    \Pi_1(x)=\begin{cases}
        \tilde{R}_{m_1}^a & \mbox{in}\quad B_{m_1}\\
        0 & \mbox{Otherwise}.
    \end{cases}
\end{equation}
Then it directly gives
\begin{eqnarray}\label{es-L40}
    L_4^0 & = & \frac{36 k}{i\pi^2}  a^3  \int_{\Omega}\int_{\Omega} \nabla\Phi_0(x, y) \,dy\,dx \times \mathbb{P}_1^{(1)} \Pi_1(y)\cdot \Lambda_1(x) \,dy \,dx \notag\\
    &-&  \frac{36 k}{i \pi^2}  a^3 d_{\rm out}^{-6} \sum_{m=1}^\aleph \int_{B_{m_1}}\int_{B_{m_1}} \nabla\Phi_0(x, y) \,dy\,dx \times {\bf P}_{0,1}^{(1)} \tilde{R}_{m_1}^a \cdot\tilde{Q}_{m_1}^a\notag\\
    |L_4^0|&\lesssim& k a^3 d_{\rm out}^{-6} \left|\int_\Omega \nabla N\left( \mathbb{P}_1^{(1)} \Pi_1 \right)(x)\cdot\Lambda_1(x)\,dx \right| + k a^3 d_{\rm out}^{-6} \sum_{m=1}^\aleph\left| \int_{B_{m_1}} \nabla N(1)\,dx \times {\bf P}_{0,1}^{(1)} \tilde{R}_{m_1}^a \cdot\tilde{Q}_{m_1}^a \right|\notag\\
    &\lesssim& k a^3 d_{\rm out}^{-6} \lVert\nabla N \rVert_{\mathcal{L}(\LL^2(\Omega); \LL^2(\Omega))}\left|{\bf P}_{0, 1}^{(1)}\right| d_{\rm out}^3 \sum_{m=1} \left|\tilde{R}_{m_1}^a\cdot\tilde{Q}_{m_1}^a\right|+ k a^3 d_{\rm out}^{-6} \cdot d_{\rm out}^3 \cdot \left|{\bf P}_{0, 1}^{(1)}\right|\sum_{m=1} \left|\tilde{R}_{m_1}^a\cdot\tilde{Q}_{m_1}^a\right|\notag\\
    &\lesssim& k \left|{\bf P}_{0, 1}^{(1)}\right| a^3 d_{\rm out}^{-3}
    \left(\sum_{m=1}^\aleph |\tilde{R}_{m_1}^a|^2 \right)^\frac{1}{2}\left(\sum_{m=1}^\aleph |\tilde{Q}_{m_1}^a|^2 \right)^\frac{1}{2},\notag\\
\end{eqnarray}
by virtue of the fact that
\begin{equation*}
    \left|\int_{B_{m_1}}\nabla N(1)(x)\,dx \right| \lesssim |B_{m_1}|\lesssim d_{\rm out}^3, 
\end{equation*}
see \cite[eq. (1.13)]{GS}. For $L_4^r$, 
% with the expression
% \begin{equation}\notag
%     (\nabla\Phi_k -\nabla\Phi_0) (z_{m_1}, z_{j_1}) = \frac{ik}{4\pi} \int_0^1 i k t\cdot e^{i k t|z_{m_1}-z_{j_1} |} \,dt \cdot\frac{z_{m_1}-z_{j_1}}{|z_{m_1}-z_{j_1}|},
% \end{equation} 
we have
\begin{eqnarray}\label{es-L4r}
     |L_4^r| & \lesssim& k a^3 \left| \sum_{m=1}^\aleph \sum_{j=1\atop j\neq m}^\aleph (\nabla\Phi_k-\nabla\Phi_0)(z_{m_1}, z_{j_1}) \times {\bf P}_{0,1}^{(1)} \tilde{R}_{j_1}^a \cdot\tilde{Q}_{m_1}^a\right|\notag\\
     &\lesssim& k a^3 d_{\rm out}^{-3} \left|{\bf P}_{0,1}^{(1)}\right| \left(\sum_{j=1}^\aleph |\tilde{R}_{j_1}^a|^2 \right)^\frac{1}{2}\left(\sum_{m=1}^\aleph |\tilde{Q}_{m_1}^a|^2 \right)^\frac{1}{2}.
\end{eqnarray}\label{es-L4r}
Combining with  \eqref{es-L40} and \eqref{es-L4r}, we can get
\begin{equation}\label{es-L4-final}
    |L_4|\lesssim |L_4^0| +|L_4^r|\lesssim k \left|{\bf P}_{0,1}^{(1)}\right| a^3 d_{\rm out}^{-3}  \left(\sum_{j=1}^\aleph |\tilde{R}_{j_1}^a|^2 \right)^\frac{1}{2}\left(\sum_{m=1}^\aleph |\tilde{Q}_{m_1}^a|^2 \right)^\frac{1}{2}.
\end{equation}

\medskip

\noindent (e) Estimate for 
\begin{equation*}
    L_5:= \sum_{m=1}^\aleph \sum_{j=1\atop j\neq m}k^4 a^5 \Upsilon_k(z_{m_1}, z_{j_2}){\bf P}_{0, 2}^{(1)}\tilde{Q}_{j_2}^a\cdot \tilde{Q}_{m_1}^a.
\end{equation*}
For the dimer $D_{m}$ and $D_j$, by taking corresponding intermediate points $z_{m_0}$ and $z_{j_0}$, it is obvious that
\begin{equation*}
    L_5\sim \sum_{m=1}^\aleph \sum_{j=1\atop j\neq m}k^4 a^5 \Upsilon_k(z_{m_0}, z_{j_0}){\bf P}_{0, 2}^{(1)}\tilde{Q}_{j_2}^a\cdot \tilde{Q}_{m_1}^a =: \tilde{L}_5.
\end{equation*}
With the similar formulation of \eqref{def-L3} for the split of $\Upsilon_k(\cdot, \cdot)$, denoting 
\begin{equation}\label{notation-integral2}
    \Omega:=\underset{m=1}{\overset{\aleph}{\cup}}B_{m_0},\quad \mathbb{P}_2^{(1)}:=\begin{cases}
        {\bf P}_{0, 2}^{(1)} & \mbox{in}\; B_{m_0}\\
        0 &\mbox{Otherwise}
    \end{cases}\quad \mbox{and}\quad \Lambda_\ell(x):=\begin{cases}
        \tilde{Q}_{m_\ell}^a & \mbox{in}\; B_{m_0}\\
        0& \mbox{Otherwise}
    \end{cases}, \ell= 1,2,
\end{equation}
with $B_{m_0}$ being the ball with radius $d_{\rm out}/2$, centered at the intermediate point $z_{m_0}$. Then following the similar argument to \eqref{es-L30-final} and \eqref{es-L3r-final}, we can derive
\begin{eqnarray}\label{es-L5-final}
    |L_5|&\lesssim& k^4 a^5 d_{\rm out}^{-6} \left| \int_\Omega \nabla M\left(\mathbb{P}_2^{(1)}\Lambda_2\right)(x) \cdot \Lambda_1(x) \,dx\right| + k^4 a^5 d_{\rm out}^{-6} \left| \sum_{m=1}^\aleph \int_{B_{m_0}}\int_{B_{m_0}} \Upsilon_0(x, y)\,dy\,dx \cdot {\bf P}_{0, 2}^{(1)} \tilde{Q}_{j_2}^a\cdot \tilde{Q}_{m_1}^a \right|\notag\\
    &+& k^4 a^5 \left| \sum_{m=1}^\aleph \sum_{j=1\atop j\neq m}^\aleph (\Upsilon_k-\Upsilon_0)(z_{m_0}, z_{j_0}) {\bf P}_{0, 2}^{(1)}\tilde{Q}_{j_2}^a \cdot\tilde{Q}_{m_1}^a \right|\notag\\
    &\lesssim& k^4 \left|{\bf P}_{0, 2}^{(1)}\right| a^5  d_{\rm out}^{-3} \left( \sum_{m=1}^\aleph |\tilde{Q}_{m_2}^a|^2 \right)^\frac{1}{2}\left( \sum_{m=1}^\aleph |\tilde{Q}_{m_1}^a|^2 \right)^\frac{1}{2}.\notag\\
\end{eqnarray}

\medskip

\noindent (f) Estimate for 
\begin{equation}\notag
    L_6:= \sum_{m=1}^\aleph \sum_{j=1\atop j\neq m}^\aleph \frac{ k \eta_2}{\pm i d_0} a^{3-h} \nabla\Phi_k(z_{m_1}, z_{j_2})\times{\bf P}_{0, 2}^{(2)}\tilde{R}_{j_2}^a\cdot\tilde{Q}_{m_1}^a.
\end{equation}
Denote 
\begin{equation*}
    \mathbb{P}_{2}^2:=\begin{cases}
         {\bf P}_{0, 2}^{(2)} & \mbox{in} \; B_{m_0}\\
         0 & \mbox{Otherwise}
    \end{cases} \quad \mbox{and}\quad 
    \Pi_2:=\begin{cases}
        \tilde{R}_{m_2}^a & \mbox{in} \; B_{m_0}\\
        0 & \mbox{Otherwise}
    \end{cases},
\end{equation*}
with $B_{m_0}$ given in \eqref{notation-integral2}. Then following the similar argument to \eqref{es-L40} and \eqref{es-L4r}, we can obtain that
\begin{eqnarray}\label{es-L6-final}
    |L_6| &\lesssim& \frac{ k \eta_2}{ d_0} a^{3-h} \left|\sum_{m=1}^\aleph \sum_{j=1\atop j\neq m}^\aleph   \nabla\Phi_k(z_{m_0}, z_{j_0})\times{\bf P}_{0, 2}^{(2)}\tilde{R}_{j_2}^a\cdot\tilde{Q}_{m_1}^a\right|\notag\\
    &\lesssim& \frac{ k \eta_2}{ d_0} a^{3-h} d_{\rm out}^{-6} \Bigg[ \left| \int_\Omega \nabla N\left(\mathbb{P}_2^{(2)}\Pi_2\right)(x)\Lambda_1(x)\,dx \right| + \sum_{m=1}^\aleph \left| \int_{B_{m_0}} \int_{B_{m_0}} \nabla\Phi_0(x, y)\,dy\,dx \times {\bf P}_{0, 2}^{(2)}\tilde{R}_{m_2}^a\cdot \tilde{Q}_{m_1}^a \right|\notag\\
    &+& \sum_{m=1}^\aleph \sum_{m=1}^\aleph \left| (\nabla\Phi_k -\nabla \Phi_0)(z_{m_0}, z_{j_0}) \times {\bf P}_{0, 2}^{(2)}\tilde{R}_{j_2}^a\cdot \tilde{Q}_{m_1}^a\right| \Bigg]\notag\\
    &\lesssim&  \frac{ k \eta_2}{ d_0} a^{3-h} \left| {\bf P}_{0, 2}^{(2)} \right| \left( d_{\rm out}^{-3} \left( \lVert \nabla N\rVert_{\mathcal{L}(\LL^2(\Omega); \LL^2(\Omega))} + 1\right)\sum_{m=1}^\aleph |\tilde{R}_{m_2}^a|\cdot |\tilde{Q}_{m_1}^a| + \aleph \left(\sum_{j=1}^\aleph|\tilde{R}_{j_2}^a|^2\right)^\frac{1}{2} \left(\sum_{m=1}^\aleph|\tilde{Q}_{m_1}^a|^2\right)^\frac{1}{2}\right)\notag\\
    &\lesssim& \frac{ k \eta_2}{ d_0} \left| {\bf P}_{0, 2}^{(2)} \right| a^{3-h}  d_{\rm out}^{-3} \left(\sum_{j=1}^\aleph|\tilde{R}_{j_2}^a|^2\right)^\frac{1}{2} \left(\sum_{m=1}^\aleph|\tilde{Q}_{m_1}^a|^2\right)^\frac{1}{2}.\notag\\
\end{eqnarray}
Combining with \eqref{es-L1}, \eqref{es-L2}, \eqref{es-L3-final}, \eqref{es-L4-final}, \eqref{es-L5-final} and \eqref{es-L6-final}, utilizing the Cauchy-Schwartz inequality, we can deduce from \eqref{Qm1-quadratic} that 
\begin{eqnarray}\label{es-Qm1a-form}
&&\left(1- \frac{k^4 \eta_0}{c_0} \left|{\bf P}_{0. 1}^{(1)}\right| a^{3-h} d_{\rm out}^{-3} \right) \left( \sum_{m=1}^\aleph |\tilde{Q}_{m_1}^a| \right)^\frac{1}{2}\notag\\
&\lesssim& k \left| {\bf P}_{0, 1}^{(1)} \right| a^3 d_{\rm out}^{-3} \left( \sum_{m=1}^\aleph |\tilde{R}_{m_1}^a|^2 \right)^\frac{1}{2} + \frac{k \eta_2}{d_0} \left|{\bf P}_{0, 2}^{(2)} \right| \left(a^{3-h} d_{\rm in}^{-2}+ a^{3-h} d_{\rm out}^{-3} \right) \left( \sum_{m=1}^\aleph |\tilde{R}_{m_2}^a|^2 \right)^\frac{1}{2}\notag\\
&+& k^4 \left|{\bf P}_{0, 2}^{(1)} \right| a^{5} d_{\rm in}^{-3} \left( \sum_{m=1}^\aleph |\tilde{Q}_{m_2}^a|^2 \right)^\frac{1}{2} + \left( \sum_{m=1}^\aleph |L_{H, m_1}|^2 \right)^\frac{1}{2},
\end{eqnarray}
by keeping the dominant terms.

With the similar arguments presented above, we can obtain the estimation respectively for $\tilde{R}_{m_1}^a, \tilde{Q}_{m_2}^a$ and $\tilde{R}_{m_2}^a$ as
\begin{eqnarray}\label{es-Rm1a-form}
    \left( \sum_{m=1}^\aleph |\tilde{R}_{m_1}^a|^2 \right)^\frac{1}{2} &\lesssim& k \left|{\bf P}_{0, 2}^{(1)}\right| \left( a^5 d_{\rm in}^{-2} + k^2 a^5 d_{\rm out}^{-3} \right)\left( \sum_{m=1}^\aleph |\tilde{Q}_{m_2}^a|^2 \right)^\frac{1}{2} \notag\\
    &+& \frac{k^2 \eta_2}{d_0} \left|{\bf P}_{0, 2}^{(2)}\right| \left( a^{3-h} d_{\rm in}^{-3} + a^{3-h} d_{\rm out}^{-3} \right)\left( \sum_{m=1}^\aleph |\tilde{R}_{m_2}^a|^2 \right)^\frac{1}{2}\notag\\
    &+& \frac{k^3 \eta_0}{c_0} \left| {\bf P}_{0, 1}^{(1)}\right| a^{3-h} d_{\rm out}^{-3} \left( \sum_{m=1}^\aleph |\tilde{Q}_{m_1}^a|^2 \right)^\frac{1}{2} + \left( |L_{E, m_1}|^2 \right)^\frac{1}{2},
\end{eqnarray}
\begin{eqnarray}\label{es-Qm2a-form}
    \left( \sum_{m=1}^\aleph |\tilde{Q}_{m_2}^a|^2 \right)^\frac{1}{2} &\lesssim& \frac{k^4 \eta_0 \eta_2}{c_0} \left|{\bf P}_{0, 1}^{(1)}\right| \left( a^{3-h} d_{\rm in}^{-3} + a^{3-h} d_{\rm out}^{-3} \right)\left( \sum_{m=1}^\aleph |\tilde{Q}_{m_1}^a|^2 \right)^\frac{1}{2} \notag\\
    &+& k \eta_2 \left|{\bf P}_{0, 1}^{(2)}\right| \left( k^2 a^{3} d_{\rm in}^{-2} + a^{3} d_{\rm out}^{-3} \right)\left( \sum_{m=1}^\aleph |\tilde{R}_{m_1}^a|^2 \right)^\frac{1}{2}\notag\\
    &+& \frac{k \eta_2^2}{d_0} \left| {\bf P}_{0, 2}^{(2)}\right| a^{3-h} d_{\rm out}^{-3} \left( \sum_{m=1}^\aleph |\tilde{R}_{m_2}^a|^2 \right)^\frac{1}{2} + \left( |L_{H, m_2}|^2 \right)^\frac{1}{2},
\end{eqnarray}
and 
\begin{eqnarray}\label{es-Rm2a-form}
    &&\left(1- \frac{k^2 \eta_2}{d_0} \left|{\bf P}_{0. 2}^{(2)}\right| a^{3-h} d_{\rm out}^{-3} \right) \left( \sum_{m=1}^\aleph |\tilde{R}_{m_2}^a| \right)^\frac{1}{2}\notag\\
&\lesssim& \frac{k^3 \eta_0}{c_0} \left| {\bf P}_{0, 1}^{(1)} \right| \left(a^{3-h} d_{\rm in}^{-2}+ a^{3-h} d_{\rm out}^{-3} \right) \left( \sum_{m=1}^\aleph |\tilde{Q}_{m_1}^a|^2 \right)^\frac{1}{2} + k^2 \left|{\bf P}_{0, 1}^{(2)} \right| \left(a^{3-h} d_{\rm in}^{-3}+ a^{3} d_{\rm out}^{-3} \right) \left( \sum_{m=1}^\aleph |\tilde{R}_{m_1}^a|^2 \right)^\frac{1}{2}\notag\\
&+& k^3 \left|{\bf P}_{0, 2}^{(1)} \right| a^{5} d_{\rm out}^{-3} \left( \sum_{m=1}^\aleph |\tilde{Q}_{m_2}^a|^2 \right)^\frac{1}{2} + \left( \sum_{m=1}^\aleph |L_{E, m_2}|^2 \right)^\frac{1}{2}.\notag\\
\end{eqnarray}

Therefore, from \eqref{es-Qm1a-form}, \eqref{es-Rm1a-form}, \eqref{es-Qm2a-form} and \eqref{es-Rm2a-form}, it indicates that under the condition
\begin{equation}\label{condition-1-es-Qm1-Rm2}
    4-h-4t_1>0\quad\mbox{and}\quad \max\left\{ \frac{k^4 \eta_0}{c_0} \left|{\bf P}_{0, 1}^{(1)}\right|, \frac{k^2 \eta_2}{d_0} \left|{\bf P}_{0, 2}^{(2)}\right| \right \} a^{3-h} d_{\rm out}^{-3} < 1,
\end{equation}
there holds
\begin{eqnarray}\label{es-Qm1a}
    \left( \sum_{m=1}^\aleph |\tilde{Q}_{m_1}^a|^2 \right)^\frac{1}{2} &\lesssim & \left( \sum_{m=1}^\aleph |L_{H, m_1}|^2 \right)^\frac{1}{2} + \frac{k \eta_2}{d_0} \left|{\bf P}_{0, 2}^{(2)}\right| a^{3-h} d_{\rm out}^{-3} \left(\sum_{m=1}^\aleph |L_{E, m_2}|^2\right)^\frac{1}{2}\notag\\
    &+& \left( \frac{k^4 \eta_2}{d_0} \left|{\bf P}_{0, 2}^{(2)}\right|\cdot \left|{\bf P}_{0, 2}^{(1)}\right| a^{8-h} d_{\rm out}^{-6} + k^4 \left|{\bf P}_{0, 2}^{(1)}\right| a^{5} d_{\rm in}^{-3} \right) \left( \sum_{m=1}^\aleph |L_{H, m_2}|^2 \right)^\frac{1}{2}\notag\\
    &+& \left( \frac{k^3 \eta_2}{d_0} \left|{\bf P}_{0, 2}^{(2)}\right|\cdot \left|{\bf P}_{0, 1}^{(2)}\right| \left(a^{6-h} d_{\rm in}^{-5} + a^{6-h} d_{\rm in}^{-3} d_{\rm out}^{-3} \right) + k \left|{\bf P}_{0, 1}^{(1)}\right| a^{3} d_{\rm out}^{-3} \right) \left( \sum_{m=1}^\aleph |L_{E, m_1}|^2 \right)^\frac{1}{2} 
\end{eqnarray}
and 
\begin{eqnarray}\label{es-Rm2a}
    \left( \sum_{m=1}^\aleph |\tilde{R}_{m_2}^a|^2 \right)^\frac{1}{2} &\lesssim & \left( \sum_{m=1}^\aleph |L_{E, m_2}|^2 \right)^\frac{1}{2} + \frac{k^3 \eta_0}{c_0} \left|{\bf P}_{0, 1}^{(1)}\right| a^{3-h} d_{\rm out}^{-3} \left(\sum_{m=1}^\aleph |L_{H, m_1}|^2\right)^\frac{1}{2}\notag\\
    &+& \left( \frac{k^7 \eta_0}{c_0} \left|{\bf P}_{0, 1}^{(1)}\right|\cdot \left|{\bf P}_{0, 2}^{(1)}\right| (a^{8-h} d_{\rm in}^{-5} + a^{8-h} d_{\rm in}^{-3} d_{\rm out}^{-3}) +k^3 \left|{\bf P}_{0, 2}^{(1)}\right| a^5 d_{\rm out}^{-3}\right) \left( \sum_{m=1}^\aleph |L_{H, m_2}|^2 \right)^\frac{1}{2}\notag\\
    &+& \left( \frac{k^4 \eta_0}{c_0} \left|{\bf P}_{0, 1}^{(1)}\right|^2  a^{6-h} d_{\rm out}^{-6} + k^2 \left|{\bf P}_{0, 1}^{(2)}\right| a^{3} d_{\rm in}^{-3} \right) \left( \sum_{m=1}^\aleph |L_{E, m_1}|^2 \right)^\frac{1}{2}.
\end{eqnarray}
Note that \eqref{condition-1-es-Qm1-Rm2} can be fulfilled for (MP1.21) and (MP1.25). Recall the notation \eqref{notation-ABU}, which directly states 
\begin{equation*}
    \tilde{Q}_{m_1}= \frac{i k \eta_0}{\pm c_0} a^{3-h} {\bf P}_{0, 1}^{(1)} \tilde{Q}_{m_1}^a\quad \mbox{and}\quad 
    \tilde{R}_{m_2} = \frac{\eta_2}{\pm d_0} a^{3-h} {\bf P}_{0, 2}^{(2)} \tilde{R}_{m_2}^a,
\end{equation*}
then we can justify the $\ell_2$-estimates for $\tilde{Q}_{m_1}$ and $\tilde{R}_{m_2}$ from \eqref{es-Qm1a} and \eqref{es-Rm2a} given by \eqref{es-Qm-main} and \eqref{es-Rm-main}, respectively, in Proposition \ref{prop-es-QR-domain}.